%% file: rewriting_affine_brauer.tex
\documentclass[11pt]{article}

\usepackage{defs/packages}
\usepackage{defs/macros}

\newcommand{\br}{1}
\newcommand{\h}{1}
\newcommand{\sca}{0.5}
\usepackage{defs/diagrammatics}
\tikzset{critshading1/.style={rounded corners=2pt,inner sep=.7mm,draw=black,dashed}}

\usepackage{defs/custom-biblatex-style}
\usepackage[a4paper, total={6in, 9in}]{geometry}

\title{Monoidal Gröbner systems\\and categories of affine Brauer type}

\author{Sigiswald Barbier \orcidlink{0000-0002-5310-882X}\thanks{Ghent University, Ghent, Belgium. \textit{Email:} \texttt{sigiswald.barbier@ugent.be}.} \and Léo Schelstraete \orcidlink{0000-0001-7167-3964}\thanks{Max Planck Institute for Mathematics, Bonn, Germany. \textit{Email:} \texttt{leoschelstraete@gmail.com}.}}

\begin{document}

\maketitle


\begin{abstract}
  We introduce an analogue of Gröbner bases, or equivalently, Bergman's diamond lemma, for linear (super)monoidal categories.
  It gives a systematic way to study monoidal ideals and to prove basis theorems.
  The central notion is that of a \emph{monoidal Gröbner system} and is based on higher linear rewriting theory.
  We then apply the theory to \emph{categories of affine Brauer type}: linear (super)monoidal categories that have the same hom-basis as the affine Brauer category, but possibly distinct composition and tensor product.
  Our main result is a criterion for a category to be of affine Brauer type, reducing the question to an explicit list of local computations, which we partially implement in the computer algebra system \textsc{Form}.
  This gives new combinatorial proofs of the basis theorems for the affine Brauer category, the nil-Brauer category and the affine VW supercategory, and yields new examples.
  In particular, we construct the \emph{odd nil-Brauer category}, a conjectural supercategorification of the split $\imath$quantum group of rank one, and the \emph{quantized affine VW supercategory}, whose non-affine part is the quantized periplectic Brauer category and which we expect to satisfy a higher quantum Schur--Weyl duality.
  Finally, we classify the categories of affine Brauer type admitting a sufficiently simple presentation.
\end{abstract}

\begin{quotation}

\end{quotation}

\tableofcontents


\newpage
\input{sections/introduction}

\input{sections/rewriting}

\input{sections/definitions}

\input{sections/confluence_analysis_general.tex}

\section{Known examples}
\label{sec:examples}

\input{sections/analysis_affine_Brauer_new}

\input{sections/analysis_nilBrauer_new}

\input{sections/analysis_affine_VW_supercategory_new}

\section{Applications}
\label{sec:applications}
\input{sections/oddnilBrauer}
\input{sections/quantization}
\input{sections/classification_nonquantum}
\appendix

\input{sections/proof_normal_forms}

\input{sections/appendix_symmetries}





\printbibliography[heading=bibintoc]

\end{document}

%% file: sections/introduction.tex
\section{Introduction}

Linear monoidal categories are ubiquitous in representation theory and related fields.
In this context, they are colloquially known as ``diagrammatic algebras'', due to the use of string diagrammatics to picture their morphisms.
In this introduction, a \emph{diagrammatic category} is a linear (strict) monoidal category defined by generators and relations and a \emph{diagram} is any morphism obtained by vertical and horizontal compositions of the generators.

As for associative algebras, bases in diagrammatic categories play a central role:
\begin{quote}
  \textsc{Problem:} given a diagrammatic category, how do we find a \emph{hom-basis}, that is, a basis for each hom-space?
\end{quote}
In practice, basis theorems for diagrammatic categories are shown by exhibiting a sufficiently faithful monoidal representation.
This representation is often quite involved; moreover, since it is specific to the case at hand, it is often not clear how it can be adapted to other similarly-presented diagrammatic categories.

\emph{Gröbner bases}\footnote{The terminology is confusing: a Gröbner basis is a basis of the \emph{ideal} presenting the algebra, satisfying some properties. Given a Gröbner basis, one can extract a basis for the algebra itself, thanks to the properties of the Gröbner basis.} first arose in Buchberger's work on the ideal-membership problem in polynomial algebras \cite{Buchberger_AlgorithmusAuffindenBasiselemente_1965}. Related ideas appeared in the independent work of Shirshov \cite{Shirshov_AlgorithmicProblemsLie_1962} (and later Bokut \cite{Bokut_ImbeddingsSimpleAssociative_1976}) and Bergman \cite{Bergman_DiamondLemmaRing_1978}; the latter is known as \emph{Bergman's diamond lemma}.
Gröbner bases have been generalized to various settings (\eg operads \cite{DK_GrobnerBasesOperads_2010}) and can be used to study various homological properties; see
\cite{%
BD_AlgebraicOperads_2016,%
BCK+_GrobnerShirshovBases_2020,
} for an overview.
These theories are instances of theories of confluent systems; as such, they fall into the more general principle of \emph{rewriting theory}, a standard framework in applied category theory and related fields.
The rewriting viewpoint for associative algebras was developed in \cite{GHM_ConvergentPresentationsPolygraphic_2019}.
Rewriting techniques were also extended to higher categories in a line of work by various authors
\cite{
  Lafont_AlgebraicTheoryBoolean_2003,
  Guiraud_TerminationOrdersThreedimensional_2006,
  GM_HigherdimensionalCategoriesFinite_2009,
  Mimram_3dimensionalRewritingTheory_2014,
  GM_PolygraphsFiniteDerivation_2018,
  Metayer_ResolutionsPolygraphs_2003,%
  Guetta_HomologyCategoriesPolygraphic_2021,%
}%
; see \cite{ABG+_PolygraphsRewritingHigher_2025} for a recent monograph.

More recently, Alleaume and Dupont suggested an approach to combine the linear and higher settings, toward applications to diagrammatic categories \cite{%
  Alleaume_RewritingHigherDimensional_2018,%
  Dupont_RewritingModuloIsotopies_2021,%
  Dupont_RewritingModuloIsotopies_2022,%
  DEL_SuperRewritingTheory_2024,%
  DN_CategorificationInfinitedimensional$mathfraksl_2$modules_2021,%
}.
In independent work, Elias studied ``Hecke-type categories'' in \cite{Elias_DiamondLemmaHecketype_2022}, building directly on Bergman's diamond lemma.
In \cite{Schelstraete_RewritingModuloDiagrammatic_2026}, the second author developed a rewriting theory suitable for diagrammatic categories, distinct from Alleaume's and Dupont's approaches.
This approach is coined \emph{higher linear rewriting theory}, or \emph{linear Gray rewriting theory} in the context of linear 2-categories; an alternative terminology could be \emph{higher Gröbner theory}.
The theory was applied in \opcit to a conjecture in categorification; at the time of writing, rewriting theory is the only known technique to solve this conjecture.
Another application of the theory appeared in \cite{JLS+_BasisSchurWeylDuality_2025}.

\medbreak

In this article, we apply higher linear rewriting theory to a certain family of diagrammatic categories modelled on the affine Brauer category \cite{RS_AffineBrauerCategory_2019}.
In a nutshell:
\begin{quote}
  \textsc{Main Theorem:} a criterion to determine if a category has the same hom-basis as the affine Brauer category.
  Such categories are called ``categories of affine Brauer type''.
\end{quote}
This gives a combinatorial and unified approach to proving basis theorems for diagrammatic categories whose presentation is modelled on the presentation of the affine Brauer category. As applications, we:
\begin{enumerate}[(i)]
  \item give a combinatorial proof of the basis theorems of the affine Brauer category \cite{RS_AffineBrauerCategory_2019}, the nil-Brauer category \cite{BWW_NilBrauerCategory_2024} and the affine VW supercategory \cite{BDE+_AffineVWSupercategory_2020};
  \item construct the \emph{odd nil-Brauer category} (\cref{defn:oddnilBr}), an odd analogue to the nil-Brauer category.
  Following the fact that the latter categorifies the split $\imath$quantum group of rank one $U^\imath_q(\mathfrak{sl}_2)$ \cite{BWW_NilBrauerCategorifiesSplit_2025}, we conjecture that the odd nil-Brauer gives a \emph{super}categorification of $U^{\imath}_q(\mathfrak{sl}_2)$ (see \cref{subsubsec:intro_odd_nilBrauer});
  \item construct the \emph{quantized affine VW supercategory} (\cref{defn:qsVW_form}), a quantization of the affine VW supercategory \cite{BDE+_AffineVWSupercategory_2020} whose non-affine part is the periplectic $q$-Brauer category \cite{RS_Periplectic$q$BrauerCategory_2025}, conjecturally related to some higher quantum Schur--Weyl duality (see \cref{subsubsec:intro-quantVW});
  \item Classify of categories of affine Brauer type with a sufficiently simple presentation in \cref{classification simple-minded} (see \cref{subsubsec:intro_classification}).
\end{enumerate}
Along the way, we further develop the higher linear rewriting toolbox, \eg with a way to classify overlapping branchings (\cref{lem:RW_overlapping_branchings_intersect}), with the notion of ``tamed preconfluence'' (\cref{subsec:tamed_preconfluence}) or with the use of symmetries (\cref{subsec:local_analysis}).
Verifying the criterion of the Main Theorem largely amounts to direct computations, which we implement in the computer algebra system \textsc{Form} \cite{FORM4,FORM5}; see \cref{subsec:implementation}.

\medbreak

\paragraph{Organization.}
\Cref{sec:RW_monoidal_grobner_system} introduces the foundation of higher Gröbner theory. \Cref{sec:definitions} gives the main definitions in the context of categories affine Brauer type and state the Main Theorem precisely, while \cref{sec:critical_branchings} gives the proof of the Main Theorem.
The reader mainly interested in the underlying theory and the proof of the Main Theorem can focus on these three sections.
\Cref{sec:examples} applies the Main Theorem to three examples already known in the literature.
The case of the affine Brauer category (\cref{subsec:affine_brauer_category}) is explained in greater details; the other cases follow a similar strategy.
The reader eager to see the Main Theorem in action can start with \cref{subsec:affine_brauer_category}, occasionally looking back at \cref{sec:RW_monoidal_grobner_system} and \cref{sec:definitions} for definitions.
Finally, \cref{sec:applications} describes novel results obtained from the Main Theorem.

The rest of the introduction gives an extended summary of the article.

\paragraph{A note on the literature.}
In \cite{Alleaume_RewritingHigherDimensional_2018}, Alleaume studied the oriented affine Brauer category using his approach of rewriting theory in linear monoidal categories.
The work of the second author in \cite{Schelstraete_RewritingModuloDiagrammatic_2026}, while inspired by \cite{Alleaume_RewritingHigherDimensional_2018} and other works, is distinct from Alleaume's work in various aspects, and so is the present work; see \cite[Section~3.7]{Schelstraete_RewritingModuloDiagrammatic_2026} for further details.
An important distinction is the use of context-dependent rewriting rules to ensure termination, while Alleaume uses a non-terminating system.
Following \cite[Remark~3.70]{Schelstraete_RewritingModuloDiagrammatic_2026}, we do not see how a basis can be extracted from a non-terminating system.

Another predecessor of our work is \cite{Elias_DiamondLemmaHecketype_2022}.
This work develops a ``diamond lemma for linear monoidal categories'', generalizing directly Bergman's diamond lemma \cite{Bergman_DiamondLemmaRing_1978}.
The main application is a criterion, similar in spirit to ours, for ``categories of Hecke type'', that is, categories having the same basis as quiver Hecke algebras \cite{KL_DiagrammaticApproachCategorification_2009,Rouquier_2KacMoodyAlgebras_2008}. Up to extanding our Main Theorem to allow different generating objects (or ``colours''), this criterion can be seen as a special case of ours.\footnote{In fact, we improve that criterion: as we show in \cref{subsubsec:horizontal_symmetry} (see in particular the first diagram in \cref{fig:horizontal_symmetries_B}), one of the ``overlaps'' from \cite{Elias_DiamondLemmaHecketype_2022} follow from other overlaps.}
The reduction algorithm for categories of Hecke type can be defined using only local rewriting rules, so that a theory in the spirit of Bergman's diamond lemma is sufficient to obtain the criterion in \cite{Elias_DiamondLemmaHecketype_2022}.
In contrast, categories of affine Brauer type require global constrains on rewriting rules. This is the main qualitative difference between the two settings, and motivate why a different theoretical foundation is needed.
See the rest of the introduction for further details.

\paragraph{A note on future directions.}
First, this work defines new examples of categories of affine Brauer type, each with a conjectural representation-theoretic interpretation; the natural next step is to further study these examples.
Second, this work only gives some applications of the Main Theorem. One could look for further examples, for instance by extending the classification result given in \cref{thm:classification_as_lin_categories} to more involved presentations.
Third, this article gives a partial implementation of a reduction-to-basis algorithm for categories of affine Brauer type.
A complete implementation could allow efficient computations in these categories.
Fourth, many diagrammatic categories have similar generating morphisms as the affine Brauer category, with similar basis theorem or basis conjecture; for instance, the basis theorem \cite{Webster_UnfurlingKhovanovLaudaRouquierAlgebras_2024} of Kac--Moody 2-categories \cite{KL_CategorificationQuantum$sln$_2010,Rouquier_2KacMoodyAlgebras_2008}, the basis conjecture of super Kac--Moody 2-categories \cite{BE_SuperKacMoody_2017} or the basis conjecture for 2-iquantum groups \cite{BWW_CategorificationQuasisplitIquantum_2025}.
We expect that the techniques developed in this article can be applied to these diagrammatic categories.

\paragraph{Acknowledgment.} L.S.\ thanks Jon Brundan for pointing out the example of the nil-Brauer category \cite{BWW_NilBrauerCategory_2024} to him, and Alexei Davydov, Yves Guiraud, Philippe Malbos and Emmanuel Wagner for discussions related to that project. The authors thank Alexis Langlois-Rémillard for pointing out the software \textsc{form} \cite{FORM4,FORM5} to them.
L.S.\ was supported by the Max Planck Institute for Mathematics. S.B.\ thanks the Max Planck Institute for Mathematics for their hospitality during his visit there.

\paragraph{Use of AI.} We used AI for tasks such as proofreading the text, literature search and assistance in writing \textsc{Form} code.
AI further contributed to the proof of the bubble-slide lemma in \cref{alt:alt-subsec:affBr_bubble_slides}, as follows.
We first wrote a complete proof without AI, using explicit formulas.
While experimenting with Claude Fable 5 in early July 2026, we were led to the idea of using generating functions to simplify the argument; a subsequent search of the literature showed the idea to be due to Savage and Webster \cite{SW_BubblesAffineBrauer_2024}, whom we credit.
AI was not used in any other mathematically relevant way.
The authors adhere to the \href{https://leidendeclaration.ai/}{Leiden Declaration on Artificial Intelligence and Mathematics}.

\subsection{Categories of affine Brauer type}
We consider a certain family of diagrammatic categories modelled on the \emph{affine Brauer category $\affBr$} \cite{RS_AffineBrauerCategory_2019}.
This category appears in the study of the representation theory of the type BCD Lie algebras and is generated by a \emph{crossing}, a \emph{cup}, a \emph{cap} and a \emph{dot}:
\begin{IEEEeqnarray}{CcCcCcCcC}
  \label{eq:intro_basic_diagrams}
  \tikzpic{\cro}
  &\quad&
  \tikzpic{\ca}
  &\quad&
  \tikzpic{\cu}
  &\quad&
  \tikzpic{\bli}
  &\mspace{50mu}&
  \tikzpic{\bubble[0][0][k]}
  \\\nonumber
  \text{\small crossing}
  &&
  \text{\small cup}
  &&
  \text{\small cap}
  &&
  \text{\small dot}
  &&
  \text{\small $k$-bubble}
\end{IEEEeqnarray}
The diagram on the right-hand side above, the \emph{$k$-bubble} for $k\in\bN$ a natural number, is obtained by composing a cup, $k$ dots and a cap.
It is shown in \cite{RS_AffineBrauerCategory_2019} that the affine Brauer category admits a hom-basis given by (a choice of representatives of) certain diagrams called \emph{reduced matchings}; see \cref{defn:end_X_reduced_matching} for a formal definition.
In a reduced matching, the only closed loops are $k$-bubbles sitting on the leftmost region of the diagram.
In the context of \cite{RS_AffineBrauerCategory_2019}, these bubbles can only carry a non-zero and even number of dots; we shall say that the affine Brauer category has bubbles $\bubbleset=2\bN_{>0}$ (here $\bN_{>0}\coloneqq\bN\setminus\{0\}$).
The affine Brauer category is the stereotypical example for the following family of diagrammatic categories:

\begin{definition}[\cref{defn:category_affine_brauer_type}]
  Let $\bubbleset$ be a subset of $\bN_{>0}$.
  Let\footnote{The pre-superscript $\bullet$ in the notation $\affC$ indicates that the discussion is about affine Brauer type, in contrast to the notation without superscript $\cC$ for a general discussion about diagrammatic categories.} $\affC$ be a linear monoidal category whose set of objects is $\bN$ and whose morphisms are generated by a crossing, a cup, a cap and a dot as in \eqref{eq:intro_basic_diagrams}.
  We say that $\affC$ is \defnemph{of affine Brauer type with $\bubbleset$-bubbles} if reduced matchings whose bubbles have dots in $\bubbleset$ define a hom-basis for $\affC$.
\end{definition}

In other words, a category of affine Brauer type is, roughly, a linear monoidal category with the same hom-spaces as the affine Brauer category, but with possibly different composition and tensor product.
This is the higher analogue of having two algebra structures on the same underlying module.
In practice, there are many unrelated algebra structures on a given module, and studying their moduli space is out of reach.
However, the notion becomes more tractable when going higher, in part thanks to the additional constraints coming from the choice of shapes (\ie domain and codomain) for the generators.

At least two other examples of categories of affine Brauer type exist in the literature: the \emph{nil-Brauer category $\nilBr[\VARev]$} \cite{BWW_NilBrauerCategory_2024} and the \emph{affine VW supercategory $\sVW$} \cite{BDE+_AffineVWSupercategory_2020}.
The nil-Brauer category $\nilBr[\VARev]$ appears as a categorification of the split $\imath$quantum group of rank one \cite{BWW_NilBrauerCategorifiesSplit_2025}; it is a category of affine Brauer type with bubbles $\bubbleset=2\bN+1$.
The affine VW supercategory $\sVW$ appears in the study of representations of the periplectic Lie superalgebra; it is a category of affine Brauer type with bubbles $\bubbleset=\emptyset$. In fact, the category $\sVW$ is not a monoidal category but a \emph{super}monoidal category \cite{BE_MonoidalSupercategories_2017}: its interchange law is twisted by a sign.
While we restrict this introduction to usual linear monoidal categories, the results of this article also apply to supermonoidal categories, and more generally to categories where the interchange is twisted by a scalar; see \cref{subsec:RW_graded_monoidal_categories}.

In each of these three papers \cite{RS_AffineBrauerCategory_2019,BWW_NilBrauerCategory_2024,BDE+_AffineVWSupercategory_2020}, the main result is the basis theorem, that is, the fact that the category is of affine Brauer type with the corresponding choice of bubbles; in each case, the proof uses a faithful monoidal representation.

\subsection{The Main Theorem}

\subsubsection{Rewriting systems of affine Brauer type}
Our main result is a criterion deciding when a diagrammatic category is of affine Brauer type.
As a preliminary, one needs to present the category in a particular way; for the affine Brauer category $\affBr$, this particular presentation is given in \cref{fig:AB_S_rewriting_steps}.
Each arrow $v\to w$ in \cref{fig:AB_S_rewriting_steps} encodes an equality $v=w$ in $\affBr$.
In general, one can consider any variation of these relations by adding ``lower-order terms'', where ``lower'' is with respect to a certain partial order $\rhd$ on diagrams; see \cref{subsec:RW_termination_order}.
For instance, one may consider a presentation where the first two relations have the following general form:
\begin{gather*}
  \tikzpic{
    \cro[0][1]
    \cro 
  } 
  \;\overset{\stRtwo}{\longrightarrow}\;
  \tikzpic{\LOTRtwo}
  \mspace{30mu}
  \an
  \mspace{30mu}
  %
  \tikzpic{
    \cro[0][2] \li[2][2]
    \li[0][1] \cro[1][1] 
    \cro \li[2][0] 
  }[scale=.7]
  \;\overset{\stRthree}{\longrightarrow}\;
  \VARRthree\;
  \tikzpic{
    \cro[1][2] \li[0][2]
    \li[2][1] \cro[0][1] 
    \cro[1][0] \li[0][0] 
  }[scale=.7]
  \;+\;
  \tikzpic{\LOTRthree}[scale=.7]\;.
\end{gather*}
Here, $\VARRthree$ is an invertible scalar and the double-lined diagrams denote a linear combination of diagrams $D$ such that
$\tikzpic{\cro\cro[0][1]}[scale=.5]\rhd D$
and
$
\tikzpic{
  \cro[0][2] \li[2][2]
  \li[0][1] \cro[1][1] 
  \cro \li[2][0] 
}[scale=.4]\rhd D$,
respectively; in these two cases, the condition boils down to $D$ having strictly less than two or three crossings, respectively.
For instance, each of the following diagrams could appear in the linear decomposition of $\tikzpic{\LOTRtwo}[scale=.6]$:
\begin{gather*}
  \tikzpic{\cro}[scale=.8]
  \;,\quad
  \tikzpic{\lili}[scale=.8]
  \;,\quad
  \tikzpic{\cc}[scale=.8]
  \;,\quad
  \tikzpic{\li\bli[1][0]}[scale=.8]
  \quad\an\quad
  \tikzpic{\bli\li[1][0]}[scale=.8]
  \;.
\end{gather*}
Note that the diagrams
$\tikzpic{
  \cro[0][2] \li[2][2]
  \li[0][1] \cro[1][1] 
  \cro \li[2][0] 
}[scale=.4]$
and
$\tikzpic{
  \cro[1][2] \li[0][2]
  \li[2][1] \cro[0][1] 
  \cro[1][0] \li[0][0] 
}[scale=.4]$
have the same number of crossings; we call the latter the ``$\rhd$-leading term'' of the (target of) the relation.
In the case of the presentation $\affS(\affBr)$ in \cref{fig:AB_S_rewriting_steps}, we have
$\tikzpic{\LOTRtwo}[scale=.6]=\tikzpic{\lili\lili[0][1]}[scale=.6]$,
$\tikzpic{\LOTRthree}[scale=.5]=0$
and $\VARRthree=1$.
The general form of the presentation we consider is given in \cref{fig:generic_grobner_of_affine_Brauer_type} in the main text.
A presentation of this form is called a \defnemph{rewriting system of affine Brauer type}.

\begin{figure}[t]
	\centering
  \renewcommand{\sca}{.4}

	\def\out{1ex}
	\def\tempspc{60mu}

	\begin{subfigure}{\textwidth}
		\begin{gather*}
				\tikzpic{
          \cro[0][1]
				  \cro 
				} 
				\;\to\;
				\tikzpic{
          \lili\lili[0][1] 
        } 
				\mspace{\tempspc}
        %
				\tikzpic{
          \cro[0][2] \li[2][2]
          \li[0][1] \cro[1][1]
          \cro \li[2][0]
        }[scale=.8]
				\;\to\;
				\tikzpic{
          \cro[1][2] \li[0][2]
          \li[2][1] \cro[0][1]
          \cro[1][0] \li[0][0]
        }[scale=.8]
				\mspace{\tempspc}
        %
				\tikzpic{
          \cu[1][0] \uli[2][0] 
				  \dli \ca 
				}
				\;\to\;
				\tikzpic{\li}
				\mspace{\tempspc} 
        %
				\tikzpic{
          \dli[2][0] \ca[1][0] 
          \uli \cu 
				}
				\;\to\;
				\tikzpic{\li}
        \\[\out]
				\tikzpic{
          \uli[0][1] \ca[1][1] 
          \cro \li[2][0] 
        }
				\;\to\;
				\tikzpic{ 
          \ca[0][1] \uli[2][1]
          \li \cro[1][0] 
        }
				\mspace{\tempspc}
        %
				\tikzpic{
          \ca \uli[2][0]
          \li[0][-1] \cro[1][-1]
          \cro[0][-2] \li[2][-2]
        }[scale=.8]
				\;\to\;
				\tikzpic{\lica}
				\mspace{\tempspc}
        %
        \tikzpic{ 
          \cro \li[2][0]
          \dli \cu[1][0] 
        }
        \;\to\;
        \tikzpic{
          \li[0][0] \cro[1][0]
          \cu[0][0] \dli[2][0]
        }
				\mspace{\tempspc}
        %
        \tikzpic{
          \cro[0][1] \li[2][1]
          \li[0][0] \cro[1][0]
          \cu \dli[2][0]
        }[scale=.8]
        \;\to\;
        \tikzpic{
          \li[0][0] \cu[1][1]
        }
			\\[\out]
				\tikzpic{
          \ca[0][1]
				  \cro 
				}
				\;\to\;
				\tikzpic{\ca}
				\mspace{\tempspc}
        %
        \tikzpic{
          \cro
          \cu
        }
        \;\to\;
        \tikzpic{
          \cu
        }
        \mspace{\tempspc}
        %
        \tikzpic{
          \ca[0][1] \uli[2][1]
          \li[0][0] \cro[1][0]
          \cu \dli[2][0]
        }[scale=.8]
        \;\to\;
        \tikzpic{
          \li
        }
  	\end{gather*}
	
		\caption{Non-affine rewriting steps.}
		\label{subfig:AB_nonaffine_rewriting_steps}
	\end{subfigure}


	\begin{subfigure}{.6\textwidth}
		\begin{IEEEeqnarray*}{rCl}
      \tikzpic{\bubble[0][0][2k+1]}[][1]
      \;&\to&\;
      -\frac{1}{2}
      \tikzpic{\bubble[0][0][2k]}[][1]
      +
      \frac{1}{2}
      \sum_{i+j=2k}
      (-1)^i
      \tikzpic{
        \bubble[0][0][2k-i]
        \bubble[1.8][0][i]
      }[][1]
      \\
      \tikzpic{\cu\ca}[][0]
      \;&\to&\;
      \VARev
    \end{IEEEeqnarray*}
	
		\caption{Bubble-evaluation rewriting steps.}
		\label{subfig:AB_bubble_evaluation_rewriting_steps}
	\end{subfigure}%
  \begin{subfigure}{.4\textwidth}
		\begin{gather*}
      \begin{IEEEeqnarraybox}{rCl}
        \tikzpic{\bubble[0][0][k]\li[-1][.5]}[][1]
        &\;\to\;&
        \tikzpic{\bubble[0][0][k]\li[1.5][.5]}[][1]
        \;+\;
        \tikzpic{\LOTbbsl}[][.5]
        \\[1ex]
        \tikzpic{\bubble[0][0][k]\li[1.5][.5]}[][1]
        &\to&
        \tikzpic{\bubble[0][0][k]\li[-1][.5]}[][1]
        \;-\;
        \tikzpic{\LOTbbsl}[][.5]
      \end{IEEEeqnarraybox}
		\end{gather*}
		\caption{Bubble-slide rewriting steps.}
		\label{subfig:AB_bubbleslide_rewriting_steps}
	\end{subfigure}

  
  \begin{subfigure}{\textwidth}
    \renewcommand{\tempspc}{40mu}
		\begin{gather*}
			\begin{IEEEeqnarraybox}{rClcrClcrClcrCl}
				\tikzpic{
          \ulcr
          \Diag[0][0][]
        }
				\;&\to&\; 
        \tikzpic{
          \drcr
          \Diag[0][0][] 
        }
				- \tikzpic{\lili} 
        + \tikzpic{\cc} 
        &
				\mspace{\tempspc}
				&
        \tikzpic{
          \drcr
          \Diag[0][0][]
        }
				\;&\to&\; 
        \tikzpic{
          \ulcr
          \Diag[0][0][]
        }
				+ \tikzpic{\lili} 
        - \tikzpic{\cc} 
        &
				\mspace{\tempspc}
				&
				\tikzpic{\rca \ca[0][0][]}
				\;&\to&\; 
				{}- \tikzpic{\lca \ca[0][0][]}
        &
				\mspace{\tempspc}
        &
        \tikzpic{\lca \ca[0][0][]}
				\;&\to&\; 
				{}- \tikzpic{\rca \ca[0][0][]}
        \\[2ex]
				\tikzpic{
          \urcr
          \AntiDiag[0][0][]
        }
				\;&\to&\; 
        \tikzpic{
          \dlcr
          \AntiDiag[0][0][]
        }
				+ \tikzpic{\lili} 
        - \tikzpic{\cc} 
        &
				\mspace{\tempspc}
				&
        \tikzpic{
          \dlcr
          \AntiDiag[0][0][]
        }
				\;&\to&\; 
        \tikzpic{
          \urcr
          \AntiDiag[0][0][]
        }
				- \tikzpic{\lili} 
        + \tikzpic{\cc} 
        &
				\mspace{\tempspc}
				&
				\tikzpic{\rcu \cu[0][0][]}
				\;&\to&\; 
				{}- \tikzpic{\lcu \cu[0][0][]}
				&
        \mspace{\tempspc}
				&
        \tikzpic{\lcu \cu[0][0][]}
				\;&\to&\; 
				{}- \tikzpic{\rcu \cu[0][0][]}
			\end{IEEEeqnarraybox}
		\end{gather*}
	
		\caption{Dot-slide rewriting steps.}
		\label{subfig:AB_dotslide_rewriting_steps}
	\end{subfigure}
	
	\caption{The linear monoidal rewriting system $\affS(\affBr)$ presenting $\affBr$. Here $k\in\bN$ and $\VARev$ is a parameter in the ground ring.
  We leave $\tikzpic{\LOTbbsl}[][.5]$ undefined in this introduction: see \cref{subsec:affine_brauer_category}.}
	\label{fig:AB_S_rewriting_steps}
\end{figure}

\subsubsection{Critical branchings and traffic rules}

\Cref{fig:critical_branchings_vertsym} gives a list of diagrams called \defnemph{critical branchings}: they encode conditions for a rewriting system of affine Brauer type to present a category of affine Brauer type.
To explain this, note that the arrows in \cref{fig:AB_S_rewriting_steps} give each relation a preferred direction; one calls them \defnemph{rewriting steps}.
A rewriting step has a source and target.
Any diagram in \cref{fig:critical_branchings_vertsym} encodes an overlap between the sources of two (possibly the same) rewriting steps.
For instance, the leftmost diagram below
\begin{gather}
  \label{eq:intro_example_confluence}
  \newcommand{\arOrNot}[1]{#1}
  \renewcommand{\sca}{.3}
  \begin{tikzcd}[ampersand replacement=\&,row sep=0,column sep=small]
    \&
    \tikzpic{
      \drcr\li[2][0][\arOrNot{->}]
      \dli\cu[1][0]
    }
    \;-\;
    \tikzpic{
      \lili\li[2][0]
      \dli\cu[1][0]
    }
    \;+\;
    \tikzpic{
      \cc\li[2][0]
      \dli\cu[1][0]
    }
    \&
    -\;
    \tikzpic{
      \cro\li[2][0][\arOrNot{->}]
      \dli\rcu[1][0]
    }
    \;-\;
    \tikzpic{
      \lili\li[2][0]
      \dli\cu[1][0]
    }
    \;+\;
    \tikzpic{
      \lili\li[2][0]
      \cu\dli[2][0]
    }
    \&
    \\
    \tikzpic{
      \ulcr\li[2][0][\arOrNot{->}]
      \dli\cu[1][0]
    }[scale=1]
    \&\&\&
    -\;
    \tikzpic{
      \li\urcr[1][0]\AntiDiag[1][0][\arOrNot{->}]
      \cu\dli[2][0]
    }
    \;-\;
    \tikzpic{
      \lili\li[2][0]
      \dli\cu[1][0]
    }
    \;+\;
    \tikzpic{
      \lili\li[2][0]
      \cu\dli[2][0]
    }
    \\
    \&
    \tikzpic{
      \bli\cro[1][0]\AntiDiag[1][0][\arOrNot{->}]
      \cu\dli[2][0]
    }
    \&
    -\;
    \tikzpic{
      \li\cro[1][0]\AntiDiag[1][0][\arOrNot{->}]
      \rcu\dli[2][0]
    }
    \&
    %
    \arrow[from=2-1,to=1-2, bend left=10pt]
    \arrow[from=1-2,to=1-3]
    \arrow[from=1-3,to=2-4, bend left=10pt]
    \arrow[from=2-1,to=3-2, bend right=10pt]
    \arrow[from=3-2,to=3-3]
    \arrow[from=3-3,to=2-4, bend right=10pt]
  \end{tikzcd}
\end{gather}
encodes an overlap between the rewriting step $\tikzpic{\ulcr}[scale=.7]\to\ldots$ and the rewriting step $\tikzpic{\cro\li[2][0]\dli\cu[1][0]}[scale=.5]\to\ldots$.
The two rewriting steps give two different ways (``branches'') of rewriting the diagram.
As above, we can continue to rewrite both branches using rewriting steps.
Eventually, each branch will reach a linear combination of ``normal forms'', that is, diagrams that cannot be rewritten by any rewriting step.
If they reach the \emph{same} linear combination (as in the example above), the condition is satisfied.
In that case, we say that the branching is \defnemph{confluent}\footnote{More precisely, one should say \emph{positively confluent}; see \cref{defn:RW_S_rewriting_step}. In this introduction, confluence always means positive confluence.}.

Some critical branchings carry extra decoration; we call these decorations \defnemph{traffic rules}:
\begin{IEEEeqnarray*}{CcCcCcC}
  \tikzpic{\draw[->] (0,0) to (1,1);}
  &\qquad&
  \tikzpic{\lorli}
  \;\text{ or }\;
  \tikzpic{\rorli}
  &\qquad&
  \tikzpic{\custop[0][0]}
  &\qquad&
  \tikzpic{\custop[0][0]\node[below right=-4pt] at (.5,-.4) {\footnotesize$\mathtt{1}$};}
  \\
  \text{\footnotesize orientation}
  &&
  \text{\footnotesize tunnels}
  &&
  \text{\footnotesize stop sign}
  &&
  \text{\footnotesize one-free-pass stop sign}
\end{IEEEeqnarray*}
For dot-slide critical branchings (\cref{subfig:critical_branchings_dot_slides_vertsym}), the strand on which the dot slides carries an orientation.
It encodes the following condition on the rewriting branches: if we apply a rewriting step on a $\rhd$-leading term, one is only allowed to slide the dot in the direction of the orientation.
For instance, in the context of the affine Brauer category:
\begin{gather*}
  \tikzpic{\rca \ca[0][0][<-]}
  \;\longrightarrow\; 
  {}- \tikzpic{\lca[0][0][<-]}
  \quad\text{but}\quad
  \tikzpic{\rca \ca[0][0][->]}
  \;\centernot\longrightarrow\; 
  {}- \tikzpic{\lca[0][0][->]}
  \;.
\end{gather*}
Bubble-slide critical branchings (\cref{subfig:critical_branchings_bubbles_vertsym}) are decorated with red arrows, called \defnemph{tunnels}; an analogous condition for sliding a bubble through a strand applies.
The first bubble-slide critical branching is an overlap between a bubble slide and a bubble evaluation, while the other two bubble-slide critical branchings are overlaps between two bubble slides.

The cup appearing in the first two exceptional critical branchings has a decoration $\tikzpic{\custop}$, called a \defnemph{stop sign}.
The terminology is self-explanatory: a dot can never cross a stop sign.
For instance, the first exceptional critical branching is an overlap between ``undoing the kink $\tikzpic{\cu\cro}[scale=.5]\to\ldots$'' and ``sliding the dot down the crossing $\tikzpic{\urcr}[scale=.5]\to\ldots$''.
Because of the stop sign, one cannot continue the second branch by sliding the dot through the cup, and must instead slide all the remaining $k$ dots down the crossing.
When $k+1$ is odd, one needs to use the bubble-evaluation rewriting step to reach confluence.
In other words, the condition encoded by the first exceptional critical branching leads one to discover bubble
evaluation and the relation it induces on bubbles in the underlying category.
Similarly, resolving the condition encoded by the second exceptional critical branching leads to the bubble-slide rewriting steps.

The last exceptional critical branching has a decoration $\tikzpic{\custop\node[below right=-4pt] at (.5,-.4) {\footnotesize$\mathtt{1}$};}[][-.3]$, called the \defnemph{one-free-pass stop sign}.
Only one dot can cross the one-free-pass stop sign: after one dot has crossed, no other dot can cross.
To resolve the condition encoded by this branching, slide one dot through the \emph{cup} from left to right, slide it back through the \emph{cap} from right to left, and if possible, simplify both the source and target using bubble evaluation:
\begin{gather}
  \label{eq:one-free-stop-sign}
  \begin{tikzcd}[ampersand replacement=\&,row sep=tiny]
    \tikzpic{\lcustop[0][0][->]\ca\node[below right=-4pt] at (.5,-.4) {\footnotesize$\mathtt{1}$};\bul\node[left=-1pt] at (0,0) {\scriptsize $k$};}[scale=1][0]
    \rar
    \arrow[dr,"\text{evaluation}"{description},dotted,bend right=18]
    \&
    (\text{scalar})
    \tikzpic{\rcustop[0][0][->]\ca\bul\node[left=-1pt] at (0,0) {\scriptsize $k$};}[scale=1][0]
    \;+\;
    (\text{other terms})
    \rar
    \&
    (\text{scalar})
    \tikzpic{\custop[0][0][->]\lca\bul\node[left=-1pt] at (0,0) {\scriptsize $k$};}[scale=1][0]
    \;+\;
    (\text{other terms})
    \arrow[dl,"\text{evaluation}"{description},dotted,bend left=15]
    \\
    \&\quad \ldots \quad\overset{?}{=}\quad \ldots \quad\&
  \end{tikzcd}
\end{gather}
Note that if $k+1\notin\bubbleset$, one cannot evaluate the source, so that the condition simply reads as ``in the rightmost term, $\text{(scalar)}=1$ and $\text{(other terms)}$ rewrites to zero''.
In this case, the diagram does not encode a branching, formally speaking, since one of the branches is just the identity.
We abuse terminology by calling it a critical branching.

\subsubsection{Statement of the Main Theorem}

We can now state the Main Theorem.
For the purpose of this introduction, we state it for categories that admit an isomorphism $\affC\cong\affC^{\op}$ ``readily visible'' in the presentation: we say that the rewriting system of affine Brauer type $\affS$ presenting $\affC$ is \emph{vertically symmetric}; see \cref{prop:vertical_symmetries}.
The version in the main text also applies when the interchange law contains scalars; here we restrict the statement to linear monoidal categories.

\begin{bigtheorem*}[simplified; see \cref{thm:affine_Grobner_system_iff_critical_branchings} and \cref{prop:vertical_symmetries}]
  Let $\affC$ be a linear monoidal category presented by a vertically symmetric rewriting system of affine Brauer type $\affS$.
  The linear monoidal category $\affC$ is a category of affine Brauer type if and only if every critical branching given in \cref{fig:critical_branchings_vertsym} is confluent.
\end{bigtheorem*}

\begin{figure}[t]
  \renewcommand{\sca}{.35}
  \newcommand{\vertsp}{.3ex}
  \newcommand{\arOrNot}[1]{#1}
  \def\tempsp{20mu}

  \begin{subfigure}{\textwidth}
    \renewcommand{\vertsp}{0ex}
    \begin{gather*}
      \tikzpic{
        \cro[0][2]
        \cro[0][1]
        \cro
      }
      \mspace{\tempsp}
      \tikzpic{
        \cro[0][3]\li[2][3]
        \cro[0][2]\li[2][2]
        \li[0][1]\cro[1][1]
        \cro[0][0]\li[2][0]
      }
      \mspace{\tempsp}
      \tikzpic{
        \cro[0][4]\li[2][4]
        \li[0][3]\cro[1][3]
        \cro[0][2]\cro[2][2]
        \li[0][1]\cro[1][1]
        \cro\li[2][0]
        \draw[] (3,0) to (3,2);\draw[] (3,3) to (3,5);
      }
      \mspace{\tempsp}
      \tikzpic{
        \ca[2][1]
        \ca\lili[2][0]
        \dli\cu[1][0]\dli[3][0]
      }
      \mspace{\tempsp}
      \tikzpic{
    		\uli[0][2]\ca[1][2]
        \cro[0][1]\li[2][1]
        \cro\li[2][0]
    	}
      \mspace{\tempsp}
      \tikzpic{
    		\uli[0][4]\ca[1][4]
        \cro[0][3]\li[2][3]
        \lili[0][2]\li[2][2]
        \li[0][1]\cro[1][1]
        \cro\li[2][0]
    	}
      \mspace{\tempsp}
      \tikzpic{
        \cro[0][2]\ca[2][2]
        \li[0][1]\cro[1][1]\li[3][1]
        \cro\lili[2][0]
      }
      \mspace{\tempsp}
      \tikzpic{
    		\uli[0][1]\ca[1][1]\uli[3][1]
        \cro\lili[2][0]
        \dli[0][0]\dli[1][0]\cu[2][0]
    	}
      \mspace{\tempsp}
      \tikzpic{ 
        \ca[2][3]
        \ca[0][2]\lili[2][2]
        \li[0][1]\cro[1][1]\li[3][1]
        \cro\lili[2][0]
      }
      \\[\vertsp]
      \tikzpic{ 
        \uli[0][1]\ca[1][1]
        \cro\li[2][0]
        \dli\cu[1][0]
      }
      \mspace{\tempsp}
      \tikzpic{ 
        \uli[0][1]\ca[1][1]
        \cro\cro[2][0]
        \dli\cu[1][0]
        \draw[] (3,-.4) to (3,0);\draw[] (3,1) to (3,1.4);
      }
      \mspace{\tempsp}
      \tikzpic{
        \ca[0][3]\uli[2][3]
        \li[0][2]\cro[1][2]
        \lili[0][1]\li[2][1]
        \cro\li[2][0]
        \dli\cu[1][0]
      }
      \mspace{\tempsp}
      \tikzpic{
        \ca[0][3]\uli[2][3]
        \li[0][2]\cro[1][2]
        \lili[0][1]\cro[2][1]
        \cro\li[2][0]
        \dli\cu[1][0]
        \draw[] (3,-.4) to (3,1);\draw[] (3,2) to (3,3.4);
      }
      \mspace{\tempsp}
      \mspace{\tempsp}
      \tikzpic{
        \ca[0][2]
    		\cro[0][1]
    		\cro[0][0]
    	}
      \mspace{\tempsp}
    	\tikzpic{
        \ca[0][3]\uli[2][3]
    		\cro[0][2]\li[2][2]
    		\cro[1][1]\li[0][1]
    		\cro[0][0]\li[2][0]
    	}
      \mspace{\tempsp}
    	\tikzpic{
    		\ca[1][1]\uli[0][1]
    		\cro[0][0]\li[2][0]
    		\cu\dli[2][0]
    	}
      \mspace{\tempsp}
    	\tikzpic{
    		\ca[0][1]
    		\cro[0][0]
    		\cu\
    	}
      \mspace{\tempsp}
    \end{gather*}
    \caption{Non-affine critical branchings.}
    \label{subfig:critical_branchings_non_affine_vertsym}
  \end{subfigure}

  \begin{subfigure}{\textwidth}
    \begin{gather*}
      \tikzpic{
        \drcr \cro[0][1]
        \AntiDiag[0][1][\arOrNot{->}]
      }
      \mspace{\tempsp}
      \tikzpic{
        \dlcr \cro[0][1]
        \Diag[0][1][\arOrNot{->}]
      }
      \mspace{\tempsp}\mspace{\tempsp}
      \tikzpic{
        \dlcr \li[2][0] \li[0][1] \cro[1][1] \cro[0][2] \li[2][2][\arOrNot{->}]
      }
      \mspace{\tempsp}
      \tikzpic{
        \drcr \li[2][0] \li[0][1] \cro[1][1] \cro[0][2] \li[2][2]
        \AntiDiag[0][2][\arOrNot{->}]
      }
      \mspace{\tempsp}
      \tikzpic{
        \cro \li[2][0] \li[0][1] \drcr[1][1] \cro[0][2] \li[2][2]
        \Diag[0][2][\arOrNot{->}]
      }
      \mspace{\tempsp} 
      \mspace{\tempsp} 
      \tikzpic{
        \dli \lca[0][0]\ca[0][0] \cu[1][0] \uli[2][0][\arOrNot{->}]
      }
      \mspace{\tempsp}
      \mspace{\tempsp}
      \tikzpic{
        \ulcr \dlcr
        \Diag[0][0][\arOrNot{<-}]
        \AntiDiag[0][0][\arOrNot{->}]
      }
      \mspace{\tempsp}
      \mspace{\tempsp}
      \tikzpic{
        \dlcr \ca[0][1]
        \Diag[0][0][\arOrNot{<-}]
      }
      \mspace{\tempsp}
      \mspace{\tempsp}
      \tikzpic{
        \dlcr \uli[0][1] \ca[1][1] \li[2][0][\arOrNot{<-}]
        \AntiDiag[0][0]
      } 
      \mspace{\tempsp}
      \tikzpic{
        \drcr \uli[0][1][\arOrNot{->}] \ca[1][1] \li[2][0]
        \Diag[0][0]
      } 
    \end{gather*}
    \caption{Dot-slide critical branchings.}
    \label{subfig:critical_branchings_dot_slides_vertsym}
  \end{subfigure}

  \begin{subfigure}{.5\textwidth}
    \begin{gather*}
      \tikzpic{
        \lorli[-1.4][0] \bul[0][.5]\ca[0][.5]\cu[0][.5]
        \node[left=-1pt] at (0,.5) {\scriptsize $k$};
        \node at (-.2,-1) {\scriptsize for ${k}\notin\bubbleset$};
      }[][.5]
      %
      \qquad
      \tikzpic{
        \lorli[-1.4][0] \bul[0][.5]\ca[0][.5]\cu[0][.5]
        \node[left=-1pt] at (0,.5) {\scriptsize $k$};
        \rorli[1.8][0]
        \draw (-1.4,1) to[out=90,in=90] (1.8,1);
      }
      \qquad
      \tikzpic{
        \cro[-2][0]
        \bul[0][.5]\ca[0][.5]\cu[0][.5]
        \node[left=-1pt] at (0,.5) {\scriptsize $k$};
        \draw[very thick,red,->] (-2+.9,.6) to (-2+.5,1);
        \draw[very thick,red,->] (-2+.9,.4) to (-2+.5,0);
      }
    \end{gather*}
    \caption{Bubble-slide critical branchings.}
    \label{subfig:critical_branchings_bubbles_vertsym}
  \end{subfigure}%
  \begin{subfigure}{.5\textwidth}
    \begin{gather*}
      \tikzpic{
        \ca[0][1]
        \urcr
        \custop
        \AntiDiag[0][0][<-]
        \bul[1][1]\node[right=-1pt] at (1,1) {\scriptsize $k$};
      }
      \qquad
      \tikzpic{
        \ca[0][2]\uli[2][2]
        \li[0][1]\cro[1][1]
        \ulcr \Diag[0][0][<-] \li[2][0]
        \dli \custop[1][0]
        \bul[0][1]\node[left=-1pt] at (0,1) {\scriptsize $k$};
      }
      \qquad
      \tikzpic{\lcustop[0][0][->]\ca\node[below right=-4pt] at (.5,-.4) {\footnotesize$\mathtt{1}$};\bul\node[left=-1pt] at (0,0) {\scriptsize $k$};}[scale=2]
    \end{gather*}
    \caption{Exceptional critical branchings.}
    \label{subfig:critical_branchings_exceptional_vertsym}
  \end{subfigure}

  \caption{Critical branchings of affine Brauer type, assuming vertical symmetry.}
  \label{fig:critical_branchings_vertsym}
\end{figure}

For most critical branchings, checking confluence is a straightforward computation; one could implement it on a computer.
In fact, we did: see \cref{subsec:implementation}.
Only the bubble-slide and exceptional critical branchings are not entirely straightforward, as they encode infinite families of critical branchings, indexed by a natural number $k$.
To compute them efficiently, we use generating functions, following an idea of Savage and Webster \cite{SW_BubblesAffineBrauer_2024}; see \cref{subsec:affine_brauer_category} for the case of the affine Brauer category.
The theorem in the main text also replaces the notion of confluence by a novel and weaker notion, called ``$\rhd$-tamed preconfluence''; see \cref{subsec:tamed_preconfluence}.
It allows greater flexibility in the use of rewriting steps.

\medbreak

As a first application, we recover the basis theorems of \cite{RS_AffineBrauerCategory_2019,BWW_NilBrauerCategory_2024,BDE+_AffineVWSupercategory_2020} in an intrinsic and combinatorial way; see \cref{sec:examples}.
To see the Main Theorem in action, the reader may wish to skip ahead to \cref{subsec:affine_brauer_category} (looking up  definitions in \cref{sec:definitions} from time to time) where the case of the affine Brauer category $\affBr$ is discussed.

\begin{proposition}
  The affine Brauer category $\affBr[\VARev]$, the nil-Brauer category $\nilBr[\VARev]$ and the affine VW supercategory $\sVW$ are categories of affine Brauer type with $\bubbleset$-bubbles, where $\bubbleset=2\bN_{>0}$, $\bubbleset=2\bN+1$ and $\bubbleset=\emptyset$ respectively.
\end{proposition}

In hindsight (and only in hindsight), one can try to interpret the conditions given by the critical branchings.
Many of them can be viewed as a naturality-like condition with respect to a braided-like or pivotal-like categorical structure.
For instance, the branching
\begin{gather*}
  \tikzpic{
    \ca[0][3]\uli[2][3]
    \cro[0][2]\li[2][2]
    \cro[1][1]\li[0][1]
    \cro[0][0]\li[2][0]
  }[scale=.7]
\end{gather*}
encodes the difference between undoing the ``kink'' and then sliding it through the strand, or sliding first everything through the strand and then undoing the kink.
As noted above, the first two exceptional critical branchings encode bubble evaluation and bubble slide, respectively.

\medbreak

We now explain the theory behind the proof of the Main Theorem.
The reader only interested in applications can skip ahead to \cref{subsec:intro_applications}.

\subsection{Monoidal Gröbner theory}
The proof of the Main Theorem uses higher linear rewriting theory \cite{Schelstraete_RewritingModuloDiagrammatic_2026}.
In this article, we specialize and rephrase the theory in the setting of linear monoidal categories.
We also further develop the rewriting toolbox; see \eg our use of symmetries as discussed in \cref{subsubsec:intro_local_analysis}.

\subsubsection{Monoidal Gröbner systems}

A presentation with preferred directions as the one given in \cref{fig:AB_S_rewriting_steps} is called a \defnemph{linear monoidal rewriting system} $\sS$ (see \cref{subsec:RW_monoidal_linear_rewriting_system}).
Directed relations are called \defnemph{rewriting steps}.
More precisely, the set $\sS$ consists of \emph{generating rewriting steps}, and a rewriting step is of the form $\Gamma[\ssf]$ for $\ssf\in\sS$ a generating rewriting step and $\Gamma[-]$ a diagram ``surrounding'' the rewriting step.
The data $\Gamma[-]$ is called a \defnemph{context}. We write $\Cont(\sS)$ the set of rewriting steps\footnote{These are really the \emph{monomial} rewriting steps. We ignore the distinction in the introduction.} in $\sS$; see \cref{subsec:RW_contexts_subsystems}.

The idea underpinning rewriting theory is as follows.
Simplifying a diagram by successively applying rewriting steps defines a ``proto'' algorithm. We say ``proto'' since at this stage, it is not clear whether the algorithm \emph{terminates} or whether it is \emph{confluent}, that is, whether applying rewriting steps in a different order would output the same answer (one says that the algorithm is \emph{non-deterministic}).
However, if the algorithm {does} terminate and is confluent, then every diagram admits a unique expression as a linear combination of \defnemph{normal diagrams}, \ie diagrams that cannot be rewritten by rewriting steps.
In other words, normal diagrams define a hom-basis.
Better still, the rewriting system describes, at least in principle, an algorithm to express any element in the hom-basis, and hence an algorithmic solution to distinguishing elements in the underlying category.
One says that the rewriting system encodes a solution to the \emph{word\footnote{In the context of diagrammatic categories, it should really be called the ``diagram problem''.} problem}.

In the context of diagrammatic categories, termination is often a problem.
For instance, if we did not impose the restriction given by stop signs, a dot could cycle forever along a closed strand.
Whether a strand is closed is a global condition, and hence so are stop signs.
More formally, traffic rules amount to a choice of subset $\lT\subset\Cont(\sS)$ of rewriting steps, called a \defnemph{linear sub-system}.
Importantly, a linear sub-system need not be compatible with contexts: if $\ssf\in\sS$ and $\Gamma[-]$ is a context, then whether $\Gamma[\ssf]$ belongs to $\lT$ typically depends on $\Gamma[-]$.

By restricting to a linear sub-system, we can reach termination.
We encode termination using a well-founded partial order $\succ$, \defnemph{compatible} with $\lT$ in the sense that each rewriting step strictly reduces $\succ$; see \cref{subsec:RW_termination_order}.\footnote{In particular, the definition of ``compatible'' is slightly stronger than what the discussion suggests.}
Our discussion leads to the following notion:

\begin{definition}[simplified; see \cref{defn:RW_monoidal_grobner_system}]
  A triple $(\sS,\lT,\succ)$ is a \defnemph{monoidal Gröbner system} if $\sS$ and $\lT$ present the same underlying category and if each diagram admits a unique normal form with respect to $\lT$.
\end{definition}

Given our discussion above, we have the following result:

\begin{theorem}[simplified; see \cref{thm:RW_monoidal_grobner_system}]
  \label{thm:intro_RW_monoidal_grobner_system}
  Let $\cC$ be a linear monoidal category presented by a linear monoidal rewriting system $\sS$.
  Let $\lT\subset\sS$ be a linear sub-system and $\succ$ a well-founded partial order compatible with $\lT$.
  Normal diagrams with respect to $\lT$ define a hom-basis for $\cC$ if and only if $(\sS,\lT,\succ)$ defines a monoidal Gröbner system.
\end{theorem}

This rephrases the motivating problem at the start of the introduction: to show that a certain hom-set $\cB$ is a hom-basis of a linear monoidal category $\cC$, find a monoidal Gröbner system $(\sS,\lT,\succ)$ presenting $\cC$ and whose normal diagrams with respect to $\lT$ coincide with $\cB$.

\subsubsection{Normal diagrams}

We now return to the setting of a rewriting system of affine Brauer type $\affS$.
We let $\affT\subset\Cont(\affS)$ be the linear sub-system consisting of rewriting steps that respect traffic rules.
A well-founded compatible order $\succ$ is defined in \cref{defn:orders_gen_brauer}; we refer to the main text for the definition.
We wish to know whether the data $(\affS,\affT,\succ)$ defines a monoidal Gröbner system.
Indeed, if that is the case, then the category presented by $\affS$ is a category of affine Brauer type, thanks to \cref{thm:intro_RW_monoidal_grobner_system} and the following result:

\begin{theorem}[simplified; see \cref{thm:normal_form_are_standard_diagrams}]
  \label{thm:intro_normal_form_are_standard_diagrams}
  Let $\affS$ be a rewriting system of affine Brauer type.
  Normal diagrams with respect to $\affT$ give a choice of representatives for reduced matchings.
\end{theorem}

\begin{wrapfigure}[6]{r}{.13\textwidth}
  \vspace*{-10pt}
  \begin{gather*}
    \tikzpic{
      \li[-1][1]\ca[0][1]\li[2][1]
      \lili[-1][0]\cro[1][0]
      \li[-1][-1]\lcustop[0][0][->]\li[2][-1]
      \begin{scope}[xscale=3]
        \cu[-1/3][-1]\ca[-1/3][2]
      \end{scope}
    }
  \end{gather*}
\end{wrapfigure}
The difficulty of this theorem lies in definitions: how do we fix the traffic rules such that the theorem holds?
The most subtle definition is the stop sign. As we noted above, we need a stop sign on each closed strand, otherwise a dot could cycle around forever, and the reduction algorithm would not terminate. But where exactly do we put the stop sign?
If we are not careful, dots may get stuck in the ``middle'' of a diagram, preventing any further simplification; see the diagram on the right.
Solving that problem involves properly understanding equivalence classes under the interchange law. For that, we rely on the work of Delpeuch and Vicary on normal forms for the interchange law \cite{DV_NormalizationPlanarString_2022}.
The proof of \cref{thm:intro_normal_form_are_standard_diagrams} (or rather \cref{thm:normal_form_are_standard_diagrams}) is given in \cref{sec:proof_normal_forms}.

\medbreak

Hence, to prove the Main Theorem, it remains to show that $(\affS,\affT,\succ)$ is a monoidal Gröbner system.
Whether $\affS$ and $\affT$ present the same category is discussed in the next subsection.
The hardest part is to show that $\affT$ is confluent. We discuss it in the remaining subsections.

\subsubsection{$\affS$ and $\affT$ present the same category}

Does $\affT(\affBr)$ present $\affBr$?
The main uncertainty lies in the stop sign, which prevents the dot from sliding through the cup in either direction.
Thus, a priori, this relation is not present in the category presented by $\affT(\affBr)$: can we recover it?
If the cup has a stop sign, then the strand must be closed.
Hence, to move the dot to the other side of the cup, we can alternatively slide it along the whole closed strand.
One then hopes that this induces the same relation as simply sliding the dot through the cup.
An example of that situation is given by the third exceptional critical branching: it essentially says that sliding a dot along the bubble should be equivalent to doing nothing.
As part of the proof of the Main Theorem, we show that this case implies all the other cases.
In other words, the third exceptional branching is a necessary and sufficient condition to ensure that $\affT(\affBr)$ and $\affS(\affBr)$ present the same underlying category.
See \cref{subsubsec:confluence_equivalence_presentation} for details.

\subsubsection{Enumerating overlaps}
\label{subsubsec:intro_enumerating_overlaps}

From now on, we focus on the hardest problem of showing confluence.
We do not discuss every part of the proof and omit many technical points; see \cref{sec:critical_branchings} for details.

In principle, one needs to show that \emph{every} diagram has a unique normal form with respect to $\lT$. In a nutshell, the Main Theorem says that it suffices to check that fact for the diagrams in \cref{fig:critical_branchings_vertsym}.
The first step consists in enumerating all the possible ways two rewriting steps can overlap.
We rely on a general statement: when the sources are connected diagrams, two rewriting steps ``interact non-trivially'' only if they actually overlap, that is, if they share at least one generator; see \cref{prop:classification_overlapping_branchings}.
As far as we are aware, this statement has not appeared before in the literature.
Then, the enumeration is a cumbersome but straightforward process: we encourage the reader to have a peek at the full list given in \cref{fig:all_overlapping_branchings}.

Note that some overlaps have an unknown diagram depicted as a box with a question mark $\tikzpic{\indexNF[0][0][?]}[scale=.7]$. These overlaps are called \defnemph{indexed branchings}.
The existence of indexed branchings is a well-known phenomenon in higher rewriting \cite{Lafont_AlgebraicTheoryBoolean_2003}: as we rewrite inside equivalence classes under the interchange law, a diagram can get ``stuck'' between the sources of two rewriting steps, while not actually overlapping with them.


\subsubsection{Local analysis}
\label{subsubsec:intro_local_analysis}

In order to reduce the number of branchings that one needs to consider, we ask: does the confluence of a given overlap imply the confluence of some other overlap?
Our main tools in this analysis are \defnemph{symmetries}, vaguely understood.
To explain this, note that the presentation of $\affBr$ induced by the rewriting system $\affS(\affBr)$ is far from being minimal.
For instance, the ``cap-pulling rewriting step'' is derived from other rewriting steps:
\begin{gather}
  \label{eq:intro_defn_cap_pulling}
  \def\scl{.7}
  \begin{tikzcd}[ampersand replacement=\&]
  	\tikzpic{ 
      \ca[0][2]\uli[2][2]
      \li[0][1]\cro[1][1]
      \cro[0][0]\li[2][0]
    }[scale=\scl]
    \&
  	\tikzpic{ 
      \uli[0][2]\ca[1][2]
      \cro[0][1]\li[2][1]
      \cro[0][0]\li[2][0]
    }[scale=\scl]
    \&
  	\tikzpic{ 
      \uli[0][2]\ca[1][2]
      \lili[0][1]\li[2][1]
      \lili[0][0]\li[2][0]
    }[scale=\scl]
  	\arrow[curve={height=-30pt}, from=1-1, to=1-3]
  	\arrow[from=1-2, to=1-1]
  	\arrow[from=1-2, to=1-3]
  \end{tikzcd}
\end{gather}
Of course, our computation is not compatible with the direction of the rewriting steps; the sequence above is not a valid sequence in  $\affS(\affBr)$, which explains why the cap-pulling rewriting step had to be manually added.
This necessity is a consequence of the arbitrary choice of direction
$
\tikzpic{
  \uli[0][1] \ca[1][1] 
  \cro \li[2][0] 
}[scale=.5]
\;\to\;
\tikzpic{ 
  \ca[0][1] \uli[2][1]
  \li \cro[1][0] 
}[scale=.5]$\;, which in some sense, ``breaks horizontal symmetry''.
If we decompose the cap-pulling rewriting step as in \eqref{eq:intro_defn_cap_pulling}, we can write certain confluences by ``patching up'' other confluences.
For instance:
\begin{gather*}
  \begingroup
  \savebox{\tempboxa}{\tikzpic{
    \uli[0][2]\ca[1][2]\uli[3][2]
    \cro[0][1]\cu[2][2]
    \cro[0][0]
    \draw[branch_overlay] (-.3,1) rectangle (3.3,2.4);
  }[scale=.4]}
  \def\scl{.6}
  \RenewDocumentCommand{\LOTRtwo}{O{0}O{0}}{%
    \lili[{#1}][{#2}]
    \lili[{#1}][{#2+1}]
  }
  \begin{tikzcd}[ampersand replacement=\&]
  	\tikzpic{ 
      \ca[0][2]\uli[2][2]\uli[3][2]
      \li[0][1]\cro[1][1]\li[3][1]
      \cro[0][0]\cu[2][1]
    }[scale=\scl*1.2]
    \&
  	\tikzpic{ 
      \uli[0][2]\ca[1][2]\uli[3][2]
      \cro[0][1]\cu[2][2]
      \cro[0][0]
    }[scale=\scl]
    \&
  	\tikzpic{ 
      \uli[0][2]\ca[1][2]\uli[3][2]
      \cu[2][2]
      \LOTRtwo[0][0]
    }[scale=\scl]
    \\
  	\tikzpic{ 
      \cro[2][4]
      \ca[0][3]\lili[2][3]
      \li[0][2]\cu[1][3]\li[3][2]
      \Diag[0][1]\AntiDiag[2][1]
      \cro[1][0]
    }[scale=\scl*.8]
  	\&
    \tikzpic{ 
      \cro[0][1]
      \cro[0][0]
    }[scale=\scl]
    \& 
    \tikzpic{ 
      \LOTRtwo[0][0]
    }[scale=\scl]
  	\arrow[curve={height=-40pt}, from=1-1, to=1-3]
  	\arrow[from=1-1, to=2-1]
  	\arrow[from=1-2, to=1-1]
  	\arrow[from=1-2, to=1-3]
  	\arrow[draw=none, from=1-2, to=2-1, "\usebox{\tempboxa}"{description}]
  	\arrow[from=1-2, to=2-2]
  	\arrow["\lrcorner"{description, pos=0}, draw=none, from=1-2, to=2-3]
  	\arrow[from=1-3, to=2-3]
  	\arrow[from=2-1, to=2-2]
  	\arrow[from=2-2, to=2-3]
  \end{tikzcd}
  \endgroup
\end{gather*}
The left square is a confluence for the branching
$\tikzpic{
  \uli[0][2]\ca[1][2]\uli[3][2]
  \cro[0][1]\cu[2][2]
}[scale=.4]$
seen in the context of an additional bottom crossing, which does not overlap with the confluence.
The right square shows a branching between the rewriting steps
$\tikzpic{\cro\cro[0][1]}[scale=.4]\to\ldots$
and
$\tikzpic{\dli\cu[1][0]\ca\uli[2][0]}[scale=.4]\to\ldots$: they do not overlap, so that confluence is automatic\footnote{These branchings are called \emph{independent branchings}, and their confluence is not quite automatic; see \cref{subsubsec:extended_MLRS_independent_branchings}.}.
This gives a confluence for the original branching.
Note that while some arrows are directed in the ``wrong direction'', they are not involved in the final confluence.

Two other symmetries can be leveraged.
The second one comes from traffic rules, where one leverages the symmetry coming from swapping the orientation of strands and tunnels.
The third one is vertical symmetry, and only occurs when the presentation has a vertical symmetry; we assumed it already in this introduction (see the discussion before the statement of the Main Theorem).

\subsubsection{Global analysis}

After considering how various ``local'' branchings are related to one another, we study local-to-global phenomena.
A prime tool in this analysis is the notion of \defnemph{tamed congruence} as defined by the second author in \cite{Schelstraete_RewritingModuloDiagrammatic_2026}; earlier incarnations of this idea can already be found in the work of Buchberger \cite{Buchberger_AlgorithmusAuffindenBasiselemente_1965}, Winkler--Buchberger \cite{WB_CriterionEliminatingUnnecessary_1986} and Bergman \cite{Bergman_DiamondLemmaRing_1978}.
In a nutshell, given a partial order $\succ$ suitably compatible with rewriting steps, tamed congruence allows one to follow a rewriting step ``in the wrong direction'', as long as one remains ``below'' the source $\bullet$ of the branching with respect to the partial order $\succ$.
Here is a schematic:
\begin{IEEEeqnarray*}{CcC}
  \begin{tikzpicture}[scale=.8,pt/.style={inner sep=1pt}]
    \node (C) at (0,0) {$\bullet$};
    \node[pt] (A1) at (2,1) {};
    \node[pt] (A2) at (2.7,.92) {};
    \node[pt] (A3) at (3.45,.63) {};
    \node[pt] (R) at (4,0) {};
    \node[pt] (B1) at (3,-1) {};
    \node[pt] (B2) at (3.35,-.87) {};
    \node[pt] (B3) at (3.68,-.6) {};
    %
    \draw[->] (A1) to[out=-5,in=162] (A2);
    \draw[
      line width=.9pt, line cap=round,
      dash pattern=on 0pt off 4pt, 
    ] (A2) to[out=-18,in=142] (A3);
    \draw[->] (A3) to[out=-38,in=120] (R);
    \draw[->] (B1) to[out=15,in=-150] (B2);
    \draw[
      line width=.9pt, line cap=round,
      dash pattern=on 0pt off 4pt, 
    ] (B2) to[out=30,in=-130] (B3);
    \draw[->] (B3) to[out=50,in=-110] (R);
    \draw[->] (C) to[out=60,in=180] node[above]{$f$} (A1);
    \draw[->] (C) to[out=-60,in=180] node[below]{$g$} (B1);
    \node at (.5,-2) {};
  \end{tikzpicture}
  &\mspace{150mu}&
  \begin{tikzpicture}[scale=.8]
    \node (C) at (0,0) {$\bullet$};
    \coordinate (A1) at (2,1);
    \coordinate (A2) at (2+.2,1-.4);
    \coordinate (A3) at (2-.6,1-2*.4);
    \coordinate (A4) at (2+1.2,1-3*.4);
    \coordinate (A5) at (2-.2,1-4*.4);
    \coordinate (A6) at (3,-1);
    \draw (A1) to (A2);
    \draw (A2) to (A3);
    \draw[
      dotted,
      dash pattern=on 3pt off 3pt, 
    ] (A3) to (A4);
    \draw (A4) to (A5);
    \draw (A5) to (A6);
    \draw[->] (C) to[out=60,in=180] node[above]{$f$} (A1);
    \draw[->] (C) to[out=-60,in=180] node[below]{$g$} (A6);
    \draw[dotted,thick] (.5,1.5) to (.5,-2.5);
    \node at (1.1,-2.2) {$\bullet\succ$};
  \end{tikzpicture}
  \\
  \text{\footnotesize confluence}
  &&
  \text{\footnotesize tamed congruence}
\end{IEEEeqnarray*}
The notion of tamed preconfluence already mentioned in this introduction is a novel intermediate notion which we introduce in this article; see \cref{subsec:tamed_preconfluence}.

Tamed congruence, in contrast to confluence, is transitive: given three rewriting steps $f$, $g$ and $h$ with the same source, if both branchings $(f,g)$ and $(g,h)$ are tamely congruent, then the branching $(f,h)$ is tamely congruent.
This property gives an efficient way to deal with the so-called ``indexed branchings'' mentioned in \cref{subsubsec:intro_enumerating_overlaps}; see \cref{subsubsec:indexed_branchings} in the main text.
Here we give another example of its use, namely to bubble slides; see \cref{subsubsec:bubble_slides} in the main text.
Consider the following schematic picturing a generic overlap between two bubble-slide rewriting steps, where the bubble slides through two ``far away'' strands bounding the same region:
\begin{gather}
  \label{eq:intro_schematic_bubble_slide}
  \tikzpic{
    \clip (0,0) circle (2cm);
    \draw[dashed] (0,0) circle (2cm);
    \draw (-2,-.5) to[out=30,in=60] (-.5,-2);
    \draw[very thick,red,->] (-.5,-.5) to (-1.2,-1.2);
    \draw (-1,1.7) to[out=-90,in=-90] (1,1.7);
    \draw[very thick,red,->] (0,.7) to (0,1.7);
    \bul[-.5][0]\ca[-.5][0]\cu[-.5][0]
    \node[left=-1pt] at (-.5,0) {\scriptsize $k$};
  }[scale=.8]
\end{gather}
Denote by $A$ and $B$ the pieces of strand through which the bubble slides.
By walking along the boundary of the region from $A$ to $B$, we can find intermediate bubble-slide rewriting steps.
Going from one rewriting step to the next gives a branching akin to the last two critical bubble-slide branchings in \cref{subfig:critical_branchings_bubbles_vertsym}; that is, a branching that encapsulates how a bubble slides around a generating diagram.
If each of these branchings is tamely congruent, then by transitivity of tamed congruence, every generic branching as in \eqref{eq:intro_schematic_bubble_slide} is tamely congruent.

\medbreak

Consider branchings $(\ssf,\ssg)$ and $\Gamma[\ssf,\ssg]$, where the latter is the same branching $(\ssf,\ssg)$ but viewed inside a bigger diagram $\Gamma$.
As part of the global analysis, we ask: if $(\ssf,\ssg)$ is confluent, is $\Gamma[\ssf,\ssg]$ also confluent?
Since traffic rules depend on the context $\Gamma$, this may fail.
For instance, consider the branching $(\ssf,\ssg)$ and its confluence $(\ssf',\ssg')$ in \eqref{eq:intro_example_confluence}.
When viewed inside a bigger diagram, it may be that the cup carries a stop sign, so that the data $\Gamma[\ssf',\ssg']$ does not provide a valid confluence for $\Gamma[\ssf,\ssg]$.
Following a careful analysis, we show that among all the cases where the cup carries a stop sign, one only needs to consider one case, namely the second exceptional critical branching in \cref{subfig:critical_branchings_exceptional_vertsym}.
The first exceptional critical branching arises similarly.
As we noted above, these two critical branchings encode the bubble-evaluation and bubble-slide rewriting steps.

It is interesting to note that these branchings arise \emph{naturally} from the confluence analysis: the rewriting approach \emph{itself} leads us to (re)discover the bubble-evaluation and bubble-slide relations. Had we attempted to apply rewriting without them---and thus tried to prove a different basis theorem---the confluence analysis would have failed, indicating that we were missing important relations and that our basis conjecture was incorrect.

\subsection{Applications}
\label{subsec:intro_applications}

\subsubsection{The odd nil-Brauer category}
\label{subsubsec:intro_odd_nilBrauer}

Quantum symmetric pairs \cite{Letzter_SymmetricPairsQuantized_1999,BW_CanonicalBasesArising_2018,Kolb_QuantumSymmetricKac_2014}, or $\imath$quantum groups, are quantum analogues of symmetric pairs.
For instance, the split $\imath$quantum group of rank one $U^\imath_q(\mathfrak{sl}_2)$ arises as a quantum analogue of the symmetric pair $(\mathrm{SL}_2,\mathrm{SO}_2)$.
Quantum groups are categorified by Kac--Moody 2-categories \cite{KL_CategorificationQuantum$sln$_2010,Rouquier_2KacMoodyAlgebras_2008}, building on Lauda's categorification of quantum $\mathfrak{sl}_2$, the quantum group of rank one \cite{Lauda_CategorificationQuantum$sl2$_2010}.
The nil-Brauer category plays the analogue of Lauda's categorification in the $\imath$world, leading to a categorification of $U^\imath_q(\mathfrak{sl}_2)$ \cite{BWW_NilBrauerCategorifiesSplit_2025}.
In this work, we give an odd analogue of the nil-Brauer category:

\begin{restatable}[label=defn:oddnilBr]{definition}{defnoddnilBr}
  \renewcommand{\sca}{.4}
  \newcommand{\vspc}{.5ex}
  \newcommand{\hspc}{50mu}
  Let $\mK$ be an integral domain in which $2$ is invertible and let $\VARev\in\{0,-1\}$.
  Let $p(\cup)\coloneqq\VARev+1\mod 2$.
  The \defnemph{odd nil-Brauer category $\oddnilBr[\VARev]$} is the $\mK$-linear monoidal category generated by a single object $\tikzpic{\li}$ and morphisms
  \begin{IEEEeqnarray*}{CcCcCcC}
    \tikzpic{\cu}
    &\qquad&
    \tikzpic{\ca}
    &\qquad&
    \tikzpic{\cro}
    &\qquad&
    \tikzpic{\bli}\;,
    \\
    p(\cup)
    &&
    p(\cup)
    &&
    \parodd
    &&
    \parodd
  \end{IEEEeqnarray*}
  up to the following relations:
  \begin{gather*}
    \tikzpic{\cro \cro[0][1]} 
    \;=\;
    0
    \mspace{\hspc}
    %
    \tikzpic{
      \cro[0][2] \li[2][2]
      \li[0][1] \cro[1][1] 
      \cro \li[2][0] 
    }[scale=.8] 
    \;=\;
    \tikzpic{
      \cro[1][2] \li[0][2]
      \li[2][1] \cro[0][1] 
      \cro[1][0] \li[0][0] 
    }[scale=.8]
    \mspace{\hspc}
    \tikzpic{\cu\ca[0][0]} 
    \;=\;\VARev
    \\[\vspc]
    \tikzpic{
      \cu[1][0] \uli[2][0] 
      \dli \ca 
    }
    \;=\;
    \tikzpic{\li}
    %
    \mspace{\hspc}
    %
    \tikzpic{
      \ca[0][1]
      \cro 
    }
    \;=\;0
    \mspace{\hspc}
    %
    \tikzpic{
      \uli[0][1] \ca[1][1] 
      \cro \li[2][0] 
    }
    \;=\;
    \tikzpic{ 
      \ca[0][1] \uli[2][1]
      \li \cro[1][0] 
    }
    \\[\vspc]
    \tikzpic{
      \ulcr
    }
    \;=\;-\;
    \tikzpic{\drcr}
    \;-\;
    \tikzpic{\lili} 
    \;-(-1)^{p(\cup)}\;
    \tikzpic{\cc} 
    \mspace{\hspc}
    \tikzpic{\rca}
    \;=\;-(-1)^{p(\cup)}\;
    \tikzpic{\lca}
  \end{gather*}
\end{restatable}

\begin{restatable}[label=thm:oddnilBr_is_affine_Brauer_type]{theorem}{oddnilBraffineBrauerthm}
  The odd nil-Brauer category $\oddnilBr[\VARev]$ (\cref{defn:oddnilBr}) is a category of affine Brauer type with $(2\bN+1)$-bubbles.
\end{restatable}

We expect that the odd analogue of that result holds as well:

\begin{conjecture}
  The odd nil-Brauer categorifies the split $\imath$quantum group of rank one.
\end{conjecture}

The proof in the even case relies on the basis theorem for the nil-Brauer category.
Given \cref{thm:oddnilBr_is_affine_Brauer_type}, we expect that the arguments of \cite{BWW_NilBrauerCategorifiesSplit_2025} can be adapted to the odd case.

Kac--Moody 2-categories admit odd analogues \cite{EL_OddCategorification$U_qmathfraksl_2$_2016,KKO_SupercategorificationQuantumKacMoody_2013,BE_SuperKacMoody_2017}, although their basis theorem remains conjectural outside rank one \cite{BK_OddGrassmannianBimodules_2026}.
At the decategorified level, their parity functor becomes a parameter $\pi$ such that $\pi^2=1$.
Quantum covering groups $U_{q,\pi}(\mathfrak{g})$ \cite{CH_QuantumSupergroupsBraid_2016,CHW_QuantumSupergroupsFoundations_2013,CHW_QuantumSupergroupsII_2014,CFL+_QuantumSupergroupsIII_2014,Clark_QuantumSupergroupsIV_2014,CSW_QuantumSupergroupsVI_2019} are associative algebras with such parameter $\pi$, so that setting $\pi=1$ recovers certain classical quantum groups, while setting $\pi=-1$ recovers certain super quantum groups.
Odd analogues to Kac--Moody 2-categories are conjectured to categorify covering quantum groups once their parity functor is properly decategorified to a formal parameter $\pi$.
Recently, Chung extended covering quantum groups to quantum symmetric pairs \cite{Chung_CanonicalBasesArising_2021,Chung_SerrePresentation$imath$quantum_2019}.
Given the above, one may expect that the odd nil-Brauer category provides a categorification of Chung's work in rank one.

\subsubsection{The quantized affine VW supercategory}
\label{subsubsec:intro-quantVW}

The periplectic Lie superalgebra $\mathfrak{p}(n)$ satisfies a Schur--Weyl duality with the periplectic Brauer algebra (also known as the odd or signed Brauer algebra). This duality has played an important role in representation theory; see, for example, \cite{Coulembier_PeriplecticBrauerAlgebra_2018}.

In \cite{AGG_QuantizedEnvelopingSuperalgebra_2021}, the authors introduced a quantized version $U_q(\mathfrak{p}(n))$ of the periplectic Lie superalgebra together with a $q$-deformation $\mathcal{B}_{q,l}$ of the periplectic Brauer algebra. These act on tensor powers of the natural representation $V=\mathbb{C}_q(n|n)$ by commuting actions:
\[
U_q(\mathfrak{p}(n))\curvearrowright V^{\otimes l} \curvearrowleft \mathcal{B}_{q,l},
\] 
leading to a quantum analogue of periplectic Schur--Weyl duality.
Subsequently, a monoidal diagrammatic category $\qpBr_q$ was constructed in \cite{RS_Periplectic$q$BrauerCategory_2025} whose endomorphism algebras are these $q$-periplectic Brauer algebras.

A different generalization of periplectic Schur--Weyl duality has been developed in \cite{BDE+_AffineVWSupercategory_2020}. There the authors introduced a higher Schur--Weyl duality by tensoring the tensor powers of the natural representation with a $\mathfrak{p}(n)$-module $M$. In this setting, the action of the periplectic Lie superalgebra commutes with the action of the affine periplectic Brauer algebra (or affine $\sVW$ superalgebra):
\[
\mathfrak{p}(n)\curvearrowright M\otimes V^{\otimes l} \curvearrowleft \sVW_l.
\]

Combining these two directions of generalization, we define the quantized affine VW supercategory.
Its full subcategory generated by $\tikzpic{\cro},\tikzpic{\ca},\tikzpic{\cu}$ is the periplectic $q$-Brauer category $\qpBr_q$ of \cite{RS_Periplectic$q$BrauerCategory_2025}, the crossing being the negative crossing of \opcit, while the specialization at $q=-1$ recovers the category $\sVW$.

We expect this quantized affine VW supercategory to satisfy a higher quantum Schur--Weyl duality with a quantized version of $\mathfrak{p}(n)$, analogous to the higher Schur--Weyl duality relating $\mathfrak{p}(n)$ and $\sVW$.

\begin{restatable}[label=defn:qsVW_form]{definition}{qsVWformdefn}
  \renewcommand{\sca}{.4}
  \newcommand{\vspc}{1.2ex}
  \newcommand{\hspc}{40mu}
  The \defnemph{quantized affine VW supercategory} $\quantsVW$ is the $\bZ[q,q^{-1}]$-linear supermonoidal category generated by a single object $\tikzpic{\li}$ and morphisms $\tikzpic{\cro}$ (even), $\tikzpic{\ca}$, $\tikzpic{\cu}$ (odd), $\tikzpic{\bli}$ (even), subject to the super interchange law and the following relations:
  \begin{gather*}
    \tikzpic{\cro \cro[0][1]}
    \;=\;
    \tikzpic{\lili}
    \;-\;(q-q^{-1})\;
    \tikzpic{\cro}
    \mspace{\hspc}
    %
    \tikzpic{
      \cro[0][2] \li[2][2]
      \li[0][1] \cro[1][1]
      \cro \li[2][0]
    }[scale=.8]
    \;=\;
    \tikzpic{
      \cro[1][2] \li[0][2]
      \li[2][1] \cro[0][1]
      \cro[1][0] \li[0][0]
    }[scale=.8]
    \\[\vspc]
    \tikzpic{
      \cu[1][0] \uli[2][0]
      \dli \ca
    }
    \;=\;
    \tikzpic{\li}
    \mspace{\hspc}
    %
    \tikzpic{
      \ca[0][1]
      \cro
    }
    \;=\;-q\;
    \tikzpic{\ca}
    \mspace{\hspc}
    %
    \tikzpic{
      \uli[0][1] \ca[1][1]
      \cro \li[2][0]
    }
    \;=\;
    \tikzpic{
      \ca[0][1] \uli[2][1]
      \li \cro[1][0]
    }
    \;-\;(q-q^{-1})\;
    \tikzpic{\li\ca[1][0]}
    \\[\vspc]
    \tikzpic{\ulcr}
    \;=\;
    \tikzpic{\drcr}
    \;+\;q\;
    \tikzpic{\lili}
    \;-\;(q-q^{-1})\;
    \tikzpic{\bli\li[1][0]}
    \;+\;q^{-1}\;
    \tikzpic{\cc}
    \;+\;(q^{-3}-q^{-1})\;
    \tikzpic{\lca\cu[0][1]}
    \\[\vspc]
    \tikzpic{\rca}
    \;=\;q^{-2}\;
    \tikzpic{\lca}
    \;+\;
    \tikzpic{\ca}
    \mspace{\hspc}
    \tikzpic{\bubble[0][0][k]}[][1]\;=\;0
  \end{gather*}
\end{restatable}

\begin{restatable}[label=thm:quantsVW_is_affine_Brauer_type]{theorem}{quantsVWaffineBrauerthm}
  The quantized affine VW supercategory $\quantsVW$ (\cref{defn:qsVW_form}) is a category of affine Brauer type with bubbles $\bubbleset=\emptyset$.
\end{restatable}

\subsubsection{Classification results}
\label{subsubsec:intro_classification}

Our final result is a classification of a certain subclass of categories of affine Brauer type, namely, those with a relatively simple presentation:

\begin{definition}\label{defn:airy_presentation_affine_Brauer_type}
  \renewcommand{\sca}{.4}
  \newcommand{\vspc}{.8ex}
  \newcommand{\hspc}{40mu}
  Fix a commutative ring $\Bbbk$.
  An \defnemph{airy presentation of affine Brauer type} is a presentation of a supermonoidal category with generators crossing, cup, cap and dot, and generating relations
  \begin{gather*}
    \tikzpic{\cro \cro[0][1]} 
    \;=\;
    \VARRtwo\;
    \tikzpic{\lili \lili[0][1]}
    \mspace{\hspc}
    %
    \tikzpic{
      \cro[0][2] \li[2][2]
      \li[0][1] \cro[1][1] 
      \cro \li[2][0] 
    }[scale=.8] 
    \;=\;\VARRthree\;
    \tikzpic{
      \cro[1][2] \li[0][2]
      \li[2][1] \cro[0][1] 
      \cro[1][0] \li[0][0] 
    }[scale=.8]
    \mspace{\hspc}
    \tikzpic{\cu\ca[0][0]} 
    \;=\;\VARev
    \\[\vspc]
    \tikzpic{
      \cu[1][0] \uli[2][0] 
      \dli \ca 
    }
    \;=\;\VARzz\;
    \tikzpic{\li}
    %
    \mspace{\hspc}
    %
    \tikzpic{
      \ca[0][1]
      \cro 
    }
    \;=\;\VARuk\;
    \tikzpic{\ca}
    \mspace{\hspc}
    %
    \tikzpic{
      \uli[0][1] \ca[1][1] 
      \cro \li[2][0] 
    }
    \;=\;\VARcasl\;
    \tikzpic{ 
      \ca[0][1] \uli[2][1]
      \li \cro[1][0] 
    }
    \\[\vspc]
    \tikzpic{
      \ulcr
    }
    \;=\;\VARdcro\;
    \tikzpic{\drcr}
    \;+\;\VARdcroID\;
    \tikzpic{\lili} 
    \;+\;\VARdcroCC\;
    \tikzpic{\cc} 
    \mspace{\hspc}
    \tikzpic{\rca}
    \;=\;\VARdcc\;
    \tikzpic{\lca}
    \;+\;\VARdccLOT\;
    \tikzpic{\ca}
  \end{gather*}
  where each symbol is a scalar in $\Bbbk$ and the scalars $\VARRthree$, $\VARzz$, $\VARcasl$, $\VARdcro$, $\VARdcroID$ and $\VARdcc$ are invertible.
\end{definition}

Perhaps the most consequential condition on parameters is the invertibility of $\VARzz$, which ensures that we can ``rotate'' the defining relations to obtain the rest of the rewriting system of affine Brauer type.
In \cref{subsec:classification}, we show the following classification result:

\begin{restatable}{theorem}{thmclassificationaslincategories}
  \label{thm:classification_as_lin_categories}
  Fix a field $\Bbbk$ where $2$ is invertible.
  Let $\cC$ be a supermonoidal category over $\Bbbk$ presented by an airy presentation of affine Brauer type (\cref{defn:airy_presentation_affine_Brauer_type}).
  Then $\cC$ is a category of affine Brauer type if and only if it is isomorphic, as a linear category, to $\affBr[\VARev]$, $\sVW$, $\nilBr[\VARev]$ or $\oddnilBr[\VARev]$.
\end{restatable}

We expect that the same techniques can be applied to more involved presentations (\eg presentations as in \cref{subsec:implementation}) and view this result, together with the results of \cite{Barbier_DiagramCategoriesBrauer_2024}, as first steps in a classification program.

\medbreak

We stress that \cref{thm:classification_as_lin_categories} is a classification \emph{as linear categories}.
Indeed, we can construct new categories of affine Brauer type out of $\affBr[\VARev]$, $\sVW$, $\nilBr[\VARev]$ and $\oddnilBr[\VARev]$ by renormalizing the generators via non-monoidal, and sometimes non-graded, normalizations:

\begin{lemma}\label{lem:isomorphism-simple-minded-not-monoidal}
  Let $\cC$ and $\widehat{\cC}$ be two supermonoidal categories presented by an airy presentation of affine Brauer type.
  We denote the parities and parameters of $\widehat{\cC}$ with a hat, \eg $\widehat{p}(\cup)$ and $\widehat{\VARuk}$.
  Denote by $U(\cC)$ and $U(\widehat{\cC})$ their underlying linear categories.
  \begin{enumerate}[(i)]
    \item Let $t$ be an element of the ground ring such that $t^2=(-1)^{p(\cup)+\widehat{p}(\cup)}$.
    If $\widehat{\VARev}= t\VARev$, $\widehat{\VARcasl}= t\VARcasl$ and $\widehat{\VARdcroCC}=t^{-1}\VARdcroCC$, and the remaining parameters are identical, then the following
    \begin{IEEEeqnarray*}{rClcrCl}
    \tikzpic{\li}[][.5]\ldots\!\tikzpic{\cu[0][1]\node[above=-2pt] at (0,1) {$\scriptstyle{i}$};}[][.5]\ldots\tikzpic{\li}[][.5]
    &\mapsto&
    t^{i}\;\tikzpic{\li}[][.5]\ldots\!\tikzpic{\cu[0][1]\node[above=-2pt] at (0,1) {$\scriptstyle{i}$};}[][.5]\ldots\tikzpic{\li}[][.5]
    &\qquad&
    \tikzpic{\li}[][.5]\ldots\!\tikzpic{\ca\node[below=-2pt] at (0,0) {$\scriptstyle{i}$};}[][.5]\ldots\tikzpic{\li}[][.5]
    &\mapsto&
    t^{-i-1}\;\tikzpic{\li}[][.5]\ldots\!\tikzpic{\ca\node[below=-2pt] at (0,0) {$\scriptstyle{i}$};}[][.5]\ldots\tikzpic{\li}[][.5]
    \end{IEEEeqnarray*}
    defines an isomorphism of linear categories between $U(\cC)$ and $U(\widehat{\cC})$ (here $i$ numbers the strands).

    \item Let $t$ be an element of the ground ring such that $t^2=1$.
    Assume $\VARuk=\widehat{\VARuk} = 0$ and $\VARdccLOT=\widehat{\VARdccLOT} = 0$.
    If $\widehat{\VARRthree}= t\VARRthree$,  $\widehat{\VARdcro}= t\VARdcro$, $\widehat{\VARdcc}= t\VARdcc$, $\widehat{\VARev}= t\VARev$ and  $\widehat{\VARdcroCC}= t\VARdcroCC$, and the remaining parameters are identical, then the following
    \begin{IEEEeqnarray*}{rClcrCl}
    \tikzpic{\li}[][.5]\ldots\!\tikzpic{\cro\node[below=-2pt] at (0,0) {$\scriptstyle{i}$};}[][.5]\ldots\tikzpic{\li}[][.5]
    &\mapsto&
    t^i\;\tikzpic{\li}[][.5]\ldots\!\tikzpic{\cro\node[below=-2pt] at (0,0) {$\scriptstyle{i}$};}[][.5]\ldots\tikzpic{\li}[][.5]
    &\qquad&
    \tikzpic{\li}[][.5]\ldots\!\tikzpic{\bli\node[below=-2pt] at (0,0) {$\scriptstyle{i}$};}[][.5]\ldots\tikzpic{\li}[][.5]
    &\mapsto&
    t^{i}\;\tikzpic{\li}[][.5]\ldots\!\tikzpic{\bli\node[below=-2pt] at (0,0) {$\scriptstyle{i}$};}[][.5]\ldots\tikzpic{\li}[][.5]
    \\
    \tikzpic{\li}[][.5]\ldots\!\tikzpic{\cu[0][1]\node[above=-2pt] at (0,1) {$\scriptstyle{i}$};}[][.5]\ldots\tikzpic{\li}[][.5]
    &\mapsto&
    t^{i}\;\tikzpic{\li}[][.5]\ldots\!\tikzpic{\cu[0][1]\node[above=-2pt] at (0,1) {$\scriptstyle{i}$};}[][.5]\ldots\tikzpic{\li}[][.5]
    &&
    \tikzpic{\li}[][.5]\ldots\!\tikzpic{\ca\node[below=-2pt] at (0,0) {$\scriptstyle{i}$};}[][.5]\ldots\tikzpic{\li}[][.5]
    &\mapsto&
    t^{i+1}\;\tikzpic{\li}[][.5]\ldots\!\tikzpic{\ca\node[below=-2pt] at (0,0) {$\scriptstyle{i}$};}[][.5]\ldots\tikzpic{\li}[][.5]
    \end{IEEEeqnarray*}
    defines an isomorphism of ($\bZ/2\bZ$-graded) linear categories between $U(\cC)$ and $U(\widehat{\cC})$.
  \end{enumerate}
\end{lemma}

In particular, if we are in situations (i) or (ii), then $\cC$ is of affine Brauer type with $\bubbleset$-bubbles if and only if $\widehat{\cC}$ is, since being of affine Brauer type is a statement about hom-spaces (and the parity of the dot).
These isomorphisms are not monoidal, and in situation (i) with $t$ a fourth root of unity, not even graded: the isomorphism flips the parity of the cup and cap.
As it turns out, applying \cref{lem:isomorphism-simple-minded-not-monoidal} to $\affBr[\VARev]$, $\sVW$, $\nilBr[\VARev]$ and $\oddnilBr[\VARev]$ gives a complete list of supermonoidal structures; see \cref{classification simple-minded} for the explicit classification.

\subsection{Implementation}
\label{subsec:implementation}

\newcommand{\VARdcroIDlb}{b_2}
\newcommand{\VARdcroIDrb}{b_1}
\newcommand{\VARdcroCCub}{c_2}
\newcommand{\VARdcroCCdb}{c_1}

We implement an algorithm in \textsc{Form} that verifies the non-affine and dot-slide critical branchings of \cref{fig:critical_branchings_vertsym} (and their vertical symmetries) automatically. The code  is available on GITHUB at \url{https://github.com/BarbierSK/AffineBrauerReductionAlgorithm}; see \cite{BS_BarbierSKAffineBrauerReductionAlgorithm_2026}. \textsc{Form} is a symbolic manipulation system that allows one to replace certain expressions by linear combinations of other expressions quickly and efficiently \cite{FORM4,FORM5}. As such, it is especially well adapted to implementing rewriting systems.

The algorithm works as follows. We start from a diagram and first apply non-affine rewriting steps until this is no longer possible. We then decorate the diagram with the traffic rules and apply dot-sliding rewriting steps accordingly. This does not yet give us a normal form, since we also need to take the interchange rules into account.
We therefore interchange some generators and then again apply the non-affine and dot-sliding rewriting steps. To make our algorithm terminate, we give a direction to the interchange rule by pushing cups up, caps down, and dots in the direction of their orientation. Unfortunately, for the crossings, there is no canonical direction and we need both directions to obtain normal forms (see the comments in the implementation for an example of how a naive implementation of the crossing interchange could lead to the algorithm getting stuck on a non-normal form). We solve this by first pushing the crossings in one direction, simplifying and then pushing in the other direction until we end up with normal forms.

We use the following names for the parameters appearing in the rewriting steps of the algorithm:
\begingroup
\renewcommand{\sca}{.4}
\newcommand{\vspc}{.8ex}
\newcommand{\hspc}{40mu}
\begin{gather*}
  \tikzpic{\cro \cro[0][1]} 
  \;=\;
  \VARRtwo\;
  \tikzpic{\lili}
  +
  e\;
  \tikzpic{\cro}
  +
  f\;
  \tikzpic{\cc}
  \mspace{\hspc}
  %
  \tikzpic{
    \cro[0][2] \li[2][2]
    \li[0][1] \cro[1][1] 
    \cro \li[2][0] 
  }[scale=.8]
  \;=\;\VARRthree\;
  \tikzpic{
    \cro[1][2] \li[0][2]
    \li[2][1] \cro[0][1] 
    \cro[1][0] \li[0][0] 
  }[scale=.8]
  \mspace{\hspc}
  \tikzpic{\cu\ca[0][0]} 
  \;=\;\VARev
  \\[\vspc]
  \tikzpic{
    \cu[1][0] \uli[2][0] 
    \dli \ca 
  }
  \;=\;\VARzz\;
  \tikzpic{\li}
  %
  %
  %
  \mspace{\hspc}
  %
  \tikzpic{
    \ca[0][1]
    \cro 
  }
  \;=\;\VARuk\;
  \tikzpic{\ca}
  \mspace{\hspc}
  %
  \tikzpic{
    \cu
    \cro 
  }
  \;=\;\tilde{\VARuk}\;
  \tikzpic{\cu}
  \mspace{\hspc}
  %
  \\[\vspc]
  \tikzpic{
    \uli[0][1] \ca[1][1] 
    \cro \li[2][0] 
  }
  \;=\;\kappa\;
  \tikzpic{\ca\li[2][0]}
  +\mu\;
  \tikzpic{ 
    \ca[0][1] \uli[2][1]
    \li \cro[1][0] 
  }
  +
  \nu\;
  \tikzpic{\li\ca[1][0]}
  \mspace{\hspc}
  %
  \tikzpic{ 
    \cro \li[2][0]
    \dli \cu[1][0] 
  }
  \;=\;\tilde{\kappa}\;
  \tikzpic{\cu[0][1]\li[2][0]}
  +\tilde{\mu}\;
  \tikzpic{
    \li[0][0] \cro[1][0]
    \cu[0][0] \dli[2][0]
  }
  +\tilde{\nu}\;
  \tikzpic{\li\cu[1][1]}
  \\[\vspc]
  \tikzpic{
    \ulcr
  }
  \;=\;\VARdcro\;
  \tikzpic{\drcr}
  \;+\;\VARdcroID\;
  \tikzpic{\lili} 
  \;+\;\VARdcroIDrb\;
  \tikzpic{\bli\li[1][0]} 
  \;+\;\VARdcroIDlb\;
  \tikzpic{\li\bli[1][0]} 
  \;+\;\VARdcroCC\;
  \tikzpic{\cc} 
  \;+\;\VARdcroCCub\;
  \tikzpic{\lca\cu[0][1]} 
  \;+\;\VARdcroCCdb\;
  \tikzpic{\ca\lcu[0][1]} 
  \mspace{\hspc}
  \tikzpic{\rca}
  \;=\;\VARdcc\;
  \tikzpic{\lca}
  \;+\;\VARdccLOT\;
  \tikzpic{\ca}
  \\[\vspc]
  \tikzpic{
    \urcr
  }
  \;=\;\tilde\VARdcro\;
  \tikzpic{\dlcr}
  \;+\;\tilde\VARdcroID\;
  \tikzpic{\lili}
  \;+\;\tilde\VARdcroIDrb\;
  \tikzpic{\bli\li[1][0]}
  \;+\;\tilde\VARdcroIDlb\;
  \tikzpic{\li\bli[1][0]}
  \;+\;\tilde\VARdcroCC\;
  \tikzpic{\cc}
  \;+\;\tilde\VARdcroCCub\;
  \tikzpic{\lca\cu[0][1]}
  \;+\;\tilde\VARdcroCCdb\;
  \tikzpic{\ca\lcu[0][1]}
  \mspace{\hspc}
  \tikzpic{\rcu}
  \;=\;\tilde\VARdcc\;
  \tikzpic{\lcu}
  \;+\;\tilde\VARdccLOT\;
  \tikzpic{\cu}
\end{gather*}
\endgroup
In the program we can choose a specific set of parameters or consider the generic parameters above.
The critical branchings are encoded as the difference $(\text{path }A)-(\text{path }B)$. For a specific set of parameters, the resulting normal forms should all be zero for the parameters to belong to a well-defined category of affine Brauer type. Starting from generic parameters, the coefficients of the resulting diagrams we obtain after applying the algorithm give us conditions on the parameters.
To compute beyond a presentation such as the one above, one can modify the code of the rewriting steps.


%% file: sections/rewriting.tex
%
%



\section{Monoidal Gröbner systems}
\label{sec:RW_monoidal_grobner_system}

We define the foundations of rewriting in linear strict (graded-)monoidal categories, based on higher linear rewriting theory \cite{Schelstraete_RewritingModuloDiagrammatic_2026}.
In particular, we define a \emph{linear monoidal rewriting system} (or \LMRS{}), see \cref{defn:MLRS}, as a special case of a higher linear rewriting system modulo.
A linear monoidal rewriting system $\sS$ comes with the extra data of a \emph{sub-system} $\lT$ (\cref{defn:RW_subsystem}) and an \emph{isp-order} $\succ$ (\cref{subsec:RW_termination_order}).
A triple $(\sS,\lT,\succ)$ encodes a reduction-to-basis algorithm if and only if it defines a \emph{monoidal Gröbner system} (\cref{defn:RW_monoidal_grobner_system}); see \cref{thm:RW_monoidal_grobner_system}.

The main purpose of rewriting is to deduce that a triple $(\sS,\lT,\succ)$ defines a monoidal Gröbner system from a reasonable set of explicit computations.
\Cref{subsec:RW_confluence_analysis} explains the generic techniques common to every linear monoidal rewriting system, reducing the set of computations to minimal overlapping branchings (\cref{subsubsec:enumeration_overlapping_brauer_type}).
In contrast with (say) Bergman's diamond lemma for associative algebras, this class of branchings often remains too large to handle directly, and further analysis is needed. We carry out this analysis in \cref{sec:critical_branchings} for categories of affine Brauer type.
The reader will find further rewriting ideas there, which may apply to other settings.

While we restrict to linear strict (graded-)monoidal categories, the theory extends verbatim to linear strict (graded)-2-categories \cite{SV_OddKhovanovHomology_2023}.

\begin{remark}[Comparison with higher linear rewriting theory]
  \label{rem:RW_comparison_with_literature}
  Unless otherwise stated, the definitions and results of this section are special cases of \cite{Schelstraete_RewritingModuloDiagrammatic_2026}, although many ideas in \opcit already appeared before in the literature; we refer to \opcit for references.
  The notion of a monoidal Gröbner system is new, and is meant to gather the relevant ideas from \cite{Schelstraete_RewritingModuloDiagrammatic_2026} in a succinct manner.
  Another exception is \cref{lem:RW_overlapping_branchings_intersect}, which as far as we know has not appeared before in the literature.
\end{remark}

\begin{notation}
  \label{not:RW_identities}
  We often abuse notation and write $i$ for the identity $\id_i$ of an object $i$ in a given category. For instance, if $f$ is a morphism, then $f\otimes i \coloneqq f\otimes\id_i$.
\end{notation}

\subsection{Graded-monoidal categories}
\label{subsec:RW_graded_monoidal_categories}

Throughout we fix a commutative ring $\ring$, an abelian group $G$ and a bicharacter map ${\bil\colon G\times G\to \ring^\times}$ (that is, $\bil$ is a bilinear map of abelian groups). For simplicity, we further assume that $\bil$ is \defnemph{skew-symmetric}, in the sense that
$\bil(a,b)^{-1} = \bil(b,a)$ for all $a,b\in G$.
Whenever $x$ is a homogeneous element of a $G$-graded set or $\ring$-module, we write $\deg_G (x)\in G$ for its degree.
Let $\ring\md\mathrm{Mod}^0_{G,\bil}$ be the symmetric monoidal $\ring$-linear category consisting of $G$-graded $\ring$-modules, degree-preserving $\ring$-linear maps, the usual monoidal structure and the symmetric structure given by
\[(v,w)\mapsto \bil(\deg_G v,\deg_G w)(w,v).\]

\begin{notation}
  \label{not:scalar_interchange}
  In \cref{subsec:classification}, we use the shorthand
  \[\inter{v}{w}\coloneqq\bil(\deg_G v,\deg_G w),\]
  which emphasize that $v$ is above $w$ right before applying the interchange. In this section we keep using the standard notation.
\end{notation}
\begin{definition}
  A  \emph{$(G,\bil)$-graded-2-category} is a category twice-enriched over $\ring\md\mathrm{Mod}^0_{G,\bil}$. A \emph{$(G,\bil)$-graded-monoidal category} (or simply \emph{graded-monoidal category}) is a one-object $(G,\bil)$-graded-2-category.
\end{definition}

In other words, a $(G,\bil)$-graded-monoidal category is almost the same as a $G$-graded linear (strict) monoidal category, except that the interchange law is twisted by $\bil$:
\begin{gather}
  \label{eq:graded_interchange_law}
  \def\scl{.6}
  \tikzpic{
    \indexNF[0][1][{x}]\lili[2][1]
    \lili\indexNF[2][0][{y}]
  }[scale=1.2][1]
  \;=\;\bil(\deg_G x,\deg_G y)\;
  \tikzpic{
    \lili[0][1]\indexNF[2][1][{y}]
    \indexNF[0][0][{x}]\lili[2][0]
  }[scale=1.2][1]
  \\
  \nonumber
  \text{\footnotesize graded interchange law}
\end{gather}
Note that a $(G,\bil)$-graded-monoidal category is \emph{not} a monoidal category with extra structure.
When $\bil$ is trivial, however, a $(G,\bil)$-graded-monoidal category is just a $G$-graded linear monoidal category.
A supermonoidal category \cite{BE_MonoidalSupercategories_2017} (also known as a ``monoidal supercategory'') is the same as a $(G,\bil)$-graded-monoidal category where $G=\bZ/2\bZ$ and $\bil(a,b) = (-1)^{a\cdot b}$.

\begin{definition}
  Let $\cC$ be a $(G,\bil)$-graded-monoidal category. A \emph{hom-basis} is a family of subsets $B(i,j)\subset\cC(i,j)$ indexed by pairs of objects $(i,j)$, such that $B(i,j)$ is a basis of $\cC(i,j)$.
\end{definition}


\subsection{Linear monoidal rewriting systems}
\label{subsec:RW_monoidal_linear_rewriting_system}

\begin{definition}
  \label{defn:MLRS}
  A \defnemph{linear monoidal rewriting system} (or \LMRS{}) is the data
  \[\sS = (\sX_0,\sX_1;\sR,\sE),\]
  where $\sX_0$, $\sX_1$, $\sR$ and $\sE$ are described in the bullet points below.
\end{definition}

In a nutshell, a \LMRS{} is a presentation of a graded-monoidal category $\cC$ where generating relations are endowed with an orientation.
In detail, $\sS$ consists of the following data:

\smallbreak

\begin{itemize}
  \item $\sX_0$ is a set, corresponding to the generating objects of $\cC$.
\end{itemize}
We write $\obj(\sS)$ for the free monoid on $\sX_0$---that is, $\obj(\sS)$ is the set of objects of~$\cC$.

\smallbreak

\begin{itemize}
  \item $\sX_1$ is a family
  \[\sX_1\coloneqq\big(\sX_1(i,j)\big)_{i,j\in\obj(\sS)}\]
  of $G$-graded sets, corresponding to the generating morphisms of $\cC$.
\end{itemize}
For $\varphi\in\sX_1(i,j)$, we write $s(\varphi)=i$ (the source) and $t(\varphi)=j$ (the target).
A generator $\varphi\in \sX_1$ can be ``extended'' by identities of objects, giving a \defnemph{whiskered generator} $\id_{k}\otimes \varphi\otimes\id_{k'}$ for each $k,k'\in\obj(\sS)$.
The source and target maps extend to whiskered generators in the natural way, and whiskered generators with matching source and target can be composed. The $G$-grading extends as
\[\deg_G(\id_{k}\otimes \varphi\otimes\id_{k'})=\deg_G(\varphi).\]
For each pair of objects $(i,j)$, we write $\sX^*(i,j)$ for the set of compositions of whiskered generators with source $i$ and target $j$. An element of $\sX^*(i,j)$ is called a \defnemph{monomial}, or when thinking in terms of diagrammatics, a \defnemph{diagram}.

For each pair of objects $(i,j)$, we write $\sX^l(i,j)\coloneqq\langle\sX^*(i,j)\rangle_\ring$, the free $G$-graded $\ring$-module generated by the set $\sX^*(i,j)$. An element of $\sX^l(i,j)$ is called a \defnemph{vector}.
We write $\sX^*$ (resp.\ $\sX^l$) for the union of all the $\sX^*(i,j)$'s (resp.\ $\sX^l(i,j)$'s).

\smallbreak

\begin{itemize}
  \item $\sR$ is a family
  \[\sR\coloneqq\big(\sR(i,j)\big)_{i,j\in\obj(\sS)}\;,\]
  such that each $r\in\sR(i,j)$ is endowed with a source $s(r) \in\sX^*(i,j)$ and target $t(r)\in\sX^l(i,j)$, where $t(r)$ is concentrated in one degree $\deg_G(t(r))$, and $\deg_G(s(r))=\deg_G(t(r))$. We write it as
  \[r\colon s(r)\to t(r).\]
  An element $r\in\sR(i,j)$ is called a \defnemph{generating rewriting step}.
\end{itemize}
Each $r$ induces a generating relation $s(r)-t(r)$. A \LMRS{} gives the extra data of an orientation $s(r)\to t(r)$. Note that while the target $v=t(r)\in \sX^l$ is (a priori) a linear combination of monomials (a vector), the source $x=s(r)\in\sX^*$ must be a monomial.
We often abuse notation and write $\sR$ for the union of the sets $\sR(i,j)$.
Note that on the set $\sR(i,j)$, we have the globular conditions $s\circ s= s\circ t=i$ and $t\circ s= t\circ t=j$.

\medbreak

We call \defnemph{$(G,\bil)$-graded interchangers} the graded interchange law \eqref{eq:graded_interchange_law} oriented from left to right.
The remaining data in the quadruple $\sS$ encapsulates the graded interchange law:
\begin{itemize}
    \item $\sE$ is the data of $(G,\bil)$-graded interchangers---given $(\sX_0,\sX_1)$, it is the same data as the data of the bilinear map $\bil$.
\end{itemize}
This ends the description of the linear monoidal rewriting system $\sS$.
\hfill$\diamond${\parfillskip0pt\par}

\subsection{Contexts and sub-systems}
\label{subsec:RW_contexts_subsystems}

A \defnemph{context} $\Gamma$ is a ``diagram with a hole'':
\begin{equation}
  \label{eq:defn_context}
  \Gamma \;=\; 
  \tikzpic{
    \def\hshift{.2}
    \def\vshift{.1}
    \draw (0,0) to (0,3);
    \draw (1,.5) to (1,2.5);
    \draw (2,.5) to (2,2.5);
    \draw (3,.5) to (3,2.5);
    \draw (4,.5) to (4,2.5);
    \draw (5,0) to (5,3);
    \draw[fill=white,rounded corners=1pt] (0-\hshift,0+\vshift) rectangle (5+\hshift,1-\vshift);
    \draw[fill=white,rounded corners=1pt] (2-\hshift,1+\vshift) rectangle (3+\hshift,2-\vshift);
    \draw[fill=white,rounded corners=1pt] (0-\hshift,2+\vshift) rectangle (5+\hshift,3-\vshift);
    \node at (2.5,.5) {\footnotesize $x$};
    \node at (2.5,2.5) {\footnotesize $y$};
    \node at (.5,1.5) {\footnotesize $i$};
    \node at (4.5,1.5) {\footnotesize $j$};
  }\;,
\end{equation}
where $i,j\in\obj(\sS)$ are objects and $x,y\in\sX^*$ are monomials (\ie diagrams), suitably composable.
Given a context $\Gamma$, we can \emph{contextualize} a generating rewriting step $r\colon s(r)\to t(r)$ as
\[{\Gamma[r]\colon \Gamma[s(r)] \to \Gamma[t(r)]},\]
``surrounding'' $r$ with $\Gamma$ (provided source and target of $r$ are compatible with the context $\Gamma$).
We call $\Gamma[r]$ a \defnemph{monomial rewriting step} and denote by $\Cont(\sR) = (\Cont(\sR)(i,j))_{i,j\in\obj(\sS)}$ the (family of) set(s) of monomial rewriting steps. We extend the source and target map to $\Cont(\sR)$ setting $s(\Gamma[r])=\Gamma[s(r)]$ and $t(\Gamma[r])=\Gamma[t(r)]$.

\medbreak

Often, it is necessary to restrict $\sS$:

\begin{definition}
  \label{defn:RW_subsystem}
  Let $\sS=(\sX_0,\sX_1,\sR,\sE)$ be a \LMRS{}.
  A \defnemph{linear sub-system $\lT$} is the data of a family of subsets $\lT(i,j)\subset\Cont(\sR)(i,j)$ for all $i,j\in\obj(\sS)$.
\end{definition}

We say that $\lT$ is \defnemph{context-agnostic} if
\[
f\in\lT\Leftrightarrow\Gamma[f]\in\lT\qquad\text{for all } f\in\lT\an\text{context }\Gamma.
\]
Otherwise, we say that $\lT$ is \defnemph{context-dependent}.
A context-dependent linear sub-system allows the encoding of ``global'' restrictions on rewriting rules, that is, restrictions that are not compatible with vertical and horizontal compositions.


\subsection{Rewriting steps}
\label{subsec:RW_rewriting_steps}

In the previous subsection, we defined contexts and went from local to global. To describe all relations induced by $\sR$, it remains to linearize.
 An \defnemph{$\sR$-rewriting step} is a rule of the form
\[\lambda\Gamma[r]+v \colon\mspace{10mu}\lambda\Gamma[s(r)] + v \mspace{10mu}\to\mspace{10mu} \lambda\Gamma[t(r)] + v,\]
where $\lambda\in \ring\setminus\{0\}$ is a non-zero scalar, $\Gamma$ is a context, $r\in\sR(i,j)$ is a generating rewriting step and $v\in\sX^l(i,j)$ is a vector.

We now take graded interchangers $\sE$ into account.
We write $\sim_\sE$ for the equivalence relation induced by graded interchangers---that is, we write $v\sim_\sE w$ for two vectors $v,w\in\sX^l$ if $v$ and $w$ are equal up to graded interchangers.
For two monomials $x,y\in\sX^*$, we write $x\projrel_\sE y$ if and only if there exists an invertible scalar $\lambda\in\ring^\times$ such that $x\sim_\sE\lambda y$.
For a vector $v\in\sX^l$, we write $\supp(v)$ for the \emph{support of $v$}---the set of monomials in its linear decomposition---and we write $\projsupp_\sE(v)$ for the set of monomials $x'\in\sX^*$ such that $x'\projrel_{\sE} x$ for some $x\in\supp(v)$.
When $\sE$ is clear from the context, we omit it and write $\sim$, $\projrel$ and $\projsupp$.

\begin{definition}
  \label{defn:RW_S_rewriting_step}
  Let $\sS$ be a \LMRS{} and $\lT\subset\sS$ a linear sub-system.
  A \defnemph{$\lT$-rewriting step} is a composition of the form
  \[u \mspace{10mu}\sim_\sE\mspace{10mu} \lambda\Gamma[s(r)] + v \mspace{10mu}\to_\lT\mspace{10mu} \lambda\Gamma[t(r)] + v \mspace{10mu}\sim_\sE\mspace{10mu} w\]
  where the middle arrow is an $\sR$-rewriting step such that $\Gamma[r]\subset\lT$.
  
  A $\lT$-rewriting step is \defnemph{positive} if
  \[\Gamma[s(r)]\notin\projsupp_\sE(v).\]
  We write $\lT^+$ the set of positive $\lT$-rewriting steps. A positive $\lT$-rewriting step is alternatively called a \defnemph{$\lT^+$-rewriting step}.
\end{definition}

In particular, an $\sS$-rewriting step consists of an arbitrary number of graded interchanges, followed by an $\sR$-rewriting step, followed by an arbitrary number of graded interchanges.
We say that this rewriting step is \emph{of type $r$}, in the notation of \cref{defn:RW_S_rewriting_step}.
We call an $\sR$-rewriting step \emph{positive} if it is positive as an $\sS$-rewriting step.
Note that positivity of an $\sR$-rewriting step depends on $\sE$.

\subsection{Termination and order}
\label{subsec:RW_termination_order}

Recall that a strict partial order is a transitive and asymmetric binary relation.
Let $\succ$ be a strict partial order on the set of monomials $\sX^*$.
We say that $\succ$ is \defnemph{$\sE$-invariant} if ($x'\projrel x$ \textsc{and} $x\succ y$ \textsc{and} $y\projrel y'$) implies ($x'\succ y'$). In other words, the strict partial order $\succ$ is $\sE$-invariant if it is independent of the interchange law.

\begin{definition}
  \label{defn:isp-order}
  An \defnemph{isp-order} is an $\sE$-invariant strict partial order on the set of monomials $\sX^*$.
\end{definition}

We say that the isp-order $\succ$ is \defnemph{context-agnostic} if
\[
x\succ y\Leftrightarrow\Gamma[x]\succ\Gamma[y]
\qquad
\text{for all }x,y\in\sX^*.
\]
Otherwise, we say that $\succ$ is \defnemph{context-dependent}.

\begin{definition}
  Let $\sS$ be a \LMRS{} and $\lT$ a linear sub-system.
  Let $\succ$ be an isp-order on the set of monomials $\sX^*$. We say that $\succ$ is \defnemph{compatible\footnote{In \cite{Schelstraete_RewritingModuloDiagrammatic_2026}, this notion is called ``strong compatibility'', while ``compatibility'' refers to a weaker notion.} with $\lT$} if
  \begin{equation}
    \label{eq:strong_compatibility}
  \Gamma[s(r)]\succ x,\quad\text{for all $r\in\sR$, context $\Gamma$ such that $\Gamma[r]\in\lT$, and $x\in\Gamma[\supp(t(r))]$}.
  \end{equation}
\end{definition}

For instance, assume $r\in\sR$ is of the form $r\colon x\to y-z$ with $y\notprojrel z$. The condition above states that, if $\Gamma$ is such that $\Gamma[r]\in\lT$, then $\Gamma[x]\succ\Gamma[y]$ and $\Gamma[x]\succ\Gamma[z]$.
Beware that one cannot replace ``$\Gamma[\supp(t(r))]$'' by ``$\supp(t(\Gamma[r]))$'' in the definition of compatibility.
Indeed, in the same example, the context $\Gamma$ could lead to an equality $\Gamma[y]=\Gamma[z]$, so that $\Gamma[r]\colon x\to 0$.
In this case, we have $\supp(t(\Gamma[r]))=\emptyset$, while $\Gamma[\supp(t(r))]=\{\Gamma[y],\Gamma[z]\}$.




\begin{definition}
  \label{defn:RW_termination}
  Let $\sS$ be a \LMRS{} and $\lT$ a linear sub-system.
  We say that $\lT^+$ \defnemph{terminates} if there is no infinite sequence of $\lT^+$-rewriting steps.
\end{definition}

Given an isp-order $\succ$ on monomials $\sX^*$, one can construct a strict partial order $\succ^+$ on vectors $\sX^l$, such that ($v'\sim v$ \textsc{and} $v\succ^+ w$ \textsc{and} $w\sim w'$) implies ($v'\succ^+ w'$); see \cite[section~3.3.4]{Schelstraete_RewritingModuloDiagrammatic_2026}.
If $x$ is a monomial and $v$ is a vector, it is such that
\[x\succ^+ v\qquad\Leftrightarrow\qquad x \succ y,\quad\forall y\in\supp(v).\]
Moreover:

\begin{lemma}
  \label{lem:RW_induced_order}
  Let $\sS$ be a \LMRS{} and $\lT$ a linear sub-system.
  Let $\succ$ be an isp-order compatible with $\lT$. Then for every $r\in\lT^+$, we have $s(r)\succ^+ t(r)$.
  In particular, if $\succ$ is well-founded, then $\lT^+$ terminates.\hfill\qed
\end{lemma}

\subsection{Tamed congruence}
\label{subsec:RW_tame_congruence}

\begin{mdframed}
  From now on we fix a linear monoidal rewriting system $\sS$, a linear sub-system $\lT$ and an isp-order $\succ$ on the set of monomials $\sX^*$.
  Moreover, we use $\aA$ or $\aB$ as placeholders for $\sS$, $\sS^+$, $\lT$ or $\lT^+$, writing e.g.\ ``$\aA$-rewriting step''.
  We omit the prefix when it is clear from the context, writing e.g.\ ``rewriting step''. In particular, we say ``$\aA$ is confluent'' to mean ``$\aA$ is $\aA$-confluent''.
\end{mdframed}

\noindent We call \defnemph{$\sE$-congruence} a zigzag of graded interchangers, and denote it $v\sim_\sE w$; this recovers the notation for the equivalence relation induced by graded interchangers introduced in \cref{subsec:RW_rewriting_steps}.
An \defnemph{$\aA$-rewriting sequence} is a composition of $\aA$-rewriting steps; by definition, a rewriting sequence consisting of an empty composition of $\aA$-rewriting steps is an $\sE$-congruence.
We write the composition of $f\colon s\to t$ with $g\colon t\to r$ as $f\starop g\colon s\to r$.
An \defnemph{$\aA$-congruence} is a zigzag of $\aA$-rewriting steps; again, an empty composition is an $\sE$-congruence.
(Note that this is the same definition as an $\sE$-congruence, replacing $\sE$ by $\aA$.)
We use the following notations for these notions (with extra subscript $\aA$ if necessary):
\begin{IEEEeqnarray*}{CcCcCcC}
  v\dashrightarrow_\aA w && v\to_\aA w && v\overset{*}{\to}_\aA w && v\sim_\aA w\\*
  \text{\footnotesize $\aA$-rewriting step}
  &\mspace{30mu}&
  \text{\footnotesize positive $\aA$-rewriting step}
  &\mspace{30mu}&
  \text{\footnotesize positive $\aA$-rewriting sequence}
  &\mspace{30mu}&
  \text{\footnotesize $\aA$-congruence}
\end{IEEEeqnarray*}

An \defnemph{$\aA$-branching} is a pair of $\aA$-rewriting sequences $(f,g)$ with the same source $s(f)=s(g)$.
An \defnemph{$\aA$-confluence} is a pair of $\aA$-rewriting sequences $(f',g')$ with the same target $t(f')=t(g')$.
An $\aA$-branching $(f,g)$ is \defnemph{$\aB$-confluent} if it admits a $\aB$-confluence $(f',g')$ such that $t(f)=s(f')$ and $t(g)=s(g')$:
\begin{IEEEeqnarray*}{CcCcC}
  \begin{tikzcd}[ampersand replacement=\&,cramped,row sep=.7em]
    \& v \\
    u \& \\
    \& {v'}
    \arrow["*"{description},"f", curve={height=-9pt}, from=2-1, to=1-2]
    \arrow["*"{description},"{g}"', curve={height=9pt}, from=2-1, to=3-2]
  \end{tikzcd}
  &&
  \begin{tikzcd}[ampersand replacement=\&,cramped,row sep=.7em]
    \& v \\
    \&\& w \\
    \& {v'}
    \arrow["*"{description},"f'", curve={height=-9pt}, from=1-2, to=2-3]
    \arrow["*"{description},"{g'}"'{shift={(-.05,.05)}}, curve={height=9pt}, from=3-2, to=2-3]
  \end{tikzcd}
  &&
  \begin{tikzcd}[ampersand replacement=\&,cramped,row sep=.7em]
    \& v \\
    u \&\& w \\
    \& {v'}
    \arrow["*"{description},"f'", curve={height=-9pt}, from=1-2, to=2-3]
    \arrow["*"{description},"f", curve={height=-9pt}, from=2-1, to=1-2]
    \arrow["*"{description},"{g}"', curve={height=9pt}, from=2-1, to=3-2]
    \arrow["*"{description},"{g'}"'{shift={(-.05,.05)}}, curve={height=9pt}, from=3-2, to=2-3]
  \end{tikzcd}
  \\*
  \text{\footnotesize a branching}
  &\mspace{70mu}&
  \text{\footnotesize a confluence}
  &\mspace{70mu}&
  \text{\footnotesize a confluent branching}
\end{IEEEeqnarray*}
A \defnemph{local $\aA$-branching} is an $\aA$-branching whose branches are rewriting \emph{steps} (rather than rewriting \emph{sequences}).
We say that $\aA$ is \defnemph{(locally) $\aB$-confluent} if every (local) $\aA$-branching is $\aB$-confluent.
An \defnemph{$\aA$-local triple} is a triple $[f,e,g]$ with $f$, $g$ are $\sR$-rewriting steps that belong to $\aA$, and $e$ is an $\sE$\nbd-congruence, such that $s(f)=s(e)$ and $t(e)=s(g)$:
\begin{IEEEeqnarray*}{CcC}
  \begin{tikzcd}[ampersand replacement=\&,cramped,row sep=.7em]
    \& y \\
    x \& \\
    \& {y'}
    \arrow["f", curve={height=-9pt}, from=2-1, to=1-2]
    \arrow["{g}"', curve={height=9pt}, from=2-1, to=3-2]
  \end{tikzcd}
  &&
  \begin{tikzcd}[ampersand replacement=\&]
    x \dar[snakecd,"e"']\arrow[r,"\sR"{subscript},"f"] \&y
    \\
    x' \rar["\sR"{subscript},"g"] \&y'
  \end{tikzcd}
  \\*
  \text{\scriptsize a local branching}
  &\mspace{100mu}&
  \text{\scriptsize a local triple}
\end{IEEEeqnarray*}
An $\aA$-local triple defines a local $\aA$-branching; in practice, the two notions are essentially identical.
For instance, one checks that $\aA$ is $\aB$-confluent if and only if every $\aA$-local triple is $\aB$-confluent. In general, local branchings are more suited for formal statements and local triples are more suited for computations.

We now introduce a weaker notion to replace confluence.
Recall the induced binary relation $\succ^+$ on vectors $\sX^l$ from \cref{subsec:RW_termination_order}.
If $W\subset \sX^l$ is a set of vectors, we write $v\succ^+ W$ if $v\succ^+ w$ for all $w\in W$.
If $f=f_n\circ\ldots\circ f_1$ is a sequence of composable arrows on vectors $\sX^l$, we write $v\succ^+ f$ to mean $v\succ^+\{s(f_1),t(f_1),\ldots,t(f_n)\}$.

\begin{definition}
  \label{defn:ARSM_tame_congruence}
  An $\aA$\nbd-branching $(f,g)$ of source $\bullet$ is said to be \defnemph{$\succ$\nbd-tamely $\aB$-congruent} (resp.\ \defnemph{$\succ$\nbd-tamely $\aB$\nbd-confluent}) if there exists a $\aB$\nbd-congruence $h$ (resp.\ $\aB$\nbd-con\-fluence $(f',g')$) such that $\bullet\succ^+ h$ (resp.\ $\bullet\succ^+ f'^{-1}\circ g'$).
\end{definition}

In particular, $\succ$\nbd-tameness implies $\bullet\succ^+ t(f)$ and $\bullet\succ^+ t(g)$.
Here is a schematic for a $\succ$\nbd-tamed congruence, where horizontal positions are used to suggest relative orderings with respect to $\succ^+$:
\begin{gather*}
  \begin{tikzpicture}[scale=.8]
    \node (C) at (0,0) {$\bullet$};
    \coordinate (A1) at (2,1);
    \coordinate (A2) at (2+.2,1-.4);
    \coordinate (A3) at (2-.6,1-2*.4);
    \coordinate (A4) at (2+1.2,1-3*.4);
    \coordinate (A5) at (2-.2,1-4*.4);
    \coordinate (A6) at (3,-1);
    \draw (A1) to (A2);
    \draw (A2) to (A3);
    \draw[
      dotted,
      dash pattern=on 3pt off 3pt, 
    ] (A3) to (A4);
    \draw (A4) to (A5);
    \draw (A5) to (A6);
    \draw[->] (C) to[out=60,in=180] node[above]{$f$} (A1);
    \draw[->] (C) to[out=-60,in=180] node[below]{$g$} (A6);
    \draw[dotted,thick] (.5,1.5) to (.5,-2.5);
    \node at (1.1,-2.2) {$\bullet\succ^+$};
  \end{tikzpicture}
\end{gather*}

Notions reminiscent of tamed congruence appeared already in the work of Buchberger \cite{Buchberger_CriterionDetectingUnnecessary_1979} (see also Winkler and Buchberger \cite{WB_CriterionEliminatingUnnecessary_1986}) and Bergman \cite{Bergman_DiamondLemmaRing_1978}.
In particular, the ``triangle lemma'' (see \eg \cite[Proposition~2.4.3.2]{BD_AlgebraicOperads_2016}) is reintepreted as a transitivity property:

\begin{lemma}[Tamed congruence is transitive]
  \label{lem:RW_tame_is_transitive}
  Let $f$, $g$ and $h$ be three $\aA$-rewriting sequences with the same source $s(f)=s(g)=s(h)$.
  If both branchings $(f,g)$ and $(g,h)$ admit a $\succ$-tamed $\aB$-congruence, then the $\aA$-branching $(f,h)$ admits a $\succ$-tamed $\aB$\nbd-congruence.\hfill\qed
\end{lemma}

As we shall see \cref{subsec:RW_confluence_analysis}, tamed congruence becomes a central tool in higher linear rewriting theory.
For now, we note that:

\begin{lemma}[Positive confluence implies tamed congruence]
  \label{lem:RW_tame_congruence_implies_confluence}
  If $\succ$ is compatible with $\lT$, then any $\lT^+$-confluence is a $\succ$-tamed $\lT$-congruence.\hfill\qed
\end{lemma}


\subsection{Monoidal Gröbner systems}
\label{subsec:RW_monoidal_grobner_system}

  
\begin{definition}
  \label{defn:RW_monoidal_grobner_system}
  Let $(\sS,\lT,\succ)$ be a triple where $\sS=(\sX_0,\sX_1;\sR,\sE)$ is a linear monoidal rewriting system, $\lT$ is a linear sub-system $\lT\subset\sS$, and $\succ$ is an isp-order.
  If:
  \begin{enumerate}[(i),itemsep=.2em]
    \item $\succ$ is well-founded and compatible with $\lT$;
    \item $\lT^+$ is confluent;
    \item for every $f\in\Cont(\sS)$, there exists a $\lT$-congruence $s(f)\overset{*}{\leftrightarrow}_\lT t(f)$;
  \end{enumerate}
  we say that $(\sS,\lT,\succ)$ is a \defnemph{monoidal Gröbner system}.
\end{definition}

\begin{remark}
  We comment on the terminology.
  In classical Gröbner theory, a Gröbner basis is a basis for the ideal of relations together with a monomial order, i.e.\ a \emph{total} and context-agnostic\footnote{If $>$ is a monomial order and $P$, $Q$, $R$ are polynomials, then $P>Q$ implies $PR>QR$.} order.
  The monomial order induces a choice of orientation on the Gröbner elements, and hence a rewriting system.
  In contrast, in rewriting theory the orientation on relations is given, and no order is assumed.
  The definition of a monoidal Gröbner system somehow fits in between the (higher analogue of the) two notions: we have a \emph{partial} order $\succ$ that guarantees termination, but which need not be context-agnostic.
  Moreover, in contrast to both the classical Gröbner approach and the classical rewriting approach, the system $\lT$ itself need not be context-agnostic.
\end{remark}

A vector $v\in\sX^l$ is said to be a \defnemph{$\lT$-normal form} if no $\lT^+$-rewriting step has source $v$; in other words, the vector $v$ cannot be rewritten by a positive $\lT$-rewriting step.
Given a pair of objects $i,j\in\obj(\sS)$, denote $\NF_\lT(i,j)$ the $\ring$-module of $\lT$-normal forms from $i$ to $j$ and $\NF_\lT(i,j) /\sE(i,j)$ the $\ring$-module obtained by quotienting the $\ring$-module of $\lT$-normal forms by the graded interchange law.

The \defnemph{graded-monoidal category presented by $\sS$} is the graded-monoidal category obtained by quotienting each $\sX^l(i,j)$ by the relations induced by $\sS(i,j)$, for each pair of objects $i,j\in\obj(\sS)$.
We denote it $[\sS]$.
Associating a normal form to its representation class in $[\sS]$ defines a canonical morphism of $\ring$-modules $\NF_\lT(i,j) /\sE(i,j)\to [\sS](i,j)$.

\begin{theorem}
  \label{thm:RW_monoidal_grobner_system}
  Let $(\sS,\lT,\succ)$ be a monoidal Gröbner system.
  For every pair of objects $i,j\in\obj(\sS)$, the canonical morphism of $\ring$-modules
  \[
  \NF_\lT(i,j) /\sE(i,j)\to [\sS](i,j)
  \]
  is an isomorphism.
  In particular, if $B = (B(i,j))_{i,j\in\obj(\sS)}\subset \sX^*$ is a family of sets of monomials such that each $B(i,j)$ is a basis of the $\ring$-module $\NF_\lT(i,j) /\sE(i,j)$, then $B$ is a hom-basis of $[\sS]$.\hfill\qed
\end{theorem}

The three conditions in \cref{defn:RW_monoidal_grobner_system} play three different roles in the above theorem.
Condition (iii) ensures that $\lT$ and $\sS$ present the same graded-monoidal category---that is, we can replace $[\sS]$ by $[\lT]$ in the statement of the theorem.
Condition (i) ensures that $\lT^+$ terminates, which in turn implies that the canonical map is surjective---that is, every vector admits a $\lT$-normal form.
Condition (ii), confluence of $\lT^+$, implies that the canonical maps are injective---that is, this $\lT$-normal form is unique up to graded-interchange.

\subsubsection{Coherence of interchangers}
\label{subsubsec:RW_coherence_interchangers}

Given \cref{thm:RW_monoidal_grobner_system}, we wish to understand equivalence classes under the graded interchange law.

We say that a monomial $y\in\sX^*$ \emph{contains} another monomial $x$ if $y\projrel_\sE\Gamma[x]$ for some context $\Gamma$, that is, ``$y$ contains $x$'' if $x$ is a sub-diagram of $y$, up to interchange.
A monomial $x\in\sX^*(\mathbb{1},\mathbb{1})$ for $\mathbb{1}$ the monoidal unit object is called a \defnemph{generalized bubble}.
In other words, a generalized bubble is a diagram without input and output.
As a consequence of the coherence of interchangers in Gray categories (see \eg \cite[Appendix~A]{Schelstraete_OddKhovanovHomology_2024}; recall that $\bil$ is assumed to be skew-symmetric), failure of coherence of graded interchangers only comes from interchanging two identical generalized bubbles:
\[x\circ x \sim_\sE \bil(\deg x,\deg x) \;x\circ x
\qquad\overset{1-\bil(\deg x,\deg x)\text{ invertible}}{\Longrightarrow}\qquad
x\circ x\sim_\sE 0.\]
We say that a diagram is \defnemph{odd-bubble-square-free} if the situation above does not occur, that is, it does not contain the diagram $x\circ x$, where $x$ is a generalized bubble with the property that $\bil(\deg x,\deg x)\neq 1$.
(The terminology comes from the super case, where that situation happens exactly when $x$ is odd.)

\begin{lemma}
  \label{lem:RW_bubble_trivial_grading_implies_coherence}
  Assume that for every generalized bubble, the scalar $1-\bil(\deg x,\deg x)$ is either zero or invertible.
  Let $\mathsf{NF}$ be a set of diagrams with the same input and output.
  Any choice of interchange representative for odd-bubble-square-free diagrams in $\mathsf{NF}$ defines a basis for the module $\langle\mathsf{NF}\rangle_\ring/\sE$, the $\ring$-module generated by $\mathsf{NF}$ and quotiented by graded interchangers.
\end{lemma}

\subsection{Confluence analysis}
\label{subsec:RW_confluence_analysis}

To use \cref{thm:RW_monoidal_grobner_system}, we wish to show that $\lT^+$ is confluent. This subsection describes general techniques to carry out this analysis, common to every linear monoidal rewriting system.
\Cref{subsubsec:branchwise_congruence} explains how to think of branchings modulo interchange.
Then, \cref{lem:RW_linear_tame_newman} reduces the confluence analysis to ``monomial branchings'', \cref{lem:RW_independent_branching_lemma} to ``overlapping branchings'' and \cref{subsubsec:RW_contextualization} gives tools to reduce to ``minimal overlapping branchings''.
\Cref{subsubsec:RW_enumeration_overlapping_branchings} explains how to enumerate minimal overlapping branchings.
Finally, \cref{subsubsec:RW_critical_branchings} gives a summary, with a view toward the confluence analysis of \cref{sec:critical_branchings}.

\subsubsection{Branchwise congruence}
\label{subsubsec:branchwise_congruence}
We discuss how to deal with the modulo data $\sE$.
Just like we think of vectors modulo graded interchangers, we would like to think of rewriting steps and branchings modulo.
Two rewriting steps are said to be \defnemph{$\sE$-congruent} if there exist $\sE$-congruences between their source and target.
When working modulo graded interchangers, we have canonical $\sE$-congruence of rewriting steps, such as the one pictured in \cref{fig:RW_interchange_E_congruence}.
Here $\bil = \bil(s(\ssf),x) = \bil(y,x)$ for each $y\in\supp(t(\ssf))$; the equality follows from the fact that $\ssf$ preserves $G$-degrees.
This ensures the existence of this canonical $\sE$-congruence.
This is encapsulated in the following lemma:

\begin{lemma}[Canonical $\sE$-congruences]
  \label{lem:RW_interchange_E_congruences}
  Let $\ssf\in\lT$ be a monomial rewriting step and $x\in\sX^*$ a monomial.
  The rewriting steps
  \[(\ssf\otimes \id_{t(x)}) \circ (\id_{s_0(\ssf)}\otimes x)
  \quad\an\quad
  (\id_{t_0(\ssf)}\otimes x) \circ (\ssf\otimes \id_{s(x)})\]
  are $\sE$-congruent, where we write $s_0(\ssf) = s(s(\ssf))$, $t_0(\ssf) = t(s(\ssf))$ and use \cref{not:RW_identities}. We have a similar statement when $x$ is on the left of $\ssf$.\hfill\qed
\end{lemma}

Two branchings $(f,g)$ and $(f',g')$ are said to be \defnemph{branchwise $\sE$-congruent} if $f$ (resp.\ $g$) is $\sE$-congruent to $f'$ (resp.\ $g'$).
We think of branchings modulo the canonical $\sE$-congruences of \cref{lem:RW_interchange_E_congruences}, which is justified by the following lemma (with $\aA$ and $\aB$ as placeholders as in \cref{subsec:RW_tame_congruence}):

\begin{lemma}[Branchwise congruence lemma]
  \label{lem:RW_branchwise_congruence_lemma}
  If $(f,g)$ and $(f',g')$ are branchwise $\sE$-congruent $\aA$-branchings, then $(f,g)$ is $\aB$-confluent (or $\succ$-tamely $\aB$-congruent) if and only if $(f',g')$ is.\hfill\qed
\end{lemma}


\begin{figure}
\centering
\begin{minipage}[b]{.5\textwidth}
  \centering  
  \begin{tikzcd}[row sep=.6cm,ampersand replacement=\&]
    \mspace{16mu}
    \tikzpic{
      \indexNF[0][1][{s(\ssf)}]\lili[2][1]
      \lili\indexNF[2][0][{x}]
    }[scale=1.2][1]
    \ar[r]\ar[d,snakecd]
    \&
    \mspace{16mu}
    \tikzpic{
      \indexNF[0][1][{t(\ssf)}]\lili[2][1]
      \lili\indexNF[2][0][{x}]
    }[scale=1.2][1]
    \ar[d,snakecd]
    \\
    \bil\;\;
    \tikzpic{
      \lili[0][1]\indexNF[2][1][{x}]
      \indexNF[0][0][{s(\ssf)}]\lili[2][0]
    }[scale=1.2][1]
    \ar[r]
    \&
    \bil\;\;
    \tikzpic{
      \lili[0][1]\indexNF[2][1][{x}]
      \indexNF[0][0][{t(\ssf)}]\lili[2][0]
    }[scale=1.2][1]
  \end{tikzcd}
  \captionof{figure}{Interchange $\sE$-congruence}
  \label{fig:RW_interchange_E_congruence}
\end{minipage}%
\begin{minipage}[b]{.5\textwidth}
  \centering
  \begin{tikzcd}[ampersand replacement=\&,row sep=-1em,column sep=3em]
    \&
    \tikzpic{
      \indexNF[0][1][{t(\ssf)}]
      \indexNF[0][0][{s(\ssg)}]
    }[scale=1.2]
    \ar[dr,bend left=20pt,"g'"{above right=2pt},"*"{description},dottedcd]
    \\
    \tikzpic{
      \indexNF[0][1][{s(\ssf)}]
      \indexNF[0][0][{s(\ssg)}]
    }[scale=1.2]
    \ar[ur,bend left=20pt,"f=\ssf\circ s(\ssg)"]
    \ar[dr,bend right=20pt,"g=s(\ssf)\circ \ssg"']
    \&\&
    \tikzpic{
      \indexNF[0][1][{t(\ssf)}]
      \indexNF[0][0][{t(\ssg)}]
    }[scale=1.2]
    \\
    \&
    \tikzpic{
      \indexNF[0][1][{s(\ssf)}]
      \indexNF[0][0][{t(\ssg)}]
    }[scale=1.2]
    \ar[ur,bend right=20pt,"f'"'{below right=2pt},"*"{description},dottedcd]
  \end{tikzcd}
  \captionof{figure}{Independent branching}
  \label{fig:RW_independent_branching}
\end{minipage}
\end{figure}

\subsubsection{Linear tamed Newman's lemma}
A local branching $(f,g)$ is \defnemph{monomial} if both branches $f$ and $g$ are monomial rewriting steps, i.e.\ their source is a monomial. Note that a monomial rewriting step is necessarily positive, so that a monomial branching is necessarily positive (meaning that both branches are positive).

\begin{lemma}[Linear tamed Newman's lemma]
  \label{lem:RW_linear_tame_newman}
  Assume that $\succ$ is compatible with $\lT$.
  If $\succ$ is well-founded and if every monomial $\lT$-branching is $\succ$-tamely $\lT$-congruent, then $\lT^+$ is confluent.\hfill\qed
\end{lemma}

The lemma reduces the study of confluence to the tamed congruence of monomial branchings.
We stress that in the above lemma, the hypothesis that $\succ$ is well-founded is necessary, even if one replaces ``$\succ$-tamely congruent'' by ``confluent'' in the statement.

\subsubsection{Independent branchings}
\label{subsubsec:RW_independent_branchings}

Consider two monomial $\sS$-rewriting steps
\[
  \ssf\colon s(\ssf)\to \lambda_1x_1+\ldots+\lambda_mx_m
  \quad\an\quad
  \ssg\colon s(\ssg)\to \mu_1y_1+\ldots+\mu_ny_n.
\]
Assuming they are suitably composable, let $f=\ssf\circ s(\ssg)$ and $g=s(\ssf)\circ\ssg$ and consider the branching $(f,g)$. Such branchings are called \defnemph{independent branchings}---if a monomial branching is not $\sE$-congruent (see \cref{subsubsec:branchwise_congruence}) to an independent branching, we say that it is \defnemph{overlapping}.
Each independent branching $(f,g)$ admits a canonical confluence $(g',f')$, where $g'$ and $f'$ are constructed from the monomial rewriting steps $x_i\circ\ssg$ and $\ssf\circ y_j$, respectively.
(More precisely, the confluence $(g',f')$ is canonical up to a choice of ordering on the $x_i$'s and the $y_j$'s.)
Diagrammatically, the branching $(f,g)$ and its canonical confluence $(g',f')$ have the form pictured in \cref{fig:RW_independent_branching}.
In general, the confluence $(g',f')$ need not be positive---however, we can expect it to be a tamed congruence.
In the context-dependent case $\lT\subset\sS$, we must further check that all the monomial rewriting steps involved belong to $\lT$.
This is the content of the following lemma:

\begin{lemma}[Independent branching lemma]
  \label{lem:RW_independent_branching_lemma}
  Assume that $\succ$ is an isp-order compatible with $\lT$. Let $(f,g)$ be an independent $\lT$\nbd-branching.
  If for each $1\leq j\leq n$ (resp.\ $1\leq i\leq m$) and with the notations above, the $\sS^+$\nbd-rewriting step $\ssf \circ y_j$ (resp.\ $x_i \circ \ssg$) is in $\lT$, then the canonical $\sS$\nbd-confluence $(g',f')$ is a $\lT$-congruence $\succ$\nbd-tamed by $s(f)=s(g)$.
\end{lemma}

\begin{proof}
  It is clear from the hypotheses that the canonical $\sS$\nbd-confluence is a $\lT$-congruence.
  Moreover, since $(f,g)$ is a $\lT$-branching, we have $s(\ssf)\circ s(\ssg)\succ x_i\circ s(\ssg)$ and $s(\ssf)\circ s(\ssg)\succ s(\ssf)\circ y_j$ for each $1\leq i\leq m$ and $1\leq j\leq n$.
  Finally, since $\ssf \circ y_j$ is in $\lT$, it follows that $s(\ssf)\circ y_j\succ x_i\circ y_j$, so that by transitivity of $\succ$, we have  $s(\ssf)\circ s(\ssg)\succ x_i\circ y_j$.
  This concludes.
\end{proof}

\subsubsection{Contextualization}
\label{subsubsec:RW_contextualization}

A branching of the form $\Gamma[f,g]$ for a non-trivial context $\Gamma$ is called a \defnemph{contextualized branching}---if a branching is not $\sE$-congruent (see \cref{subsubsec:branchwise_congruence}) to a contextualized branching, we say that it is \defnemph{minimal}.

In general, one wishes to deduce the tamed congruence of contextualized branchings from the tamed congruence of minimal branchings.
This works in the context-agnostic setting:

\begin{lemma}[Contextualization lemma; context-agnostic order]
  \label{lem:RW_contextualization_lemma_agnostic}
  Assume that $\strongsucc$ is context-agnostic and that $\strongsucc$ is compatible with $\lT$.
  Let $(f,g)$ be a monomial $\sS$-branching which admits a $\strongsucc$-tamed $\lT$-congruence $(f',g')$.
  If $\Gamma$ is a context such that both $\Gamma[f,g]$ and $\Gamma[f',g']$ belong to $\lT$, then $\Gamma[f',g']$ is a $\strongsucc$-tamed $\lT$-congruence for $\Gamma[f,g]$.\hfill\qed
\end{lemma}

However, the analogue of that lemma is \emph{not} true when the order is not context-agnostic.
Instead, one relies on the following lemma:

\begin{lemma}\label{lem:RW_contextualization-sequence}
  Assume $\succ$ is compatible with $\lT$.
  Consider a composition $f\starop g$, where $f$ is a monomial $\sS$\nbd-rewriting step and $g$ is a positive $\sS$\nbd-rewriting sequence.
  If $\Gamma$ is a context such that both $\Gamma[f]$ and $\Gamma[g]$ are in $\lT$,
  then $\Gamma[g]$ is $\succ$\nbd-tamed by $\Gamma[s(f)]$.
\end{lemma}

In \cref{subsec:tamed_preconfluence}, we introduce the novel notion of ``tamed preconfluence'' that interpolates between confluence and tamed congruence.
In \cref{sec:critical_branchings}, we use \cref{lem:RW_contextualization-sequence} to contextualized tamed preconfluence.



\subsubsection{Enumeration of overlapping branchings}
\label{subsubsec:RW_enumeration_overlapping_branchings}
Given a diagram $D$, its \emph{multiset of generators $\mathcal{G}en(D)$} is the multiset of generating diagrams that constitute $D$ under whiskering and composition.
Interchange preserves this multiset; more precisely, it induces a canonical bijection between the multisets of generators.
Let $[f,e,g]$ be a monomial local triple, with $f=\Gamma_f[\ssf]$ and $g=\Gamma_g[\ssg]$ for some contexts $\Gamma_f$ and $\Gamma_g$ and some generating rewriting steps $\ssf$ and $\ssg$.
If $D$ is the source of $[f,e,g]$, the multiset $\mathcal{G}en(s(\ssf))$ is canonically a sub-multiset of $\mathcal{G}en(D)$; and similarly for $\ssg$.
We write $\mathcal{G}en_D(s(\ssf))=\mathcal{G}en(s(\ssf))$ to emphasize that point.

\begin{lemma}
  \label{lem:RW_overlapping_branchings_intersect}
  Assume that $s(\ssf)$ and $s(\ssg)$ are connected diagrams with at least one generator.
  Then $[f,e,g]$ is an overlapping branching if and only if $\mathcal{G}en_D(s(\ssf))\cap\mathcal{G}en_D(s(\ssg))\neq\emptyset$ as sub-multisets of $\mathcal{G}en(D)$.
\end{lemma}

In other words, under certain hypotheses, a branching is overlapping precisely when there is an ``overlap'' between the sources of the two branches; that is, when they ``share'' a generator.

\begin{proof}
  The direction $\Leftarrow$ is clear.

  Assume then that $\mathcal{G}en_D(s(\ssf))\cap\mathcal{G}en_D(s(\ssg))=\emptyset$ as sub-multisets of $\mathcal{G}en(D)$.
  While the elements of $\mathcal{G}en_D(s(\ssf))$ sit as $s(\ssf)$ in the diagram $D$, the elements of $\mathcal{G}en_D(s(\ssg))$ are a priori scattered across $D$.
  We wish to show that there exists a diagram $D'\sim D$ such that both $\mathcal{G}en_D(s(\ssf))$ and $\mathcal{G}en_D(s(\ssg))$ sit as $s(\ssf)$ and $s(\ssg)$ in the diagram $D'$, respectively.

There exist two elements of $\mathcal{G}en_D(s(\ssg))$ that can be interchanged to lie one above the other and, since $s(\ssg)$ is connected, that share some input and output strands. Viewed as sitting inside $D$, these two generators must be in one of the following four cases:
  \begin{gather*}
    \includegraphics[width=3cm]{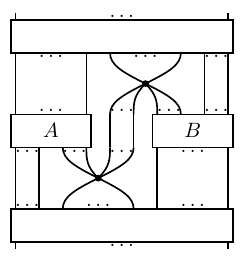}
    \qquad
    \includegraphics[width=3cm]{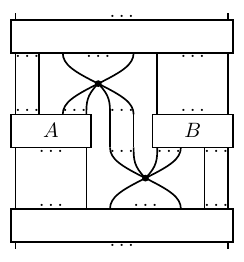}
    \qquad
    \includegraphics[width=3cm]{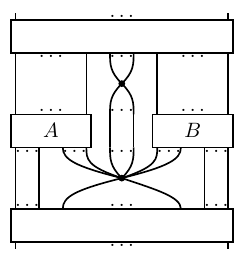}
    \qquad
    \includegraphics[width=3cm]{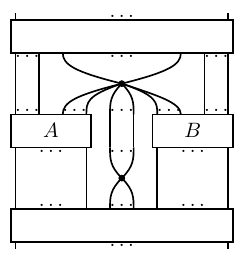}
  \end{gather*}
  To move the two generators close to one another, we wish to interchange the diagram $A$ upward (or downward) and the diagram $B$ downward (or upward).
  We can do so without affecting the fact that the elements of $\mathcal{G}en_D(s(\ssf))$ sit as $s(\ssf)$. Indeed, as $s(\ssf)$ is a connected diagram with at least one generator, we have that if $\mathcal{G}en_D(s(\ssf))\cap\big(\mathcal{G}en_D(A)\cup\mathcal{G}en_D(B)\big)\neq\emptyset$, then either $\mathcal{G}en_D(s(\ssf))\subset\mathcal{G}en_D(A)$ or $\mathcal{G}en_D(s(\ssf))\subset\mathcal{G}en_D(B)$.
  Once the two generators are close to one another, we can view them as one generator (or ``contract'' them into one generator), and conclude by induction on the number of generators of $s(\ssg)$.
\end{proof}

\begin{remark}
  \newcommand{\tempdiag}[2]{
    \draw ({#1-.1},{#2+.15}) rectangle ({#1+1.1},{#2+.85});
    \draw ({#1},{#2}) to ({#1},{#2+.15});
    \draw ({#1+1},{#2}) to ({#1+1},{#2+.15});
    \draw ({#1},{#2+.85}) to ({#1},{#2+1});
    \draw ({#1+1},{#2+.85}) to ({#1+1},{#2+1});
  }
  \newcommand{\tempdiagbis}[2]{
    \draw[pattern=north east lines] ({#1-.1},{#2+.15}) rectangle ({#1+1.1},{#2+.85});
    \draw ({#1},{#2}) to ({#1},{#2+.15});
    \draw ({#1+1},{#2}) to ({#1+1},{#2+.15});
    \draw ({#1},{#2+.85}) to ({#1},{#2+1});
    \draw ({#1+1},{#2+.85}) to ({#1+1},{#2+1});
  }
  In the above lemma, it is necessary to assume that at least one of $s(\ssf)$ and $s(\ssg)$ is a connected diagram.
  Indeed, if $s(\ssf)$ and $s(\ssg)$ are respectively
  \begin{gather*}
  \tikzpic{
    \tempdiag{0}{1}\li[2][1]
    \li\tempdiag{1}{0}
  }[scale=.6]
  \qquad\an\qquad
  \tikzpic{
    \lili[0][1]\li[2][1]\tempdiagbis{3}{1}
    \tempdiagbis{0}{0}\lili[2][0]\li[4][0]
  }[scale=.6]\;,
  \end{gather*}
  then the following is an overlapping branching:
  \begin{gather*}
    \tikzpic{
      \li[0][2]\tempdiag{1}{2}\lili[3][2]
      \lili[0][1]\li[2][1]\tempdiagbis{3}{1}
      \tempdiagbis{0}{0}\lili[2][0]\li[4][0]
      \lili[0][-1]\tempdiag{2}{-1}\li[4][-1]
    }[scale=.6]
    \;\sim_{\sE}\;
    \tikzpic{
      \lili[0][2]\li[2][2]\tempdiagbis{3}{2}
      \li[0][1]\tempdiag{1}{1}\lili[3][1]
      \lili[0][0]\tempdiag{2}{0}\li[4][0]
      \tempdiagbis{0}{-1}\lili[2][-1]\li[4][-1]
    }[scale=.6]
    \;.
  \end{gather*}
  We expect that there exists a generalization of the above lemma to include non-connected diagrams and another combinatorial data than $\mathcal{G}en(-)$, although we do not pursue that direction further.
\end{remark}

\subsubsection{Critical branchings}
\label{subsubsec:RW_critical_branchings}

In practice, how do we go about showing confluence of $\lT^+$?
In the context-agnostic setting, we can give a summary as follows:

\begin{theorem}
  \label{thm:RW_confluence_analysis_agnostic}
  Let $\sS$ be a \LMRS{} and $\lT$ a context-agnostic linear sub-system.
  Let $\succ$ be an isp-order compatible with $\lT$, which is context-agnostic and well-founded.
  Then $\lT^+$ is confluent if and only if all minimal overlapping branchings\footnote{Or more precisely, for each choice of representative under branchwise $\sE$-congruence.} are positively $\lT$-confluent.
\end{theorem}

\begin{proof}
  Note that since $\succ$ is compatible with $\lT$, positive confluence implies $\succ$-tamed $\lT$-congruence (\cref{lem:RW_tame_congruence_implies_confluence}).
  By hypothesis and the context-agnostic contextualization lemma (\cref{lem:RW_contextualization_lemma_agnostic}), every overlapping branching is $\succ$-tamely $\lT$-congruent.
  Up to branchwise $\sE$-congruence, every monomial branching is either overlapping or independent:
  by the branchwise congruence lemma (\cref{lem:RW_branchwise_congruence_lemma}) and the independent branching lemma (\cref{lem:RW_independent_branching_lemma}), every monomial branching is $\succ$-tamely $\lT$-congruent.
  Finally, by the tamed linear Newman lemma (\cref{lem:RW_linear_tame_newman}), $\lT^+$ is confluent.
\end{proof}

In the context-dependent setting, additional arguments are usually needed: for instance, whether independent branchings are $\succ$-tamely congruent is not automatic.
One must use \cref{lem:RW_independent_branching_lemma} and \cref{lem:RW_contextualization-sequence}, or depending on the specifics of the \LMRS{}, variations thereof.

\medbreak

In practice, we look for the smallest\footnote{Here we mean ``smallest'' in an informal way, as this set is often infinite anyway.} set that can replace ``minimal overlapping branchings'' in \cref{thm:RW_confluence_analysis_agnostic}.
We call any set that makes \cref{thm:RW_confluence_analysis_agnostic} hold a \defnemph{set of critical branchings}.
In that regard, transitivity of tamed congruence (\cref{lem:RW_tame_is_transitive}) often turns out to be a powerful tool.
A typical situation is as follows.
Consider a monomial branching $(f,g)$ and assume that there exists a third rewriting step $h$ on the source of $(f,g)$, such that $h$ is independent of both $f$ and $g$. If independent branchings are tamely congruent, then $(f,h)$ and $(h,g)$ are tamely congruent, and by transitivity of tamed congruence, so is $(f,g)$.
This kind of argument is called \defnemph{independent rewriting}.


%% file: sections/definitions.tex
\section{Rewriting systems of affine Brauer type}
\label{sec:definitions}

\subsection{Terminology}

In this article, a \defnemph{diagram} is any string diagram obtained from whiskering (horizontal composition with identities) and vertical composition of the generators $\tikzpic{\ca}$, $\tikzpic{\cu}$, $\tikzpic{\cro}$ and $\tikzpic{\bli}$.
Topologically, a diagram with $m$ bottom points and $n$ top points consists of $(m+n)/2$ intervals immersed in the strip $\bR\times[0,1]$, called \defnemph{open strands}, and a number of immersed closed curves, called \defnemph{closed strands}; strands are further decorated with dots.
We regard diagrams up to rectilinear isotopies preserving the relative vertical positions---in particular, two generators always have distinct vertical positions.
When instead we wish to regard a diagram up to all rectilinear isotopies, we say ``up to interchange''.
When we wish to emphasize the presence or absence of dots, we say \emph{dotted} or \emph{undotted} diagram.

Fix a diagram $D$.
A \defnemph{generalized cup} (resp.\ \defnemph{generalized cap}) is a strand connecting two top (resp.\ bottom) points. A \defnemph{propagating line} is a strand connecting a top point to a bottom point.
A \defnemph{straight line} is a piece of strand that has no critical point and does not cross other strands.
A \defnemph{bubble} is a closed strand consisting only of a cup and a cap, possibly with some dots---in this case, we assume dots sit on the left.
Given a subset $\bubbleset\subset\bN_{>0}$, an \emph{$\bubbleset$-bubble} is a bubble whose number of dots is in $\bubbleset$.

As a set of immersed curves, the diagram $D$ divides the strip $\bR\times[0,1]$ into connected components that we call \defnemph{regions}. The region containing all $(x,y)$ for $x\ll0$ is called the \defnemph{leftmost region}.
Any connected component $c$ without boundary delimits a union of compact regions. We call the unique region adjacent to this union the \defnemph{outer region} of $c$.

\subsection{Traffic rules}

\subsubsection{Reduced matchings}


Let $m,n\in\bN$ such that $m+n$ is even.
A \emph{pairing} is a partition of $\{1,\ldots,n+m\}$ into size-two subsets.
Any diagram induces a pairing on the set of its endpoints.
Given a pairing, a \defnemph{choice of dead ends} is the choice $\trafficend$ of one element (the \emph{dead end}) in each pair, and a \emph{dot map} is a function $d\colon \trafficend\to\bN$.
If, for each open strand of a diagram, all dots sit close to the same endpoint, then these endpoints induce a choice $\trafficend$ of dead ends.

\begin{definition}
  \label{defn:reduced_matching}
  Fix a pairing and a choice $\trafficend$ of dead ends.
  An \defnemph{$(\trafficend,\emptyset)$-reduced matching} is a diagram where each strand has at most one critical point, there are no self-inter\-sec\-tions of strands, distinct strands intersect at most once, and for each strand, dots sit close to their dead ends.
  Two reduced matchings are \emph{equivalent} if they induce the same pairing, the same choice of dead ends and the same dot map.
\end{definition}

A \defnemph{choice of bubbles} is a subset $\bubbleset\subset\bN_{> 0}$.
A \emph{bubble multiset} is a multiset with elements in $\bubbleset$.
Any diagram whose bubbles are $\bubbleset$-bubbles induces a bubble multiset, consisting of the numbers of dots on its bubbles.

\begin{definition}
  \label{defn:end_X_reduced_matching}
  An \defnemph{$(\trafficend,\bubbleset)$-reduced matching} is a horizontal juxtaposition of $\bubbleset$-bubbles on the leftmost region and an $(\trafficend,\emptyset)$-reduced matching.
  Two $(\trafficend,\bubbleset)$-reduced matchings are \emph{equivalent} if their underlying $(\trafficend,\emptyset)$-reduced matchings are equivalent and if they induce the same bubble multiset.
\end{definition}

\begin{definition}
  \label{defn:category_affine_brauer_type}
  Let $\bubbleset$ be a subset of $\bN_{>0}$.
  Let $\affC$ be a graded-monoidal category whose set of objects is $\bN$ and whose morphisms are generated by the following homogeneous morphisms:
  \begin{IEEEeqnarray}{CcCcCcC}
    \label{eq:generators_affine_Brauer}
    \tikzpic{\cro}\;,\quad
    \tikzpic{\ca}\;,\quad
    \tikzpic{\cu}\quad\an\quad
    \tikzpic{\bli}\;.
  \end{IEEEeqnarray}
  We say that $\affC$ is \defnemph{of affine Brauer type with $\bubbleset$-bubbles} if for any $(m,n)$, any choice $\trafficend$ of dead ends and any choice $\affB(m,n)$ of representatives of $(\trafficend,\emptyset)$-reduced matchings, the $\ring$-module $\affC(m,n)$ is a free left module over the algebra of $\bubbleset$-bubbles, with basis given by $\affB(m,n)$.
\end{definition}

Here the algebra of bubbles is the quotient of the free associative algebra over the set $\bubbleset$ by relations $nm=\mu(\deg n,\deg m)mn$, where we identify $n\in\bubbleset$ with the bubble with $n$ dots; in particular, $deg(n)$ has a well-defined meaning.
In the context of monoidal categories, the algebra of bubbles is the symmetric algebra over $\bubbleset$.
In the context of supermonoidal categories, if $\bubbleset$ only consists of bubbles of odd parity, the algebra of bubbles is the exterior algebra over $\bubbleset$.





\subsubsection{Tunnels}

Let $D$ be a diagram.
The \emph{dual graph} of $D$ is the graph whose vertices are the regions delimited by $D$, and whose edges connect regions sharing a strand.
We endow this graph with an orientation as follows:
\begin{itemize}
  \item each edge joining the leftmost region is oriented toward the leftmost region;
  \item each remaining edge incident to an already-oriented vertex is oriented toward that vertex.
\end{itemize}
In other words, each edge is oriented toward its endpoint closer to the leftmost region. This is well-defined, as $D$ admits a checkerboard colouring\footnote{a black and white colouring of its regions, such that adjacent regions have distinct colours}, and hence the dual graph has only even cycles.
\defnemph{Tunnels} record this orientation, drawn as red arrows crossing strands:
\begin{gather*}
  \tikzpic{\lorli}
  \qquad\an\qquad
  \tikzpic{\rorli}\;.
\end{gather*}
Every vertex $v$ of the dual graph has a well-defined quantity $\disttunnel^D(v)$ given by the length (that is, the number of tunnels) of any oriented path from $v$ to the leftmost region.

Recall that a \emph{simple closed curve} is a closed strand without any crossing.

\begin{definition}
  \label{defn:distance_tunnel}
  Let $D$ be a diagram.
  For each vertex $v$ of the dual graph, we write $\nbscc^D(v;k)$ for the number of simple closed curves with $k$ dots whose outer region is $v$, and denote by $\nbscc^D(v)\colon\bN\to\bN$ the associated function.
  For each $t\in\bN$, we write
  \[\disttunnel^D(t) \coloneqq \sum_{\disttunnel^D(v)=t}\nbscc^D(v),\]
  and denote by $\disttunnel^D\colon\bN\to(\bN)^\bN$ the associated function.
\end{definition}

\subsubsection{Stop signs}

Orient the interchange in the following direction:
\begin{gather*}
  \tikzpic{
    \indexNF[0][1][{x}]\lili[2][1]
    \lili\indexNF[2][0][{y}]
  }[scale=1.2][1]
  \;\to_{\textsc{lint}}\;
  \tikzpic{
    \lili[0][1]\indexNF[2][1][{y}]
    \indexNF[0][0][{x}]\lili[2][0]
  }[scale=1.2][1]
\end{gather*}
In \cite{DV_NormalizationPlanarString_2022}, Delpeuch and Vicary showed that this process terminates and is confluent on connected diagrams $D$, leading to a unique normal form $N_\textsc{lint}(D)$ for equivalence classes under the interchange law.

\begin{definition}
  Let $D$ be a diagram and $S$ a closed strand in $D$.
  Write $\ov{S}$ for the connected component of $D$ that contains $S$.
  The \defnemph{\textsc{lint}-minimal cup} (of $S$ in $D$) is the cup in $S$ which is the bottom generator of $S$ in $N_\textsc{lint}(\ov{S})$.
\end{definition}

The \textsc{lint}-minimal cup depends both on $S$ and on how $S$ sits inside $D$.
We decorate the \textsc{lint}-minimal cups of a diagram (one for each closed strand) with a \defnemph{stop sign}:
\begin{equation*}\label{eq:stop_sign}
  \tikzpic{\custop[0][0][->]}\;.
\end{equation*}

Here are two examples:
\[
\tikzpic{
  \ca[3][6]\uli[5][6]
  \AntiDiag[2][5]\cro[4][5]
  \li[2][4]\AntiDiag[3][4]\Diag[5][4]
  \ca[0][3]\cro[2][3]\ca[4][3]\li[6][3]
  \li[0][2]\cro[1][2]\cu[3][3]\lili[5][2]
  \cro[0][1]\begin{scope}[xscale=3]
    \cu[2/3][2]
  \end{scope}\li[6][1]
  \dli[0][1]\begin{scope}[xscale=5]
    \custop[1/5][1]
  \end{scope}
}[scale=.6]
\qquad\an\qquad
\tikzpic{
  \begin{scope}[xscale=3]
    \ca[3/3][4]
  \end{scope}
  \ca[1][3]\li[3][3]\ca[4][3]\li[6][3]
  \lili[1][2]\custop[3][3]\lili[5][2]
  \li[1][1]\begin{scope}[xscale=3]
    \cu[2/3][2]
  \end{scope}\li[6][1]
  \begin{scope}[xscale=5]
    \cu[1/5][1]
  \end{scope}
}[scale=.8]\;.
\]
Assume a choice $\trafficend$ of dead ends is fixed.
We extend the notion of dead ends to closed strands, where the dead end refers to the small piece of strand directly on the right of the stop sign. This also induces an orientation on closed strands, going from the stop sign to the dead end. Open strands are oriented toward the end point with the dead end.

\begin{definition}
  \label{defn:distance_traffic}
  Fix a choice $\trafficend$ of dead ends.
  If $d$ is a dot in a diagram $D$, we write $\disttraffic^D(d)$ for the number of crossings, cups and caps encountered when travelling along the strand from the dot to its dead end.
  We write $\disttraffic^D$ for the sum of $\disttraffic^D(d)$ over all dots in $D$.
\end{definition}

\subsubsection{Traffic rules}

\begin{definition}
  A \defnemph{choice of traffic rules} is a pair $\trafficrule=(\trafficend,\bubbleset)$ where $\trafficend$ is a choice of dead ends for each pairing and $\bubbleset$ is a choice of bubbles.
\end{definition}

Informally, we think of traffic rules as all the extra decorations on diagrams (tunnels, stop signs, orientations, dead ends), either induced from the choice of traffic rules $\trafficrule$ or canonically defined.

\subsection{Orders}

\subsubsection{Review}

We review the order on compactly supported functions. Fix a set $\Lambda$ and a well-founded total order $>$ on $\Lambda$, and let $\bot$ denote its unique minimal element.
Recall that a function $f\colon\bN\to\Lambda$ is \defnemph{compactly supported} if it eventually ``vanishes'', that is, if there exists a natural number $N$ such that $f(x)=\bot$ for all $x>N$.
Denote by $(\Lambda^\bN)_{\mathrm{cpt}}\subset\Lambda^\bN$ the set of compactly supported functions from $\bN$ to $\Lambda$.
There is a canonical well-founded total order on $(\Lambda^\bN)_{\mathrm{cpt}}$ induced by $>$, defined as $f>g$ if and only if $f\neq g$ and $f(N)>g(N)$ for $N$ the largest index with $f(N)\neq g(N)$.
In particular, there are canonical well-founded total orders on the set $(\bN^\bN)_{\mathrm{cpt}}$ and on the set $(((\bN^{\bN})_{\mathrm{cpt}})^\bN)_{\mathrm{cpt}}$.

Recall that given a set $X$ and partial orders $>_1,\ldots,>_n$ on $X$, the \defnemph{lexicographic order for $(>_1,\ldots,>_n)$} is inductively defined as $x>_{1,\ldots,n}y$ if and only if either $x>_1y$, or ($x\not>_1y$ and $x\not<_1y$) and $x>_{2,\ldots,n}y$, where $>_{2,\ldots,n}$ is the lexicographic order for $(>_2,\ldots,>_n)$.
When $\#\colon X\to (\Lambda,>)$ is a function to a poset $(\Lambda,>)$, we sometimes abuse notation and refer to $\#$ as a partial order on $X$.

\subsubsection{Orders of affine Brauer type}

For a diagram $D$, we set:
\begin{description}[align=left, labelsep=1em, labelindent=\parindent,
    labelwidth=\widthof{$\nbbubble^D$},
    leftmargin=\widthof{$\nbbubble^D$}+1em+\parindent]
  \item[$\nbcro^D$] the number of crossings in $D$;
  \item[$\nbdot^D$] the number of dots in $D$;
  \item[$\nbcc^D$] the number of cups and caps in $D$;
  \item[$\distcc^D$] the number of crossings encountered when following the left strand of a cup or a cap before encountering an endpoint or a critical point, summed over all cups and caps;
  \item[$\distcro^D$] the number of crossings encountered when following the bottom left strand of a crossing before encountering an endpoint or a critical point, summed over all crossings;
  \item[$\nbbubble^D$] the function $\nbbubble^D\colon\bN\to\bN$ where $\nbbubble^D(k)$ is the total number of bubbles with $k$ dots in $D$;
\end{description}
In what follows, we often omit $D$ if it is clear from the context.

Recall $\disttunnel\colon\bN\to\bN^\bN$ from \cref{defn:distance_tunnel} and $\disttraffic$ from \cref{defn:distance_traffic}, where the latter depends on a choice of traffic rules.
Note that $\nbbubble\in(\bN^\bN)_{\mathrm{cpt}}$ and $\disttunnel\in(((\bN^{\bN})_{\mathrm{cpt}})^\bN)_{\mathrm{cpt}}$; each induces a partial order on diagrams via the canonical total orders on $(\bN^\bN)_{\mathrm{cpt}}$ and $(((\bN^{\bN})_{\mathrm{cpt}})^\bN)_{\mathrm{cpt}}$ given above.

\begin{definition}
  \label{defn:orders_gen_brauer}
  On the set of diagrams, we define:
  \begin{description}[align=left, labelsep=1em, labelindent=\parindent,
    labelwidth=\widthof{$\strongsucc$},
    leftmargin=\widthof{$\strongsucc$}+1em+\parindent]
    \item[$\strongsucc$] the partial order given by the lexicographic order for $(\nbcro,\nbdot)$;
    \item[$\asymp$] the equivalence relation where $x\asymp y$ if and only if $x$ and $y$ are identical diagrams up to different dots on open strands and different dotted bubbles in regions;
    \item[$\rhd$] the partial order where $x\rhd y$ if $x\strongsucc y$ or if $(\nbcro^x,\nbdot^x)=(\nbcro^y,\nbdot^y)$ and $x\asymp y$ and $\nbbubble^x>\nbbubble^y$;
    \item[${\succ}$] for a choice of traffic rules $\trafficrule$, the partial order ${\succ_\trafficrule}$ (written $\succ$ when $\trafficrule$ is clear from the context) given by the lexicographic order for
    \[(\nbcro,\nbdot,\disttunnel,\nbcc,\distcc,\distcro,\disttraffic).\]
  \end{description}
\end{definition}




\subsection{Rewriting}
\label{subsec:defn_rewriting_affine_brauer}
\subsubsection{Monoidal rewriting system}
\label{subsubsec:defn_rewriting_system}

Fix a ring $\ring$, an abelian group $G$, and a symmetric bilinear map $\bil\colon G\times G\to\ring^\times$, as in \cref{subsec:RW_graded_monoidal_categories}.
We denote by $\affE$ the set of $(G,\bil)$-graded interchangers.
Let also $\affX_0= \{\;\tikzpic{\li}\;\}$ and 
\begin{gather*}
  \affX_1
  = \left\{\;
    \tikzpic{\cro}\;,\;
    \tikzpic{\cu}\;,\;
    \tikzpic{\ca}\;,\;
    \tikzpic{\bli}
  \;\right\}
\end{gather*}
viewed as a $G$-graded set and with source and target maps given by reading the diagrams of $\affX_1$ from bottom to top.
Given this data, we consider the set $\affR$ of rewriting steps given in \cref{fig:generic_grobner_of_affine_Brauer_type}, together with the following extra data:
\begin{enumerate}[(i)]
  \item A set $\bubbleset\subset\bN_{>0}$, which gives a choice of bubbles.
  \item Scalars $\VARRthree$, $\VARcasl$, $\VARcusl$, $\VARbbsl$, $\VARdcro$, $\widetilde{\VARdcro}$, $\VARdcc$ and $\VARdccbis$, assumed to be invertible, and scalars $\VARzz$, $\VARzzbis$, $\VARdccLOT$ and $\VARdccLOTbis$ (which need not be invertible).
  \item Vectors pictured with double-lined borders:
  \begin{gather*}
    \renewcommand{\sca}{.4}
    \newcommand{\hspc}{10mu}
    \underbrace{
    \tikzpic{\LOTRtwo}
    ,\mspace{\hspc}
    \tikzpic{\clip (-.3,.5) rectangle (2.3,2.5);\LOTRthree}
    ,\mspace{\hspc}
    \tikzpic{\LOTcasl}
    ,\mspace{\hspc}
    \tikzpic{\LOTcusl}
    ,\mspace{\hspc}
    \tikzpic{\LOTuk}
    ,\mspace{\hspc}
    \tikzpic{\LOTdk}
    ,\mspace{\hspc}
    \tikzpic{\LOTdcro}
    ,\mspace{\hspc}
    \tikzpic{\LOTdcrobis}
    }_{\text{$\strongsucc$-lower-order terms}}
    ,\mspace{\hspc}
    \tikzpic{\LOTev}
    \mspace{\hspc}\an\mspace{\hspc}
    \tikzpic{\LOTbbsl}
    \;.
  \end{gather*}
  We assume that each of them is a linear combination of diagrams that are lower, with respect to the order $\rhd$, than the source of their respective rewriting step.
  We call them the \defnemph{$\rhd$-lower-order terms} and the remaining diagram in the target, if it exists, the \defnemph{$\rhd$-leading term}.
\end{enumerate}
As indicated, many $\rhd$-lower-order terms are actually \defnemph{$\strongsucc$-lower-order terms}, that is, lower with respect to the order $\strongsucc$; indeed, notice that the source $x$ of their rewriting step does not contain any bubble, so that if $y$ is such that $x\rhd y$, then we must have $x\strongsucc y$.
In practice, these diagrams are often even ``$\nbcro$-lower-order terms'', in the sense that they have strictly fewer crossings than the source of their rewriting step.
Note that the diagrams $\tikzpic{\ca}[scale=.5]$ and $\tikzpic{\cu}[scale=.5]$ in the target of the dot-slide-through-cup-and-cap rewriting steps are also $\strongsucc$-lower-order terms.

A diagram appearing in one of the last two $\rhd$-lower-order terms has either strictly fewer than $k$ dots or exactly $k$ dots; in the latter case, it must be a union of bubbles with at most $k-1$ dots (case $\tikzpic{\LOTev}$) or a straight line possibly carrying dots, together with unions of bubbles with at most $k-1$ dots on its sides (case $\tikzpic{\LOTbbsl}$).

\begin{definition}
  \label{defn:main_text_affine_Brauer_rewriting_system}
  A linear monoidal \defnemph{rewriting system} is \defnemph{of affine Brauer type} with $\bubbleset$-bubbles if it is of the form $\affS=(\affX_0,\affX_1;\affR,\affE)$ with $\bubbleset$ as the choice of bubbles.
\end{definition}

\begin{figure}[p]
	\centering

	\def\spaceout{2ex}
	\def\spacein{.5ex}
  \renewcommand{\sca}{.4}
  \newcommand{\hspc}{60mu}

	\begin{subfigure}{\textwidth}
		\begin{gather*}
			\begin{IEEEeqnarraybox}{CcC}
				\tikzpic{
          \cro[0][1]
				  \cro 
				} 
				\;\overset{\stRtwo}{\longrightarrow}\;
        \tikzpic{\LOTRtwo}
				\mspace{30mu}
        %
				\tikzpic{
          \cro[0][2] \li[2][2]
          \li[0][1] \cro[1][1] 
          \cro \li[2][0] 
        } 
				\;\overset{\stRthree}{\longrightarrow}\;
        \VARRthree\;
				\tikzpic{
          \cro[1][2] \li[0][2]
          \li[2][1] \cro[0][1] 
          \cro[1][0] \li[0][0] 
        }
        \;+\;
        \tikzpic{\LOTRthree}
        &
				\mspace{60mu}
        &
				\tikzpic{
          \cu[1][0] \uli[2][0] 
				  \dli \ca 
				}
				\;\overset{\stzz}{\longrightarrow}\;
        \VARzz\;
        \tikzpic{\uli\dli}
				\mspace{30mu}
        %
				\tikzpic{
          \dli[2][0] \ca[1][0] 
          \uli \cu 
				}
				\;\overset{\stzzbis}{\longrightarrow}\;
        \VARzzbis\;
        \tikzpic{\uli\dli}
        \\[\spacein]
        \text{\footnotesize braiding}
        &&
        \text{\footnotesize zigzags}
      \end{IEEEeqnarraybox}
      \\[\spaceout]
      \begin{IEEEeqnarraybox}{rClcrCl}
				\tikzpic{
          \uli[0][1] \ca[1][1] 
          \cro \li[2][0] 
        }
				\;&\overset{\stcasl}{\longrightarrow}&\;
        \VARcasl\;
				\tikzpic{ 
          \ca[0][1] \uli[2][1]
          \li \cro[1][0] 
        }
        \;+\;
        \tikzpic{\LOTcasl}
				&\mspace{100mu}&
        %
				\tikzpic{
          \ca \uli[2][0] 
          \li[0][-1] \cro[1][-1]
          \cro[0][-2] \li[2][-2]
        }
				\;&\overset{\stcapu}{\longrightarrow}&\;
        \VARcasl^{-1}\;
        \tikzpic{
          \uli[0][2]\ca[1][2]
          \LOTRtwo\li[2][0]\li[2][1]
        }
        -
        \VARcasl^{-1}\;
        \tikzpic{
          \LOTcasl[0][1]
          \cro[0][0] \li[2][0]
        }
        \\[\spacein]
        \IEEEeqnarraymulticol{3}{c}{
          \text{\footnotesize cap sliding}
        }
				&&
        \IEEEeqnarraymulticol{3}{c}{
          \text{\footnotesize cap pulling}
        }
        \\[\spaceout]
        \tikzpic{ 
          \cro \li[2][0]
          \dli \cu[1][0] 
        }
        \;&\overset{\stcusl}{\longrightarrow}&\;
        \VARcusl\;
        \tikzpic{
          \li[0][0] \cro[1][0]
          \cu[0][0] \dli[2][0]
        }
        \;+\;
        \tikzpic{\LOTcusl}
        &&
        \tikzpic{
          \cro[0][1] \li[2][1] 
          \li[0][0] \cro[1][0]
          \cu \dli[2][0]
        }
        \;&\overset{\stcupu}{\longrightarrow}&\;
        \VARcusl^{-1}\;
        \tikzpic{
          \LOTRtwo[0][1]\li[2][1]\li[2][2]
          \dli[0][1]\cu[1][1]
        }
        -
        \VARcusl^{-1}\;
        \tikzpic{
          \cro[0][2] \li[2][2]
          \LOTcusl[0][0]
        }
        \\[\spacein]
        \IEEEeqnarraymulticol{3}{c}{
          \text{\footnotesize cup sliding}
        }
				&&
        \IEEEeqnarraymulticol{3}{c}{
          \text{\footnotesize cup pulling}
        }
			\end{IEEEeqnarraybox}
			\\[\spaceout]
      \begin{IEEEeqnarraybox}{CcCcC}
				\tikzpic{
          \ca[0][1]
				  \cro 
				}
        \;\overset{\stuk}{\longrightarrow}\;
        \tikzpic{\LOTuk}
				&\mspace{100mu}&
        %
        \tikzpic{
          \cro
          \cu
        }
        \;\overset{\stdk}{\longrightarrow}\;
        \tikzpic{\LOTdk}
        &\mspace{100mu}&
        %
        \tikzpic{
          \ca[0][1] \uli[2][1]
          \li[0][0] \cro[1][0]
          \cu \dli[2][0]
        }
        \;\overset{\stlk}{\longrightarrow}\;
        \VARcusl^{-1}\;
        \tikzpic{\LOTuk[0][0]\dli[0][0]\cu[1][0]\li[2][0]\uli[2][1]}
        -
        \VARcusl^{-1}\;
        \tikzpic{\LOTcusl\ca[0][2]\uli[2][2]}
        \\[\spacein]
        \text{\footnotesize upward kink}
        &&
        \text{\footnotesize downward kink}
        &&
        \text{\footnotesize leftward kink}
      \end{IEEEeqnarraybox}
  	\end{gather*}
	
		\caption{Non-affine rewriting steps.}
		\label{subfig:generic_nonaffine_rewriting_steps}
	\end{subfigure}

	\vspace*{.5cm}

  \begin{subfigure}{\textwidth}
    \begin{gather*}
      \tikzpic{\bubble[0][0][k]}[][1]
      \;\overset{\stev_k}{\longrightarrow}\;
      \tikzpic{\LOTev}
    \end{gather*}
		\caption{Bubble-evaluation rewriting step, where $k\notin \bubbleset$.}
		\label{subfig:generic_bubble_evaluation_rewriting_steps}
  \end{subfigure}

  \vspace*{.5cm}

  \begin{subfigure}{\textwidth}
    \begin{gather*}
      \tikzpic{\bubble[0][0][k]\li[-1.2][.5]}[][1]
      \;\overset{\stbbsl_k}{\longrightarrow}\;
      \VARbbsl\;
      \tikzpic{\bubble[0][0][k]\li[1.5][.5]}[][1]
      \;+\;
      \tikzpic{\LOTbbsl[0][0][k]}
      \mspace{50mu}
      \tikzpic{\bubble[0][0][k]\li[1.5][.5]}[][1]
      \;\overset{(\stbbsl_k)^*}{\longrightarrow}\;
      \VARbbsl^{-1}\;
      \tikzpic{\bubble[0][0][k]\li[-1.2][.5]}[][1]
      \;-\;\VARbbsl^{-1}\;
      \tikzpic{\LOTbbsl[0][0][k]}
    \end{gather*}
		\caption{Bubble-slide rewriting steps.}
		\label{subfig:generic_bubble_slide_rewriting_steps}
  \end{subfigure}

  \vspace*{.5cm}

	\begin{subfigure}{\textwidth}
		\begin{gather*}
      \begin{IEEEeqnarraybox}{rClcrCl}
        \tikzpic{
          \ulcr
          \Diag[0][0][]
        }
        \;&\overset{\stdcro}{\longrightarrow}&\;
        \VARdcro\;
        \tikzpic{
          \drcr
          \Diag[0][0][]
        }
        \;+\;
        \tikzpic{\LOTdcro}
        &\mspace{50mu}&
        \tikzpic{
          \drcr
          \Diag[0][0][]
        }
        \;&\overset{(\stdcro)^*}{\longrightarrow}&\;
        \VARdcro^{-1}\;
        \tikzpic{
          \ulcr
          \Diag[0][0][]
        }
        \;-\;\VARdcro^{-1}\;
        \tikzpic{\LOTdcro}
        \\[\spacein]
        \tikzpic{
          \urcr
          \AntiDiag[0][0][]
        }
        \;&\overset{\stdcrobis}{\longrightarrow}&\;
        \VARdcrobis\;
        \tikzpic{
          \dlcr
          \AntiDiag[0][0][]
        }
        \;+\;
        \tikzpic{\LOTdcrobis}
        &\mspace{50mu}&
        \tikzpic{
          \dlcr
          \AntiDiag[0][0][]
        }
        \;&\overset{(\stdcrobis)^*}{\longrightarrow}&\;
        \VARdcrobis^{-1}\;
        \tikzpic{
          \urcr
          \AntiDiag[0][0][]
        }
        \;-\;\VARdcrobis^{-1}\;
        \tikzpic{\LOTdcrobis}
        \\[\spaceout]
        \tikzpic{\rca \ca[0][0][]}
        \;&\overset{\stdca}{\longrightarrow}&\;
        \VARdcc\;
        \tikzpic{\lca \ca[0][0][]}
        \;+\;
        \VARdccLOT\;\tikzpic{\ca}
        &\mspace{100mu}&
        \tikzpic{\lca \ca[0][0][]}
        \;&\overset{(\stdca)^*}{\longrightarrow}&\;
        \VARdcc^{-1}\;
        \tikzpic{\rca \ca[0][0][]}
        \;-\;\VARdcc^{-1}
        \VARdccLOT\;\tikzpic{\ca}
        \\[\spacein]
        \tikzpic{\rcu \cu[0][0][]}
        \;&\overset{(\stdcu)^*}{\longrightarrow}&\;
        \VARdccbis^{-1}\;
        \tikzpic{\lcu \cu[0][0][]}
        \;-\;
        \VARdccbis^{-1}\;
        {\VARdccLOTbis}\;\tikzpic{\cu}
        &\mspace{100mu}&
        \tikzpic{\lcu \cu[0][0][]}
        \;&\overset{\stdcu}{\longrightarrow}&\;
        \VARdccbis\;
        \tikzpic{\rcu \cu[0][0][]}
        \;+\;
        {\VARdccLOTbis}\;\tikzpic{\cu}
      \end{IEEEeqnarraybox}
		\end{gather*}
	
		\caption{Dot-slide rewriting steps.}
		\label{subfig:generic_dot_slide_rewriting_steps}
	\end{subfigure}
	
	\caption{The general form of a rewriting system of affine Brauer type with $\bubbleset$-bubbles.
  We assume that the scalars $\VARRthree$, $\VARcasl$, $\VARcusl$, $\VARbbsl$, $\VARdcro$, $\widetilde{\VARdcro}$, $\VARdcc$ and $\VARdccbis$ are invertible.
  The schematics with double-lined borders are either $\strongsucc$-lower-order terms or $\rhd$-lower-order terms, consisting of diagrams that are lower, with respect to $\strongsucc$ or $\rhd$ respectively, than the source of the rewriting step.}
	\label{fig:generic_grobner_of_affine_Brauer_type}
\end{figure}

\subsubsection{Sub-system}

\begin{definition}
  \label{defn:ABT}
  Fix a choice of traffic rules $\trafficrule$.
  We define $\affT=\affT_\trafficrule$ as the linear sub-system of $\affS$ such that for each $\ssr\in\affR$, the rewriting $\Gamma[\ssr]$ belongs to $\affT$ if and only if:
  \begin{enumerate}[(i)]
    \item when $\ssr$ is a non-affine rewriting step (\cref{subfig:generic_nonaffine_rewriting_steps}), $\Gamma$ is any context;

    \item when $\ssr$ is a bubble evaluation (see \cref{subfig:generic_bubble_evaluation_rewriting_steps}), $\Gamma$ is any context;

    \item when $\ssr$ is a bubble slide (see \cref{subfig:generic_bubble_slide_rewriting_steps}), $\Gamma$ is such that, according to tunnels in $\Gamma[s(\ssr)]$, the bubble follows the orientation of the tunnel;
    
    \item when $\ssr$ is a dot slide (see \cref{subfig:generic_dot_slide_rewriting_steps}), $\Gamma$ is such that, according to the traffic rules $\trafficrule$ on $\Gamma[s(\ssr)]$, the dot follows the orientation of the strand, and in the case of $s(\ssr) = \tikzpic{\lcu \cu[0][0][->]}$, there is no stop sign on the cup.
  \end{enumerate}
\end{definition}

We let $\affR_{\nonaff}\subset\affR$ be the set of non-affine rewriting steps (i) and set $\affT_{\nonaff}=\Cont(\affR_{\nonaff})$. Likewise, we let $\affR_{\aff} = \affR\setminus\affR_{\nonaff}$ and $\affT_\aff=\affT\cap\Cont(\affR_\aff)$; their rewriting steps constitute the \emph{affine rewriting steps}.
Note that $\affR_{\nonaff}$ is context-agnostic while $\affT_\aff$ is context-dependent. (Although $\affT_\aff$ does contain context-agnostic rewriting steps, namely bubble-evaluation rewriting steps.)

\subsection{Normal forms}
\label{subsec:normal_brauer_diagrams}

By definition, a monomial $\affT_\trafficrule$-normal form is a diagram whose dots sit close to their dead ends, whose only bubbles are $\bubbleset$-bubbles on the leftmost region, and which, up to interchange, avoids the following patterns:
\begin{gather}
  \label{eq:forbidden_patterns}
  \tikzpic{
    \cro[0][1]
    \cro 
  }
  \quad
  \tikzpic{
    \cro[0][2] \li[2][2]
    \li[0][1] \cro[1][1] 
    \cro \li[2][0] 
  }[scale=.8]
  \quad
  \tikzpic{
    \cu[1][0] \uli[2][0] 
    \dli \ca 
  }
  \quad
  \tikzpic{
    \dli[2][0] \ca[1][0] 
    \uli \cu 
  }
  \quad
  \tikzpic{
    \uli[0][1] \ca[1][1] 
    \cro \li[2][0] 
  }
  \quad
  \tikzpic{
    \ca \uli[2][0] 
    \li[0][-1] \cro[1][-1]
    \cro[0][-2] \li[2][-2]
  }[scale=.8]
  \quad
  \tikzpic{ 
    \cro \li[2][0]
    \dli \cu[1][0] 
  }
  \quad
  \tikzpic{
    \cro[0][1] \li[2][1] 
    \li[0][0] \cro[1][0]
    \cu \dli[2][0]
  }[scale=.8]
  \quad
  \tikzpic{
    \ca[0][1]
    \cro 
  }
  \quad
  \tikzpic{
    \cro
    \cu
  }
  \quad
  \tikzpic{
    \ca[0][1] \uli[2][1]
    \li[0][0] \cro[1][0]
    \cu \dli[2][0]
  }[scale=.8]\;.
\end{gather}
We call them the \defnemph{forbidden patterns}.
We now give two further descriptions of normal forms: one combinatorial, one topological.
\Cref{fig:normal_diagram_example} shows an example of a monomial $\affT_\trafficrule$-normal form.

A diagram is \defnemph{left-cc-justified} if the left strand of any cup (resp. cap) on an open strand is a straight line (possibly with dots) to a top (resp. bottom) point.
A diagram is \defnemph{braid-justified} if, up to interchange, it avoids the pattern $\tikzpic{
    \cro[0][2] \li[2][2]
    \li[0][1] \cro[1][1] 
    \cro \li[2][0] 
  }[scale=.5]$. Standard diagrams are defined below.

\begin{restatable}{theorem}{thmnormalformstandarddiagrams}
  \label{thm:normal_form_are_standard_diagrams}
  Let $\trafficrule=(\trafficend,\bubbleset)$ be a choice of traffic rules. For a diagram $D$, the following are equivalent:
  \begin{enumerate}[(i)]
    \item $D$ is a monomial $\affT_\trafficrule$-normal form;
    \item $D$ is, up to interchange, a $(\trafficend,\bubbleset)$-standard diagram;
    \item $D$ is a left-cc-justified and braid-justified $(\trafficend,\bubbleset)$-reduced matching.
  \end{enumerate}
  Moreover, the set of such diagrams is, up to choosing an interchange representative, a representative set of $(\trafficend,\bubbleset)$-reduced matchings.
\end{restatable}

We now define the $(\trafficend,\bubbleset)$-standard diagrams appearing in \cref{thm:normal_form_are_standard_diagrams}(ii).

A \emph{cap-standard diagram} is defined inductively as a diagram which is either the identity or of the form $\Gamma[D]$, where $\Gamma=\Gamma_{\mathrm{cap\text{-}std}}$ is the \emph{cap-standard context} from \cref{subfig:cap_standard_context} and $D$ is a cap-standard diagram.
\emph{Braid-standard diagrams} are defined in the same way, using the \emph{braid-standard context} $\Gamma_{\mathrm{braid\text{-}std}}$ from \cref{subfig:braid_standard_context}.
Cup-standard diagrams are defined analogously, with the ``cup-standard context'' being the vertical symmetry of the cap-standard context.
A \emph{dotted identity} is an identity with extra dots on its strands.
A \emph{$(\trafficend,\emptyset)$-standard diagram} is a diagram of the form:
\[
(\text{dotted identity})\circ(\text{cup-standard})\circ(\text{braid-standard})\circ(\text{cap-standard})\circ(\text{dotted identity}),
\]
where dots sit close to their dead ends.
Finally, a \defnemph{$(\trafficend,\bubbleset)$-standard diagram} is a horizontal juxtaposition of $\bubbleset$-bubbles on the leftmost region and of a $(\trafficend,\emptyset)$-standard diagram.

\begin{figure}
  \centering
  \begin{minipage}[b]{.4\textwidth}
    \centering
    \tikzpic{
      \lili[-1][4] \bul[-1][5-.25][2] \li[1][4]\bul[1][5-.25] \cu[1][4] \cro[2][4]
      \Diag[-1][3]\Diag[0][3]\AntiDiag[2][3]
      \li[0][2] \cro[1][2]
      \cc[0][1] \li[2][1]
      \li\bul[0][.25][3] \dlcr[1][0]
      \bubble[-2][1][4]
    }
    \vspace{2em}
    \captionof{figure}{A normal diagram ($4\in\bubbleset$).}
    \label{fig:normal_diagram_example}
  \end{minipage}%
  \hfill
  \begin{minipage}[b]{.6\textwidth}
    \centering
  \begin{subfigure}[b]{.49\linewidth}
    \begin{gather*}
      \newcommand{\vsh}{.1}
      \newcommand{\hsh}{.15}
      \tikzpic{
        \draw (0,0) to (0,4+\vsh);
        \draw (1,0) to (1,4+\vsh);
        \draw (2,1-\vsh) to (2,3);
        \draw (3,1-\vsh) to (3,2) to (4,3) to (4,4+\vsh);
        \draw (4,1-\vsh) to (5,2) to (5,4+\vsh);
        \draw (6,1-\vsh) to (6,4+\vsh);
        \draw (7,1-\vsh) to (7,4+\vsh);
        \ca[2][3]\draw (3,3) to (5,1) to (5,1-\vsh);
        \draw (2-\hsh,0+\vsh) rectangle (7+\hsh,1-\vsh);
        \node[] at (.5,0) {\scriptsize\ldots};
        \node[] at (.5,4) {\scriptsize\ldots};
        \node[] at (3.5,1) {\scriptsize\ldots};
        \node[] at (4.5,4) {\scriptsize\ldots};
        \node[] at (6.5,1) {\scriptsize\ldots};
        \node[] at (6.5,4) {\scriptsize\ldots};
        \draw[<->,thick] (3.5-.2,2.5-.2) to node[draw=none,circle,fill=white,inner sep=.2pt] {\scriptsize $n$} (4.5-.2,1.5-.2);
        %
      }
    \end{gather*}
    \caption{Cap-standard context.}
    \label{subfig:cap_standard_context}
  \end{subfigure}
  \begin{subfigure}[b]{.49\linewidth}
    \begin{gather*}
      \newcommand{\vsh}{.1}
      \newcommand{\hsh}{.15}
      \tikzpic{
        \draw (0,0) to (0,2) to (3,5);
        \draw (1,0) to (1,1) to (2,2) to (2,2+\vsh);
          \draw (2,3-\vsh) to (2,3) to (1,4) to (1,5); 
        \draw (2,0) to (3,1) to (3,2) to (3,2+\vsh);
          \draw (3,3-\vsh) to (3,3) to (3,4) to (2,5); 
        \draw (3,0) to (0,3) to (0,5);
        \draw (3.5,0) to (3.5,5);
        \draw (4.5,0) to (4.5,5);
        \draw[fill=white] (2-\hsh,2+\vsh) rectangle (4.5+\hsh,3-\vsh);
        \node[] at (1.5,0) {\scriptsize\ldots};
        \node[] at (2.5,2) {\scriptsize\ldots};
        \node[] at (2.5,3) {\scriptsize\ldots};
        \node[] at (1.5,5) {\scriptsize\ldots};
        \node[] at (4,0) {\scriptsize\ldots};
        \node[] at (4,2) {\scriptsize\ldots};
        \node[] at (4,3) {\scriptsize\ldots};
        \node[] at (4,5) {\scriptsize\ldots};
        \draw[<->,thick] (1.5-.2,1.5-.2) to node[draw=none,circle,fill=white,inner sep=0pt] {\scriptsize $m$} (2.5-.2,.5-.2);
        \draw[<->,thick] (1.5-.2,3.5+.2) to node[draw=none,circle,fill=white,inner sep=0pt] {\scriptsize $n$} (2.5-.2,4.5+.2);
      }
    \end{gather*}
    \caption{Braid-standard context.}
    \label{subfig:braid_standard_context}
  \end{subfigure}
  \caption{Standard contexts.}
  \label{fig:standard_contexts}
  \end{minipage}
\end{figure}

\medbreak
The proof of \cref{thm:normal_form_are_standard_diagrams} is postponed to \cref{sec:proof_normal_forms}, which also contains the proof of the following lemma:

\begin{restatable}{lemma}{lembasisinduced}
  \label{lem:basis_induced_matching_works_for_all}
  Let $(\trafficend,\bubbleset)$ be a choice of traffic rules and $\affC$ be a category presented by a rewriting system of affine Brauer type $\affS$.
  If some representative set of $(\trafficend,\bubbleset)$-reduced matchings defines a hom-basis of $\affC$, then every such set does.
\end{restatable}

\subsection{Tamed preconfluence}
\label{subsec:tamed_preconfluence}

\newcommand{\temporder}{>}

For a partial order $\temporder$ and a monomial $x\in\sX^*$, we denote
\[\smallO_{\temporder}(x)=\{y\in\sX^*\mid x\temporder y\}\]
and call it the \defnemph{set of $\temporder$-lower-order terms}.

Assume that a linear rewriting system is given with an order $\temporder$, such that rewriting steps are either of the form $\ssr\colon s\to \lambda t + v$ or of the form $\ssr\colon s\to v$, where $\smallO_{\temporder}(s)=\smallO_{\temporder}(t)$ and $s\temporder v$.
A $\temporder$-tamed preconfluence can be understood as follows:
\begin{quote}
  \textsc{Tamed preconfluence (idea):} given a branching, rewrite the $\temporder$-leading term on each branch as much as possible until they either become equal or disappear; then, relate the remaining $\temporder$-lower-order terms using a $\temporder$-tamed congruence.
\end{quote}
This is formalized in the following definition.

\begin{definition}
  Let $\sS$ be a linear monoidal rewriting system and $\temporder$ an isp-order. Consider sub-systems $\aA$ and $\aB$ of $\sS$, as in \cref{subsec:RW_tame_congruence}.
  Let $(f,g)$ be an $\aA$-branching with source $\bullet$.
  A \defnemph{$\temporder$-tamed $\aB$-preconfluence} for $(f,g)$ is a triple $(f',h,g')$, such that:
  \begin{enumerate}[(i)]
    \item $f'=f_n'\starop\ldots\starop f_1'$ with $s(f_1')=t(f)$, and
    \[f'_i = \lambda_i \ssf'_i + v_i,\qquad\smallO_{\temporder}(s(\ssf'_i))=\smallO_{\temporder}(\bullet)\quad\an\quad \bullet\temporder v_i.\]
    Here $\ssf'_i$ is a monomial $\aB$-rewriting step, $\lambda_i$ is an invertible scalar and $v_i$ is a vector.
    \item Similarly, $g'=g_m'\starop\ldots\starop g_1'$ with $s(g_1')=t(g)$, and
    \[g'_j = \mu_j \ssg'_j + w_j,\qquad\smallO_{\temporder}(s(\ssg'_j))=\smallO_{\temporder}(\bullet)\quad\an\quad \bullet\temporder w_j.\]
    Here $\ssg'_j$ is a monomial $\aB$-rewriting step, $\mu_j$ is an invertible scalar and $w_j$ is a vector.
    \item 
    We have $t(f')=\lambda x + v$ and $t(g')=\mu y + w$, with $\lambda,\mu$ not necessarily invertible scalars, $x,y\in\sX^*$ monomials such that $\smallO_{\temporder}(x)=\smallO_{\temporder}(y)=\smallO_{\temporder}(\bullet)$ and $v,w$ vectors such that $\bullet\temporder v,w$.
    Moreover:
    \[\widehat{\bullet} \coloneqq\lambda x = \mu y\quad\an\quad h = \widehat{\bullet} + h_{\temporder}\]
    where $h_{\temporder}$ is an $\sS$-congruence $\temporder$-tamed by $\bullet$.
  \end{enumerate}
\end{definition}

Here is a schematic for a $\temporder$-tamed $\aB$-preconfluence.
\begin{gather*}
\begin{tikzcd}[ampersand replacement=\&,row sep=small]
	\& {\lambda_1x_1+v_1} \& {\lambda_2x_2+v_2} \& \ldots \& \lambda x + v \& \widehat{\bullet} + v\\
	\bullet \\
	\& {\mu_1y_1+w_1} \& {\mu_2y_2+w_2} \& \ldots \& \mu y + w \& \widehat{\bullet} + w
	\arrow["{f'_1}","\aB"{subscript}, from=1-2, to=1-3]
	\arrow["{f'_2}","\aB"{subscript}, from=1-3, to=1-4]
	\arrow["{f'_n}","\aB"{subscript}, from=1-4, to=1-5]
	\arrow["{h_{\temporder}}", shift left=4, tail reversed, from=1-6, to=3-6]
	\arrow[shift right=4, equals, from=1-6, to=3-6]
	\arrow["f", bend left=25, from=2-1, to=1-2]
	\arrow["g"', bend right=25, from=2-1, to=3-2]
	\arrow["{g'_1}","\aB"{subscript}, from=3-2, to=3-3]
	\arrow["{g'_2}","\aB"{subscript}, from=3-3, to=3-4]
	\arrow["{g'_m}","\aB"{subscript}, from=3-4, to=3-5]
  \arrow[draw=none,"="{description},from=1-5,to=1-6]
  \arrow[draw=none,"="{description},from=3-5,to=3-6]
\end{tikzcd}
\end{gather*}
The conditions $\smallO_{\temporder}(s(\ssf'_i))=\smallO_{\temporder}(\bullet)$ and $\bullet\temporder v_i$ imply that $s(\ssf'_i)\temporder v_i$. In particular, the $\aB$-rewriting step $f_i'$ is necessarily positive; and similarly for $g_j'$.
The vector $\widehat{\bullet}$ may be zero.

\subsubsection{Tamed preconfluence for non-affine branchings}
\label{subsubsec:example_tamed_preconfluence_nonaffine}

Below is an example of a $\strongsucc$-tamed $\affT$-preconfluence for a non-affine branching:
\begin{gather*}
  \renewcommand{\sca}{.25}
  \begin{tikzcd}[ampersand replacement=\&,row sep=0pt]
    \&
    \tikzpic{\li\cro[1][0]\LOTRthree[0][1]}
    \;+\;
    \VARRthree\;\tikzpic{\cro[0][3]\li[2][3]\li[0][2]\cro[1][2]\cro[0][1]\li[2][1]\li[0][0]\cro[1][0]}
    \&
    \VARRthree^2\;\tikzpic{
      \li[0][3]\cro[1][3]  
      \cro[0][2]\li[2][2]
      \cro[1][0]\cro[1][1]\li[0][0]\li[0][1]
    }
    \;+\;
    \VARRthree\;\tikzpic{\li\cro[1][0]\LOTRthree[0][1]}
    \;+\;
    \tikzpic{\LOTRthree\cro[0][3]\li[2][3]}
    \&
    \VARRthree^2\;\tikzpic{
      \li[0][3]\cro[1][3]  
      \cro[0][2]\li[2][2]
      \LOTRtwo[1][0]\li[0][0]\li[0][1]
    }
    \;+\;
    \VARRthree\;\tikzpic{\li\cro[1][0]\LOTRthree[0][1]}
    \;+\;
    \tikzpic{\LOTRthree\cro[0][3]\li[2][3]}
    \\
    \tikzpic{
      \cro[0][3]\li[2][3]
      \cro[0][2]\li[2][2]
      \li[0][1]\cro[1][1]
      \cro[0][0]\li[2][0]
    }
    \\
    \&
    \tikzpic{
      \LOTRtwo[0][2]\li[2][2]\li[2][3]
      \li[0][1]\cro[1][1]  
      \cro[0][0]\li[2][0]
    }
    \&\&
    \tikzpic{
      \LOTRtwo[0][2]\li[2][2]\li[2][3]
      \li[0][1]\cro[1][1]  
      \cro[0][0]\li[2][0]
    }
    \arrow[from=2-1,to=1-2, bend left=20pt]
    \arrow[from=1-2,to=1-3]
    \arrow[from=1-3,to=1-4]
    \arrow[from=2-1,to=3-2, curve={height=12pt}]
    \arrow[from=3-2,to=3-4,equals]
    \arrow[from=1-4,to=3-4, "{h_{\strongsucc}}", tail reversed, dottedcd]
  \end{tikzcd}
\end{gather*}
Here, the vector $\widehat{\bullet}$ is zero.
The branching is $\strongsucc$-tamely $\affT$-preconfluent if and only if there exists a $\affS$-congruence
\begin{gather*}
  \renewcommand{\sca}{.25}
  h_{\strongsucc}\colon\qquad
  \VARRthree^2\;\tikzpic{
    \li[0][3]\cro[1][3]  
    \cro[0][2]\li[2][2]
    \LOTRtwo[1][0]\li[0][0]\li[0][1]
  }
  \;+\;
  \VARRthree\;\tikzpic{\li\cro[1][0]\LOTRthree[0][1]}
  \;+\;
  \tikzpic{\LOTRthree\cro[0][3]\li[2][3]}
  \;\overset{*}{\longleftrightarrow}\;
  \tikzpic{
    \LOTRtwo[0][2]\li[2][2]\li[2][3]
    \li[0][1]\cro[1][1]  
    \cro[0][0]\li[2][0]
  }
\end{gather*}
which passes only through diagrams with at most three crossings.

\subsubsection{Tamed preconfluence for dot-slide branchings}
\label{subsubsec:example_tamed_preconfluence_dotslide}

We give two examples of $\strongsucc$-tamed $\affT$-preconfluence for dot-slide branchings; we restrict to the context of the affine Brauer category (see \cref{fig:AB_S_rewriting_steps}).
The first example is a branching between an R2 move and a dot-slide:
\begin{gather*}
  \newcommand{\arOrNot}[1]{#1}
  \renewcommand{\sca}{.3}
  \begin{tikzcd}[ampersand replacement=\&,row sep=tiny]
    \& 
    \tikzpic{\cro[0][1]\cro\bul[1][1]\AntiDiag[0][0][<-]}
    \;-\;
    \tikzpic{\lili[0][1]\cro}
    \;+\;
    \tikzpic{\cc[0][1]\cro}
    \& 
    \tikzpic{\cro[0][1]\dlcr\AntiDiag[0][0][<-]}
    \;+\;
    \tikzpic{\cro[0][1]\lili}
    \;-\;
    \tikzpic{\cro[0][1]\cc}
    \;-\;
    \tikzpic{\lili[0][1]\cro}
    \;+\;
    \tikzpic{\cc[0][1]\cro}
    \& 
    \tikzpic{\lili\lili[0][1]\bul[0][1]}
    \;-\;
    \tikzpic{\cro[0][1]\cc}
    \;+\;
    \tikzpic{\cc[0][1]\cro}
    \\
    \tikzpic{\ulcr[0][1]\cro\AntiDiag[0][0][<-]}
    \\
    \&
    \tikzpic{\lili\lili[0][1]\bul[0][1]}
    \&\&
    \tikzpic{\lili\lili[0][1]\bul[0][1]}
    \arrow[from=2-1,to=1-2, bend left=20pt]
    \arrow[from=1-2,to=1-3]
    \arrow[from=1-3,to=1-4]
    \arrow[from=2-1,to=3-2, curve={height=12pt}]
    \arrow[from=3-2,to=3-4,equals]
    \arrow[from=1-4,to=3-4, "{h_{\strongsucc}}", tail reversed, dottedcd]
  \end{tikzcd}
\end{gather*}
Here, the vector $\widehat{\bullet}$ is zero.
The branching is $\strongsucc$-tamely $\affT$-preconfluent if and only if there exists an $\affS$-congruence
\begin{gather*}
  \renewcommand{\sca}{.3}
  h_{\strongsucc}\colon\quad
  \tikzpic{\lili\lili[0][1]\bul[0][1]}
  \;-\;
  \tikzpic{\cro[0][1]\cc}
  \;+\;
  \tikzpic{\cc[0][1]\cro}
  \;\overset{*}{\longleftrightarrow}\;
  \tikzpic{\lili\lili[0][1]\bul[0][1]}
\end{gather*}
which passes only through diagrams with either at most one crossing, or with two crossings and no dot. (In this case, we can use the upward and downward kink rewriting steps.)

The second example is a branching between a cap-slide and a dot-slide:
\begin{gather*}
  \newcommand{\arOrNot}[1]{#1}
  \renewcommand{\sca}{.3}
  \begin{tikzcd}[ampersand replacement=\&,row sep=tiny]
    \&
    -\;
    \tikzpic{
      \cro \uli[0][1] \lca[1][1] \li[2][0]
      \AntiDiag[0][0][\arOrNot{<-}]
    }
    \&
    -\;
    \tikzpic{
      \dlcr \uli[0][1] \ca[1][1] \li[2][0]
      \AntiDiag[0][0][\arOrNot{<-}]
    }
    \;-\;
    \tikzpic{
      \lili \uli[0][1] \ca[1][1] \li[2][0]
    }
    \;+\;
    \tikzpic{
      \cc \uli[0][1] \ca[1][1] \li[2][0]
    }
    \&
    -\;
    \tikzpic{
      \ca[0][1] \uli[2][1]
      \bli[0][0] \cro[1][0]
    }
    \;-\;
    \tikzpic{
      \lili \uli[0][1] \ca[1][1] \li[2][0]
    }
    \;+\;
    \tikzpic{
      \cc \uli[0][1] \ca[1][1] \li[2][0]
    }
    \\
    \tikzpic{
      \cro \uli[0][1] \rca[1][1] \li[2][0]
      \AntiDiag[0][0][\arOrNot{<-}]
    }
    \\
    \&
    \tikzpic{
      \ca[0][1] \uli[2][1]
      \li[0][0][<-] \drcr[1][0]
    }
    \&
    \tikzpic{
      \ca[0][1] \uli[2][1]
      \li[0][0][<-] \cro[1][0] \bul[1][1]
    }
    \;+\;
    \tikzpic{
      \ca[0][1] \uli[2][1]
      \li[0][0] \lili[1][0]
    }
    \;-\;
    \tikzpic{
      \ca[0][1] \uli[2][1]
      \li[0][0] \cc[1][0]
    }
    \&
    -\;
    \tikzpic{
      \ca[0][1] \uli[2][1]
      \bli[0][0] \cro[1][0]
    }
    \;+\;
    \tikzpic{
      \ca[0][1] \uli[2][1]
      \li[0][0] \lili[1][0]
    }
    \;-\;
    \tikzpic{
      \ca[0][1] \uli[2][1]
      \li[0][0] \cc[1][0]
    }
    %
    \arrow[from=2-1,to=1-2, bend left=20pt]
    \arrow[from=1-2,to=1-3]
    \arrow[from=1-3,to=1-4]
    \arrow[from=2-1,to=3-2, bend right=20pt]
    \arrow[from=3-2,to=3-3]
    \arrow[from=3-3,to=3-4]
    \arrow[from=1-4,to=3-4, "{h_{\strongsucc}}", shift left=10.5, tail reversed, dottedcd]
    \arrow[from=1-4,to=3-4, shift right=14, equals]
  \end{tikzcd}
\end{gather*}
Here, the vector $\widehat{\bullet}$ is non-zero.
The branching is $\strongsucc$-tamely $\affT$-preconfluent if and only if the coefficients of the $\strongsucc$-leading term agree (they do) and if there exists an $\affS$-congruence
\begin{gather*}
  \renewcommand{\sca}{.3}
  h_{\strongsucc}\colon\qquad
  \tikzpic{
    \lili \uli[0][1] \ca[1][1] \li[2][0]
  }
  \;-\;
  \tikzpic{
    \cc \uli[0][1] \ca[1][1] \li[2][0]
  }
  \;\tamearrow
  \;-\;
  \tikzpic{
    \ca[0][1] \uli[2][1]
    \li[0][0] \lili[1][0]
  }
  \;+\;
  \tikzpic{
    \ca[0][1] \uli[2][1]
    \li[0][0] \cc[1][0]
  }
\end{gather*}
which passes only through diagrams that either have no crossing, or have one crossing and no dot. (In this case, it suffices to use the zigzag rewriting steps.)


\subsubsection{Tamed preconfluence for bubble-slide branchings}

As the examples above illustrate, tamed preconfluence boils down to the existence of a tamed congruence (and possibly, an equality between two scalars).
For certain branchings, called ``bubble-slide'' in what follows, it will be useful to spell out this tamed congruence.

\begin{notation}
  \label{not:tame_arrow}
  We denote by
  \[x\overset{k}{\tamearrow}y\]
  an $\affS$-congruence between $x$ and $y$ that is $\rhd$-tamed by a bubble with $k$ dots.
\end{notation}

\begin{lemma}
  \label{lem:bubble-slide-as-tame}
  Let $\affS$ be a rewriting system of affine Brauer type and $\trafficrule$ a choice of traffic rules.
  The branchings
  \begin{gather*}
    \renewcommand{\sca}{.4}
    \tikzpic{
      \lorli[-1.4][0] \bul[0][.5]\ca[0][.5]\custop[0][.5]
      \node[left=-1pt] at (0,.5) {\scriptsize $k$};
    }
    \qquad
    \tikzpic{
      \cro[-2][0]
      \bul[0][.5]\ca[0][.5]\custop[0][.5]
      \node[left=-1pt] at (0,.5) {\scriptsize $k$};
      \draw[very thick,red,->] (-2+.9,.6) to (-2+.5,1);
      \draw[very thick,red,->] (-2+.9,.4) to (-2+.5,0);
    }
    \qquad
    \tikzpic{
      \lorli[-1.4][0] \bul[0][.5]\ca[0][.5]\custop[0][.5]
      \node[left=-1pt] at (0,.5) {\scriptsize $k$};
      \rorli[1.8][0]
      \draw (-1.4,1) to[out=90,in=90] (1.8,1);
    }
    \qquad
    \tikzpic{
      \lorli[-1.4][0] \bul[0][.5]\ca[0][.5]\custop[0][.5]
      \node[left=-1pt] at (0,.5) {\scriptsize $k$};
      \rorli[1.8][0]
      \draw (-1.4,0) to[out=-90,in=-90] (1.8,0);
    }
  \end{gather*}
  are $\rhd$-tamely $\affT$-preconfluent if and only if $\VARbbsl^2=1$ and there exist $\affS$-congruences as follows:
  \begin{gather*}    
    \tikzpic{\li[-.5][0]\LOTev}
    -\VARbbsl\;
    \tikzpic{\li[1.5][0]\LOTev}
    \;\overset{k}{\tamearrow}\;
    \tikzpic{\LOTbbsl}
    \qquad
    \tikzpic{\ulcrbbsl}[][.5]
    \;+\;
    \VARbbsl\;\tikzpic{\urcrbbsl}[][.5]
    \;\overset{k}{\tamearrow}\;
    \tikzpic{\dlcrbbsl}[][.5]
    \;+\;
    \VARbbsl\;\tikzpic{\drcrbbsl}[][.5]
    \\[1ex]
    \VARbbsl\tikzpic{\LOTbbsl[0][0][k][left]\li[1][0]\ca[0][1]}[][.5]
    \;+\;
    \tikzpic{\li\LOTbbsl[1][0]\ca[0][1]}[][.5]
    \;\overset{k}{\tamearrow}\;
    0
    \qquad
    \VARbbsl\tikzpic{\LOTbbsl[0][0][k][left]\li[1][0]\cu}[][.5]
    \;+\;
    \tikzpic{\li\LOTbbsl[1][0]\cu}[][.5]
    \;\overset{k}{\tamearrow}\;
    0
  \end{gather*}
\end{lemma}

We leave the proof to the reader.

\subsection{Statement of the Main Theorem}

The following is the main result of this paper; \cref{sec:critical_branchings} is devoted to its proof.

\begin{restatable}[Main Theorem]{theorem}{thmAffineGrobnerSystem}
  \label{thm:affine_Grobner_system_iff_critical_branchings}
  Let $\affC$ be a category presented by a rewriting system of affine Brauer type $\affS$. Fix a choice $\trafficrule$ of traffic rules. Assume that every critical branching from \cref{fig:critical_branchings} is $\rhd$-tamely $\affT$-preconfluent.
  Then $\affC$ is a category of affine Brauer type with $\bubbleset$-bubbles (\cref{defn:category_affine_brauer_type}). Moreover, the triple $(\affS,\affT,\succ)$ defines a monoidal Gröbner system for $\affC$.
\end{restatable}

\begin{figure}[t]
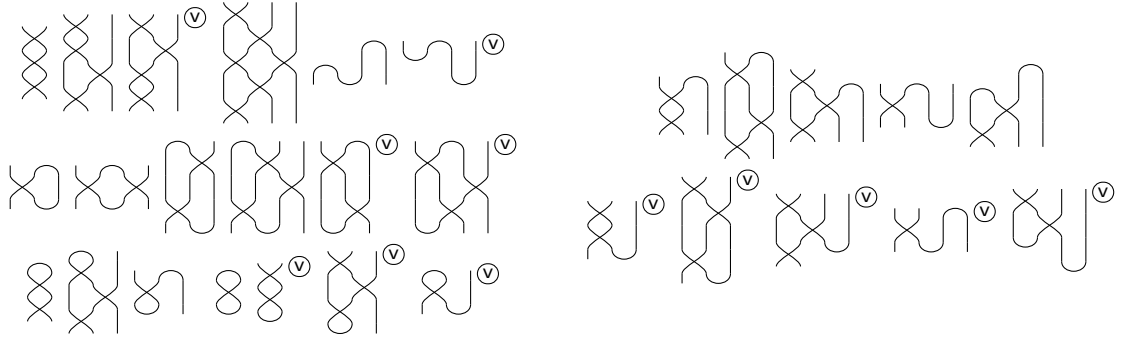
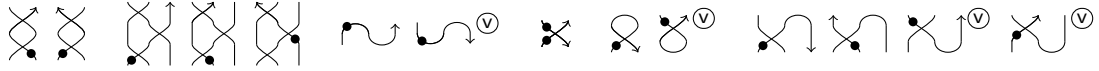
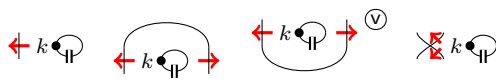
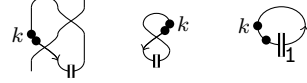

  \renewcommand{\sca}{.3}
  \newcommand{\vertsp}{.3ex}
  \def\tempsp{50mu}

  \begin{subfigure}{\textwidth}
    \renewcommand{\vertsp}{0ex}
    \renewcommand{\sca}{.32}
    \def\tempsp{10mu}
    \centering
    \begin{minipage}[c]{0.52\textwidth}
    \centering
    \begin{gather*}
      \tikzpic{
        \cro[0][2]
        \cro[0][1]
        \cro
      }
      \mspace{\tempsp}
      \tikzpic{
        \cro[0][3]\li[2][3]
        \cro[0][2]\li[2][2]
        \li[0][1]\cro[1][1]
        \cro[0][0]\li[2][0]
      }
      \mspace{\tempsp}
      \tikzpic{
        \cro[0][3]\li[2][3]
        \li[0][2]\cro[1][2]
        \cro[0][1]\li[2][1]
        \cro[0][0]\li[2][0]
      }^{\;\vertsym}
      \mspace{\tempsp}
      \tikzpic{
        \cro[0][4]\li[2][4]
        \li[0][3]\cro[1][3]
        \cro[0][2]\cro[2][2]
        \li[0][1]\cro[1][1]
        \cro\li[2][0]
        \draw[] (3,0) to (3,2);\draw[] (3,3) to (3,5);
      }
      \mspace{\tempsp}
      \tikzpic{
        \ca[2][1]
        \ca\lili[2][0]
        \dli\cu[1][0]\dli[3][0]
      }
      \mspace{\tempsp}
      \tikzpic{
        \uli[0][1]\ca[1][1]\uli[3][1]
        \cu[0][1]\lili[2][0]
        \cu[2][0]
      }^{\;\vertsym}
      \\[\vertsp]
      \tikzpic{ 
        \uli[0][1]\ca[1][1]
        \cro\li[2][0]
        \dli\cu[1][0]
      }
      \mspace{\tempsp}
      \tikzpic{ 
        \uli[0][1]\ca[1][1]
        \cro\cro[2][0]
        \dli\cu[1][0]
        \draw[] (3,-.4) to (3,0);\draw[] (3,1) to (3,1.4);
      }
      \mspace{\tempsp}
      \tikzpic{
        \ca[0][3]\uli[2][3]
        \li[0][2]\cro[1][2]
        \lili[0][1]\li[2][1]
        \cro\li[2][0]
        \dli\cu[1][0]
      }
      \mspace{\tempsp}
      \tikzpic{
        \ca[0][3]\uli[2][3]
        \li[0][2]\cro[1][2]
        \lili[0][1]\cro[2][1]
        \cro\li[2][0]
        \dli\cu[1][0]
        \draw[] (3,-.4) to (3,1);\draw[] (3,2) to (3,3.4);
      }
      \mspace{\tempsp}
      \tikzpic{
        \uli[0][3]\ca[1][3] 
        \cro[0][2]\li[2][2]
        \lili[0][1]\li[2][1]
        \li[0][0]\cro[1][0]
        \cu \dli[2][0]
      }^{\;\vertsym}
      \mspace{\tempsp}
      \tikzpic{
        \uli[0][3]\ca[1][3] 
        \cro[0][2]\li[2][2]
        \lili[0][1]\cro[2][1]
        \li[0][0]\cro[1][0]
        \cu \dli[2][0]
        \draw[] (3,-.4) to (3,1);\draw[] (3,2) to (3,3.4);
      }^{\;\vertsym}
      \\[\vertsp]
      \tikzpic{
        \ca[0][2]
    		\cro[0][1]
    		\cro[0][0]
    	}
      \mspace{\tempsp}
    	\tikzpic{
        \ca[0][3]\uli[2][3]
    		\cro[0][2]\li[2][2]
    		\cro[1][1]\li[0][1]
    		\cro[0][0]\li[2][0]
    	}
      \mspace{\tempsp}
    	\tikzpic{
    		\ca[1][1]\uli[0][1]
    		\cro[0][0]\li[2][0]
    		\cu\dli[2][0]
    	}
      \mspace{\tempsp}
      \mspace{\tempsp}
    	\tikzpic{
    		\ca[0][1]
    		\cro[0][0]
    		\cu\
    	}
      \mspace{\tempsp}
      \tikzpic{
        \cro[0][1]
    		\cro[0][0]
        \cu[0][0]
    	}^{\;\vertsym}
      \mspace{\tempsp}
    	\tikzpic{
        \cro[0][2]\li[2][2]
    		\cro[1][1]\li[0][1]
    		\cro[0][0]\li[2][0]
        \cu[0][0]\dli[2][0]
    	}^{\;\vertsym}
      \mspace{\tempsp}
    	\tikzpic{
    		\ca[0][1]\uli[2][1]
    		\cro[0][0]\li[2][0]
    		\dli \cu[1][0]
    	}^{\;\vertsym}
    \end{gather*}
    \end{minipage}%
    \hfill
    \begin{minipage}[c]{0.46\textwidth}
    \centering
    \begin{gather*}
      \tikzpic{
    		\uli[0][2]\ca[1][2]
        \cro[0][1]\li[2][1]
        \cro\li[2][0]
    	}
      \mspace{\tempsp}
      \tikzpic{
    		\uli[0][4]\ca[1][4]
        \cro[0][3]\li[2][3]
        \lili[0][2]\li[2][2]
        \li[0][1]\cro[1][1]
        \cro\li[2][0]
    	}
      \mspace{\tempsp}
      \tikzpic{
        \cro[0][2]\ca[2][2]
        \li[0][1]\cro[1][1]\li[3][1]
        \cro\lili[2][0]
      }
      \mspace{\tempsp}
      \tikzpic{
    		\uli[0][1]\ca[1][1]\uli[3][1]
        \cro\lili[2][0]
        \dli[0][0]\dli[1][0]\cu[2][0]
    	}
      \mspace{\tempsp}
      \tikzpic{ 
        \ca[2][3]
        \ca[0][2]\lili[2][2]
        \li[0][1]\cro[1][1]\li[3][1]
        \cro\lili[2][0]
      }
      \\[\vertsp]
      \tikzpic{
        \cro[0][1]\li[2][1]
        \cro\li[2][0]
    		\dli[0][0]\cu[1][0]
    	}^{\;\vertsym}
      \mspace{\tempsp}
      \tikzpic{
        \cro[0][3]\li[2][3]
        \li[0][2]\cro[1][2]
        \lili[0][1]\li[2][1]
        \cro\li[2][0]
    		\dli\cu[1][0]
    	}^{\;\vertsym}
      \mspace{\tempsp}
      \tikzpic{
        \cro[0][2]\lili[2][2]
        \li[0][1]\cro[1][1]\li[3][1]
        \cro\cu[2][1]
      }^{\;\vertsym}
      \mspace{\tempsp}
      \tikzpic{
        \uli[0][1]\uli[1][1]\ca[2][1]
        \cro\lili[2][0]
    		\dli[0][0]\cu[1][0]\dli[3][0]
    	}^{\;\vertsym}
      \mspace{\tempsp}
      \tikzpic{ 
        \cro[0][2]\lili[2][2]
        \li[0][1]\cro[1][1]\li[3][1]
        \cu[0][1]\lili[2][0]
        \cu[2][0]
      }^{\;\vertsym}
    \end{gather*}
    \end{minipage}
    \caption{Non-affine critical branchings.}
    \label{subfig:critical_branchings_non_affine}
  \end{subfigure}

  \begin{subfigure}{\textwidth}
    \renewcommand{\sca}{.35}
    \newcommand{\arOrNot}[1]{#1}
    \def\tempsp{13mu}
    \begin{gather*}
      \tikzpic{
        \drcr \cro[0][1]
        \AntiDiag[0][1][\arOrNot{->}]
      }
      \mspace{\tempsp}
      \tikzpic{
        \dlcr \cro[0][1]
        \Diag[0][1][\arOrNot{->}]
      }
      \mspace{\tempsp}\mspace{\tempsp}
      \tikzpic{
        \dlcr \li[2][0] \li[0][1] \cro[1][1] \cro[0][2] \li[2][2][\arOrNot{->}]
      }[scale=.8]
      \mspace{\tempsp}
      \tikzpic{
        \drcr \li[2][0] \li[0][1] \cro[1][1] \cro[0][2] \li[2][2]
        \AntiDiag[0][2][\arOrNot{->}]
      }[scale=.8]
      \mspace{\tempsp}
      \tikzpic{
        \cro \li[2][0] \li[0][1] \drcr[1][1] \cro[0][2] \li[2][2]
        \Diag[0][2][\arOrNot{->}]
      }[scale=.8]
      \mspace{\tempsp} 
      \mspace{\tempsp} 
      \tikzpic{
        \dli \lca[0][0]\ca[0][0] \cu[1][0] \uli[2][0][\arOrNot{->}]
      }
      \mspace{\tempsp} 
      \tikzpic{
        \uli \lcu[0][0]\cu[0][0] \ca[1][0] \dli[2][0][\arOrNot{->}]
      }^{\;\vertsym}
      \mspace{\tempsp}\mspace{\tempsp}
      \tikzpic{
        \ulcr \dlcr
        \Diag[0][0][\arOrNot{<-}]
        \AntiDiag[0][0][\arOrNot{->}]
      }
      \mspace{\tempsp}
      \mspace{\tempsp}
      \tikzpic{
        \dlcr \ca[0][1]
        \Diag[0][0][\arOrNot{<-}]
      }
      \mspace{\tempsp}
      \tikzpic{
        \ulcr
        \cu
        \AntiDiag[0][0][\arOrNot{->}]
      }^{\;\vertsym}
      \mspace{\tempsp}
      \mspace{\tempsp} 
      \tikzpic{
        \dlcr \uli[0][1] \ca[1][1] \li[2][0][\arOrNot{<-}]
      } 
      \mspace{\tempsp}
      \tikzpic{
        \drcr \uli[0][1][\arOrNot{->}] \ca[1][1] \li[2][0]
      } 
      \mspace{\tempsp} 
      \tikzpic{ 
        \ulcr \dli \cu[1][0] \li[2][0][\arOrNot{->}] 
      }^{\;\vertsym}
      \mspace{\tempsp}
      \tikzpic{ 
        \dlcr \dli \cu[1][0] \li[2][0]
        \AntiDiag[0][0][\arOrNot{->}]
      }^{\;\vertsym}
      \mspace{\tempsp}
    \end{gather*}
    \caption{Dot-slide critical branchings.}
    \label{subfig:critical_branchings_dot_slides}
  \end{subfigure}

  \begin{subfigure}[b]{0.5\textwidth}
    \centering
    \begin{gather*}
      \renewcommand{\sca}{.35}
      \tikzpic{
        \lorli[-1.4][0] \bul[0][.5]\ca[0][.5]\custop[0][.5]
        \node[left=-1pt] at (0,.5) {\scriptsize $k$};
      }
      \quad
      \tikzpic{
        \lorli[-1.4][0] \bul[0][.5]\ca[0][.5]\custop[0][.5]
        \node[left=-1pt] at (0,.5) {\scriptsize $k$};
        \rorli[1.8][0]
        \draw (-1.4,1) to[out=90,in=90] (1.8,1);
      }
      \quad
      \tikzpic{
        \lorli[-1.4][0] \bul[0][.5]\ca[0][.5]\custop[0][.5]
        \node[left=-1pt] at (0,.5) {\scriptsize $k$};
        \rorli[1.8][0]
        \draw (-1.4,0) to[out=-90,in=-90] (1.8,0);
      }^{\;\vertsym}
      \quad
      \tikzpic{
        \cro[-2][0]
        \bul[0][.5]\ca[0][.5]\custop[0][.5]
        \node[left=-1pt] at (0,.5) {\scriptsize $k$};
        \draw[very thick,red,->] (-2+.9,.6) to (-2+.5,1);
        \draw[very thick,red,->] (-2+.9,.4) to (-2+.5,0);
      }
    \end{gather*}
    \caption{Bubble-slide critical branchings.}
    \label{subfig:critical_branchings_bubbles}
  \end{subfigure}%
  \begin{subfigure}[b]{0.5\textwidth}
    \centering
    \begin{gather*}
      \renewcommand{\sca}{.35}
      \tikzpic{
        \ca[0][2]\uli[2][2]
        \li[0][1]\cro[1][1]
        \ulcr \Diag[0][0][<-] \li[2][0]
        \dli \custop[1][0]
        \bul[0][1]\node[left=-1pt] at (0,1) {\scriptsize $k$};
      }
      \qquad
      \tikzpic{
        \ca[0][1]
        \urcr
        \custop
        \AntiDiag[0][0][<-]
        \bul[1][1]\node[right=-1pt] at (1,1) {\scriptsize $k$};
      }
      \quad
      \tikzpic{\lcustop[0][0][->]\ca\node[below right=-4pt] at (.5,-.4) {\footnotesize$\mathtt{1}$};\bul\node[left=-1pt] at (0,0) {\scriptsize $k$};}[scale=1.8]
    \end{gather*}
    \caption{Exceptional critical branchings.}
    \label{subfig:critical_branchings_exceptional}
  \end{subfigure}

  \caption{Critical branchings of affine Brauer type. The symbol $\vertsym$ indicates vertical symmetries: these critical branchings can be ignored under certain conditions.%
  }
  \label{fig:critical_branchings}
\end{figure}

Each diagram in \cref{fig:critical_branchings} encodes a critical branching, as exactly two rewriting steps can be applied to the diagram.
The third exceptional critical branching is the sole exception: as we explained in the introduction around \cref{eq:one-free-stop-sign}, it encodes the rewriting step consisting in sliding the dot on the bottom left through the cup, from left to right.
After sliding the dot, the bubble carries a stop sign, as prescribed by the traffic rules.
While a rewriting step is not a branching, it still makes sense to ask that it is tamely preconfluent.
We abuse terminology by calling it a ``critical branching''.

\medbreak

In many examples, one has an isomorphism between the category $\affC$ and its opposite category $\affC^{\op}$.
This vertical symmetry can be leveraged at the level of the rewriting system $(\affS,\affT,\succ)$ to reduce the number of critical branchings:

\begin{proposition}
  \label{prop:vertical_symmetries}
  Assume that $(\affS,\affT,\succ)$ has vertical symmetry in the sense of \cref{subsec:symmetries}.
  Then \cref{thm:affine_Grobner_system_iff_critical_branchings} still holds if we ignore every critical branching labelled $\vertsym$.
\end{proposition}

We prove the proposition in \cref{subsec:symmetries}.

%% file: sections/confluence_analysis_general.tex
\section{Critical branchings of affine Brauer type}
\label{sec:critical_branchings}

This section is devoted to the proof of \cref{thm:affine_Grobner_system_iff_critical_branchings} and \cref{prop:vertical_symmetries}, building on \cref{thm:RW_monoidal_grobner_system}.
The main difficulty is to show that $\affT^+$ is confluent.
The section is organized as follows:
\begin{itemize}
  \item \Cref{subsubsec:enumeration_overlapping_brauer_type} enumerates minimal overlapping branchings. It also gives some terminology for the branchings we shall encounter throughout the proof.
  \item \Cref{subsec:confluence_preliminaries} gathers various preliminary considerations, including the proof that $\succ$ is compatible with $\affT$ (\cref{lem:order_strongly_compatible}), tamed congruence of independent branchings (\cref{subsubsec:extended_MLRS_independent_branchings}), contextualization (\cref{subsubsec:confluence_contextualization}) and the fact that $\affS$ and $\affT$ present the same category (\cref{lem:affine_Brauer_equivalence_of_presentations}).
  The necessity of the third exceptional branching will come out naturally in the proof of this latter statement.
  \item \Cref{subsec:local_analysis} contains the ``local analysis'', that is, how the tamed preconfluence of many minimal branchings follows from the tamed preconfluence of critical branchings.
  The main idea is to leverage various symmetries, namely horizontal symmetries, symmetries of traffic rules, and vertical symmetries.
  In particular, we prove \cref{prop:vertical_symmetries} in this subsection.
  \item \Cref{subsec:global_analysis} contains the ``global analysis'', that is, how the tamed preconfluence of a minimal branching implies the tamed congruence of a contextualized branching.
  The main difficulty comes from understanding how traffic rules behave with respect to context.
  The necessity of the first two exceptional branchings (\cref{subfig:critical_branchings_exceptional}) will come out naturally from this global analysis.
  Here, transitivity of tamed congruence plays an important role.
  \item Finally, \cref{subsec:proof_main_thm} gathers all previous results and gives the proof of \cref{thm:affine_Grobner_system_iff_critical_branchings}.
\end{itemize}

\subsection{Enumeration of overlapping branchings}
\label{subsubsec:enumeration_overlapping_brauer_type}

\begin{proposition}
  \label{prop:classification_overlapping_branchings}
  \Cref{fig:all_overlapping_branchings} gives a complete list of minimal overlapping $\affS$-branchings.
\end{proposition}

Branchings in \cref{subfig:all_overlapping_branchings_non_affine} and \cref{subfig:all_overlapping_branchings_affine} are called \emph{minimal non-affine branchings} and \emph{minimal dot-slide branchings}, respectively; we reserve the terminology ``minimal bubble-slide branchings'' for the branchings in \cref{fig:minimal_bubble_slide} (in \cref{subsubsec:twin_branchings} below).
Contextualizations of the branchings in \cref{subfig:all_overlapping_branchings_non_affine}, \cref{subfig:all_overlapping_branchings_affine} and \cref{subfig:all_overlapping_branchings_bubbles} are called \defnemph{non-affine branchings}, \defnemph{dot-slide branchings} and \defnemph{bubble-slide branchings}, respectively.
The non-affine branchings that depend on the choice of an extra diagram, depicted with a box marked ``?'', are called \defnemph{indexed branchings}.

\begin{proof}
  Let $(f,g)=\Gamma[\ssf,\ssg]$ be a minimal branching, and denote by $D$ its source. Assume that $s(\ssf)$ and $s(\ssg)$, as diagrams in $D$, have generators in common.
  Going over all the possibilities gives all the branchings in \cref{subfig:all_overlapping_branchings_non_affine,subfig:all_overlapping_branchings_affine}, as well as the first, second and fourth (family of) branchings in \cref{subfig:all_overlapping_branchings_bubbles}.

  Assume instead that $s(\ssf)$ and $s(\ssg)$, as diagrams in $D$, do not share any generator.
  Following \cref{lem:RW_overlapping_branchings_intersect}, at least one of the two must be disconnected, that is, at least one of the two (say $\ssg$) must be a bubble-slide rewriting step.
  Then we can always find a sequence of interchanges such that $D$ remains of the form $D=\Gamma_f[s(\ssf)]$ for some context $\Gamma_f$, while also being of the form $D=\Gamma_g[s(\ssg)]$ for some other context $\Gamma_g$.
  The branching $(f,g)$ is then readily seen to be either of the form the third schematic in \cref{subfig:all_overlapping_branchings_bubbles} or an independent branching.
\end{proof}

\begin{figure}
  \renewcommand{\sca}{.32}
  \newcommand{\vertsp}{.3ex}
  \newcommand{\arOrNot}[1]{#1}

  \begin{subfigure}{\textwidth}
    \renewcommand{\vertsp}{0ex}
    \begin{gather*}
      \tikzpic{
        \cro[0][2]
        \cro[0][1]
        \cro
      }
      \qquad
      \tikzpic{
        \cro[0][3]\li[2][3]
        \cro[0][2]\li[2][2]
        \li[0][1]\cro[1][1]
        \cro[0][0]\li[2][0]
      }
      \qquad
      \tikzpic{
        \cro[0][3]\li[2][3]
        \li[0][2]\cro[1][2]
        \cro[0][1]\li[2][1]
        \cro[0][0]\li[2][0]
      }
      \qquad
      \tikzpic{
        \cro[0][4]\li[2][4]
        \li[0][3]\cro[1][3]
        \cro[0][2]\indexNF[2][2][?]
        \li[0][1]\cro[1][1]
        \cro\li[2][0]
        %
      }
      \qquad
      \tikzpic{
        \ca[2][1]
        \ca\lili[2][0]
        \dli\cu[1][0]\dli[3][0]
      }
      \qquad
      \tikzpic{
        \uli[0][1]\ca[1][1]\uli[3][1]
        \cu[0][1]\lili[2][0]
        \cu[2][0]
      }
      \\[\vertsp]
      \tikzpic{
    		\uli[0][2]\ca[1][2]
        \cro[0][1]\li[2][1]
        \cro\li[2][0]
    	}
      \qquad
      \tikzpic{
    		\uli[0][4]\ca[1][4]
        \cro[0][3]\li[2][3]
        \lili[0][2]\indexNF[2][2][?]
        \li[0][1]\cro[1][1]
        \cro\li[2][0]
    	}
      \qquad
      \tikzpic{
        \cro[0][2]\ca[2][2]
        \li[0][1]\cro[1][1]\li[3][1]
        \cro\lili[2][0]
      }
      \qquad
      \tikzpic{
    		\uli[0][1]\ca[1][1]\uli[3][1]
        \cro\lili[2][0]
        \dli[0][0]\dli[1][0]\cu[2][0]
    	}
      \qquad
      \tikzpic{ 
        \ca[0][3]\uli[2][3]
        \li[0][2]\cro[1][2]
        \cro[0][1]\li[2][1]
        \cro\li[2][0]
      }
      \qquad
      \tikzpic{
        \ca[0][5]\uli[2][5]
        \li[0][4]\cro[1][4]
        \cro[0][3]\li[2][3]
        \lili[0][2]\indexNF[2][2][?]
        \li[0][1]\cro[1][1]
        \cro\li[2][0]
      }
      \qquad
      \tikzpic{ 
        \ca[2][3]
        \ca[0][2]\lili[2][2]
        \li[0][1]\cro[1][1]\li[3][1]
        \cro\lili[2][0]
      }
      \\[\vertsp]
      \tikzpic{
        \cro[0][1]\li[2][1]
        \cro\li[2][0]
    		\dli[0][0]\cu[1][0]
    	}
      \qquad
      \tikzpic{
        \cro[0][3]\li[2][3]
        \li[0][2]\cro[1][2]
        \lili[0][1]\indexNF[2][1][?]
        \cro\li[2][0]
    		\dli\cu[1][0]
    	}
      \qquad
      \tikzpic{
        \cro[0][2]\lili[2][2]
        \li[0][1]\cro[1][1]\li[3][1]
        \cro\cu[2][1]
      }
      \qquad
      \tikzpic{
        \uli[0][1]\uli[1][1]\ca[2][1]
        \cro\lili[2][0]
    		\dli[0][0]\cu[1][0]\dli[3][0]
    	}
      \qquad
      \tikzpic{ 
        \cro[0][2]\li[2][2]
        \cro[0][1]\li[2][1]
        \li\cro[1][0]
        \cu\dli[2][0]
      }
      \qquad
      \tikzpic{
        \cro[0][4]\li[2][4]
        \li[0][3]\cro[1][3]
        \lili[0][2]\indexNF[2][2][?]
        \cro[0][1]\li[2][1]
        \li[0][0]\cro[1][0]
        \cu\dli[2][0]
      }
      \qquad
      \tikzpic{ 
        \cro[0][2]\lili[2][2]
        \li[0][1]\cro[1][1]\li[3][1]
        \cu[0][1]\lili[2][0]
        \cu[2][0]
      }
      \\[\vertsp]
      \tikzpic{ 
        \ca[0][2]\uli[2][2]\uli[3][2]
        \li[0][1]\cro[1][1]\li[3][1]
        \cro\cu[2][1]
      }
      \qquad
      \tikzpic{ 
        \uli[0][1]\ca[1][1]
        \cro\indexNF[2][0][?]
        \dli\cu[1][0]
      }
      \qquad
      \tikzpic{
        \ca[0][3]\uli[2][3]
        \li[0][2]\cro[1][2]
        \lili[0][1]\indexNF[2][1][?]
        \cro\li[2][0]
        \dli\cu[1][0]
      }
      \qquad
      \tikzpic{
        \ca[0][4]\uli[2][4]
        \li[0][3]\cro[1][3] 
        \lili[0][2]\indexNF[2][2][?]
        \cro[0][1]\li[2][1] 
        \li[0][0]\cro[1][0]
        \dli[2][0]\cu
      }
      \qquad\qquad
      \tikzpic{ 
        \cro[0][1]\ca[2][1]
        \li[0][0]\cro[1][0]\li[3][0]
        \cu[0][0]\dli[2][0]\dli[3][0]
      }
      \qquad
      \tikzpic{
        \uli[0][3]\ca[1][3] 
        \cro[0][2]\li[2][2]
        \lili[0][1]\indexNF[2][1][?]
        \li[0][0]\cro[1][0]
        \cu \dli[2][0]
      }
      \\[\vertsp]
      \def\tempsp{20mu}
      \tikzpic{
        \ca[0][2]
    		\cro[0][1]
    		\cro[0][0]
    	}
      \mspace{\tempsp}
    	\tikzpic{
        \ca[0][3]\uli[2][3]
    		\cro[0][2]\li[2][2]
    		\cro[1][1]\li[0][1]
    		\cro[0][0]\li[2][0]
    	}
      \mspace{\tempsp}
    	\tikzpic{
    		\ca[1][1]\uli[0][1]
    		\cro[0][0]\li[2][0]
    		\cu\dli[2][0]
    	}
      \mspace{\tempsp}
    	\tikzpic{
    		\ca[0][2]\uli[2][2]
    		\cro[0][0]\li[2][0]
    		\li[0][1]\cro[1][1]
    		\cu\dli[2][0]
    	}
      \mspace{\tempsp}
    	\tikzpic{
    		\ca[0][1]
    		\cro[0][0]
    		\cu\
    	}
      \mspace{\tempsp}
    	\tikzpic{
        \ca[0][2]
    		\li[0][1]\li[1][1]\ca[2][1]
    		\li\cro[1][0]\li[3][0]
    		\cu\dli[2][0]\dli[3][0]
    	}
      \qquad
      \tikzpic{
        \cro[0][1]
    		\cro[0][0]
        \cu[0][0]
    	}
      \mspace{\tempsp}
    	\tikzpic{
        \cro[0][2]\li[2][2]
    		\cro[1][1]\li[0][1]
    		\cro[0][0]\li[2][0]
        \cu[0][0]\dli[2][0]
    	}
      \mspace{\tempsp}
    	\tikzpic{
    		\ca[0][1]\uli[2][1]
    		\cro[0][0]\li[2][0]
    		\dli \cu[1][0]
    	}
      \mspace{\tempsp}
    	\tikzpic{
    		\ca[0][2]\uli[2][2]
    		\li[0][0]\cro[1][0]
    		\cro[0][1]\li[2][1]
    		\cu\dli[2][0]
    	}
      \mspace{\tempsp}
    	\tikzpic{
        \ca[0][2]\uli[2][2]\uli[3][2]
    		\li[0][1]\cro[1][1]\li[3][1]
    		\lili[0][0]\cu[2][1]
    		\cu
    	}
    \end{gather*}
    \caption{Minimal overlapping branchings involving non-affine rewriting steps; we call them \emph{minimal non-affine branchings} and contextualizations thereof \defnemph{non-affine branchings}. A box marked ``?'' denotes any diagram with non-empty input and output; branchings with such a box are called \emph{indexed branchings}.}
    \label{subfig:all_overlapping_branchings_non_affine}
  \end{subfigure}

  \begin{subfigure}{\textwidth}
    \begin{gather*}
      \tikzpic{
        \drcr \cro[0][1]
        \AntiDiag[0][1][\arOrNot{->}]
      }
      \qquad
      \tikzpic{
        \dlcr \cro[0][1]
        \Diag[0][1][\arOrNot{->}]
      }
      \qquad
      \tikzpic{
        \cro \ulcr[0][1]
        \AntiDiag[0][0][\arOrNot{<-}]
      }
      \qquad
      \tikzpic{
        \cro \urcr[0][1]
        \Diag[0][0][\arOrNot{<-}]
      }
      \qquad\qquad
      \tikzpic{
        \dlcr \li[2][0] \li[0][1] \cro[1][1] \cro[0][2] \li[2][2][\arOrNot{->}]
        \AntiDiag[0][0]
      }
      \qquad
      \tikzpic{
        \drcr \li[2][0] \li[0][1] \cro[1][1] \cro[0][2] \li[2][2]
        \AntiDiag[0][2][\arOrNot{->}]
      }
      \qquad
      \tikzpic{
        \cro \li[2][0] \li[0][1] \drcr[1][1] \cro[0][2] \li[2][2]
        \Diag[0][2][\arOrNot{->}]
      }
      \qquad
      \tikzpic{
        \cro \li[2][0] \li[0][1] \urcr[1][1] \cro[0][2] \li[2][2]
        \AntiDiag[0][0][\arOrNot{<-}]
      }
      \qquad
      \tikzpic{
        \cro \li[2][0] \li[0][1] \cro[1][1] \urcr[0][2] \li[2][2]
        \Diag[0][0][\arOrNot{<-}]
      }
      \qquad
      \tikzpic{
        \cro \li[2][0][\arOrNot{<-}] \li[0][1] \cro[1][1] \ulcr[0][2] \li[2][2]
        \Diag[0][2]
      }
      \\[\vertsp]
      \tikzpic{
        \dli \lca[0][0]\ca[0][0] \cu[1][0] \uli[2][0][\arOrNot{->}]
      }
      \qquad
      \tikzpic{
        \dli[0][0][\arOrNot{->}] \ca \rcu[1][0]\cu[1][0] \uli[2][0]
      }
      \qquad
      \tikzpic{
        \uli \lcu[0][0]\cu[0][0] \ca[1][0] \dli[2][0][\arOrNot{->}]
      }
      \qquad
      \tikzpic{
        \uli[0][0][\arOrNot{->}] \cu \rca[1][0]\ca[1][0] \dli[2][0]
      }
      \qquad
      \qquad
      \tikzpic{
        \ulcr \dlcr
        \Diag[0][0][\arOrNot{<-}]
        \AntiDiag[0][0][\arOrNot{->}]
      }
      \qquad
      \tikzpic{
        \ulcr \urcr
        \Diag[0][0][\arOrNot{<-}]
        \AntiDiag[0][0][\arOrNot{<-}]
      }
      \qquad
      \tikzpic{
        \urcr \drcr
        \Diag[0][0][\arOrNot{->}]
        \AntiDiag[0][0][\arOrNot{<-}]
      }
      \qquad
      \tikzpic{
        \dlcr \drcr
        \Diag[0][0][\arOrNot{->}]
        \AntiDiag[0][0][\arOrNot{->}]
      }
      \\[\vertsp]
      \tikzpic{
        \dlcr \uli[0][1] \ca[1][1] \li[2][0][\arOrNot{<-}]
        \AntiDiag[0][0]
      } 
      \qquad
      \tikzpic{
        \drcr \uli[0][1][\arOrNot{->}] \ca[1][1] \li[2][0]
        \Diag[0][0]
      } 
      \qquad
      \tikzpic{
        \ulcr \uli[0][1] \ca[1][1] \li[2][0]
        \Diag[0][0][\arOrNot{<-}]
      } 
      \qquad
      \tikzpic{
        \cro \uli[0][1] \rca[1][1]\ca[1][1] \li[2][0]
        \AntiDiag[0][0][\arOrNot{<-}]
      }
      \qquad\qquad
      \tikzpic{ 
        \ca \urcr[1][-1] \li[0][-1] \li[2][-2] \cro[0][-2] \uli[2][0]
        \AntiDiag[0][-2][\arOrNot{<-}]
      } 
      \qquad
      \tikzpic{ 
        \ca \drcr[1][-1] \li[0][-1] \li[2][-2] \cro[0][-2] \uli[2][0]
        \Diag[0][-2][\arOrNot{<-}]
      } 
      \qquad
      \tikzpic{ 
        \ca \cro[1][-1] \li[0][-1] \li[2][-2] \dlcr[0][-2] \uli[2][0][\arOrNot{->}]
        \AntiDiag[0][-2]
      } 
      \qquad
      \tikzpic{ 
        \ca \cro[1][-1] \li[0][-1] \li[2][-2][\arOrNot{<-}] \drcr[0][-2] \uli[2][0]
        \Diag[0][-2]
      } 
      \\[\vertsp]
      \tikzpic{ 
        \ulcr \dli \cu[1][0] \li[2][0][\arOrNot{->}]
        \Diag[0][0]
      }
      \qquad
      \tikzpic{ 
        \dlcr \dli \cu[1][0] \li[2][0]
        \AntiDiag[0][0][\arOrNot{->}]
      } 
      \qquad
      \tikzpic{ 
        \urcr \dli[0][0][\arOrNot{->}] \cu[1][0] \li[2][0]
        \AntiDiag[0][0]
      }
      \qquad
      \tikzpic{ 
        \cro \dli \rcu[1][0]\cu[1][0] \li[2][0]
        \Diag[0][0][\arOrNot{->}]
      }
      \qquad\qquad
      \tikzpic{
        \dli[2][0] \cu \urcr[1][0] \li \li[2][1] \cro[0][1]
        \AntiDiag[0][1][\arOrNot{->}]
      } 
      \qquad
      \tikzpic{
        \dli[2][0] \cu \drcr[1][0] \li \li[2][1] \cro[0][1]
        \Diag[0][1][\arOrNot{->}]
      }
      \qquad
      \tikzpic{
        \dli[2][0][\arOrNot{->}] \cu \cro[1][0] \li \li[2][1] \ulcr[0][1]
        \Diag[0][1]
      }
      \qquad
      \tikzpic{
        \dli[2][0] \cu \cro[1][0] \li \li[2][1][\arOrNot{->}] \urcr[0][1]
        \AntiDiag[0][1]
      }
      \\[\vertsp]
      \tikzpic{
        \dlcr \ca[0][1]
        \AntiDiag[0][0]
        \Diag[0][0][\arOrNot{<-}]
      }
      \qquad
      \tikzpic{
        \drcr \ca[0][1]
        \Diag[0][0]
        \AntiDiag[0][0][\arOrNot{<-}]
      }
      \qquad
      \tikzpic{
        \ca[0][1]\uli[2][1]
        \li[0][0] \urcr[1][0]
        \cu\dli[2][0][\arOrNot{->}]
        \AntiDiag[1][0]
      }
      \qquad
      \tikzpic{
        \ca[0][1]\uli[2][1][\arOrNot{->}]
        \li[0][0] \drcr[1][0]
        \cu\dli[2][0]
        \Diag[1][0]
      }
      \qquad
      \tikzpic{
        \ulcr
        \cu
        \Diag[0][0]
        \AntiDiag[0][0][\arOrNot{->}]
      }
      \qquad 
      \tikzpic{
        \urcr
        \cu
        \AntiDiag[0][0]
        \Diag[0][0][\arOrNot{->}]
      }
    \end{gather*}
    \caption{Minimal overlapping branchings involving dot slides; we call them minimal dot-slide branchings and contextualizations thereof \defnemph{dot-slide branchings}.}
    \label{subfig:all_overlapping_branchings_affine}
  \end{subfigure}

  \begin{subfigure}{\textwidth}
    \begin{gather*}
      \tikzpic{
        \lorli[-1.4][0] 
        \bul[0][.5]\ca[0][.5]\cu[0][.5]
        \node[left=-1pt] at (0,.5) {\scriptsize $k$};
      }
      \qquad
      \tikzpic{
        \rorli[1.6][0] 
        \bul[0][.5]\ca[0][.5]\cu[0][.5]
        \node[left=-1pt] at (0,.5) {\scriptsize $k$};
      }
      \qquad\qquad
      \tikzpic{
        \clip (0,0) circle (2cm);
        \draw[dashed] (0,0) circle (2cm);
        \draw[out=45,in=-90] (-2,-2) to (0,2);
        \draw (-.8,2) to[out=-90,in=180] (0,1) to[out=0,in=-90] (.8,2);        
        \draw (-.5,-2) to[out=90,in=0] (-2,-.5);        
        \begin{scope}[shift={(1,-.4)}]
          \bul[-.5][0]\ca[-.5][0]\cu[-.5][0]
          \node[left=-1pt] at (-.5,0) {\scriptsize $k$};
        \end{scope}
        \node[rotate=60] at (-.8,.4) {\ldots};
      \draw[very thick,red,->] (0,0) to (-.7,.4);
      }
      \qquad\qquad
      \tikzpic{
        \clip (0,0) circle (2cm);
        \draw[dashed] (0,0) circle (2cm);
        \draw (-2,-.5) to[out=30,in=60] (-.5,-2);
        \draw[very thick,red,->] (-.5,-.5) to (-1.2,-1.2);
        \draw (-1,1.7) to[out=-90,in=-90] (1,1.7);
        \draw[very thick,red,->] (0,.7) to (0,1.7);
        \bul[-.5][0]\ca[-.5][0]\cu[-.5][0]
        \node[left=-1pt] at (-.5,0) {\scriptsize $k$};
      }
    \end{gather*}
    \caption{Minimal overlapping branchings involving bubble steps; we call contextualizations thereof \defnemph{bubble-slide branchings} (we use ``minimal bubble-slide branchings'' for the branchings in \cref{fig:minimal_bubble_slide}). The first two branchings are overlaps between bubble evaluation and bubble slide. The third picture is a schematic for overlaps between a rewriting step with an inner region and a bubble slide. The last picture is a schematic for overlaps between two bubble slides.
    See \cref{subsubsec:bubble_slides} for details.
    }
    \label{subfig:all_overlapping_branchings_bubbles}
  \end{subfigure}

  \caption{Minimal overlapping $\affS$-branchings.}
  \label{fig:all_overlapping_branchings}
\end{figure}

\subsection{Preliminaries}
\label{subsec:confluence_preliminaries}


\subsubsection{Extended linear monoidal rewriting systems}
\label{subsec:extended_MLRS}

Following \cref{subsec:defn_rewriting_affine_brauer}, we have extra data attached to a rewriting system of affine Brauer type $(\affS,\affT,\succ)$, namely:
\begin{itemize}
  \item a splitting $\affR=\affR_\agn\sqcup\affR_{\dpt}$, where
  \[\affR_\agn=\affR_\nonaff\quad\an\quad\affR_\dpt=\affR_\aff,\]
  respectively context-agnostic and context-dependent;
  \item a context-agnostic order $\strongsucc$, stronger than $\succ$;
  \item a context-agnostic equivalence relation $\asymp$.
\end{itemize}
Moreover, we have that for each $\ssr\in\affR_\dpt$:
\begin{enumerate}[(i)]
  \item for each element $x\in\supp(t(\ssr))$, either $s(\ssr)\asymp x$ or $s(\ssr)\strongsucc x$;
  \item whether $\Gamma[\ssr]\in\affT_\dpt$ only depends on $\ssr$ and on the $\asymp$-equivalence class of $\Gamma[s(\ssr)]$.
\end{enumerate}
In our setting, the combination of these properties can be understood as ``up to $\strongsucc$-lower-order terms, context-dependent rewriting steps do not affect traffic rules''.
We will often use that property in the proofs below.
In general, we call any data $(\sS,\lT,\sR_\agn,\sR_\dpt;\succ,\strongsucc,\asymp)$ as above an \defnemph{extended linear monoidal rewriting system} (or extended \LMRS{}).

\subsubsection{Compatible order}

We show that:

\begin{lemma}
  \label{lem:order_strongly_compatible}
  Let $(\affS,\affT,\succ)$ be a rewriting system of affine Brauer type and $\trafficrule$ a choice of traffic rules.
  The order $\succ$ is compatible with $\affT$.
\end{lemma}

\begin{proof}
  The $\strongsucc$-lower-order terms are compatible, since $\strongsucc$ is context-agnostic and stronger than $\succ$.
  The fact that $\rhd$-lower-order terms are compatible follows from the first three cases in \cref{lem:core_lemma_contextualization_bubbles} below.
  Finally, we consider $\strongsucc$-leading and $\rhd$-leading terms:
  \begin{itemize}[]
    \setlength\itemsep{0em}
    \item \emph{braiding R3}: strictly reduces $\distcro$, while preserving $(\nbcro,\nbdot,\disttunnel,\nbcc,\distcc)$.
    \item \emph{zigzags}: strictly reduces $\nbcc$, while preserving $(\nbcro,\nbdot,\disttunnel)$.
    \item \emph{cap- and cup-slides}: strictly reduces $\distcc$, while preserving $(\nbcro,\nbdot,\disttunnel,\nbcc)$.
    \item \emph{dot-slide}: strictly reduces $\disttraffic$ while preserving $(\nbcro,\nbdot,\disttunnel,\nbcc,\distcc,\distcro)$.
    \item \emph{bubble-slide}: reduces $\disttunnel$, while preserving $(\nbcro,\nbdot)$.
  \end{itemize}
  This concludes.
\end{proof}

\begin{lemma}
  \label{lem:core_lemma_contextualization_bubbles}
  Assume that $\mathsf{d}$ is one of the following diagrams and $\Gamma$ is a context such that tunnels in $\Gamma[\mathsf{d}]$ are oriented as depicted:
  \def\tpscl{1}
  \begin{gather*}
    \tikzpic{\bubble[0][0][k]}
    \;,\qquad
    \tikzpic{\lorli\bubble[.8][-.5][k]}
    \;,\qquad
    \tikzpic{\rorli\bubble[-1.8][-.5][k]}
    \\
    \tikzpic{
      \cro[-2][0]
      \draw[very thick,red,->] (-2+1,.6) to (-2+.6,1);
      \draw[very thick,red,->] (-2+1,.4) to (-2+.6,0);
      \bubble[-.5][-.5][k]
    }[scale=\tpscl]
    \;,\qquad
    \tikzpic{
      \cro[-2][0]
      \draw[very thick,red,<-] (-2+1,.6) to (-2+.6,1);
      \draw[very thick,red,->] (-2+.4,1) to (-2,.6);
      \bubble[-2][.5][k]
    }[scale=\tpscl]
    \;,\qquad
    \tikzpic{
      \cro[-2][0]
      \draw[very thick,red,<-] (-2+1,.4) to (-2+.6,0);
      \draw[very thick,red,->] (-2+.4,0) to (-2,.4);
      \bubble[-2][-1.5][k]
    }[scale=\tpscl]
    \;,\qquad
    \tikzpic{
      \cro[-2][0]
      \draw[very thick,red,<-] (-2+.4,1) to (-2,.6);
      \draw[very thick,red,<-] (-2+.4,0) to (-2,.4);
      \bubble[-3.5][-.5][k]
    }[scale=\tpscl]
    \\
    \tikzpic{
      \lorli\rorli[3][0]
      \begin{scope}[xscale=3]
        \ca[0][1]
      \end{scope}
      \bubble[1][-.5][k]
    }
    \;,\qquad
    \tikzpic{\rorli\lorli[1][0]\ca[0][1]
      \bubble[-1.8][-.5][k]
    }
    \;,\qquad
    \tikzpic{\lorli\rorli[3][0]
      \begin{scope}[xscale=3]
        \cu[0][0]
      \end{scope}
      \bubble[1][-.5][k]
    }
    \qquad\an\qquad
    \tikzpic{\rorli\lorli[1][0]\cu[0][0]
      \bubble[-1.8][-.5][k]
    }
    \;.
  \end{gather*}
  If $\mathsf{d}\rhd x$ for a diagram $x$, then $\Gamma[\mathsf{d}]\succ \Gamma[x]$.
\end{lemma}

\begin{proof}
  If $\mathsf{d}\strongsucc x$, then $\Gamma[\mathsf{d}]\succ \Gamma[x]$ since $\strongsucc$ is context-agnostic and stronger than $\succ$.
  Otherwise, the diagram $x$ is identical to $\mathsf{d}$ but with the $k$-dotted bubble removed, additional dots on open strands and additional bubbles with at most $k-1$ dots.
  In each case, the $k$-dotted bubble in $\mathsf{d}$ lies in a region $v$ at maximal distance $t$ with respect to tunnels.
  If we remove the $k$-dotted bubble, then $\disttunnel^{\Gamma[\mathsf{d}]}(t)$ strictly reduces; adding bubbles with at most $k-1$ dots to $\mathsf{d}$ or adding dots on the pieces of strands that are open in $\mathsf{d}$ do not affect that fact.
  Note that for the latter, the open pieces of strands in $\mathsf{d}$ may belong to a bubble in $\Gamma[\mathsf{d}]$, so that, compared to $\Gamma[\mathsf{d}]$, the diagram $\Gamma[x]$ has an extra bubble with possibly more than $k$ dots; but in this case, the outer region of that bubble has tunnel distance at most $t-1$, so that $\disttunnel^{\Gamma[\mathsf{d}]}(t)$ is not affected.
\end{proof}

\subsubsection{Independent branchings}
\label{subsubsec:extended_MLRS_independent_branchings}

Tamed congruence of independent branchings holds in any extended \LMRS{}:

\begin{lemma}[Independent branching lemma for extended \LMRS{}]
  \label{lem:extended_MLRS_independent_branching}
  Let $(\sS,\lT,\sR_\agn,\sR_\dpt;\succ,\strongsucc,\asymp)$ be an extended \LMRS{} and let $(f,g)=(\ssf\circ s(\ssg),s(\ssf)\circ\ssg)$ be an independent branching, following the notations of \cref{subsubsec:RW_independent_branchings}.
  Then $(f,g)$ is $\succ$-tamely $\sS$-congruent.
\end{lemma}

\begin{proof}
  We wish to show that:
  \begin{gather*}
    s(\ssf)\circ s(\ssg)\succ x\circ y\qquad\text{ for all } x\in\supp(t(\ssf))\an y\in\supp(t(\ssg)).
  \end{gather*}
  Fix such $x$ and $y$. Since $f$ and $g$ are positive and $\succ$ is compatible with $\lT$, we always have $s(\ssf)\circ s(\ssg)\succ x\circ s(\ssg)$ and $s(\ssf)\circ s(\ssg)\succ s(\ssf)\circ y$.
  Hence, it suffices to show that either $x\circ s(\ssg)\succ x\circ y$ or $s(\ssf)\circ y\succ x\circ y$, which follows if either $x\circ\ssg$ or $\ssf\circ y$ belongs to $\lT$.
  In particular, this holds if either $\ssf$ or $\ssg$ belongs to $\lT_{\agn}=\Cont(\sR_{\agn})$, as by definition $\lT_{\agn}$ is context-agnostic.
  
  Assume instead that $\ssf$ and $\ssg$ belong to $\lT_{\dpt}$ and let $x\in\supp(t(\ssf))$.
  By property (i), either $s(\ssf)\asymp x$ or $s(\ssf)\strongsucc x$.
  In the former case, we find that $s(\ssf)\circ s(\ssg)\asymp x\circ s(\ssg)$ since $\asymp$ is context-agnostic.
  Since $\ssf\circ\ssg\in\lT_\dpt$, it follows from property (ii) that $x\circ \ssg\in\lT$; the claim follows.
  In the latter case, we find that $s(\ssf)\circ y\strongsucc x\circ y$ since $\strongsucc$ is context-agnostic, and hence $s(\ssf)\circ y\succ x\circ y$ since $\strongsucc$ is stronger than $\succ$.
  This concludes.
\end{proof}

When the rewriting steps only have $\strongsucc$-lower-order terms, independent branchings are better behaved:

\begin{lemma}[Independent branching lemma for $\strongsucc$]
  \label{lem:independent_branching_for_strongsucc}
  Let $(f,g)=(\ssf\circ s(\ssg),s(\ssf)\circ\ssg)$ be an independent branching, where $\ssf$ is either a non-affine or a dot-slide rewriting step; and similarly for $\ssg$.
  Denote by $t_{\strongsucc}(\ssf)$ (resp.\ $t_{\strongsucc}(\ssg)$) the $\strongsucc$-leading term of $t(\ssf)$ (resp.\ $t(\ssg)$), if it exists.
  If $t_{\strongsucc}(\ssf)\circ\ssg$ and $\ssf\circ t_{\strongsucc}(\ssg)$ belong to $\lT$, then the canonical $\affS$-confluence of $(f,g)$ is a $\strongsucc$-tamed $\affT$-preconfluence.
\end{lemma}

\begin{proof}
  Write $t(\ssf) = \lambda_\ssf t_{\strongsucc}(\ssf) + v_f$ and $t(\ssg) = \lambda_\ssg t_{\strongsucc}(\ssg) +v_g$, where we have $s(\ssf)\strongsucc x$ (resp.\ $s(\ssg)\strongsucc y$) for each $x\in\supp(v_f)$ (resp.\ $y\in\supp(v_g)$). (Here $\lambda_\ssf$ or $\lambda_\ssg$ could be zero.)
  Recall that the canonical $\affS$-confluence of $(f,g)$ is canonical only up to an ordering of the terms of $t(f)$ and $t(g)$. Choosing the $\strongsucc$-leading term to be first in each case defines a triple $(f',h,g')$, where $f'=\lambda_\ssf t_{\strongsucc}(\ssf)\circ\ssg + v_f\circ s(\ssg)$ and $g'=\lambda_\ssg \ssf\circ t_{\strongsucc}(\ssg) + s(\ssf)\circ v_g$, and $h$ is defined using the rewriting steps of the form $x\circ s(\ssg)$ and $s(\ssf)\circ y$ for $x\in\supp(v_f)$ and $y\in\supp(v_g)$. (Here $f'$ or $g'$ could be zero.)

  By hypothesis, the rewriting step $f'$, if non-trivial, belongs to $\affT$.
  Since $\strongsucc$ is context-agnostic, we have $t_{\strongsucc}(\ssf)\circ s(\ssg)\strongsucc\supp(v_f\circ s(\ssg))$, and in particular $t_{\strongsucc}(\ssf)\circ s(\ssg)\notin\supp(v_f\circ s(\ssg))$; it follows that $f'$ is positive. Similarly, the rewriting step $g'$ is a positive $\affT$-rewriting step.
  On the other hand, the terms in the sources and targets of $x\circ \ssg$ and $\ssf\circ y$ are all $\strongsucc$-tamed by $s(\ssf)\circ s(\ssg)$, again thanks to $\strongsucc$ being context-agnostic; it follows that $h$ has the desired form for a $\strongsucc$-tamed $\affT$-preconfluence.
\end{proof}

\subsubsection{Contextualization}
\label{subsubsec:confluence_contextualization}

As for independent branchings, when the rewriting steps only have $\strongsucc$-lower-order terms, contextualized branchings are better behaved:

\begin{lemma}[Contextualization lemma for $\strongsucc$]
  \label{lem:contextualization_for_strongsucc}
  Let $(\ssf,\ssg)$ be a monomial $\sS$-branching which admits a $\strongsucc$-tamed $\sS$-preconfluence $(\ssf',\ssh,\ssg')$, where $\ssf$ is either a non-affine or a dot-slide rewriting step; and similarly for $\ssg$.
  If $\Gamma$ is a context such that $\Gamma[\ssf]$, $\Gamma[\ssg]$, $\Gamma[\ssf']$ and $\Gamma[\ssg']$ belong to $\lT$, then $(\Gamma[\ssf'],\Gamma[\ssh],\Gamma[\ssg'])$ is a $\strongsucc$-tamed $\affT$-preconfluence for $(\Gamma[\ssf],\Gamma[\ssg])$.
\end{lemma}

\begin{proof}
  Since $\strongsucc$ is context-agnostic and $\ssf'$ and $\ssg'$ are rewriting sequences that always apply to the $\strongsucc$-leading term, so are $\Gamma[\ssf']$ and $\Gamma[\ssg']$; hence they must be positive rewriting sequences.
  It also follows from $\strongsucc$ being context-agnostic that since $\ssh$ is $\strongsucc$-tamed by $s(\ssf)=s(\ssg)$, we have that $\Gamma[\ssh]$ is $\strongsucc$-tamed by $s(\Gamma[\ssf])=s(\Gamma[\ssg])$.
  This concludes.
\end{proof}

Below is the analogous statement for bubble-slide branchings:
\begin{lemma}
  \label{lem:contextualization_bubbles}
  Let $(\ssf,\ssg)$ be a minimal bubble-slide branching.
  If $(\ssf,\ssg)$ is $\rhd$-tamely $\affT$-preconfluent, then for every context $\Gamma$, the branching $(\Gamma[\ssf],\Gamma[\ssg])$ is $\succ$-tamely $\affS$-congruent.
\end{lemma}

\begin{proof}
  Let $(\ssf',\ssh,\ssg')$ be the $\rhd$-tamed $\affT$-preconfluence given by hypothesis.
  Since sliding bubbles does not affect traffic rules, if $\Gamma[\ssf]$ and $\Gamma[\ssg]$ belong to $\affT$, then so do $\Gamma[\ssf']$ and $\Gamma[\ssg']$; in particular, they are $\succ$-tamed by the source of $(\Gamma[\ssf],\Gamma[\ssg])$.
  On the other hand, \cref{lem:core_lemma_contextualization_bubbles} shows that $\Gamma[\ssh]$ is $\succ$-tamed by the source of $(\Gamma[\ssf],\Gamma[\ssg])$.
\end{proof}

\subsubsection{Equivalence of presentations}
\label{subsubsec:confluence_equivalence_presentation}

In the statement of \cref{lem:extended_MLRS_independent_branching}, the conclusion is ``$\succ$-tamed $\sS$-congruence'', rather than ``$\succ$-tamed $\lT$-congruence''.
Indeed, we cannot expect that the canonical congruence associated to an independent branching will belong to $\lT$.
For instance (ignoring scalars and $\rhd$-lower-order terms):
\begin{gather*}
\begin{tikzcd}[ampersand replacement=\&,row sep=0pt]
	\&
  \tikzpic{
    \ca[1][3]
    \AntiDiag[0][2]\Diag[2][2]
    \li[0][1]\ca[1][1]\li[3][1]
    \custop[0][1]\lili[2][0]
    \rcu[2][0][->]
  }[scale=.5]
  \&
  \\
  \tikzpic{
    \ca[1][3]
    \AntiDiag[0][2]\Diag[2][2]
    \li[0][1]\ca[1][1]\li[3][1]
    \custop[0][1]\lili[2][0]
    \lcu[2][0][->]
  }[scale=.5]
  \&\&
  \tikzpic{
    \rcustop[0][0][->]\ca
  }[scale=1]
  \\
	\&
  \tikzpic{
    \lcustop[0][0][->]\ca
  }[scale=1]
	\arrow[from=1-2, to=2-3,bend left=15]
	\arrow[from=2-1, to=1-2,bend left=15]
	\arrow[from=2-1, to=3-2,bend right=15]
	\arrow["\times"{description}, dotted, leftrightarrow,from=3-2, to=2-3,bend right=15]
\end{tikzcd}
\end{gather*}
The arrow marked ``$\times$'' is not a $\affT$-congruence, as neither reading it forward nor backward is a valid $\affT$-rewriting step.
However, the next lemma shows that we can always replace an $\affS$-congruence by a $\affT$-congruence, and in a way that preserves tameness.

For an order $\succ$ and monomials $x,y\in\sX^*$, we say that $y$ is \defnemph{relatively $\succ$-tamed by $x$} if for any monomial $z\in\sX^*$, we have that $z\succ x$ implies $z\succ y$.
Similarly to our previous convention for orders, we say that a congruence $h$ is \emph{relatively $\succ$-tamed by $x$} if every term appearing in $h$ is relatively $\succ$-tamed by $x$.

Recall \cref{not:tame_arrow}:

\begin{lemma}
  \label{lem:affine_Brauer_equivalence_of_presentations}
  Assume that the third exceptional critical branching (\cref{subfig:critical_branchings_exceptional}) is $\rhd$-tamely $\affT$-preconfluent, that is, assume that the following $\affS$-rewriting step
  \begin{gather*}
  \begin{tikzcd}[ampersand replacement=\&,row sep=tiny]
    \tikzpic{\lcustop[0][0][->]\ca\bul\node[left=-1pt] at (0,0) {\scriptsize $k$};}[scale=1][0]
    \rar
    \&
    \VARdccbis
    \tikzpic{\rcustop[0][0][->]\ca\bul\node[left=-1pt] at (0,0) {\scriptsize $k$};}[scale=1][0]
    \;+\;
    \VARdccLOTbis
    \tikzpic{\custop[0][0][->]\ca\bul\node[left=-1pt] at (0,0) {\scriptsize $k$};}[scale=1][0]
  \end{tikzcd}
  \end{gather*}
  is $\rhd$-tamely $\affT$-preconfluent.
  Then for every $f\in\Cont(\affS)$, there exists a $\affT$-congruence $h\colon s(f)\overset{*}{\longleftrightarrow}_{\affT} t(f)$ relatively $\succ$-tamed by $f$.
\end{lemma}

\begin{proof}
  The dot-slide rewriting steps are defined in such a way that sliding in one direction, and then sliding back, comes back to the original diagram.
  For instance:
  \begin{gather*}
    \tikzpic{
      \ulcr
      \Diag[0][0][]
    }
    \;\overset{\stdcro}{\longrightarrow}\;
    \VARdcro\;
    \tikzpic{
      \drcr
      \Diag[0][0][]
    }
    \;+\;
    \tikzpic{\LOTdcro}
    \;\overset{\VARdcro(\stdcro)^*+\tikzpic{\LOTdcro}[scale=.7]}{\longrightarrow}\;
    \VARdcro\;
    \left(
    \VARdcro^{-1}\;
    \tikzpic{
      \ulcr
      \Diag[0][0][]
    }
    \;-\;\VARdcro^{-1}\;
    \tikzpic{\LOTdcro}
    \right)
    \;+\;
    \tikzpic{\LOTdcro}
    \;=\;
    \tikzpic{
      \ulcr
      \Diag[0][0][]
    }
  \end{gather*}
  As long as there is no stop sign, one of these two rewriting steps belongs to $\lT$. A similar argument applies to sliding a bubble through tunnels.

  It remains to consider the case of a dot sliding through a stop sign.
  By definition of the stop sign, there is a sequence of interchanges that pushes the cup down, so that the rest of the connected component lies above.
  Moreover, when doing an interchange with the cup, we can always combine it with an interchange with the dot.
  Hence, the sequence of interchanges can be seen as an $\affE$-congruence of the rewriting step, and we can assume that the dotted cup is at the bottom of its connected component.
  
  Next, we can rewrite the connected component above the dotted cup and assume it is in normal form.
  Indeed, if $s(f)=\Gamma[\;\tikzpic{\lcustop}\;]$ for some context $\Gamma$ and $g\colon\Gamma\overset{*}{\to}\Gamma'$ is an $\affS$-rewriting sequence such that $g[\;\tikzpic{\lcustop}\;]\colon\Gamma[\;\tikzpic{\lcustop}\;]\overset{*}{\to}\Gamma'[\;\tikzpic{\lcustop}\;]$ belongs to $\affT$, then the rewriting sequences $g[\;\tikzpic{\rcustop}\;]$ and $g[\;\tikzpic{\custop}\;]$ also belong to $\affT$, since the position of dots does not affect the traffic rules.
  Since $\succ$ is weakly compatible with $\affT$, each of the three instances of $g$ is relatively $\succ$-tamed by $f$.
  Denote by $f'$ the dot-slide through the cup in $\Gamma'[\;\tikzpic{\lcuquestion}\;]$; it is relatively $\succ$-tamed by $f$.
  \begin{gather*}
    \savebox{\tempboxa}{\tikzpic{\lcu}[scale=.7]}
    \savebox{\tempboxb}{\tikzpic{\rcu}[scale=.7]}
    \savebox{\tempboxc}{\tikzpic{\cu}[scale=.7]}
    \begin{tikzcd}[ampersand replacement=\&,row sep=large]
      \Gamma[\;\tikzpic{\lcustop}\;]
      \&
      \VARdccbis\;\Gamma[\;\tikzpic{\rcustop}\;] + \VARdccLOTbis\;\Gamma[\;\tikzpic{\custop}\;]
      \\
      \Gamma'[\;\tikzpic{\lcuquestion}\;]
      \&
      \VARdccbis\;\Gamma'[\;\tikzpic{\rcuquestion}\;] + \VARdccLOTbis\;\Gamma'[\;\tikzpic{\cuquestion}\;]
      \arrow[from=1-1,to=1-2,"f"]
      \arrow[from=2-1,to=2-2,"f'"]
      \arrow[from=1-1,to=2-1,"*"{description},"{g[\;\usebox{\tempboxa}\;]}"'{xshift=-2mm}]
      \arrow[from=1-2,to=2-2,"*"{description},"{\VARdccbis\;g\Gamma[\;\usebox{\tempboxb}\;] + \VARdccLOTbis\;g\Gamma[\;\usebox{\tempboxc}\;]}"{xshift=2mm}]
    \end{tikzcd}
  \end{gather*}
  If the cup does not have a stop sign $\Gamma'[\;\tikzpic{\lcuquestion}\;]=\Gamma'[\;\tikzpic{\lcu}\;]$, then $f'$ is a $\affT$-congruence of length one; this concludes.
  On the other hand, if the cup has a stop sign $\Gamma'[\;\tikzpic{\lcuquestion}\;]=\Gamma'[\;\tikzpic{\lcustop}\;]$, then the strand is a closed strand.
  Since the connected component above the dotted cup is now assumed to be a normal diagram, we have, up to contextualization, that $f'=\Gamma''[\ssf']$ where $\ssf'$ is the $\affS$-rewriting step given in the statement of the lemma.
  Following 
  Following \cref{lem:core_lemma_contextualization_bubbles},The contextualization of $\ssf'$ is $\succ$-tamed by the source of $f'$, and hence by the source of $f$.
  Moreover, it still belongs to $\affT$ after contextualization.
  This concludes.
\end{proof}

\begin{corollary}
  If a branching $(f,g)$ is $\succ$-tamely $\sS$-congruent, then it is $\succ$-tamely $\lT$-cong\-ruent.\hfill\qed
\end{corollary}

\subsection{Local analysis}
\label{subsec:local_analysis}

In this subsection, we prove that if critical branchings are tamely preconfluent, then so are most minimal overlapping branchings.
\begin{proposition}
  \label{prop:local_analysis}
  Assume that every critical branching (\cref{fig:critical_branchings}) is $\rhd$-tamely $\affT$-preconfluent. Then:
  \begin{enumerate}[(i)]
    \item every minimal non-affine branching (\cref{subfig:all_overlapping_branchings_non_affine}), with either $\tikzpic{\indexNF[][][?]}[scale=.6] = \tikzpic{\cro}[scale=.6]$ or $\tikzpic{\indexNF[][][?]}[scale=.6] = \tikzpic{\li}[scale=.6]$ for indexed branchings,
    \item every minimal dot-slide branching (\cref{subfig:all_overlapping_branchings_affine}),
    \item every minimal bubble-slide branching (\cref{fig:minimal_bubble_slide}),
  \end{enumerate}
  is $\rhd$-tamely $\affT$-pre\-conf\-luent.
\end{proposition}

The primary rewriting tools are ``symmetries'', used in different ways.
The proof of the proposition consists of four lemmas below: (i) follows from \cref{lem:reduction_horizontal_symmetries_non-affine}, (ii) follows from \cref{lem:twin_branchings_dot_slides} and \cref{lem:reduction_horizontal_symmetries_dot-slide} and (iii) follows from \cref{lem:twin_branchings_bubble_slides}.
We also define vertical symmetries (\cref{defn:vertical_symmetries}) and prove \cref{prop:vertical_symmetries} (assuming \cref{thm:affine_Grobner_system_iff_critical_branchings}).

\subsubsection{Symmetries of traffic rules}
\label{subsubsec:twin_branchings}

Below, we make precise the idea that the following branchings, and others, are ``twins'':
\begin{equation*}
  (f,g)=\tikzpic{
    \dli \lca[0][0]\ca[0][0][] \cu[1][0] \uli[2][0][->]
  }
  \qquad\an\qquad
  (\widehat{f},\widehat{g})=\tikzpic{
    \dli[0][0][->] \ca \rcu[1][0]\cu[1][0] \uli[2][0] 
  }\;.
\end{equation*}
where $f,\widehat{f}$ are dot-slides and $g,\widehat{g}$ are zigzags.
This will imply that one is $\strongsucc$-tamely $\affT$-preconfluent if and only if the other is.
In some sense, imposing traffic rules broke symmetry of the category of affine Brauer type; twins are a way to recover that symmetry.

Note already that $g$ and $\widehat{g}$ are two contextualizations of the same rewriting step $\ssg$, the zigzag.
A $\rhd$-tamed $\affT$-preconfluence for $(f,g)$ has the following form, where for simplicity, we assume that the coefficients of the interchange and the zigzag $\ssg$ are trivial:
\begin{equation*}
  \begin{tikzcd}[ampersand replacement=\&]
    \tikzpic{
      \dli \lca[0][0]\ca[0][0] \cu[1][0] \uli[2][0][->]
    }
    \&
    \VARdcc^{-1}\;
    \tikzpic{
      \dli \rca[0][0]\ca[0][0] \cu[1][0] \uli[2][0][->]
    }
    -
    \VARdcc^{-1}\VARdccLOT\;
    \tikzpic{
      \dli \ca[0][0]\ca[0][0] \cu[1][0] \uli[2][0][->]
    }
    \&
    \VARdcc^{-1}\VARdccbis\;
    \tikzpic{
      \dli \ca[0][0]\ca[0][0] \rcu[1][0] \uli[2][0][->]
    }
    +
    \VARdcc^{-1}(\VARdccLOTbis-\VARdccLOT)\;
    \tikzpic{
      \dli \ca[0][0]\ca[0][0] \cu[1][0] \uli[2][0][->]
    }
    \\
    \tikzpic{\bli}
    \&
    \VARdcc^{-1}\VARdccbis\;
    \tikzpic{\bli}
    +
    \VARdcc^{-1}(\VARdccLOTbis-\VARdccLOT)\;
    \tikzpic{\li}
    \&
    \VARdcc^{-1}\VARdccbis\;
    \tikzpic{\bli}
    +
    \VARdcc^{-1}(\VARdccLOTbis-\VARdccLOT)\;
    \tikzpic{\dli \ca[0][0]\ca[0][0] \cu[1][0] \uli[2][0][->]}
    \arrow[from=1-1,to=1-2,"f"]
    \arrow[from=1-1,to=2-1,"g"']
    \arrow[from=1-2,to=1-3,"*"{description},"f'"]
    \arrow[from=1-3,to=2-3,"\VARdcc^{-1}\VARdccbis\widehat{g}+v"]
    \arrow[from=2-3,to=2-2,dottedcd,leftrightarrow]
    \arrow[from=2-1,to=2-2,dottedcd,leftrightarrow,"?"]
  \end{tikzcd}
\end{equation*}
Here
$v=\VARdcc^{-1}(\VARdccLOTbis-\VARdccLOT)\;\;
\tikzpic{
  \dli \ca[0][0]\ca[0][0] \cu[1][0] \uli[2][0][->]
}$.
Tamed congruence reduces to the conditions $\VARdcc=\VARdccbis$ and  $\VARdccLOT=\VARdccLOTbis$. Note that the top part, the positive $\affT$-rewriting steps, is forced: we need to slide the dot fully before we can simplify the zigzag.
At the end, we read the source of the branching $(\widehat{f},\widehat{g})$ and the next step is $\widehat{g}$, up to scalar multiplication and vector addition.

On the other hand, a $\rhd$-tamed $\affT$-preconfluence for $(\widehat{f},\widehat{g})$ has the following form:
\begin{equation*}
  \begin{tikzcd}[ampersand replacement=\&]
    \VARdccbis^{-1}\VARdcc\;
    \tikzpic{
      \dli[0][0][->] \ca[0][0]\ca[0][0] \rcu[1][0] \uli[2][0]
    }
    +
    \VARdccbis^{-1}(\VARdccLOT-\VARdccLOTbis)\;
    \tikzpic{
      \dli[0][0][->] \ca[0][0]\ca[0][0] \cu[1][0] \uli[2][0]
    }
    \&
    \VARdccbis^{-1}\;
    \tikzpic{
      \dli[0][0][->] \ca[0][0]\ca[0][0] \lcu[1][0] \uli[2][0]
    }
    -
    \VARdccbis^{-1}
    \VARdccLOTbis\;
    \tikzpic{
      \dli[0][0][->] \ca[0][0]\ca[0][0] \cu[1][0] \uli[2][0]
    }
    \&
    \tikzpic{
      \dli[0][0][->] \ca[0][0]\ca[0][0] \rcu[1][0] \uli[2][0]
    }
    \\    
    \VARdccbis^{-1}\VARdcc\;
    \tikzpic{\bli}
    +
    \VARdccbis^{-1}(\VARdccLOT-\VARdccLOTbis)\;
    \tikzpic{\dli[0][0][->] \ca[0][0]\ca[0][0] \cu[1][0] \uli[2][0]}
    \&
    \VARdccbis^{-1}\VARdcc\;
    \tikzpic{\bli}
    +
    \VARdccbis^{-1}(\VARdccLOT-\VARdccLOTbis)\;
    \tikzpic{\li}
    \&
    \tikzpic{\bli}
    \arrow[from=1-3,to=1-2,"\widehat{f}"']
    \arrow[from=1-1,to=2-1,"\VARdccbis^{-1}\VARdcc g-\VARdccbis\VARdcc^{-1} v"']
    \arrow[from=1-2,to=1-1,"*"{description},"\widehat{f}'"']
    \arrow[from=1-3,to=2-3,"\widehat{g}"]
    \arrow[from=2-3,to=2-2,dottedcd,leftrightarrow,"?"']
    \arrow[from=2-1,to=2-2,dottedcd,leftrightarrow]
  \end{tikzcd}
\end{equation*}
This leads to the same conditions $\VARdcc=\VARdccbis$ and $\VARdccLOT=\VARdccLOTbis$.
That is not surprising: given the form of the dot-slide rewriting steps, sliding a dot in one direction is the opposite of sliding the dot in the other direction.
In the example above, we have that $f'=\VARdcc^{-1}\widehat{f}^{*}-\VARdcc^{-1}\VARdccLOT\;\tikzpic{\dli \ca[0][0]\ca[0][0] \cu[1][0] \uli[2][0][->]}[scale=.6]$ and $\widehat{f}'=\VARdccbis^{-1} f^{*}-\VARdccbis^{-1}\VARdccLOTbis\;\tikzpic{\dli \ca[0][0]\ca[0][0] \cu[1][0] \uli[2][0][->]}[scale=.6]$.

\begin{definition}
  Let $\trafficrule$ and $\widehat{\trafficrule}$ be two choices of traffic rules.
  Two monomial branchings $(f,g)$ and $(\widehat{f},\widehat{g})$, respectively belonging to $\affT_{\trafficrule}$ and $\affT_{\widehat{\trafficrule}}$, are said to be \defnemph{twins} if every $\rhd$-tamed $\affT_{\trafficrule}$-preconfluence (resp.\ $\affT_{\widehat{\trafficrule}}$-preconfluence) for $(f,g)$ (resp.\ $(\widehat{f},\widehat{g})$) has the following form, where $s$ and $\widehat{s}$ denote the sources of the branchings, $\lambda$ is an invertible scalar, and $v$ is a vector with $s\rhd v$ (resp.\ $\widehat{s}\rhd v$):
  \begin{equation*}
    \begin{tikzcd}
      s & t(f) & \lambda \widehat{s} + v\\
      t(g) && \lambda t(\widehat{g}) + v\\
      t(g') && \lambda t(\widehat{g}') + v
      \arrow[from=1-1,to=1-2,"f"]
      \arrow[from=1-1,to=2-1,"g"']
      \arrow[from=2-1,to=3-1,"*"{description},"g'"']
      \arrow[from=2-3,to=3-3,"*"{description},"\lambda \widehat{g}' + v"]
      \arrow[from=1-2,to=1-3,"*"{description},"f'"]
      \arrow[from=1-3,to=2-3,"\lambda\widehat{g}+v"]
      \arrow[from=3-1,to=3-3,leftrightarrow,dottedcd,"h"]
    \end{tikzcd}
    \qquad\left(\text{resp.\ }
    \begin{tikzcd}
      \lambda^{-1} s-\lambda^{-1}v & t(\widehat{f}) & \widehat{s}\\
      \lambda^{-1} t(g)-\lambda^{-1}v && t(\widehat{g})
      \\
      \lambda^{-1}t(g')-\lambda^{-1}v && t(\widehat{g}')
      \arrow[from=1-3,to=1-2,"\widehat{f}"']
      \arrow[from=2-1,to=3-1,"*"{description},"\lambda^{-1}g'-\lambda^{-1}v"']
      \arrow[from=2-3,to=3-3,"*"{description},"\lambda^{-1}\widehat{g}'-\lambda^{-1}v"']
      \arrow[from=1-1,to=2-1,"\lambda^{-1} g-\lambda^{-1}v"']
      \arrow[from=1-2,to=1-1,"*"{description},"\widehat{f}'"']
      \arrow[from=1-3,to=2-3,"\widehat{g}"]
      \arrow[from=3-1,to=3-3,leftrightarrow,dottedcd,"\lambda^{-1} h-\lambda^{-1}v"]
    \end{tikzcd}\right).
  \end{equation*}
\end{definition}


In particular, if the branchings $(f,g)$ and $(\widehat{f},\widehat{g})$ are twins, then $(f,g)$ is $\rhd$-tamely $\affT_{\trafficrule}$-preconfluent if and only if $(\widehat{f},\widehat{g})$ is $\rhd$-tamely $\affT_{\widehat{\trafficrule}}$-preconfluent.
Dot-slide and bubble-slide branchings all have a twin: the branching $\tikzpic{
  \ulcr \dlcr
  \Diag[0][0][<-]
  \AntiDiag[0][0][->]
}[scale=.6]$ even has two twins.
This leads to the following lemmas:

\begin{lemma}
  \label{lem:twin_branchings_dot_slides}
  Assume that in each of the following pairs (and in one case, quadruple) of branchings, at least one of them is $\rhd$-tamely $\affT$-preconfluent:
  \renewcommand{\sca}{.35}
  \newcommand{\vertsp}{.3ex}
  \newcommand{\arOrNot}[1]{#1}
  \begin{gather*}
    \tikzpic{
      \drcr \cro[0][1]
      \AntiDiag[0][1][\arOrNot{->}]
    }
    \mid
    \tikzpic{
      \cro \urcr[0][1]
      \Diag[0][0][\arOrNot{<-}]
    }
    \qquad
    \tikzpic{
      \dlcr \cro[0][1]
      \Diag[0][1][\arOrNot{->}]
    }
    \mid
    \tikzpic{
      \cro \ulcr[0][1]
      \AntiDiag[0][0][\arOrNot{<-}]
    }
    \quad\quad
    \tikzpic{
      \dlcr \li[2][0] \li[0][1] \cro[1][1] \cro[0][2] \li[2][2][\arOrNot{->}]
      \AntiDiag[0][0]
    }
    \mid
    \tikzpic{
      \cro \li[2][0] \li[0][1] \urcr[1][1] \cro[0][2] \li[2][2]
      \AntiDiag[0][0][\arOrNot{<-}]
    }
    \qquad
    \tikzpic{
      \drcr \li[2][0] \li[0][1] \cro[1][1] \cro[0][2] \li[2][2]
      \AntiDiag[0][2][\arOrNot{->}]
    }
    \mid
    \tikzpic{
      \cro \li[2][0] \li[0][1] \cro[1][1] \urcr[0][2] \li[2][2]
      \Diag[0][0][\arOrNot{<-}]
    }
    \qquad
    \tikzpic{
      \cro \li[2][0] \li[0][1] \drcr[1][1] \cro[0][2] \li[2][2]
      \Diag[0][2][\arOrNot{->}]
    }
    \mid
    \tikzpic{
      \cro \li[2][0][\arOrNot{<-}] \li[0][1] \cro[1][1] \ulcr[0][2] \li[2][2]
      \Diag[0][2]
    }
    \\[\vertsp]
    \tikzpic{
      \dli \lca[0][0]\ca[0][0] \cu[1][0] \uli[2][0][\arOrNot{->}]
    }
    \mid
    \tikzpic{
      \dli[0][0][\arOrNot{->}] \ca \rcu[1][0]\cu[1][0] \uli[2][0]
    }
    \qquad 
    \tikzpic{
      \uli \lcu[0][0]\cu[0][0] \ca[1][0] \dli[2][0][\arOrNot{->}]
    }
    \mid
    \tikzpic{
      \uli[0][0][\arOrNot{->}] \cu \rca[1][0]\ca[1][0] \dli[2][0]
    }
    \quad
    \quad
    \tikzpic{
      \ulcr \dlcr
      \Diag[0][0][\arOrNot{<-}]
      \AntiDiag[0][0][\arOrNot{->}]
    }
    \mid
    \tikzpic{
      \ulcr \urcr
      \Diag[0][0][\arOrNot{<-}]
      \AntiDiag[0][0][\arOrNot{<-}]
    }
    \mid
    \tikzpic{
      \urcr \drcr
      \Diag[0][0][\arOrNot{->}]
      \AntiDiag[0][0][\arOrNot{<-}]
    }
    \mid
    \tikzpic{
      \dlcr \drcr
      \Diag[0][0][\arOrNot{->}]
      \AntiDiag[0][0][\arOrNot{->}]
    }
    \\[\vertsp]
    \tikzpic{
      \dlcr \uli[0][1] \ca[1][1] \li[2][0][\arOrNot{<-}]
      \AntiDiag[0][0]
    } 
    \mid
    \tikzpic{
      \cro \uli[0][1] \rca[1][1]\ca[1][1] \li[2][0]
      \AntiDiag[0][0][\arOrNot{<-}]
    }
    \qquad
    \tikzpic{
      \drcr \uli[0][1][\arOrNot{->}] \ca[1][1] \li[2][0]
      \Diag[0][0]
    } 
    \mid
    \tikzpic{
      \ulcr \uli[0][1] \ca[1][1] \li[2][0]
      \Diag[0][0][\arOrNot{<-}]
    } 
    \quad\quad
    \tikzpic{ 
      \ca \urcr[1][-1] \li[0][-1] \li[2][-2] \cro[0][-2] \uli[2][0]
      \AntiDiag[0][-2][\arOrNot{<-}]
    } 
    \mid
    \tikzpic{ 
      \ca \cro[1][-1] \li[0][-1] \li[2][-2] \dlcr[0][-2] \uli[2][0][\arOrNot{->}]
      \AntiDiag[0][-2]
    } 
    \qquad
    \tikzpic{ 
      \ca \drcr[1][-1] \li[0][-1] \li[2][-2] \cro[0][-2] \uli[2][0]
      \Diag[0][-2][\arOrNot{<-}]
    } 
    \mid
    \tikzpic{ 
      \ca \cro[1][-1] \li[0][-1] \li[2][-2][\arOrNot{<-}] \drcr[0][-2] \uli[2][0]
      \Diag[0][-2]
    } 
    \\[\vertsp]
    \tikzpic{ 
      \ulcr \dli \cu[1][0] \li[2][0][\arOrNot{->}]
      \Diag[0][0]
    }
    \mid
    \tikzpic{ 
      \cro \dli \rcu[1][0]\cu[1][0] \li[2][0]
      \Diag[0][0][\arOrNot{->}]
    }
    \qquad
    \tikzpic{ 
      \urcr \dli[0][0][\arOrNot{->}] \cu[1][0] \li[2][0]
      \AntiDiag[0][0]
    }
    \mid
    \tikzpic{ 
      \dlcr \dli \cu[1][0] \li[2][0]
      \AntiDiag[0][0][\arOrNot{->}]
    } 
    \quad\quad
    \tikzpic{
      \dli[2][0] \cu \urcr[1][0] \li \li[2][1] \cro[0][1]
      \AntiDiag[0][1][\arOrNot{->}]
    } 
    \mid
    \tikzpic{
      \dli[2][0] \cu \cro[1][0] \li \li[2][1][\arOrNot{->}] \urcr[0][1]
      \AntiDiag[0][1]
    }
    \qquad
    \tikzpic{
      \dli[2][0][\arOrNot{->}] \cu \cro[1][0] \li \li[2][1] \ulcr[0][1]
      \Diag[0][1]
    }
    \mid
    \tikzpic{
      \dli[2][0] \cu \drcr[1][0] \li \li[2][1] \cro[0][1]
      \Diag[0][1][\arOrNot{->}]
    }
    \\[\vertsp]
    \tikzpic{
      \dlcr \ca[0][1]
      \AntiDiag[0][0]
      \Diag[0][0][\arOrNot{<-}]
    }
    \mid
    \tikzpic{
      \drcr \ca[0][1]
      \Diag[0][0]
      \AntiDiag[0][0][\arOrNot{<-}]
    }
    \qquad
    \tikzpic{
      \ca[0][1]\uli[2][1]
      \li[0][0] \urcr[1][0] 
      \cu\dli[2][0][\arOrNot{->}]
      \AntiDiag[1][0]
    }
    \mid
    \tikzpic{
      \ca[0][1]\uli[2][1][\arOrNot{->}]
      \li[0][0] \drcr[1][0]
      \cu\dli[2][0]
      \Diag[1][0]
    }
    \qquad
    \tikzpic{
      \ulcr
      \cu
      \Diag[0][0]
      \AntiDiag[0][0][\arOrNot{->}]
    }
    \mid 
    \tikzpic{
      \urcr
      \cu
      \AntiDiag[0][0]
      \Diag[0][0][\arOrNot{->}]
    }
  \end{gather*}
  Then every minimal dot-slide branching in \cref{subfig:all_overlapping_branchings_affine} is $\rhd$-tamely $\affT$-preconfluent.\hfill\qed
\end{lemma}

\begin{figure}
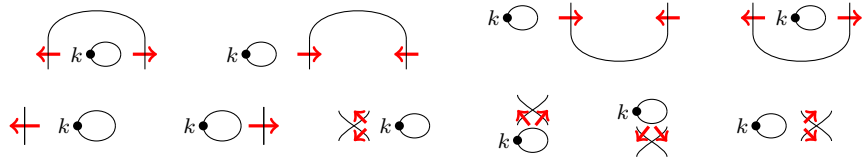

    \centering
    \def\tpscl{.8}
    \begin{gather*}
      \qquad
      \tikzpic{
        \lorli[-1.4][0] \bul[0][.5]\ca[0][.5]\cu[0][.5]
        \node[left=-1pt] at (0,.5) {\scriptsize $k$};
        \rorli[1.8][0]
        \draw (-1.4,1) to[out=90,in=90] (1.8,1);
      }[scale=\tpscl]
      \qquad
      \tikzpic{
        \begin{scope}[shift={(-3.5,0)}]
          \bul[0][.5]\ca[0][.5]\cu[0][.5]
          \node[left=-1pt] at (0,.5) {\scriptsize $k$};
        \end{scope}
        \rorli[-1.4][0]
        \lorli[1.8][0]
        \draw (-1.4,1) to[out=90,in=90] (1.8,1);
      }[scale=\tpscl]
      \qquad
      \tikzpic{
        \begin{scope}[shift={(-3.5,0)}]
          \bul[0][.5]\ca[0][.5]\cu[0][.5]
          \node[left=-1pt] at (0,.5) {\scriptsize $k$};
        \end{scope}
        \rorli[-1.4][0]
        \lorli[1.8][0]
        \draw (-1.4,0) to[out=-90,in=-90] (1.8,0);
      }[scale=\tpscl]
      \qquad
      \tikzpic{
        \lorli[-1.4][0] \bul[0][.5]\ca[0][.5]\cu[0][.5]
        \node[left=-1pt] at (0,.5) {\scriptsize $k$};
        \rorli[1.8][0]
        \draw (-1.4,0) to[out=-90,in=-90] (1.8,0);
      }[scale=\tpscl]
      \\[1ex]
      \tikzpic{
        \lorli[-1.4][0] 
        \li[-1.4][0] 
        \bul[0][.5]\ca[0][.5]\cu[0][.5]
        \node[left=-1pt] at (0,.5) {\scriptsize $k$};
      }
      \qquad
      \tikzpic{
        \rorli[1.6][0] 
        \li[1.6][0] 
        \bul[0][.5]\ca[0][.5]\cu[0][.5]
        \node[left=-1pt] at (0,.5) {\scriptsize $k$};
      }
      \qquad
      \tikzpic{
        \cro[-2][0]
        \bul[0][.5]\ca[0][.5]\cu[0][.5]
        \node[left=-1pt] at (0,.5) {\scriptsize $k$};
        \draw[very thick,red,->] (-2+.9,.6) to (-2+.5,1);
        \draw[very thick,red,->] (-2+.9,.4) to (-2+.5,0);
      }[scale=\tpscl]
      \qquad
      \tikzpic{
        \cro[-2][0]
        \begin{scope}[shift={(-2,-1)}]
          \bul[0][.5]\ca[0][.5]\cu[0][.5]
          \node[left=-1pt] at (0,.5) {\scriptsize $k$};
        \end{scope}
        \draw[very thick,red,<-] (-2+0,.5) to (-2+.4,0);
        \draw[very thick,red,<-] (-2+1,.5) to (-2+.6,0);
      }[scale=\tpscl]
      \qquad
      \tikzpic{
        \cro[-2][0]
        \begin{scope}[shift={(-2,1)}]
          \bul[0][.5]\ca[0][.5]\cu[0][.5]
          \node[left=-1pt] at (0,.5) {\scriptsize $k$};
        \end{scope}
        \draw[very thick,red,<-] (-2+1,.5) to (-2+.6,1);
        \draw[very thick,red,<-] (-2+0,.5) to (-2+.4,1);
      }[scale=\tpscl]
      \qquad
      \tikzpic{
        \cro[-2][0]
        \begin{scope}[shift={(-3.5,0)}]
          \bul[0][.5]\ca[0][.5]\cu[0][.5]
          \node[left=-1pt] at (0,.5) {\scriptsize $k$};
        \end{scope}
        \draw[very thick,red,->] (-2+.1,.6) to (-2+.5,1);
        \draw[very thick,red,->] (-2+.1,.4) to (-2+.5,0);
      }[scale=\tpscl]
    \end{gather*}
    \caption{Minimal bubble-slide branchings.}
    \label{fig:minimal_bubble_slide}
\end{figure}

A similar result holds for bubble-slide branchings:

\begin{lemma}
  \label{lem:twin_branchings_bubble_slides}
  Assume that every $\affT$-branching from \cref{subfig:critical_branchings_bubbles} is $\rhd$-tamely $\affT$-preconfluent.
  Then every minimal bubble-slide branching from \cref{fig:minimal_bubble_slide} is $\rhd$-tamely $\affT$-preconfluent.\hfill\qed
\end{lemma}

\begin{remark}
  \label{rem:bubble-slide-more-tunnels}
  \def\tpscl{1}
  The tunnels in the branching
  $
  \tikzpic{
    \cro[-2][0]
    \bul[0][.5]\ca[0][.5]\cu[0][.5]
    \node[left=-1pt] at (0,.5) {\scriptsize $k$};
    \draw[very thick,red,->] (-2+1,.6) to (-2+.6,1);
    \draw[very thick,red,->] (-2+1,.4) to (-2+.6,0);
  }[scale=\tpscl]
  $ are oriented as $
  \tikzpic{
    \cro[-2][0]
    \draw[very thick,red,->] (-2+1,.6) to (-2+.6,1);
    \draw[very thick,red,->] (-2+1,.4) to (-2+.6,0);
    \draw[very thick,red,->] (-2+.4,1) to (-2,.6);
    \draw[very thick,red,->] (-2+.4,0) to (-2,.4);
  }[scale=\tpscl]
  $.
  When we contextualize, they may be oriented as
  $
  \tikzpic{
    \cro[-2][0]
    \draw[very thick,red,->] (-2+1,.6) to (-2+.6,1);
    \draw[very thick,red,->] (-2+1,.4) to (-2+.6,0);
    \draw[very thick,red,<-] (-2+.4,1) to (-2,.6);
    \draw[very thick,red,->] (-2+.4,0) to (-2,.4);
  }[scale=\tpscl]
  $
  or
  $
  \tikzpic{
    \cro[-2][0]
    \draw[very thick,red,->] (-2+1,.6) to (-2+.6,1);
    \draw[very thick,red,->] (-2+1,.4) to (-2+.6,0);
    \draw[very thick,red,->] (-2+.4,1) to (-2,.6);
    \draw[very thick,red,<-] (-2+.4,0) to (-2,.4);
  }[scale=\tpscl]
  $.
  These alternative cases are analogous to ``twins''; in particular, if one is $\rhd$-tamely $\affT$-preconfluent, then so are the other two.
\end{remark}

\subsubsection{Horizontal symmetries}
\label{subsubsec:horizontal_symmetry}

The cap pulling rewriting step can be understood as the composition of a cap sliding in the opposite direction followed by an R2 rewriting step:
\begin{gather*}
  \tikzpic{
    \ca[0][2] \li[2][2]
    \li[0][1] \cro[1][1]
    \cro[0][0] \li[2][0]
  }
  \;\overset{\stcasl}{\longleftarrow}\;
  \VARcasl^{-1}\;
  \tikzpic{
    \ca[1][2] \li[0][2]
    \li[2][1] \cro[0][1]
    \cro[0][0] \li[2][0]
  }
  \;-\;\VARcasl^{-1}\;
  \tikzpic{
    \LOTcasl[0][1]
    \cro[0][0] \li[2][0]
  }
  \;\overset{\stRtwo}{\longrightarrow}\;
  \VARcasl^{-1}\;
  \tikzpic{
    \li[0][2]\ca[1][2]
    \LOTRtwo\li[2][0]\li[2][1]
  }
  \;-\;\VARcasl^{-1}\;
  \tikzpic{
    \LOTcasl[0][1]
    \cro[0][0] \li[2][0]
  }
\end{gather*}
Given that fact, one may expect that some branchings involving a cap pulling are consequences of other branchings. Indeed, consider the following branching between a cap pulling and an R2 rewriting step, where we decomposed the former as above:
\begin{gather}
  \label{eq:horizontal_symmetries_example}
  \begin{tikzcd}[ampersand replacement=\&]
    \tikzpic{
      \ca[0][2] \li[2][2]
      \li[0][1] \cro[1][1]
      \cro[0][0] \li[2][0]
      \cro[0][-1] \li[2][-1]
    }[scale=.7]
    \&
    \VARcasl^{-1}\;
    \tikzpic{
      \ca[1][2] \li[0][2]
      \li[2][1] \cro[0][1]
      \cro[0][0] \li[2][0]
      \cro[0][-1] \li[2][-1]
    }[scale=.7]
    \;-\;\VARcasl^{-1}\;
    \tikzpic{
      \LOTcasl[0][1]
      \cro[0][0] \li[2][0]
      \cro[0][-1] \li[2][-1]
    }[scale=.7]
    \&
    \VARcasl^{-1}\;
    \tikzpic{
      \li[0][2]\ca[1][2]
      \LOTRtwo\li[2][0]\li[2][1]
      \cro[0][-1] \li[2][-1]
    }[scale=.7]
    \;-\;\VARcasl^{-1}\;
    \tikzpic{
      \LOTcasl[0][1]
      \cro[0][0] \li[2][0]
      \cro[0][-1] \li[2][-1]
    }[scale=.7]
    \\
    \tikzpic{
      \ca[0][3]\li[2][3]
      \cro[1][2]\li[0][2] 
      \LOTRtwo\li[2][0]\li[2][1]
    }[scale=.7]
    \&
    \VARcasl^{-1}\;
    \tikzpic{
      \li[0][3]\ca[1][3]
      \cro[0][2]\li[2][2] 
      \LOTRtwo\li[2][0]\li[2][1]
    }[scale=.7]
    \;-\;\VARcasl^{-1}\;
    \tikzpic{
      \LOTcasl[0][1]
      \cro[0][0] \li[2][0]
      \cro[0][-1] \li[2][-1]
    }[scale=.7]
    \arrow[from=1-2,to=1-1,"\stcasl"']
    \arrow[from=1-2,to=2-2,"\stRtwo"']
    \arrow[from=1-2,to=1-3,"\stRtwo"]
    \arrow[from=1-1,to=1-3,"\stcapu",bend left=20pt]
    \arrow[from=1-1,to=2-1,"\stRtwo"']
    \arrow[from=2-2,to=2-1,"\stcasl"']
  \end{tikzcd}
\end{gather}
We added the bottom right vector; note that this vector is $\strongsucc$-tamed by the source of the branching.
To achieve $\strongsucc$-tamed $\affT$-preconfluence, it remains to check that we have a congruence, $\strongsucc$-tamed by the source of the branching (for which $(\nbcro,\nbdot)=(3,0)$), as follows:
\begin{gather*}
  \tikzpic{
    \li[0][2]\ca[1][2]
    \LOTRtwo\li[2][0]\li[2][1]
    \cro[0][-1] \li[2][-1]
  }[scale=.7]
  \;\tamearrow\;
  \tikzpic{
    \li[0][3]\ca[1][3]
    \cro[0][2]\li[2][2] 
    \LOTRtwo\li[2][0]\li[2][1]
  }[scale=.7]
  \qquad\Leftrightarrow\qquad
  \tikzpic{
    \li[0][2]\ca[1][2]
    \cro[0][1]\li[2][1]
    \cro[0][0]\li[2][0]
    \cro[0][-1] \li[2][-1]
    \draw[branch_overlay] (-.3,-1) rectangle (1+.3,2);
  }[scale=.7]\quad\text{is $\strongsucc$-tamely $\affT$-preconfluent}.
\end{gather*}
As pictured above, this congruence follows from $\strongsucc$-tamed $\affT$-preconfluence of (a contextualization of) a $\stRtwo$-$\stRtwo$ branching.
We conclude that if this $\stRtwo$-$\stRtwo$ branching is $\strongsucc$-tamely $\affT$-preconfluent, so is the branching we started with.

More conceptually, the ``cap sliding in the opposite direction'' is the cap sliding of the reverse rewriting system $\affR^\rev$:

\begin{definition}
  The set of rewriting steps $\affR^\rev$ is defined as in \cref{fig:nonaffine_rewriting_steps_rev} for non-affine rewriting steps and as in \cref{defn:main_text_affine_Brauer_rewriting_system} for affine rewriting steps.
  We write $(\affS^\rev,\affT^\rev)$ for the associated \LMRS{}.
\end{definition}

\begin{figure}
	\centering

	\def\spaceout{2ex}
	\def\spacein{.5ex}
  \renewcommand{\sca}{.4}

	\begin{gather*}
		\begin{IEEEeqnarraybox}{CcC}
			\tikzpic{
        \cro[0][1]
			  \cro
			}
			\;\overset{\stRtwo}{\longrightarrow}\;
      \tikzpic{\LOTRtwo}
			\mspace{30mu}
      %
			\tikzpic{
        \cro[1][2] \li[0][2]
        \li[2][1] \cro[0][1]
        \cro[1][0] \li[0][0]
      }
			\;\overset{\stRthree}{\longrightarrow}\;
      \VARRthree^{-1}\;
			\tikzpic{
        \cro[0][2] \li[2][2]
        \li[0][1] \cro[1][1]
        \cro \li[2][0]
      }
      \;-\;\VARRthree^{-1}\;
      \tikzpic{\LOTRthree}
      &
			\mspace{60mu}
      &
			\tikzpic{
        \cu[1][0] \uli[2][0]
			  \dli \ca
			}
			\;\overset{\stzz}{\longrightarrow}\;
      \VARzz\;
      \tikzpic{\uli\dli}
			\mspace{30mu}
      %
			\tikzpic{
        \dli[2][0] \ca[1][0]
        \uli \cu
			}
			\;\overset{\stzzbis}{\longrightarrow}\;
      \VARzzbis\;
      \tikzpic{\uli\dli}
    \end{IEEEeqnarraybox}
    \\[\spaceout]
    \begin{IEEEeqnarraybox}{rClcrCl}
			\tikzpic{
        \ca[0][1] \uli[2][1]
        \li \cro[1][0]
      }
			&\overset{\stcasl}{\longrightarrow}&
      \VARcasl^{-1}\;
			\tikzpic{
        \uli[0][1] \ca[1][1]
        \cro \li[2][0]
      }
      \;-\;\VARcasl^{-1}\;
      \tikzpic{\LOTcasl}
			&\mspace{100mu}&
      %
			\tikzpic{
        \ca \uli[2][0]
        \li[0][-1] \cro[1][-1]
        \cro[0][-2] \li[2][-2]
      }[xscale=-1]
			\;&\overset{\stcapu}{\longrightarrow}&\;
      \VARcasl\;
      \tikzpic{
        \uli[0][2]\ca[1][2]
        \LOTRtwo\li[2][0]\li[2][1]
      }[xscale=-1]
      \;+\;
      \tikzpic{
        \LOTcasl[0][1]
        \cro[0][0] \li[2][0]
      }[xscale=-1]
      \\[\spaceout]
      \tikzpic{
        \li[0][0] \cro[1][0]
        \cu[0][0] \dli[2][0]
      }
      \;&\overset{\stcusl}{\longrightarrow}&\;
      \VARcusl^{-1}\;
      \tikzpic{
        \cro \li[2][0]
        \dli \cu[1][0]
      }
      \;-\;\VARcusl^{-1}\;
      \tikzpic{\LOTcusl}
      &&
      \tikzpic{
        \cro[0][1] \li[2][1]
        \li[0][0] \cro[1][0]
        \cu \dli[2][0]
      }[xscale=-1]
      \;&\overset{\stcupu}{\longrightarrow}&\;
      \VARcusl\;
      \tikzpic{
        \LOTRtwo[0][1]\li[2][1]\li[2][2]
        \dli[0][1]\cu[1][1]
      }[xscale=-1]
      \;+\;
      \tikzpic{
        \cro[0][2] \li[2][2]
        \LOTcusl[0][0]
      }[xscale=-1]
    \end{IEEEeqnarraybox}
		\\[\spaceout]
    \begin{IEEEeqnarraybox}{CcCcC}
			\tikzpic{
        \ca[0][1]
			  \cro
			}
      \;\overset{\stuk}{\longrightarrow}\;
      \tikzpic{\LOTuk}
			&\mspace{100mu}&
      %
      \tikzpic{
        \cro
        \cu
      }
      \;\overset{\stdk}{\longrightarrow}\;
      \tikzpic{\LOTdk}
      &\mspace{100mu}&
      %
      \tikzpic{
        \ca[0][1] \uli[2][1]
        \li[0][0] \cro[1][0]
        \cu \dli[2][0]
      }[xscale=-1]
      \;\overset{\stlk}{\longrightarrow}\;
      \VARcusl\;
      \tikzpic{\LOTuk[0][0]\dli[0][0]\cu[1][0]\li[2][0]\uli[2][1]}[xscale=-1]
      \;+\;
      \tikzpic{\LOTcusl\ca[0][2]\uli[2][2]}[xscale=-1]
    \end{IEEEeqnarraybox}
  \end{gather*}

	\caption{Non-affine rewriting steps of $\affR^\rev$.}
	\label{fig:nonaffine_rewriting_steps_rev}
\end{figure}

When we defined $\affR$, we broke horizontal symmetry and chose a preferred direction for certain rewriting steps, namely the R3 move ($\stRthree$), the cap-slide ($\stcasl$) and the cup-slide ($\stcusl$).
Every argument made for $\affR$ could have been made for $\affR^\rev$.
In particular:

\begin{lemma}
  Let $(f,g)$ be a branching in $\affR$ and denote by $(f^\rev,g^\rev)$ the associated branching in $\affR^\rev$.
  Then $(f,g)$ is $\strongsucc$-tamely $\affT$-preconfluent if and only if $(f^\rev,g^\rev)$ is $\strongsucc$-tamely $\affT^\rev$-preconfluent.\hfill\qed
\end{lemma}

In hindsight, the square in \cref{eq:horizontal_symmetries_example} is an independent branching in $\affR^\rev$; this is why it leads to a $\strongsucc$-tamed $\affT$-preconfluence, using the independent branching for $\strongsucc$ given above (\cref{lem:independent_branching_for_strongsucc}). In the schematic below, we use the symbol $\llcorner$ to denote independent branchings:
\begingroup
\savebox{\tempboxa}{\tikzpic{
  \li[0][2]\ca[1][2]
  \cro[0][1]\li[2][1]
  \cro[0][0]\li[2][0]
  \cro[0][-1] \li[2][-1]
  \draw[branch_overlay] (-.3,-1) rectangle (1+.3,2);
}[scale=.5]}
\def\scl{.6}
\begin{gather*}
  \begin{tikzcd}[ampersand replacement=\&]
    \tikzpic{
      \ca[0][2] \li[2][2]
      \li[0][1] \cro[1][1]
      \cro[0][0] \li[2][0]
      \cro[0][-1] \li[2][-1]
    }[scale=\scl]
    \&
    \tikzpic{
      \ca[1][2] \li[0][2]
      \li[2][1] \cro[0][1]
      \cro[0][0] \li[2][0]
      \cro[0][-1] \li[2][-1]
    }[scale=\scl]
    \&
    \tikzpic{
      \li[0][2]\ca[1][2]
      \LOTRtwo\li[2][0]\li[2][1]
      \cro[0][-1] \li[2][-1]
    }[scale=\scl]
    \\
    \tikzpic{
      \ca[0][3]\li[2][3]
      \cro[1][2]\li[0][2] 
      \LOTRtwo\li[2][0]\li[2][1]
    }[scale=\scl]
    \&
    \tikzpic{
      \li[0][3]\ca[1][3]
      \cro[0][2]\li[2][2] 
      \LOTRtwo\li[2][0]\li[2][1]
    }[scale=\scl]
    \arrow[from=1-2,to=1-1,"\stcasl"']
    \arrow[from=1-2,to=2-2,"\stRtwo"']
    \arrow[from=1-2,to=1-3,"\stRtwo"]
    \arrow[from=1-1,to=2-1,"\stRtwo"']
    \arrow[from=2-2,to=2-1,"\stcasl"']
    \arrow[from=1-2,to=2-1, phantom, "\llcorner", very near start]
    \arrow[from=2-2,to=1-3,""{name=0, anchor=center, inner sep=0},dottedcd,bend right=50pt,leftrightarrow]
    \arrow[from=1-2,phantom,to=0,"{\usebox{\tempboxa}}"]
  \end{tikzcd}
\end{gather*}
\endgroup
If instead $\strongsucc$-tamed $\affT$-preconfluence is given by another branching, we picture that branching inside the commuting diagram.
Note that the information of the leading scalar and the $\strongsucc$-lower-order terms is superfluous, as both scalars and tamed congruences are transitive when gluing patches.

Going through every branching, one finds that many of them can be recovered by ``patching up'' the tamed preconfluences of other branchings. In each case, the proof consists in providing the appropriate patches; the reader can find them in \cref{fig:horizontal_symmetries_A,fig:horizontal_symmetries_B,fig:horizontal_symmetries_C} in \cref{sec:horizontal_symmetries}.
Formally, we use the independent-branching lemma for $\strongsucc$ (\cref{lem:independent_branching_for_strongsucc}; part (ii)) and the contextualization lemma for $\strongsucc$ (\cref{lem:contextualization_for_strongsucc}).
Note that if a schematic exists for a certain branching, the vertical symmetry of that branching also has a schematic, by taking the vertical symmetry of the schematic of the original branching.

We sum up the results in the following two lemmas:

\begin{lemma}
  \label{lem:reduction_horizontal_symmetries_non-affine}
  \def\tempsp{20mu}
  \renewcommand{\sca}{.3}
  Assume that every $\affT$-branching below is $\strongsucc$-tamely $\affT$-preconfluent:
  \begin{gather*}
    \tikzpic{
      \cro[0][2]
      \cro[0][1]
      \cro
    }
    \mspace{\tempsp}
    \tikzpic{
      \cro[0][3]\li[2][3]
      \cro[0][2]\li[2][2]
      \li[0][1]\cro[1][1]
      \cro[0][0]\li[2][0]
    }
    \mspace{\tempsp}
    \tikzpic{
      \cro[0][3]\li[2][3]
      \li[0][2]\cro[1][2]
      \cro[0][1]\li[2][1]
      \cro[0][0]\li[2][0]
    }
    \mspace{\tempsp}
    \tikzpic{
      \uli[0][2]\ca[1][2]
      \cro[0][1]\li[2][1]
      \cro\li[2][0]
    }
    \mspace{\tempsp}
    \tikzpic{
      \uli[0][4]\ca[1][4]
      \cro[0][3]\li[2][3]
      \lili[0][2]\li[2][2]
      \li[0][1]\cro[1][1]
      \cro\li[2][0]
    }
    \mspace{\tempsp}
    \tikzpic{
      \cro[0][2]\ca[2][2]
      \li[0][1]\cro[1][1]\li[3][1]
      \cro\lili[2][0]
    }
    \mspace{\tempsp}
    \tikzpic{
      \uli[0][1]\uli[1][1]\ca[2][1]
      \cro\lili[2][0]
      \dli[0][0]\cu[1][0]\dli[3][0]
    }
    \mspace{\tempsp}
    \tikzpic{
      \cro[0][1]
      \cro[0][0]
      \cu[0][0]
    }
    \mspace{\tempsp}
    \\
    \tikzpic{
      \cro[0][1]\li[2][1]
      \cro\li[2][0]
      \dli[0][0]\cu[1][0]
    }
    \mspace{\tempsp}
    \tikzpic{
      \cro[0][3]\li[2][3]
      \li[0][2]\cro[1][2]
      \lili[0][1]\li[2][1]
      \cro\li[2][0]
      \dli\cu[1][0]
    }
    \mspace{\tempsp}
    \tikzpic{
      \cro[0][2]\lili[2][2]
      \li[0][1]\cro[1][1]\li[3][1]
      \cro\cu[2][1]
    }
    \mspace{\tempsp}
    \tikzpic{
      \uli[0][1]\ca[1][1]\uli[3][1]
      \cro\lili[2][0]
      \dli[0][0]\dli[1][0]\cu[2][0]
    }
    \mspace{\tempsp}
    \tikzpic{
      \ca[0][2]
      \cro[0][1]
      \cro[0][0]
    }
    \mspace{\tempsp}
  \end{gather*}
  Then each of the following $\affT$-branchings
  \begin{gather*}
    \tikzpic{
      \cro[0][4]\li[2][4]
      \li[0][3]\cro[1][3]
      \cro[0][2]\li[2][2]
      \li[0][1]\cro[1][1]
      \cro[0][0]\li[2][0]
    }
    \mspace{\tempsp}
    \mspace{\tempsp}
    \tikzpic{
      \uli[0][4]\ca[1][4]
      \cro[0][3]\li[2][3]
      \lili[0][2]\cro[2][2]
      \li[0][1]\cro[1][1]
      \cro\li[2][0]
      \draw[] (3,0) to (3,2);\draw[] (3,3) to (3,4.4);
    }
    \mspace{\tempsp}
    \tikzpic{ 
      \ca[0][3]\uli[2][3]
      \li[0][2]\cro[1][2]
      \cro[0][1]\li[2][1]
      \cro\li[2][0]
    }
    \mspace{\tempsp}
    \tikzpic{
      \ca[0][5]\uli[2][5]
      \li[0][4]\cro[1][4]
      \cro[0][3]\li[2][3]
      \lili[0][2]\li[2][2]
      \li[0][1]\cro[1][1]
      \cro\li[2][0]
    }
    \mspace{\tempsp}
    \tikzpic{
      \ca[0][5]\uli[2][5]
      \li[0][4]\cro[1][4]
      \cro[0][3]\li[2][3]
      \lili[0][2]\cro[2][2]
      \li[0][1]\cro[1][1]
      \cro\li[2][0]
      \draw[] (3,0) to (3,2);\draw[] (3,3) to (3,5.4);
    }
    \mspace{\tempsp}
    \mspace{\tempsp}
    \tikzpic{ 
      \ca[0][2]\uli[2][2]\uli[3][2]
      \li[0][1]\cro[1][1]\li[3][1]
      \cro\cu[2][1]
    }
    \mspace{\tempsp}
    \tikzpic{
      \ca[0][4]\uli[2][4]
      \li[0][3]\cro[1][3] 
      \lili[0][2]\li[2][2]
      \cro[0][1]\li[2][1] 
      \li[0][0]\cro[1][0]
      \dli[2][0]\cu
    }
    \mspace{\tempsp}
    \tikzpic{
      \ca[0][4]\uli[2][4]
      \li[0][3]\cro[1][3] 
      \lili[0][2]\cro[2][2]
      \cro[0][1]\li[2][1] 
      \li[0][0]\cro[1][0]
      \dli[2][0]\cu
      \draw[] (3,-.4) to (3,2);\draw[] (3,3) to (3,4.4);
    }
    \mspace{\tempsp}
    \mspace{\tempsp}
    \tikzpic{
      \ca[0][2]\uli[2][2]
      \cro[0][0]\li[2][0]
      \li[0][1]\cro[1][1]
      \cu\dli[2][0]
    }
    \mspace{\tempsp}
    \tikzpic{
      \ca[0][2]
      \li[0][1]\li[1][1]\ca[2][1]
      \li\cro[1][0]\li[3][0]
      \cu\dli[2][0]\dli[3][0]
    }
    \\
    \tikzpic{
      \cro[0][3]\li[2][3]
      \li[0][2]\cro[1][2]
      \lili[0][1]\cro[2][1]
      \cro\li[2][0]
      \dli\cu[1][0]
      \draw[] (3,-.4) to (3,1);\draw[] (3,2) to (3,4);
    }
    \mspace{\tempsp}
    \tikzpic{
      \cro[0][2]\li[2][2]
      \cro[0][1]\li[2][1]
      \li\cro[1][0]
      \cu\dli[2][0]
    }
    \mspace{\tempsp}
    \tikzpic{
      \cro[0][4]\li[2][4]
      \li[0][3]\cro[1][3]
      \lili[0][2]\li[2][2]
      \cro[0][1]\li[2][1]
      \li\cro[1][0]
      \cu\dli[2][0]
    }
    \mspace{\tempsp}
    \tikzpic{
      \cro[0][4]\li[2][4]
      \li[0][3]\cro[1][3]
      \lili[0][2]\cro[2][2]
      \cro[0][1]\li[2][1]
      \li\cro[1][0]
      \cu\dli[2][0]
      \draw[] (3,-.4) to (3,2);\draw[] (3,3) to (3,5);
    }
    \mspace{\tempsp}
    \mspace{\tempsp}
    \tikzpic{
      \cro[0][1]\ca[2][1]
      \li[0][0]\cro[1][0]\li[3][0]
      \cu\dli[2][0]\dli[3][0]
    }
    \mspace{\tempsp}
    \tikzpic{
      \ca[0][2]\uli[2][2]
      \cro[0][1]\li[2][1]
      \li\cro[1][0]
      \cu\dli[2][0]
    }
    \mspace{\tempsp}
    \tikzpic{
      \ca[0][2]\uli[2][2]\uli[3][2]
      \li[0][1]\cro[1][1]\li[3][1]
      \lili[0][0]\cu[2][1]
      \cu
    }
  \end{gather*}
  is $\strongsucc$-tamely $\affT$-preconfluent.\hfill\qed
\end{lemma}

\begin{lemma}
  \label{lem:reduction_horizontal_symmetries_dot-slide}
  \def\tempsp{20mu}
  \renewcommand{\sca}{.3}
  \newcommand{\arOrNot}[1]{#1}
  Assume that every $\affT$-branching below is $\strongsucc$-tamely $\affT$-preconfluent:
  \begin{gather*}
    \tikzpic{
      \uli[0][2]\ca[1][2]
      \cro[0][1]\li[2][1]
      \cro\li[2][0]
    }
    \mspace{\tempsp}
    \tikzpic{
      \ca[0][1]\uli[2][1]
      \cro[0][0]\li[2][0]
      \dli \cu[1][0]
    }
    \mspace{\tempsp}\mspace{\tempsp}
    \tikzpic{
      \drcr \cro[0][1]
      \AntiDiag[0][1][\arOrNot{->}]
    }
    \mspace{\tempsp}
    \tikzpic{
      \dlcr \cro[0][1]
      \Diag[0][1][\arOrNot{->}]
    }
    \mspace{\tempsp}\mspace{\tempsp}
    \tikzpic{
      \dlcr \uli[0][1] \ca[1][1] \li[2][0][\arOrNot{<-}]
    } 
    \mspace{\tempsp}
    \tikzpic{
      \drcr \uli[0][1][\arOrNot{->}] \ca[1][1] \li[2][0]
    } 
    \mspace{\tempsp}
    \tikzpic{ 
      \cro \dli \cu[1][0] \bli[2][0]
      \Diag[0][0][\arOrNot{->}]
    }
    \mspace{\tempsp}\mspace{\tempsp}
    \tikzpic{
      \drcr \ca[0][1]
      \AntiDiag[0][0][\arOrNot{<-}]
    }
    \\
    \tikzpic{
      \cro[0][1]\li[2][1]
      \cro\li[2][0]
      \dli[0][0]\cu[1][0]
    }
    \mspace{\tempsp}
    \tikzpic{
      \uli[0][1] \ca[1][1]
      \cro[0][0]\li[2][0]
      \cu[0][0]\uli[2][0]
    }
    \mspace{\tempsp}\mspace{\tempsp}
    \tikzpic{
      \cro \urcr[0][1]
      \Diag[0][0][\arOrNot{<-}]
    }
    \mspace{\tempsp}
    \tikzpic{
      \cro \ulcr[0][1]
      \AntiDiag[0][0][\arOrNot{<-}]
    }
    \mspace{\tempsp}\mspace{\tempsp}
    \tikzpic{
      \ulcr \li[2][0][\arOrNot{->}]
      \dli[0][0] \cu[1][0]
    } 
    \mspace{\tempsp}
    \tikzpic{
      \urcr \li[2][0]
      \dli[0][0][\arOrNot{->}] \cu[1][0]
    } 
    \mspace{\tempsp}
    \tikzpic{ 
      \cro \uli[0][1] \ca[1][1] \bli[2][0]
      \AntiDiag[0][0][\arOrNot{<-}]
    }
    \mspace{\tempsp}\mspace{\tempsp}
    \tikzpic{
      \urcr \cu[0][0]
      \Diag[0][0][\arOrNot{->}]
    }
  \end{gather*}
  Then each of the following $\affT$-branchings
  \begin{gather*}
    \tikzpic{
      \ca[0][1]\uli[2][1]
      \li[0][0] \urcr[1][0]
      \cu\dli[2][0][\arOrNot{->}]
    }
    \mspace{\tempsp} 
    \tikzpic{ 
      \ca \cro[1][-1] \li[0][-1] \li[2][-2] \dlcr[0][-2] \uli[2][0][\arOrNot{->}]
    } 
    \mspace{\tempsp}
    \tikzpic{ 
      \ca \cro[1][-1] \li[0][-1] \li[2][-2][\arOrNot{<-}] \drcr[0][-2] \uli[2][0]
    } 
    \mspace{\tempsp}
    \mspace{\tempsp}
    \tikzpic{
      \ca[0][1]\uli[2][1][\arOrNot{->}]
      \li[0][0] \drcr[1][0]
      \cu\dli[2][0]
    }
    \mspace{\tempsp} 
    \tikzpic{ 
      \ulcr[0][1]\li[2][1]
      \li\cro[1][0]
      \cu\dli[2][0][\arOrNot{->}]
    } 
    \mspace{\tempsp}
    \tikzpic{ 
      \urcr[0][1]\li[2][1][\arOrNot{->}]
      \li\cro[1][0]
      \cu\dli[2][0]
    } 
  \end{gather*}
  is $\strongsucc$-tamely $\affT$-preconfluent.\hfill\qed
\end{lemma}

\subsubsection{Vertical symmetries}
\label{subsec:symmetries}


Often, the category of affine Brauer type, as well as the \LMRS{} presenting it, exhibit a vertical symmetry.
Note that the previous two symmetries we discussed, symmetries of traffic rules and horizontal symmetries, were symmetries ``internal'' to the rewriting system, and were independent of any further assumption.
In contrast, the vertical symmetries we discuss here are ``external'', in the sense that they are induced on the rewriting system by a symmetry already assumed in the category of affine Brauer type presented by the rewriting system.

First, we explain what we mean by vertical symmetry.
Let $\cC$ be a $(G,\bil)$-graded-monoidal category.

Recall the following:

\begin{definition}
  The \emph{opposite category $\cC^{\op}$} is the $(G,\bil)$-graded-monoidal category with the same objects as $\cC$, hom-spaces given by $\Hom_{\cC^{\op}}(i,j) = \Hom_{\cC}(j,i)$, the same monoidal structure, but with composition given by
  \[u\circ_{\cC^{\op}} v\coloneqq v\circ_{\cC} u.\]  
\end{definition}

Informally, the opposite category is the vertical symmetry of $\cC$, obtained by flipping the diagrams around a horizontal axis.





\medbreak

We first formalize vertical symmetry generically.
Fix a linear monoidal rewriting system $\sS = (\sX_0,\sX_1;\sR,\sE)$, a linear sub-system $\lT\subset\sS$ and an isp-order $\succ$ on the set of monomials $\sX^*$.
Denote by $\sX^l/\sE$ the free $(G,\bil)$-graded-monoidal category generated by $(\sX_0,\sX_1)$. Let $\ring^\times\sX^*$ be the set of monomials together with the data of an invertible scalar.
There are canonical ways of viewing $\sX^*$ and $\ring^\times\sX^*$ as sitting in $\sX^l$. We (momentarily) abuse notation and write $\sX^*$ and $\ring^\times\sX^*$ for their image in $\sX^l/\sE$.

A \defnemph{vertical symmetry of $(\sX_0,\sX_1)$} is the data of a graded-monoidal isomorphism
\[\iota\colon \sX^l/\sE\to (\sX^l/\sE)^{\op},\]
such that $\iota(\sX^*)\subset\ring^\times\sX^*$. In other words, $\iota$ sends diagrams to diagrams, up to invertible scalar.

Let $r_1,r_2\in\sR$.
We say that \defnemph{$r_1$ and $r_2$ are $\iota$-symmetric} if $\iota(s(r_1)) = \lambda s(r_2)$ and $\iota(t(r_1)) = \lambda t(r_2)$ for some invertible scalar $\lambda\in\ring^\times$.
In this case, we write $\iota(r_1)\coloneqq\lambda r_2$.
If $r_1$ and $r_2$ are $\iota$-symmetric, we can extend the symmetry to any $\sS$-rewriting step of type $r_1$, setting $\iota(\tau\Gamma[r_1]+v) = \tau\iota(\Gamma)[\iota(r_1)]+\iota(v)$.

We now turn to vertical symmetry in the context of categories of affine Brauer type.

\begin{definition}
  \label{defn:vertical_symmetries}
  Let $(\affS,\affT)$ be a rewriting system of affine Brauer type and $\iota$ be a vertical symmetry on $\affX$.
  We say that \defnemph{the vertical symmetry $\iota$ extends to $\affS$} if:
  \begin{itemize}
    \item R2 move, R3 move, leftward kink, bubble evaluation and bubble slide are $\iota$-self-symmetric;
    \item zigzags, cap and cup slides, cap and cup pulls, upward and downward kinks, and appropriate pairings of dot slides, are each other $\iota$-symmetric.
  \end{itemize}
\end{definition}
 
\begin{proof}[Proof of \cref{prop:vertical_symmetries}, assuming \cref{thm:affine_Grobner_system_iff_critical_branchings}]
  If $(\ssf,\ssg)$ is a critical branching and $(\ssf',\ssh,\ssg')$ is a $\rhd$-tamed $\affT$-preconfluence for $(\ssf,\ssg)$, then $(\iota(\ssf'),\iota(\ssh),\iota(\ssg'))$ is a $\rhd$-tamed $\affT$-preconfluence for the branching $(\iota(\ssf),\iota(\ssg))$. (Note that $\rhd$ is invariant under $\iota$.)
  The only subtlety is the exceptional branching: the vertical symmetry of the first exceptional branching is
  $\tikzpic{
    \uli[0][2] \ca[1][2]
    \dlcr[0][1]\AntiDiag[0][1][<-] \li[2][1]
    \li \cro[1][0]
    \custop \dli[2][0]
    \bul[0][1]\node[below left=-2pt] at (0,1) {\scriptsize $k$};
  }[scale=.5]$,
  which is a twin (see \cref{subsubsec:twin_branchings}) of the second exceptional branching.
\end{proof}

\subsection{Global analysis}
\label{subsec:global_analysis}

In this subsection, we prove that if most minimal branchings are $\rhd$-tamely $\affT$-preconfluent, then every overlapping branching is $\succ$-tamely $\affS$-congruent.

\begin{proposition}
  \label{prop:global_analysis}
  Assume that:
  \begin{enumerate}[(i)]
    \item every minimal non-affine branching (\cref{subfig:all_overlapping_branchings_non_affine}), with either $\tikzpic{\indexNF[][][?]}[scale=.6] = \tikzpic{\cro}[scale=.6]$ or $\tikzpic{\indexNF[][][?]}[scale=.6] = \tikzpic{\li}[scale=.6]$ for indexed branchings,
    \item every minimal dot-slide branching (\cref{subfig:all_overlapping_branchings_affine}),
    \item every minimal bubble-slide branching (\cref{fig:minimal_bubble_slide}),
    \item and every exceptional branching (\cref{subfig:critical_branchings_exceptional}),
  \end{enumerate}
  is $\rhd$-tamely $\affT$-preconfluent.
  Then every non-affine branching (\cref{subfig:all_overlapping_branchings_non_affine}), dot-slide branching (\cref{subfig:all_overlapping_branchings_affine}) and bubble-slide branching (\cref{subfig:all_overlapping_branchings_bubbles}) is $\succ$-tamely $\affS$-congruent.
\end{proposition}

The primary rewriting tool is the transitivity of tamed congruence; the core combinatorial subtlety is to understand how contextualization interacts with traffic rules.

The proof of the proposition consists in three lemmas.
As a preliminary, note that by the context-agnostic contextualization lemma (\cref{lem:RW_contextualization_lemma_agnostic}), every non-affine branching, with either $\tikzpic{\indexNF[][][?]}[scale=.6] = \tikzpic{\cro}[scale=.6]$ or $\tikzpic{\indexNF[][][?]}[scale=.6] = \tikzpic{\li}[scale=.6]$ for indexed branchings, is $\succ$-tamely $\affS$-congruent.
In \cref{lem:dot-slide_branching_tamed_congruent}, we first show that every dot-slide $\affT$-branching is $\succ$-tamely $\affS$-congruent.
We then show every non-affine branching is $\succ$-tamely $\affS$-congruent in \cref{lem:indexed_branching_affine_Brauer}, and finally show that every bubble-slide branching is $\succ$-tamely $\affS$-congruent in \cref{lem:bubble-slide_branching_tamed_congruent}.

\subsubsection{Dot-slide branchings}
\label{subsubsec:dot-slide_branchings}

\begin{lemma}
  \label{lem:dot-slide_branching_tamed_congruent}
  Assume that every contextualization of the following branchings
  \begin{gather*}
  \tikzpic{
    \uli[0][1]\uli[1][1]\ca[2][1]
    \cro\lili[2][0]
    \dli[0][0]\cu[1][0]\dli[3][0]
  }[scale=.7]
  \qquad
  \tikzpic{ 
    \cro[0][1]\ca[2][1]
    \li[0][0]\cro[1][0]\li[3][0]
    \cu[0][0]\dli[2][0]\dli[3][0]
  }[scale=.7]
  \qquad
  \tikzpic{
    \ca[0][2]
    \li[0][1]\li[1][1]\ca[2][1]
    \li\cro[1][0]\li[3][0]
    \cu\dli[2][0]\dli[3][0]
  }[scale=.7]
\end{gather*}
is $\succ$-tamely $\affS$-congruent and that every minimal dot-slide branching and every exceptional branching is $\rhd$-tamely $\affT$-pre\-con\-fluent.
  Then every dot-slide $\affT$-branching is $\succ$-tamely $\affS$-congruent.
\end{lemma}

\begin{proof}
  Let $(f,g)=\Gamma[\ssf,\ssg]$ be a dot-slide branching where $(\ssf,\ssg)$ is minimal and $\ssf$ is a dot-slide.
  Denote by $\mathsf{d}$ the source of $(\ssf,\ssg)$ and set $D=\Gamma[\mathsf{d}]$.
  By hypothesis, the branching $(\ssf,\ssg)$, for the traffic rules induced by $\Gamma$, admits a $\rhd$-tamed $\affT$-preconfluence $(\ssf',\ssh,\ssg')$; in fact, it must be a $\strongsucc$-tamed $\affT$-preconfluence since the source $\mathsf{d}$ does not have any bubble.
  Tamed preconfluences of minimal dot-slide branchings have the following schematics (see \cref{subsubsec:example_tamed_preconfluence_dotslide} for some examples):
  \begin{gather*}
  \begin{tikzcd}[ampersand replacement=\&,row sep=small]
    \& \bullet \& \bullet \& t_{\strongsucc}(\ssf_2') + \smallO_{\strongsucc}(\mathsf{d}) \\
    \mathsf{d} \\
    \& \bullet \&\& t_{\strongsucc}(\ssg') + \smallO_{\strongsucc}(\mathsf{d})
    \arrow["{*}"{description},"\ssf_1'", from=1-2, to=1-3]
    \arrow["\ssf_2'", from=1-3, to=1-4,thickcd]
    \arrow["{\ssh_{\strongsucc}}", shift left=8, tail reversed, dottedcd, from=1-4, to=3-4]
    \arrow[shift right=8, equals, from=1-4, to=3-4]
    \arrow["\ssf", curve={height=-12pt}, from=2-1, to=1-2]
    \arrow["\ssg"', curve={height=12pt}, from=2-1, to=3-2,thickcd]
    \arrow["{*}"{description},"\ssg'"', from=3-2, to=3-4]
  \end{tikzcd}
  \\
  \text{dot-slide \& $\{\stRthree,\stzz,\stcasl,\stcusl\}$}
  \\[2ex]
  \begin{IEEEeqnarraybox}{CcC}
    \begin{tikzcd}[ampersand replacement=\&,row sep=small]
      \& \bullet \& \bullet \& \smallO_{\strongsucc}(\mathsf{d}) \\
      \mathsf{d} \\
      \& \bullet \&\& \smallO_{\strongsucc}(\mathsf{d})
      \arrow["{*}"{description},"\ssf_1'", from=1-2, to=1-3]
      \arrow["\ssf_2'", from=1-3, to=1-4,thickcd]
      \arrow["{\ssh_{\strongsucc}}", tail reversed, dottedcd, from=1-4, to=3-4]
      \arrow["\ssf", curve={height=-12pt}, from=2-1, to=1-2]
      \arrow["\ssg"', curve={height=12pt}, from=2-1, to=3-2,thickcd]
      \arrow[equals, from=3-2, to=3-4]
    \end{tikzcd}
    &\qquad&
    \begin{tikzcd}[ampersand replacement=\&,row sep=small]
      \& \bullet \& t_{\strongsucc}(\ssf') + \smallO_{\strongsucc}(\mathsf{d}) \\
      \mathsf{d} \\
      \& \bullet \& t_{\strongsucc}(\ssg') + \smallO_{\strongsucc}(\mathsf{d})
      \arrow["\ssf'", from=1-2, to=1-3]
      \arrow["{\ssh_{\strongsucc}}", shift left=8, tail reversed, dottedcd, from=1-3, to=3-3]
      \arrow[shift right=8, equals, from=1-3, to=3-3]
      \arrow["\ssf", curve={height=-12pt}, from=2-1, to=1-2]
      \arrow["\ssg"', curve={height=12pt}, from=2-1, to=3-2]
      \arrow["\ssg'"', from=3-2, to=3-3]
    \end{tikzcd}
    \\
    \text{dot-slide \& $\{\stRtwo,\stcapu,\stcupu,\stuk,\stlk,\stdk\}$}
    &&
    \text{dot-slide \& dot-slide}
  \end{IEEEeqnarraybox}
  \end{gather*}
  Here a thick arrow denotes a context-agnostic rewriting step.
  Recall from \cref{lem:contextualization_for_strongsucc} that if $\Gamma[\ssf']$ and $\Gamma[\ssg']$ belong to $\affT$, then $(\Gamma[\ssf'],\Gamma[\ssh],\Gamma[\ssg'])$ is a $\strongsucc$-tamed $\affT$-preconfluence for $(\Gamma[\ssf],\Gamma[\ssg])$.
  In particular, if $(\ssf,\ssg)$ is as in the third schematic, then $(\Gamma[\ssf],\Gamma[\ssg])$ is $\strongsucc$-tamely $\affT$-preconfluent for any context $\Gamma$, since $\Gamma[\ssf]$ and $\Gamma[\ssg]$ are dot slides in $\affT$ that do no affect traffic rules.

  Consider the bottom branch of the first schematic.
  Recall that $\affT_\nonaff$ denotes non-affine rewriting steps.
  Denote by $\succ_\nonaff$ the lexicographic order for $(\nbcro,\nbcc,\distcc,\distcro)$ and by $\asymp_{\mathrm{dot}}$ the equivalence relation such that $x\asymp_{\mathrm{dot}}y$ if and only if $x$ and $y$ only differ by the position of dots.
  Note that:
  \begin{enumerate}[(i)]
    \item if $\ssf$ is a dot-slide, then $s(\ssf)\asymp_{\mathrm{dot}} t_{\strongsucc}(\ssf)$;
    \item $\succ_\nonaff$ is compatible with $\affT_\nonaff$;
    \item if $x\asymp_{\mathrm{dot}}y$ and $z\succ_\nonaff x$, then $z\succ_\nonaff y$.
  \end{enumerate} 
  Now, the rewriting sequence $g'=\Gamma[\ssg']$ only consists of dot-slides, so fact (i) gives that $t_{\strongsucc}(g)\asymp_{\mathrm{dot}} x$ for each $x$ the $\strongsucc$-leading term at a given step in $g'$.  
  Then, it follows from fact (ii) that $D\succ_\nonaff t_{\strongsucc}(g)$, and from fact (iii) that $D\succ_\nonaff x$ for each such $x$, so that the branch $g'$ is $\succ_\nonaff$-tamed by $D$.
  Finally, the diagram $D$ does not have any bubble, so that ``$\succ_\nonaff$-tamed by $D$'' is equivalent to ``$\succ$-tamed by $D$''.

  It remains to consider the top branch of the first and second schematics.
  If $\Gamma[\ssf_1']$ belongs to $\affT$, we can conclude; if not, that means the traffic rules induced by $\Gamma$ are not compatible with some dot-slides.
  Since dot-slides do not affect traffic rules and $f$ is also a dot-slide, this can only happen if a dot slides through a stop sign, which must lie on a cup in $\mathsf{d}$.
  This can happen only in one of the following situations:
  \begin{gather}
    \label{eq:situations_stop_sign}
    \newcommand{\arOrNot}[1]{#1}
    \tikzpic{
      \dli \lca[0][0]\ca[0][0] \custop[1][0] \uli[2][0][\arOrNot{->}]
    }[scale=.8]
    \qquad
    \tikzpic{
      \dli[2][0] \custop \cro[1][0] \li \li[2][1][\arOrNot{->}] \urcr[0][1] 
      \AntiDiag[0][1]
    }[scale=.8]
    \qquad
    \tikzpic{ 
      \ulcr \dli \custop[1][0] \li[2][0][\arOrNot{->}]
      \Diag[0][0]
    }[scale=.8]
    \qquad
    \tikzpic{
      \ca[0][1]\uli[2][1][\arOrNot{->}]
      \li[0][0] \drcr[1][0] 
      \custop\dli[2][0]
      \Diag[1][0]
    }[scale=.8]
    \qquad 
    \tikzpic{
      \urcr
      \custop
      \Diag[0][0][\arOrNot{->}]
    }[scale=.8]
  \end{gather}
  Denote by $S$ the strand where the dot slides.

  Consider the first situation in \eqref{eq:situations_stop_sign}. Actually, the diagram cannot have a stop sign, as otherwise we can find another generator of $S$ continuing downward ($S$ must be closed), and this would prevent the cup from being \textsc{lint}-minimal.
  Note that this argument was used before in the proof of \cref{lem:defn_normal_form_avoid_stop_sign_patterns}; we will continue to use similar arguments.
  We adopt the same notation $\Gamma=T\circ(\id_a\otimes -\otimes \id_b)\circ B$, where we call $T$ the \emph{top diagram} and $B$ the \emph{bottom diagram}. 
  We can assume that $B$ and $T$ are normal forms by rewriting them independently; moreover, we can interchange generators between $B$ and $T$ if they do not intersect with $D$, by branchwise $\affE$-congruence.

  Consider the second, third and fourth situations in \eqref{eq:situations_stop_sign}.
  Follow the rightmost top endpoint of $\mathsf{d}$ upward in $T$.
  If we reach a cap from its left strand, this strand is a straight line from $\mathsf{d}$ to the cap, and we can pull the cap down.
  This leads to a non-affine forbidden pattern, and hence a third $\affT$-rewriting step $h$, which creates an independent branching $(f,h)$ (assuming $f$ is the dot slide) and a non-affine branching $(h,g)$; the latter is one of the three non-affine branchings in the hypotheses of the lemma.
  Following $\succ$-tamed congruence of independent branchings (\cref{lem:extended_MLRS_independent_branching}) and of these non-affine branchings, we conclude with transitivity of tamed congruence.

  Assume instead we reach a cap from its right strand.
  If we find that it continues downward in $B$ to the left of $\mathsf{d}$, the cup cannot be \textsc{lint}-minimal; hence this cap connects to the top of $\mathsf{d}$ in a straight line.
  This eliminates the fourth situation, and up to contextualization, gives the following four possibilities, using that $T$ is a normal form:
  \begin{gather*}
    \tikzpic{
      \dli[2][0] \custop \cro[1][0] \li \li[2][1] \urcr[0][1] 
      \AntiDiag[0][1][<-]
      \uli[0][2]\ca[1][2]\bul[1][2]
      \node[right=-1pt] at (1,2) {\scriptsize $k$};
    }[scale=.8]
    \qquad
    \tikzpic{
      \dli[2][0] \custop \cro[1][0] \li \li[2][1] \urcr[0][1] 
      \AntiDiag[0][1][<-]
      \ca[0][3]\uli[2][3]
      \li[0][2]\cro[1][2]
      \bul[1][2]\node[right=-1pt] at (1,2) {\scriptsize $k$};
      \bul[0][2]\node[left=-1pt] at (0,2) {\scriptsize $l$};
    }[scale=.8]
    \qquad
    \tikzpic{ 
      \ulcr \dli \custop[1][0] \li[2][0] 
      \Diag[0][0][<-]
      \uli[0][1]\ca[1][1]
      \bul[1][1]\node[right=-1pt] at (1,1) {\scriptsize $k$};
    }[scale=.8]
    \qquad
    \tikzpic{ 
      \ulcr \dli \custop[1][0] \li[2][0] 
      \Diag[0][0][<-]
      \ca[0][2]\uli[2][2]
      \li[0][1]\cro[1][1]
      \bul[1][1]\node[right=-1pt] at (1,1) {\scriptsize $l$};
      \bul[0][1]\node[left=-1pt] at (0,1) {\scriptsize $k$};
    }[scale=.8]
  \end{gather*}
  Using horizontal symmetry (see \cref{subsubsec:horizontal_symmetry}), the tamed congruence of the first branching above follows from the tamed congruence of the second exceptional critical branching; see the schematics \cref{fig:horizontal_symmetries_D} in \cref{sec:horizontal_symmetries}.
  Consider the second case; here we can argue as in the previous paragraph. Namely, we follow the sole output of the diagram in $T$: since it belongs to a closed strand, it must meet a cap, which either leads to a third branch $h$ which is independent of both $f$ and $g$, or it leads to a contradiction with the \textsc{lint}-minimality of the cup.
  The third case can actually never occur, as the cup cannot be \textsc{lint}-minimal in this situation.
  Finally, in the fourth case, sliding one of the $l$ dots gives a rewriting step $h$ such that the branchings $(f,h)$ and $(h,g)$ are $\succ$-tamely congruent: via independent branchings if the dot slides up, and via two dot-slide branchings, whose tamed congruence was shown earlier in this proof, if it slides down.
  Hence, we may assume that $l=0$, and this gives the second exceptional branching.
  Using arguments similar to those above, one checks that if an exceptional branching is $\rhd$-tamely $\affT$-preconfluent, then any of its contextualizations is $\succ$-tamely $\affS$-congruent.

  Finally, consider the last situation in \eqref{eq:situations_stop_sign}.
  Since the cup is \textsc{lint}-minimal, there is a sequence of left-directed interchanges that makes it the bottom generator of $S$.
  In fact, the remaining generators in $\mathsf{d}$ can also be assumed to be at the bottom, and we can assume that we always interchange the diagram $\mathsf{d}$ as a whole. (We used the same argument in \cref{lem:affine_Brauer_equivalence_of_presentations}.)
  We can then view this sequence of interchanges as a sequence of branchwise $\affE$-congruences, and assume that $\mathsf{d}$ as a whole lies at the bottom of $S$.
  It follows that the rest of the strand $S$ lies in $T$ and hence is a dotted cap (using that $T$ is a normal form). This gives the last exceptional branching.
\end{proof}

\subsubsection{Indexed branchings}
\label{subsubsec:indexed_branchings}

Consider the branchings from \cref{subfig:all_overlapping_branchings_non_affine} with a box $\tikzpic{\indexNF[][][?]}[scale=.6]$: recall from \cref{subsubsec:enumeration_overlapping_brauer_type} that we call them \emph{indexed branchings}.
A priori, any choice of diagram for $\tikzpic{\indexNF[][][?]}[scale=.6]$ gives a branching. In this subsection, we show how to drastically reduce the number of such branchings.

\begin{lemma}
  \label{lem:indexed_branching_affine_Brauer}
  Assume that non-indexed non-affine branchings, dot-slide branchings and indexed branchings where $\tikzpic{\indexNF[][][?]}[scale=.6] = \tikzpic{\cro}[scale=.6]$ or $\tikzpic{\indexNF[][][?]}[scale=.6] = \tikzpic{\li}[scale=.6]$ are $\succ$-tamely congruent.
  Then every indexed non-affine branching is $\succ$-tamely congruent.
\end{lemma}

To illustrate the main idea of the proof, consider the indexed branching
$(f,g)=\tikzpic{\uli[0][1]\ca[1][1]
\cro \indexNF[2][0][?]
\dli[0][0][->] \cu[1][0]}[scale=.6]$
and the situations where $\tikzpic{\indexNF[][][?]}[scale=.6] = \tikzpic{\bli}[scale=.6]$ and $\tikzpic{\indexNF[][][?]}[scale=.6] = \tikzpic{\ca\li[2][0]}[scale=.6]$:
\begin{gather*}
  \tikzpic{ 
    \uli[0][1]\ca[1][1]
    \cro \li[2][0] \bul[2][1]
    \dli[0][0][->] \cu[1][0]
    \draw[branch_overlay]
    (1-.3,1-.5) rectangle (2+.3,1.5);
  }[scale=.8]
\qquad\an\qquad
  \tikzpic{ 
    \uli[0][1]\draw (1,1) to[out=90,in=90] (4,1) to (4,-.5);
    \cro \ca[2][0]
    \dli \cu[1][0] \dli[3][0]
    \draw[branch_overlay]
    (1-.3,-.5) rectangle (3+.3,+.5);
  }[scale=.8]\;.
\end{gather*}
In each situation, there exists a rewriting step $h$ which is independent of one of the branches of $(f,g)$, and defines with the other branch either a non-indexed branching from \cref{subfig:all_overlapping_branchings_non_affine} or a branching from \cref{subfig:all_overlapping_branchings_affine} (in the diagrams, the source of $h$ is shaded).
Transitivity of tamed congruence allows us to conclude.

\begin{proof}
  Let $(f,g)$ be an indexed branching, with $f$ on top and $g$ on the bottom, and $\Gamma$ a context.
  Throughout we use independent rewriting to freely rewrite $\tikzpic{\indexNF[][][?]}[scale=.6]$; in particular, we can assume $\tikzpic{\indexNF[][][?]}[scale=.6]$ is a normal form, and following \cref{thm:normal_form_are_standard_diagrams}, a standard diagram.
  Moreover, if the first input (resp.\ output) is not an endpoint of the top (resp.\ bottom) generator of $\tikzpic{\indexNF[][][?]}[scale=.6]$, we can interchange this generator with the rest of the diagram; formally, we are using the canonical branchwise $\sE$-congruence induced by the interchange.

  The situation where the first output connects with the first input as a straight line in $\tikzpic{\indexNF[][][?]}[scale=.6]$ is part of the hypotheses.

  Assume the first output of $\tikzpic{\indexNF[][][?]}[scale=.6]$ is oriented outward, with a dot in $\tikzpic{\indexNF[][][?]}[scale=.6]$. Pushing the dot outside $\tikzpic{\indexNF[][][?]}[scale=.6]$ defines a rewriting step $h$ which overlaps with $f$ and is independent of $g$.
  Combining the hypothesis and $\succ$-tamed congruence of independent branchings (\cref{lem:extended_MLRS_independent_branching}), we conclude that $(f,g)$ is $\succ$-tamely congruent by transitivity of $\succ$-tamed congruence.
  Similar arguments apply to other situations where the top or bottom generator of $\tikzpic{\indexNF[][][?]}[scale=.6]$ is a dot.

  Assume the top generator of $\tikzpic{\indexNF[][][?]}[scale=.6]$ is a cup, with its left strand being the first output of $\tikzpic{\indexNF[][][?]}[scale=.6]$. Since the left strand is not a straight line, the diagram is not reduced, and there exists a rewriting step $h$ that overlaps with $f$. We argue as in the previous paragraph, and similar arguments apply to related situations, including the ones with a cap.

  Finally, assume the top generator of $\tikzpic{\indexNF[][][?]}[scale=.6]$ is a crossing. We can assume $\tikzpic{\indexNF[][][?]}[scale=.6]$ is a braid-standard diagram. Given the form of the braid-standard context (\cref{subfig:braid_standard_context}), we can reduce to $\tikzpic{\indexNF[][][?]}[scale=.6]$ being a single crossing, which is part of the hypotheses.
\end{proof}

\subsubsection{Bubble-slide branchings}
\label{subsubsec:bubble_slides}

In \cref{subfig:all_overlapping_branchings_bubbles}, we gave two schematics for branchings involving a bubble slide:
\begin{gather*}
  \tikzpic{
    \clip (0,0) circle (2cm);
    \draw[dashed] (0,0) circle (2cm);
    \draw[out=45,in=-90] (-2,-2) to (0,2);
    \draw (-.8,2) to[out=-90,in=180] (0,1) to[out=0,in=-90] (.8,2);        
    \draw (-.5,-2) to[out=90,in=0] (-2,-.5);        
    \begin{scope}[shift={(1,-.4)}]
      \bul[-.5][0]\ca[-.5][0]\cu[-.5][0]
      \node[left=-1pt] at (-.5,0) {\scriptsize $k$};
    \end{scope}
    \node[rotate=60] at (-.8,.4) {\ldots};
    \draw[very thick,red,->] (0,0) to (-.7,.4);
  }[scale=.8]
  \qquad\an\qquad
  \tikzpic{
    \clip (0,0) circle (2cm);
    \draw[dashed] (0,0) circle (2cm);
    \draw (-2,-.5) to[out=30,in=60] (-.5,-2);
    \draw[very thick,red,->] (-.5,-.5) to (-1.2,-1.2);
    \draw (-1,1.7) to[out=-90,in=-90] (1,1.7);
    \draw[very thick,red,->] (0,.7) to (0,1.7);
    \bul[-.5][0]\ca[-.5][0]\cu[-.5][0]
    \node[left=-1pt] at (-.5,0) {\scriptsize $k$};
  }[scale=.8]
  \;.
\end{gather*}
The first picture is a schematic for overlaps between a rewriting step with an inner region and a bubble-slide. The second picture is a schematic for overlaps between two bubble-slides.

\begin{lemma}
  \label{lem:bubble-slide_branching_tamed_congruent}
  Assume that every minimal bubble-slide branching (\cref{fig:minimal_bubble_slide}) is $\rhd$-tamely $\affT$-preconfluent. Then every bubble-slide branching (\cref{subfig:all_overlapping_branchings_bubbles}) is $\succ$-tamely $\affS$-congruent.
\end{lemma}

\begin{proof}
  \def\tpscl{1}
  \Cref{lem:contextualization_bubbles} shows that the contextualization of any minimal bubble-slide branching is $\succ$-tamely $\affS$-congruent.

  Let $(f,g)=\Gamma[\ssf,\ssg]$ be a branching as in the first schematic above, with $f$ being the bubble-slide.
  Then there exists a third branch $h$, consisting in sliding the bubble through an adjacent piece of strand and landing in a region which is not an inner region of the source of $\ssg$. This third branch $h$ is independent of $g$ and defines a contextualization of a minimal bubble-slide branching with $f$.
  Both branchings are $\succ$-tamed $\affS$-congruent by the independent-branching lemma (\cref{lem:extended_MLRS_independent_branching}) and the previous paragraph. Transitivity of tamed congruence concludes.

  Let $(f,g)$ be a bubble-slide branching where both $f$ and $g$ are bubble-slides, as in the second schematic. Consider the outer boundary of the region where the bubble is, and denote by $P_f$ and $P_g$ the points where the bubble slides according to $f$ and $g$, respectively.
  Move along that boundary from $P_f$ to $P_g$. For at least one direction, we only move along strands, and not along the boundary of the diagram itself; this would otherwise contradict the properties of tunnels.
  Therefore, that path can be decomposed as going past a cap, a cup or a crossing. At each stage, there exists another bubble-slide rewriting step, leading to a contextualization of a minimal bubble-slide branching.
  We again conclude by transitivity of tamed congruence.
\end{proof}

\subsection{Proof of the Main Theorem}
\label{subsec:proof_main_thm}

We can finally give the proof of \cref{thm:affine_Grobner_system_iff_critical_branchings}, which we reproduce here for convenience:

\thmAffineGrobnerSystem*

We first gather the results of this section in the following proposition:

\begin{proposition}
  \label{prop:critical_branchings_iff_grobner}
  Let $(\affS,\affT,\succ)$ be a rewriting system of affine Brauer type. Assume that $\bubbleset$-bubbles have trivial grading,
  that for every $k\in\bN$ there exists a $\rhd$-tamed $\affT$-congruence
  \begin{gather*}
    (1-\VARdcc\VARdccbis\inter{\bullet}{\bullet}^k)\;
    \tikzpic{\LOTev[0][0][k{+}1]}
    \quad
    \overset{k{+}1}{\longtamearrow}
    \quad
    (\VARdccbis\VARdccLOT\inter{\bullet}{\bullet}^k\;+\;\VARdccLOTbis)\;
    \tikzpic{\LOTev}\;,
  \end{gather*}
  and that every critical branching from \cref{fig:critical_branchings} is $\rhd$-tamely $\affT$-preconfluent.
  Then $(\affS,\affT,\succ)$ is a monoidal Gröbner system.
\end{proposition}

\begin{proof}
  We prove the four properties of a monoidal Gröbner system (\cref{defn:RW_monoidal_grobner_system}).
  For property (i), the isp-order $\succ$ is well-founded and compatible with $\affT$ by \cref{lem:order_strongly_compatible}.
  For property (iii), since $\bubbleset$-bubbles have trivial grading,
  it follows from \cref{lem:RW_bubble_trivial_grading_implies_coherence} that $\sE$ is coherent on monomial $\affT$-normal forms.
  Property (iv) was proved in \cref{lem:affine_Brauer_equivalence_of_presentations}.

  Finally, we prove property (ii). Following \cref{prop:classification_overlapping_branchings}, \cref{prop:local_analysis} and \cref{prop:global_analysis}, every overlapping branching is
  $\succ$-tamely $\affS$-congruent.
  Since independent branchings are also $\succ$-tamely $\affS$-congruent (\cref{lem:extended_MLRS_independent_branching}), every monomial branching is $\succ$-tamely congruent.
  The tamed linear Newman's lemma (\cref{lem:RW_linear_tame_newman}) implies that $\affT^+$ is confluent.
\end{proof}

We can now prove \cref{thm:affine_Grobner_system_iff_critical_branchings}.

\begin{proof}[Proof of \cref{thm:affine_Grobner_system_iff_critical_branchings}]
  By \cref{prop:critical_branchings_iff_grobner}, the data $(\affS,\affT,\succ)$ defines a monoidal Gröbner system.
  Following \cref{thm:RW_monoidal_grobner_system}, any set of interchange representatives of monomial $\affT$-normal forms defines a hom-basis of the category presented by $\affS$.
  Thanks to \cref{thm:normal_form_are_standard_diagrams}, this gives a hom-basis consisting of a representative set of $(\trafficend,\bubbleset)$-reduced matchings.
  Finally, it follows from \cref{lem:basis_induced_matching_works_for_all} that every representative set of $(\trafficend,\bubbleset)$-reduced matchings defines a hom-basis.
\end{proof}

%% file: sections/analysis_affine_brauer_new.tex
\subsection{The affine Brauer category}
\label{subsec:affine_brauer_category}

Throughout we assume that $2$ is invertible in the communtative ring $\ring$.

\begin{definition}\label{alt:def:affBr}
  \renewcommand{\sca}{.4}
  \newcommand{\vspc}{.5ex}
  \newcommand{\hspc}{60mu}
  The \emph{affine Brauer category $\affBr$} \cite[Definition~1.2]{RS_AffineBrauerCategory_2019} is the linear monoidal category generated by a single object $\tikzpic{\li}$ and morphisms $\tikzpic{\cu}$, $\tikzpic{\ca}$, $\tikzpic{\cro}$, $\tikzpic{\bli}$ subject to the following relations:
  \begin{gather*}
    \tikzpic{
      \cro[0][1]
      \cro 
    } 
    \;=\;
    \tikzpic{
      \lili[0][1] 
      \lili[0][0] 
    } 
    \mspace{\hspc}
    %
    \tikzpic{
      \cro[0][2] \li[2][2]
      \li[0][1] \cro[1][1] 
      \cro \li[2][0] 
    }[scale=.8]
    \;=\;
    \tikzpic{
      \cro[1][2] \li[0][2]
      \li[2][1] \cro[0][1] 
      \cro[1][0] \li[0][0] 
    }[scale=.8]
    \mspace{\hspc}
    \tikzpic{
      \cu\ca[0][0]
    } 
    \;=\;\VARev
    \\[\vspc]
    \tikzpic{
      \cu[1][0] \uli[2][0] 
      \dli \ca 
    }
    \;=\;
    \tikzpic{\li}
    %
    \mspace{\hspc}
    %
    \tikzpic{
      \ca[0][1]
      \cro 
    }
    \;=\;
    \tikzpic{\ca}
    \mspace{\hspc}
    %
    \tikzpic{
      \uli[0][1] \ca[1][1] 
      \cro \li[2][0] 
    }
    \;=\;
    \tikzpic{ 
      \ca[0][1] \uli[2][1]
      \li \cro[1][0] 
    }
    \\[\vspc]
    \tikzpic{
      \ulcr
    }
    \;=\;
    \tikzpic{\drcr}
    \;-\;
    \tikzpic{\lili} 
    \;+\;
    \tikzpic{\cc} 
    \mspace{\hspc}
    \tikzpic{\rca}
    \;=\;-\;
    \tikzpic{\lca}
  \end{gather*}
\end{definition}

In this section, we use the Main Theorem (\cref{thm:affine_Grobner_system_iff_critical_branchings}) to prove the following:

\begin{proposition}\label{alt:prop:affBr_main}
  The affine Brauer category $\affBr$ is a category of affine Brauer type with $2\bN_{>0}$-bubbles (\cref{defn:category_affine_brauer_type}).
\end{proposition}

The proof also serves as a blueprint for other categories of affine Brauer type.

\begin{lemma}\label{alt:lem:relations}
  \renewcommand{\sca}{.4}
  \newcommand{\vspc}{.5ex}
  \newcommand{\hspc}{60mu}
  In $\affBr$, the following relations hold:
  \begin{gather*}
    \tikzpic{
      \dli[2][0] \ca[1][0] 
      \uli \cu 
    }
    \; = \;
    \tikzpic{\li}
    \mspace{\hspc}
    %
    \tikzpic{
      \ca \uli[2][0] 
      \li[0][-1] \cro[1][-1]
      \cro[0][-2] \li[2][-2]
    }[scale=.8]
    =
    \tikzpic{\li\ca[1][0]}
    \mspace{\hspc}
    %
    \tikzpic{ 
      \cro \li[2][0]
      \dli \cu[1][0] 
    }
    =
    \tikzpic{
      \li[0][0] \cro[1][0]
      \cu[0][0] \dli[2][0]
    }
    \mspace{\hspc}
    %
    \tikzpic{
      \cro[0][1] \li[2][1] 
      \li[0][0] \cro[1][0]
      \cu \dli[2][0]
    }[scale=.8]
    =
    \tikzpic{\li\cu[1][1]}
    \\[\vspc]
    \tikzpic{
      \cro
      \cu
    }
    \;=\;
    \tikzpic{\cu}
    \mspace{\hspc}
    %
    \tikzpic{
      \ca[0][1] \uli[2][1]
      \li[0][0] \cro[1][0]
      \cu \dli[2][0]
    }[scale=.8]
    \;=\;
    \tikzpic{\li}
    \mspace{\hspc}
    \tikzpic{\urcr}
    \;=\; 
    \tikzpic{\dlcr}
    +\tikzpic{\lili} 
    -\tikzpic{\cc} 
    \mspace{\hspc}
    \tikzpic{\rcu \cu[0][0]}
    \;=\;-\;
    \tikzpic{\lcu \cu[0][0]}.
  \end{gather*}
\end{lemma}
To define a rewriting system of affine Brauer type as in \cref{fig:generic_grobner_of_affine_Brauer_type} associated to $\affBr$, it remains to define the bubble evaluation and the bubble slide.
It will be convenient to package these rewriting steps using generating functions, as we now explain.

\subsubsection{Diagrams and their generating functions}
\label{alt:subsubsec:generating_functions}
 
Following an idea of Savage and Webster \cite{SW_BubblesAffineBrauer_2024}, we use generating functions to compute with bubbles.
For any rewriting system of affine Brauer type, we set:
\begin{gather*}
\Delta_k \coloneqq\tikzpic{\bubble[0][0][k]}\qquad\text{for all }k\geq 0
\end{gather*}
and
\begin{gather*}
\Delta(x)\coloneqq\sum_{k\geq0}x^k\,\Delta_k,
\qquad
\Delta_\subeven(x)\coloneqq\sum_{n\geq0}x^{2n}\,\Delta_{2n},
\qquad
\Delta_\subodd(x)\coloneqq\sum_{n\geq0}x^{2n+1}\,\Delta_{2n+1},
\\
\tikzpic{\GFbli} \coloneqq\sum_{k\geq0}x^k\,\tikzpic{\bli[0][0][][k]},
\qquad
\tikzpic{\GFblirev} \coloneqq\sum_{k\geq0}(-1)^k x^k\,\tikzpic{\bli[0][0][][k]}.
\end{gather*}
Here $x$ is a formal variable; formally, we are working over the ring $\ring[[x]]$ of formal power series (we sometimes further extend to formal Laurent series $\ring((x))$).
Note that by definition
\[\Delta(x) = \tikzpic{\GFbli \cu \ca[0][1] \li[1][0]}
\quad\an\quad\Delta(-x) = \tikzpic{\GFblirev \cu \ca[0][1] \li[1][0]}\;.\]
An identity in $\affBr$ given over $\ring[[x]]$ amounts to a family of identities for each natural number $k\geq 0$, obtained by extracting the coefficient of $x^k$.
For instance, in the case of $\affBr$, we have the following simple but useful identity:
\begin{gather}
  \label{alt:eq:identity_GF_dot}
  \tikzpic{
    \GFblirev[0][0]
    \GFbli[0][1]
  }
  =
  \frac{1}{2}\left(
  \tikzpic{
    \GFbli[0][0]
  }
  +
  \tikzpic{
    \GFblirev[0][0]
  }\right)
\end{gather}
To prove it, note that the generating functions of $\tikzpic{\GFbli}$ and $\tikzpic{\GFblirev}$ are respectively $\frac{1}{1-x\bft}$ and $\frac{1}{1+x\bft}$, where $\bft$ is another formal variable keeping track of the label of the dot.
The identity \eqref{alt:eq:identity_GF_dot} is a consequence of the equality $\frac{1}{1-x\bft}\frac{1}{1+x\bft} = \frac{1}{2}\left(\frac{1}{1-x\bft}+\frac{1}{1+x\bft}\right)$.

We express the dot-slide identities in terms of generating functions:

\begin{lemma}
  In $\affBr$, we have the following identities:
  \begin{gather*}
    \tikzpic{\uli[0][1]\uli[1][1]\cro\GFbul[0][1]\dli[0][0]\dli[1][0]}
    =
    \tikzpic{\uli[0][1]\uli[1][1]\cro\GFbul[1][0]\dli[0][0]\dli[1][0]}
    -x\;
    \tikzpic{\uli[0][1]\uli[1][1]\GFbul[0][1]\lili\dli[0][0]\dli[1][0]\GFbul[1][0]}
    +x\;
    \tikzpic{\uli[0][1]\uli[1][1]\GFbul[0][1]\cc\dli[0][0]\dli[1][0]\GFbul[1][0]}
    \mspace{60mu}
    \tikzpic{\uli[0][1]\uli[1][1]\cro\GFbul[1][1]\dli[0][0]\dli[1][0]}
    =
    \tikzpic{\uli[0][1]\uli[1][1]\cro\GFbul[0][0]\dli[0][0]\dli[1][0]}
    +x\;
    \tikzpic{\uli[0][1]\uli[1][1]\GFbul[1][1]\lili\dli[0][0]\dli[1][0]\GFbul[0][0]}
    -x\;
    \tikzpic{\uli[0][1]\uli[1][1]\GFbul[1][1]\cc\dli[0][0]\dli[1][0]\GFbul[0][0]}
    \\
    \tikzpic{\cu\GFbul\uli\uli[1][0]}
    =
    \tikzpic{\cu\uli\GFbulrev[1][0]\uli[1][0]}
    \qquad
    \tikzpic{\cu\GFbulrev\uli\uli[1][0]}
    =
    \tikzpic{\cu\uli\GFbul[1][0]\uli[1][0]}
    \qquad
    \tikzpic{\ca\GFbul\dli\dli[1][0]}
    =
    \tikzpic{\ca\dli\GFbulrev[1][0]\dli[1][0]}
    \qquad
    \tikzpic{\ca\GFbulrev\dli\dli[1][0]}
    =
    \tikzpic{\ca\dli\GFbul[1][0]\dli[1][0]}
  \end{gather*}
\end{lemma}

Using these identities, we compute that:

\begin{lemma}\label{alt:lem:bubble_identities}
  In $\affBr$, we have the following identities:
  \begin{gather*}
  \Delta_\subodd(x)
  =
  -
  \frac{x}{4}(\Delta(x)+\Delta(-x))
  +
  \frac{x}{2} \Delta(x)\Delta(-x)
  \\
  \tikzpic{\li}\;\Delta(x)
  =
  \Delta(x)\;\tikzpic{\li}
  +
  x\big(1-\frac{x}{2}\big)
  \left(
  \tikzpic{
    \GFbli[0][0]
    \GFbli[0][1]
  }
  -
  \tikzpic{
    \GFblirev[0][0]
    \GFblirev[0][1]
  }
  \right)
  +
  x^2
  \left(
  \tikzpic{
    \GFbli[0][0]
    \GFbli[0][1]
  }\;\Delta(x)
  -
  \Delta(x)\;
  \tikzpic{
    \GFblirev[0][0]
    \GFblirev[0][1]
  }
  \right)
  \end{gather*}
\end{lemma}

\begin{proof}
\renewcommand{\sca}{.4}
Indeed, using the identity \eqref{alt:eq:identity_GF_dot}:
\begin{IEEEeqnarray*}{rCl}
  \Delta(-x)
  &=&
  \tikzpic{\ca[0][1]\cro\GFbul[1][1]\cu}
  =
  \tikzpic{\ca[0][1]\cro\GFbul\cu}
  +
  x\;
  \tikzpic{\ca[0][1]\lili\GFbul[1][1]\GFbul\cu}
  -
  x\;
  \tikzpic{\ca[0][1]\cc\GFbul[1][1]\GFbul\cu}
  =
  \Delta(x)
  +
  \frac{x}{2}(\Delta(x)+\Delta(-x))
  -
  x\Delta(x)\Delta(-x)
\end{IEEEeqnarray*}
Note that the last two terms are odd functions: the coefficient of $x^k$ for $k$ even is always zero. If instead we restrict to odd monomials and use $\Delta_\subodd(-x)=-\Delta_\subodd(x)$, we get the statement.

On the other hand, we compute, again using  \eqref{alt:eq:identity_GF_dot}:
\begin{IEEEeqnarray*}{rCl}
  \Delta(x)\;\tikzpic{\li}
  &=&
  \tikzpic{
    \ca[0][2]\uli[2][2]
    \li[0][1]\cro[1][1]
    \cro \li[2][0]
    \dli \cu[1][0]
    \GFbul[0][1]
  }[scale=.8]
  =
  \tikzpic{
    \ca[0][2]\uli[2][2]
    \li[0][1]\cro[1][1]
    \cro \li[2][0]
    \dli \cu[1][0]
    \GFbul[1][0]
  }[scale=.8]
  -x\;
  \tikzpic{
    \ca[0][2]\uli[2][2]
    \li[0][1] \cro[1][1]
    \lili \li[2][0]
    \dli\cu[1][0]
    \GFbul[0][1]
    \GFbul[1][0]
  }[scale=.8]
  +x\;
  \tikzpic{
    \ca[0][2]\uli[2][2]
    \li[0][1] \cro[1][1]
    \cc \li[2][0]
    \dli\cu[1][0]
    \GFbul[0][1]
    \GFbul[1][0]
  }[scale=.8]
  =
  \tikzpic{\li}\;\Delta(x)
  -x
  \tikzpic{
    \ca[0][2]\uli[2][2]
    \li[0][1] \cro[1][1]
    \lili \li[2][0]
    \dli\cu[1][0]
    \GFbul[0][0]
    \GFbul[1][1]
  }[scale=.8]
  +x
  \tikzpic{
    \ca[0][2]\uli[2][2]
    \cro[0][1] \li[2][1]
    \li \cu[1][1]
    \dli
    \GFbul[0][2]
    \GFbulrev[0][0]
  }[scale=.8]
  \\
  &=&
  \tikzpic{\li}\;\Delta(x)
  -x
  \left(
  \tikzpic{
    \ca[0][2]\uli[2][2]
    \li[0][1] \cro[1][1]
    \lili \li[2][0]
    \dli\cu[1][0]
    \GFbul[0][0]
    \GFbul[2][2]
  }[scale=.8]
  -x
  \tikzpic{
    \ca[0][2]\uli[2][2]
    \li[0][1] \lili[1][1]
    \lili \li[2][0]
    \dli\cu[1][0]
    \GFbul[0][0]
    \GFbul[1][1]
    \GFbul[2][2]
  }[scale=.8]
  +x
  \tikzpic{
    \ca[0][2]\uli[2][2]
    \li[0][1] \cc[1][1]
    \lili \li[2][0]
    \dli\cu[1][0]
    \GFbul[0][0]
    \GFbul[1][1]
    \GFbul[2][2]
  }[scale=.8]
  \right)
  +x
  \left(
  \tikzpic{
    \ca[0][2]\uli[2][2]
    \cro[0][1] \li[2][1]
    \li \cu[1][1]
    \dli
    \GFbul[1][1]
    \GFbulrev[0][0]
  }[scale=.8]
  -x
  \tikzpic{
    \ca[0][2]\uli[2][2]
    \lili[0][1] \li[2][1]
    \li \cu[1][1]
    \dli
    \GFbul[0][2]
    \GFbul[1][1]
    \GFbulrev[0][0]
  }[scale=.8]
  +x
  \tikzpic{
    \ca[0][2]\uli[2][2]
    \cc[0][1] \li[2][1]
    \li \cu[1][1]
    \dli
    \GFbul[0][2]
    \GFbul[1][1]
    \GFbulrev[0][0]
  }[scale=.8]
  \right)
  \\
  &=&
  \tikzpic{\li}\;\Delta(x)
  -x
  \left(
  \tikzpic{
    \GFbli[0][0]
    \GFbli[0][1]
  }
  -
  \tikzpic{
    \GFblirev[0][0]
    \GFblirev[0][1]
  }
  \right)
  +x^2
  \left(
  \tikzpic{
    \GFbli[0][0]
    \GFblirev[0][1]
    \GFbli[0][2]
  }
  -
  \tikzpic{
    \GFblirev[0][0]
    \GFbli[0][1]
    \GFblirev[0][2]
  }
  \right)
  -x^2
  \left(
  \tikzpic{
    \GFbli[0][0]
    \GFbli[0][1]
  }\;\Delta(x)
  -
  \Delta(x)\;
  \tikzpic{
    \GFblirev[0][0]
    \GFblirev[0][1]
  }
  \right)
  \\
  &=&
  \tikzpic{\li}\;\Delta(x)
  -x\big(1-\frac{x}{2}\big)
  \left(
  \tikzpic{
    \GFbli[0][0]
    \GFbli[0][1]
  }
  -
  \tikzpic{
    \GFblirev[0][0]
    \GFblirev[0][1]
  }
  \right)
  -x^2
  \left(
  \tikzpic{
    \GFbli[0][0]
    \GFbli[0][1]
  }\;\Delta(x)
  -
  \Delta(x)\;
  \tikzpic{
    \GFblirev[0][0]
    \GFblirev[0][1]
  }
  \right)\;.
\end{IEEEeqnarray*}
This gives the statement.
\end{proof}

\subsubsection{The affine Brauer rewriting system}
We associate an affine Brauer rewriting system $\affS(\affBr)$ to the affine Brauer category $\affBr$.
Following the previous subsection, we use generating functions and avoid expanding formulas.
For that purpose, we define for every rewriting system of affine Brauer type:
\begin{gather}\label{alt:eq:defn_bbsl_bbev}
\bbslgenfun(x)\coloneqq\sum_{k\geq0}\,x^k\,\tikzpic{\LOTbbsl[0][0][k]}
\qquad\an\qquad
\evgenfun(x) \coloneqq \sum_{k\notin\bubbleset}x^{k}\,\tikzpic{\LOTev[0][0][k]}.
\end{gather}
For $\affBr$, we have $\bubbleset=2\bN_{>0}$, so that the latter sum is over odd natural numbers and zero.

\begin{definition}\label{alt:def:S_affBr}
Let $\affS(\affBr)$ be the rewriting system of affine Brauer type (\cref{fig:generic_grobner_of_affine_Brauer_type}) with set of bubbles $\bubbleset=2\bN_{>0}$ and the following data: all scalars equal $1$, except the dot slides through cups and caps, whose scalar is $-1$; we have
\begin{align*}
  \renewcommand{\sca}{.4}
  \tikzpic{\LOTRtwo} = \tikzpic{\lili\lili[0][1]}\;,\quad
  \tikzpic{\LOTuk} = \tikzpic{\ca}\;,\quad
  \tikzpic{\LOTdk} = \tikzpic{\cu}\;,
  \quad\an\quad
\tikzpic{\LOTdcro}\;=\;-\,\tikzpic{\lili}+\tikzpic{\cc}\;=\;-\,\tikzpic{\LOTdcrobis},
\end{align*}
the lower-order term of the bubble evaluation is given by
\begin{gather}\label{alt:eq:bbev_formula}
  \evgenfun(x)
  =
  \VARev
  -
  \frac{x}{4}(\Delta(x)+\Delta(-x))
  +
  \frac{x}{2} \Delta(x)\Delta(-x)\;;
\end{gather}
the lower-order term of the bubble slide is recursively defined by
\begin{gather}\label{alt:eq:bbsl_formula}
  \bbslgenfun(x)
  =
  x\big(1-\frac{x}{2}+x\Delta(x)\big)
  \left(
  \tikzpic{
    \GFbli[0][0]
    \GFbli[0][1]
  }
  -
  \tikzpic{
    \GFblirev[0][0]
    \GFblirev[0][1]
  }
  \right)
  +
  x^2
  \tikzpic{
    \GFbli[0][0]
    \GFbbsl[0][2]
    \GFbli[0][1]
  }\;;
\end{gather}
and the remaining lower-order terms vanish.
\end{definition}

This recovers \cref{fig:AB_S_rewriting_steps} in the introduction, with the addition of the $\rhd$-lower order term of the bubble-slide rewriting step.
Each $\tikzpic{\LOTev[0][0][2n{+}1]}$ has the shape prescribed in \cref{fig:generic_grobner_of_affine_Brauer_type}: extracting the coefficient of $x^{2n+1}$ in \eqref{alt:eq:bbev_formula}, each bubble in the right-hand side has at most $2n$ dots, since each term has a leading factor $x$.
By the same argument, each $\tikzpic{\LOTbbsl[0][0][k]}$ has the shape prescribed in \cref{fig:generic_grobner_of_affine_Brauer_type}, as the right-hand side of \eqref{alt:eq:bbsl_formula} is, by induction, a linear combination of straight lines carrying at most $k-1$ dots together with bubbles carrying at most $k-2$ dots.

\begin{corollary}
The linear monoidal rewriting system $\affS(\affBr)$ presents $\affBr$.
\end{corollary}
\begin{proof}
This follows from \cref{alt:lem:relations} and \cref{alt:lem:bubble_identities}.
\end{proof}

We define the following vertical symmetry on the generators of $\affBr$:
  \begin{gather}\label{alt:eq:vertsymm}
    \left(\tikzpic{\cu}\right)^{\op}
    \;=\;
    \tikzpic{\ca}
    \qquad
    \left(\tikzpic{\ca}\right)^{\op}
    \;=\;
    \tikzpic{\cu}
    \qquad
    \left(\tikzpic{\cro}\right)^{\op}
    \;=\;
    \tikzpic{\cro}
    \qquad
    \left(\tikzpic{\bli}\right)^{\op}
    \;=\;
    \tikzpic{\bli}
  \end{gather}
 One can easily verify that it extends to a vertical symmetry on $\affS(\affBr)$ in the sense of \cref{defn:vertical_symmetries}.

\begin{corollary}\label{alt:cor:reduction_vertsym}
  If the critical branchings of \cref{fig:critical_branchings_vertsym} are $\rhd$-tamely $\affT$-preconfluent, then $\affBr$ is a category of affine Brauer type with $2\bN_{>0}$-bubbles.
\end{corollary}
\begin{proof}
The Main Theorem (\cref{thm:affine_Grobner_system_iff_critical_branchings}) asserts that  $\affBr$ will be a category of affine Brauer type with $2\bN_{>0}$-bubbles if all critical branchings given in \cref{fig:critical_branchings}  are $\rhd$-tamely $\affT$-preconfluent; the other conditions are trivial to verify. 
Since $\affBr$ has a vertical symmetry, it is sufficient by \cref{prop:vertical_symmetries} to only check the critical branchings given in \cref{fig:critical_branchings_vertsym}.
\end{proof}

In the following sections, we will show that the four types of critical branching in \cref{fig:critical_branchings_vertsym} are indeed $\rhd$-tamely $\affT$-preconfluent.
By \cref{alt:cor:reduction_vertsym}, this gives a proof of \cref{alt:prop:affBr_main}.

\subsubsection{Non-affine critical branchings}
\label{alt:subsec:affBr_nonaffine_critical_branchings}

The following is easy to check:

\begin{lemma}\label{alt:lem:nonaffine_branchings}
  For $\affS(\affBr)$, the non-affine critical branchings from \cref{subfig:critical_branchings_non_affine_vertsym} are $\rhd$-tamely $\affT$-congruent.
\end{lemma}

\begin{proof}
  In fact, one can argue purely formally, as follows.
  Since the non-affine rewriting steps are of the form $D_s\to D_t$ for diagrams $D_s$ and $D_t$ with $D_s$ and $D_t$ inducing the same pairing on their endpoints, each branch rewrites into the unique normal Brauer diagram associated to the given pairing.
  Hence, every critical branching is even positively $\affT$-confluent.
\end{proof}

\subsubsection{Dot-slide critical branchings}
\label{alt:subsec:affBr_dotslide_without_stop}

\begin{lemma}\label{alt:lem:dotslide_branchings}
  For $\affS(\affBr)$, the dot-slide critical branchings from \cref{subfig:critical_branchings_dot_slides_vertsym} are $\rhd$-tamely $\affT^+$-preconfluent.
\end{lemma}
\begin{proof}
        For $   \tikzpic{
        \ulcr \dlcr
      }$, the branch where we first push the top left dot down and then the down left dot up is given by
      \begin{align*}
         \tikzpic{
        \ulcr \dlcr
      } \;\toT \;  
        \tikzpic{
        \urcr \drcr
      }+
      \tikzpic{\lca \cu[0][1]}-\tikzpic{\bli \li[1][0]} +
      \tikzpic{\rca \cu[0][1]}-\tikzpic{\li \bli[1][0]}
      \toT \;
      \tikzpic{
        \urcr \drcr
      }
      -\tikzpic{\bli \li[1][0]}-\tikzpic{\li \bli[1][0]},
      \end{align*}
     while first  pushing the down left dot up and then the up left dot down leads to
 \begin{align*}
         \tikzpic{
        \ulcr \dlcr
      } \;\toT \;  
        \tikzpic{
        \urcr \drcr
      }+
      \tikzpic{\ca \lcu[0][1]}-\tikzpic{\bli \li[1][0]} +
      \tikzpic{\ca \rcu[0][1]}-\tikzpic{\li \bli[1][0]}
      \toT \;
      \tikzpic{
        \urcr \drcr
      }
      -\tikzpic{\bli \li[1][0]}-\tikzpic{\li \bli[1][0]},
      \end{align*}
      giving us a  $\rhd$-tamely $\affT^+$-preconfluence.
      
For the other critical branchings with one dot, one branch entails sliding the dot as far as we can, followed by a non-affine rewriting step on the term containing a dot, while for the other branch we start with a non-affine rewriting step followed by a dot-slide. The term still containing a dot can be clearly seen to be the same after rewriting for both branches. So we only need to check that the lower order terms we obtained during a dot-slide are tamed congruent. 
We will consider the case 
$\tikzpic{
  \dlcr \li[2][0] \li[0][1] \cro[1][1] \cro[0][2] \li[2][2]
}[scale=.5]$. The other cases are similar.
After dot-sliding, the lower order terms of both branches are given by
\begin{align*}
  \renewcommand{\sca}{.4}
-\tikzpic{
  \LOTdcrobis \li[2][0] \li[0][1] \cro[1][1] \cro[0][2] \li[2][2]
}
-
\tikzpic{
  \cro \li[2][0] \li[0][1] \LOTdcrobis[1][1] \cro[0][2] \li[2][2]
}
\text{ and }
- \tikzpic{
  \li[0][0] \cro[1][0] \li[2][1] \LOTdcrobis[0][1] \cro[1][2] \li[0][2]
}
-\tikzpic{
  \li[0][0] \cro[1][0] \li[2][1] \cro[0][1] \LOTdcrobis[1][2] \li[0][2]
}.
\end{align*}
Filling in the definition of $\tikzpic{
\LOTdcrobis}=\tikzpic{
\lili}-\tikzpic{\cc}  $ leads to an expression which is clearly tamed congruent.
\end{proof}

\subsubsection{Exceptional critical branchings}
\label{alt:subsec:affBr_dotslide_with_stop}

\begin{lemma}
  \label{alt:lem:affBr_dotslide_with_stop}
  For $\affS(\affBr)$, the exceptional critical branchings from \cref{subfig:critical_branchings_exceptional_vertsym} are $\rhd$-tamely $\affT^+$-preconfluent.
\end{lemma}

Since the formulas for the bubble evaluation and the bubble slide were found precisely by computing the first two exceptional branchings (see the proof of \cref{alt:lem:bubble_identities}), the proof follows almost entirely by definition of $\evgenfun$ and $\bbslgenfun$.

\begin{proof}
\renewcommand{\sca}{.4}
The third exceptional critical branching is straightforward to verify.
Following the computation given in the proof of \cref{alt:lem:bubble_identities}, we have:
\begin{IEEEeqnarray*}{rCl}
  \Delta(-x)
  \;\overset{*}{\longleftarrow}_{\affT}\;
  \tikzpic{\ca[0][1]\cro\GFbul[1][1]\cu}
  \;&\overset{*}{\longrightarrow}_{\affT}&\;
  \tikzpic{\ca[0][1]\cro\GFbul\cu}
  +
  x\;
  \tikzpic{\ca[0][1]\lili\GFbul[1][1]\GFbul\cu}
  -
  x\;
  \tikzpic{\ca[0][1]\cc\GFbul[1][1]\GFbul\cu}
  \\
  \;&\overset{*}{\longrightarrow}_{\affT}&\;
  \Delta(x)
  +
  \frac{x}{2}(\Delta(x)+\Delta(-x))
  -
  x\Delta(x)\Delta(-x).
\end{IEEEeqnarray*}
Extracting the coefficient of $x^{2n}$, the two bubbles give $\Delta_{2n}$ and the lower-order terms are zero, since both terms are odd functions in $x$.
Extracting instead the coefficient of $x^{2n+1}$, one finds a tamed congruence by applying the bubble evaluation on both sides.

Following again the computation given in the proof of \cref{alt:lem:bubble_identities}, we have:
\begin{IEEEeqnarray*}{rCl}
  \Delta(x)\;\tikzpic{\li}
  \;\overset{*}{\longleftarrow}_{\affT}\;
  \tikzpic{
    \ca[0][2]\uli[2][2]
    \li[0][1]\cro[1][1]
    \cro \li[2][0]
    \dli \cu[1][0]
    \GFbul[0][1]
  }[scale=.8]
  \;&\overset{*}{\longrightarrow}_{\affT}&\;
  \tikzpic{\li}\;\Delta(x)
  -x\big(1-\frac{x}{2}\big)
  \left(
  \tikzpic{
    \GFbli[0][0]
    \GFbli[0][1]
  }
  -
  \tikzpic{
    \GFblirev[0][0]
    \GFblirev[0][1]
  }
  \right)
  -x^2
  \left(
  \tikzpic{
    \GFbli[0][0]
    \GFbli[0][1]
  }\;\Delta(x)
  -
  \Delta(x)\;
  \tikzpic{
    \GFblirev[0][0]
    \GFblirev[0][1]
  }
  \right)\;.
\end{IEEEeqnarray*}
Sliding the leading bubble left gives a tamed congruence, by definition of $\bbslgenfun$.
\end{proof}

\subsubsection{Bubble-slide critical branchings}
\label{alt:alt-subsec:affBr_bubble_slides}

It remains to treat the critical branchings involving the bubble-slide rewriting steps (\cref{subfig:critical_branchings_bubbles_vertsym}). By \cref{lem:bubble-slide-as-tame} (here $\VARbbsl=1$), it suffices to construct the tamed congruences listed there.

\medbreak

For any rewriting system of affine Brauer type, given a power series
\[H(x;\bft) = \sum_{k=0}^{\infty}\sum_{l=0}^{\infty} H(k,l)x^k\bft^l,\]
we denote:
\begin{align*}
\tikzpic{\bigbli[0][0][H(x;\bft)]}\;\coloneqq
= \sum_{k=0}^\infty\sum_{l=0}^\infty H(k,l)\,x^k\,\tikzpic{\bli[0][0][][l]}.
\end{align*}
The addition of dots becomes multiplicative when using generating functions:
\begin{align*}
\tikzpic{\bigbli[0][0][H(x;\bft)]}\circ \tikzpic{\bigbli[0][0][H'(x;\bft)]}
=
\tikzpic{\bigbli[0][0][H(x;\bft)H'(x;\bft)]}
\qquad\text{for all } H,H'\in\ring[[x,\bft]].
\end{align*}
Recall from \cref{not:tame_arrow} that the tamed congruences required for the bubble-slide branchings are $\affS$-congruences $\rhd$-tamed by a bubble with $k$ dots: every diagram appearing along the congruence must be strictly $\rhd$-smaller than a bubble carrying $k$ dots. In the case of the affine Brauer category, it will be sufficient to check that such diagrams have at most $k-1$ dots.

We extend taming at the level of generating functions. Let $g(x)$ be a generating function in $x$ whose coefficient of $x^k$, written $[x^k]\,g(x)$, is a diagram. An $\affS$-congruence between generating functions in $x$ is \emph{$\rhd$-tamed by $g(x)$} if, for every $k$, the congruence between the coefficients of $x^k$ is $\rhd$-tamed by the diagram $[x^k]\,g(x)$. Extending \cref{not:tame_arrow}, we write
\begin{align*}
f\ \overset{g(x)}{\tamearrow}\ f'
\end{align*}
for an $\affS$-congruence between the generating functions $f$ and $f'$ that is $\rhd$-tamed by $g(x)$. In all the statements below, the taming generating function is $\Delta(x)$.

\subsubsection*{The bubble-slide as a formal power series}

We solve the recursive definition of the bubble-slide and identify its relevant properties.
Recall that for $\affS(\affBr)$, we have $\VARbbsl=1$; we express everything with $\VARbbsl$ generic below, for latter use in other examples.

\begin{lemma}\label{alt:lem:GF_dot_pushing} 
  In $\affS(\affBr)$, we have the following identity:
\begin{align*}
\bbslgenfun(x)\ =\ x\,\;\bbslgenfun(x)\circ \tikzpic{\bigbli[0][0][\GFbil(x;\bft)]}\;+\VARbbsl\;\big(\GFnaba(x)+x\,\Delta(x)\big)\,
\big(\tikzpic{\bigbli[0][0][\GFbil(x;\bft)]}-\tikzpic{\bigbli[0][0][\GFbil(x;-\bft)]}\big),
\end{align*}
where $\GFnaba(x)$ is formal power series in $x$ and where $\GFbil(x;\bft)$ is a formal power series in $x$ and $\bft$, such that $\GFbil(-x;\bft)=-\GFbil(x;-\bft)$.
\end{lemma}

\begin{proof}
This follows from \cref{alt:lem:bubble_identities},
with
$
\GFbil(x;\bft) \coloneqq 
\frac{x}{(1-x\bft)^2}
$
and
$
\GFnaba(x) \coloneqq 1-\frac{x}{2}$.
\end{proof}

For any rewriting system of affine Brauer type such that the above data is defined, we set the following shorthands:
\begin{align}\label{alt:eq:definition_nabla_T}
  \GFnabla(x)\coloneqq \VARbbsl (\GFnaba(x)+x\,\Delta(x) )
  \quad\an\quad
  \GFbbsltemp(x;\bft)\coloneqq\frac{\GFbil(x;\bft)-\GFbil(x;-\bft)}{1-x\,\GFbil(x;\bft)} 
.
\end{align}
Note that $\GFbbsltemp(x;\bft)$ is well-defined, since the formal power series $1-x\,\GFbil(x;\bft)$ has constant term $1$ and hence is invertible.
We gather some basic properties of $\GFbbsltemp(x;\bft)$:

\begin{lemma}\label{alt:lem:properties_of_F}
Let $\GFbil(x;\bft)\in\ring[[x,\bft]]$ be a formal power series and let $\GFbbsltemp(x;\bft)$ be defined as in \eqref{alt:eq:definition_nabla_T}.
If $\GFbil(-x;\bft)=-\GFbil(x;-\bft)$, then:
\begin{enumerate}[(i)]
\item\label{alt:item:involution} $(1+x\GFbbsltemp(x;\bft))(1+ x\GFbbsltemp(x;-\bft))=1$;
\item\label{alt:item:support} $\GFbbsltemp(-x;\bft)=-\GFbbsltemp(x;-\bft)$.
\end{enumerate}
\end{lemma}
\begin{proof}
Straightforward.
\end{proof}

\begin{lemma}\label{alt:lem:single_slide}
  Let $(\affS,\affT)$ be a rewriting system of affine Brauer type such that \cref{alt:lem:GF_dot_pushing} holds.
  We have the following identity:
\begin{align*}
\bbslgenfun(x)\ =\ \GFnabla(x)\; \tikzpic{\bigbli[0][0][\GFbbsltemp(x;\bft)]},
\end{align*}
where $\GFbbsltemp(x;\bft)$ is defined as in \eqref{alt:eq:definition_nabla_T}.
\end{lemma}
\begin{proof}
The identity of \cref{alt:lem:GF_dot_pushing} rewrites as
\begin{align*}
\bbslgenfun(x)\circ \tikzpic{\bigbli[0][0][1-x\,\GFbil(x;\bft)]}
\ =\ \GFnabla(x)\,
\big(\tikzpic{\bigbli[0][0][\GFbil(x;\bft)]}-\tikzpic{\bigbli[0][0][\GFbil(x;-\bft)]}\big).
\end{align*}
The statement follows.
\end{proof}

\subsubsection*{Sliding two bubbles past a strand}

\begin{lemma}\label{alt:lem:sliding_two_bubbles}
Let $(\affS,\affT)$ be a rewriting system of affine Brauer type such that \cref{alt:lem:GF_dot_pushing} holds and such that $\VARbbsl=1$.
There exists a $\rhd$-tamed congruence
\begin{gather*}
  \tikzpic{\li}\;\Delta(x)\Delta(-x)-\Delta(x)\Delta(-x)\;\tikzpic{\li}
  \;\overset{\Delta(x)}{\tamearrow}\;
  \frac{\GFnaba(-x)}{x}\,\bbslgenfun(x)\ -\ \frac{\GFnaba(x)}{x}\,\bbslgenfun(-x).
\end{gather*}
\end{lemma}
\begin{proof}
Sliding the product $\Delta_i\,\Delta_j$ costs
\begin{align*}
\Delta_i\,\tikzpic{\LOTbbsl[0][0][j]}
+\Delta_j\,\tikzpic{\LOTbbsl[0][0][i]}
+\tikzpic{\LOTbbsl[0][0][j] \LOTbbsl[0][1][i]}.
\end{align*}
Summing these against $(-1)^i x^{i+j}$, the factor $\Delta_i$ coming from $\Delta(-x)$, gives
\begin{IEEEeqnarray*}{rCl}
\Delta(-x)\,\bbslgenfun(x)\ +\ \Delta(x)\,\bbslgenfun(-x)\ +\ \bbslgenfun(-x)\circ\bbslgenfun(x).
\end{IEEEeqnarray*}
Using \cref{alt:lem:properties_of_F,alt:lem:single_slide}, the last term simplifies as
\[\frac{1}{x}\GFnabla(-x)\,\bbslgenfun(x)\ -\ \frac{1}{x}\GFnabla(x)\,\bbslgenfun(-x).\]
Recalling that $\GFnabla(x)\coloneqq\GFnaba(x)+x\Delta(x)$, we find the statement.
\end{proof}

\subsubsection*{Bubble slide and bubble evaluation}

Recall that in $\affS(\affBr)$, we have $\evgenfun(x)
  =
  \VARev
  -
  \frac{x}{4}(\Delta(x)+\Delta(-x))
  +\frac{x}{2} \Delta(x)\Delta(-x)$ and $\GFnaba(x)=1-\frac{x}{2}$.
Thanks to the following lemma, the branching between bubble evaluation and bubble slide is $\rhd$-tamely $\affT$-preconfluent in $\affS(\affBr)$.

\begin{lemma}\label{alt:lem:general_bubble_evaluation}
Let $(\affS,\affT)$ be a rewriting system of affine Brauer type with bubble set $2\bN_{>0}$ such that \cref{alt:lem:GF_dot_pushing} holds and such that $\VARbbsl=1$.
Assume that $p(x)$ is of the form $p(x)=p_0+p_1 x$ and that $\evgenfun(x)$ has the shape
\begin{align}\label{alt:eq:shape_of_evgenfun}
\evgenfun(x)\ =a(x) + b(x)\,\big(\Delta(x)+\Delta(-x)\big)\ +\ c(x)\,\Delta(x)\,\Delta(-x),
\end{align}
such that $p_0c(x)=\frac{1}{2}x$ and $b(x)=p_1 c(x)$.
Then the branching
\begin{gather*}
  \tikzpic{
    \lorli[-1.4][0] \bul[0][.5]\ca[0][.5]\custop[0][.5]
    \node[left=-1pt] at (0,.5) {\scriptsize $k$};
  }
\end{gather*}
is $\rhd$-tamely $\affT$-preconfluent.
\end{lemma}
\begin{proof}
Following \cref{alt:lem:sliding_two_bubbles}, we have a tamed congruence
\begin{IEEEeqnarray*}{rCl}
\tikzpic{\li}\;\evgenfun(x)\ -\ \evgenfun(x)\;\tikzpic{\li}
&\overset{\Delta(x)}{\tamearrow}&
b(x)\big(\bbslgenfun(x)+\bbslgenfun(-x)\big)+c(x)\,\left(
\frac{\GFnaba(-x)}{x}\,\bbslgenfun(x)\ -\ \frac{\GFnaba(x)}{x}\,\bbslgenfun(-x)\right)
\\
&=&
p_0\frac{c(x)}{x}\Big(\bbslgenfun(x)- \bbslgenfun(-x)\Big)
+
\left(b(x)-p_1 c(x)\right)\Big(\bbslgenfun(x)+ \bbslgenfun(-x)\Big).
\end{IEEEeqnarray*}
By assumption, the leading factor of the first term is $\frac{1}{2}$, while the leading factor of the second term vanishes.
Restricting to odd coefficients, we find $\bbslgenfun(x)$ and conclude with \cref{lem:bubble-slide-as-tame}.
\end{proof}

\subsubsection*{Sliding a bubble past two strands}

\begin{lemma}\label{alt:lem:double_slide}
  Let $(\affS,\affT)$ be a rewriting system of affine Brauer type such that \cref{alt:lem:GF_dot_pushing} holds and such that $\VARbbsl=1$ or $p(x)=0$.
  There exists a tamed congruence
\begin{align*}
  \tikzpic{\GFbbsl\li[.5][0]}\,(x)
  +
   \VARbbsl\; \tikzpic{\GFbbsl\li[-.5][0]}\,(x)
\ \overset{\Delta(x)}{\tamearrow}\
\GFnabla(x)
\Big(\tikzpic{\bigbli[0][0][\GFbbsltemp(x;\bft)][left]\li[.5][0]}
+
\tikzpic{\li[-.5][0]\bigbli[0][0][\GFbbsltemp(x;\bft)]}
+x
\tikzpic{\bigbli[0][0][\GFbbsltemp(x;\bft)][left]\bigbli[.5][0][\GFbbsltemp(x;\bft)]}
\Big).
\end{align*}
\end{lemma}
\begin{proof}
Apply \cref{alt:lem:single_slide} to both terms on the left-hand side, slide the bubble that lands in between the two strands, and apply \cref{alt:lem:single_slide} again:
\begin{IEEEeqnarray*}{rCl}
  \tikzpic{\GFbbsl\li[.5][0]}\,(x)
  +
 \VARbbsl\; \tikzpic{\GFbbsl\li[-.5][0]}\,(x)
  &=&
\GFnabla(x)\tikzpic{\bigbli[0][0][\GFbbsltemp(x;\bft)][left]\li[.5][0]}
+
\GFnaba(x)\;\tikzpic{\li[-.5][0]\bigbli[0][0][\GFbbsltemp(x;\bft)]}
+x\,\;
\tikzpic{\li}\;\Delta(x)\; \tikzpic{\bigbli[0][0][\GFbbsltemp(x;\bft)]}
\\
&\overset{\Delta(x)}{\tamearrow}&
\GFnabla(x)
\Big(\tikzpic{\bigbli[0][0][\GFbbsltemp(x;\bft)][left]\li[.5][0]}
+
\tikzpic{\li[-.5][0]\bigbli[0][0][\GFbbsltemp(x;\bft)]}
\Big)
+x\,\;
\tikzpic{\GFbbsl[-.7][0]\bigbli[0][0][\GFbbsltemp(x;\bft)]}
\\
&=&
\GFnabla(x)
\Big(\tikzpic{\bigbli[0][0][\GFbbsltemp(x;\bft)][left]\li[.5][0]}
+
\tikzpic{\li[-.5][0]\bigbli[0][0][\GFbbsltemp(x;\bft)]}
+x
\tikzpic{\bigbli[-.5][0][\GFbbsltemp(x;\bft)][left]\bigbli[0][0][\GFbbsltemp(x;\bft)]}
\Big).
\end{IEEEeqnarray*}
This concludes.
\end{proof}

\subsubsection*{Bubble slide through a cap}

\begin{lemma}\label{alt:lem:cap_slide}
  Let $(\affS,\affT)$ be a rewriting system of affine Brauer type such that \cref{alt:lem:GF_dot_pushing} holds, such that $\VARbbsl = 1$ or $p(x)=0$,
  and such that $\tikzpic{\lca}[scale=.7]\to-\tikzpic{\rca}[scale=.7]$.
  Then the bubble-slide critical branching
  \begin{gather*}
    \tikzpic{
      \lorli[-1.4][0] \bul[0][.5]\ca[0][.5]\custop[0][.5]
      \node[left=-1pt] at (0,.5) {\scriptsize $k$};
      \rorli[1.8][0]
      \draw (-1.4,1) to[out=90,in=90] (1.8,1);
    }
  \end{gather*}
  is $\rhd$-tamely $\affT$-preconfluent.
\end{lemma}
\begin{proof}
  Following \cref{lem:bubble-slide-as-tame}, the statement is equivalent to the existence of a tamed congruence as follows, which we can also express in terms of formal power series:
  \begin{gather*}
  \tikzpic{
  \ca[0][1]
  \LOTbbsl[0][0][k] \li[1][0] }
  + \VARbbsl \;
  \tikzpic{
  \ca[0][1]
  \li \LOTbbsl[1][0][k] }
  \ \overset{k}{\tamearrow}\ 0
  \quad\Leftrightarrow\quad
  \tikzpic{
  \ca[0][1]
  \GFbbsl \li[1][0] }\,(x)
  +\VARbbsl\;
  \tikzpic{
  \ca[0][1]
  \li \GFbbsl[1][0] }\,(x)
  \ \overset{\Delta(x)}{\tamearrow}\ 0.
  \end{gather*}
Cap the tamed congruence of \cref{alt:lem:double_slide}. Sliding all dots across the cap onto the left strand introduces one sign per dot, i.e.\ it substitutes $\bft$ by $-\bft$.
The relation $x\,\GFbbsltemp(x;\bft)\GFbbsltemp(x;-\bft)+\GFbbsltemp(x;\bft)+\GFbbsltemp(x;-\bft)=0$ from \cref{alt:lem:properties_of_F}\ref{alt:item:involution} concludes.
\end{proof}

\subsubsection*{Bubble slide through a crossing}

We say that a linear combination of diagrams
\NewDocumentCommand{\tempdiag}{O{0}O{0}}{
  \draw (#1,#2) to (#1,#2+1);
  \draw (#1+1,#2) to (#1+1,#2+1);
  \draw[fill=white] (#1-\tikzLOTh,#2+.5-\tikzLOTv) rectangle (#1+1+\tikzLOTh,#2+.5+\tikzLOTv);
  \node at (#1+.5,#2+.5) {\tiny $D$};
}
\tikzpic{
  \tempdiag
} \defnemph{commutes with the crossing} if
\[\tikzpic{\cro\tempdiag[0][1]}\ \tamearrow\ \tikzpic{\cro\tempdiag[0][-1]}\;,\]
where the taming is understood from the context. 
This property is preserved under linear combinations and under multiplication by bubbles.

\begin{lemma}\label{alt:lem:diagrams_property_P}
\renewcommand{\sca}{.3}
In $\affS(\affBr)$,
for all $i,j\geq0$, the element
\renewcommand{\sca}{.5}
\begin{align*}
  \tikzpic{\bli[0][0][][i][left]\bli[0.5][0][][j][right]}
  +
  \tikzpic{\bli[0][0][][j][left]\bli[0.5][0][][i][right]}
  -(-1)^{\evenindicator{i}\evenindicator{j}+1}
  \big(
  \tikzpic{\bli[0][0][][i+j][left]\li[0.5][0]}
  +
  \tikzpic{\li[0][0]\bli[0.5][0][][i+j][right]}
  \big)
\end{align*}
commutes with the crossing.
\end{lemma}
\begin{proof}
  \renewcommand{\sca}{.4}
  We only use the form of dot slide through crossing in $\affS(\affBr)$:
  \[\tikzpic{\ulcr}\to\tikzpic{\drcr}-\tikzpic{\lili}+\tikzpic{\cc}
  \qquad\an\qquad
  \tikzpic{\urcr}\to\tikzpic{\dlcr}+\tikzpic{\lili}-\tikzpic{\cc}.\]
  For the purpose of the proof, we use two formal variables $x$ and $y$, and denote
  \begin{gather*}
    \tikzpic{\GFbli[0][0][x]} \coloneqq\sum_{k\geq0}x^k\,\tikzpic{\bli[0][0][][k]}
    \quad\an\quad
    \tikzpic{\GFbli[0][0][y]} \coloneqq\sum_{k\geq0}y^k\,\tikzpic{\bli[0][0][][k]}\;.
  \end{gather*}
  Using the identity 
  $
  \tikzpic{
    \GFblirev[0][0][y]
    \GFbli[0][1][x]
  }
  =
  \frac{1}{x+y}\left(
  x\;
  \tikzpic{
    \GFbli[0][0][x]
  }
  +
  y\;
  \tikzpic{
    \GFblirev[0][0][y]
  }\right)
  $, a generalization of \eqref{alt:eq:identity_GF_dot}, we find:
  \begin{IEEEeqnarray*}{rCl}
    \tikzpic{\uli[0][1]\uli[1][1]\cro\GFbul[0][1][x][left]\GFbul[1][1][y]\dli[0][0]\dli[1][0]}
    &\overset{*}{\longrightarrow}_{\affT}&
    \tikzpic{\uli[0][1]\uli[1][1]\cro\GFbul[0][0][y][left]\GFbul[1][0][x]\dli[0][0]\dli[1][0]}
    -x\;
    \tikzpic{\uli[0][1]\uli[1][1]\GFbul[0][1][x][left]\lili\dli[0][0]\dli[1][0]\GFbul[1][0][x]\GFbul[1][1][y]}
    +x\;
    \tikzpic{\uli[0][1]\uli[1][1]\GFbul[0][1][x][left]\cc\dli[0][0]\dli[1][0]\GFbul[1][0][x]\GFbul[1][1][y]}
    +y\;
    \tikzpic{\uli[0][1]\uli[1][1]\GFbul[1][1][y]\lili\dli[0][0]\dli[1][0]\GFbul[0][0][y][left]\GFbul[1][0][x]}
    -y\;
    \tikzpic{\uli[0][1]\uli[1][1]\GFbul[1][1][y]\cc\dli[0][0]\dli[1][0]\GFbul[0][0][y][left]\GFbul[1][0][x]}\;.
    \\
    &=&
    \tikzpic{\uli[0][1]\uli[1][1]\cro\GFbul[0][0][y][left]\GFbul[1][0][x]\dli[0][0]\dli[1][0]}
    -x\;
    \tikzpic{\uli[0][1]\uli[1][1]\GFbul[0][1][x][left]\lili\dli[0][0]\dli[1][0]\GFbul[1][0][x]\GFbul[1][1][y]}
    +y\;
    \tikzpic{\uli[0][1]\uli[1][1]\GFbul[1][1][y]\lili\dli[0][0]\dli[1][0]\GFbul[0][0][y][left]\GFbul[1][0][x]}
    +
    \frac{1}{x+y}
    \left(x^2
    \tikzpic{\uli[0][1]\uli[1][1]\GFbul[0][1][x][left]\cc\dli[0][0]\dli[1][0]\GFbulrev[0][0][x][left]}
    -y^2
    \tikzpic{\uli[0][1]\uli[1][1]\GFbulrev[0][1][y][left]\cc\dli[0][0]\dli[1][0]\GFbul[0][0][y][left]}\;
    \right).
  \end{IEEEeqnarray*}
  Setting
  $B(x)\coloneqq x^2\left(
  \tikzpic{\uli[0][1]\uli[1][1]\GFbul[0][1][x][left]\cc\dli[0][0]\dli[1][0]\GFbulrev[0][0][x][left]}[scale=.7]
  -
  \tikzpic{\uli[0][1]\uli[1][1]\GFbulrev[0][1][x][left]\cc\dli[0][0]\dli[1][0]\GFbul[0][0][x][left]}[scale=.7]\;
  \right)$,
  we then compute that:
  \begin{IEEEeqnarray*}{rCl}
    \tikzpic{\uli[0][1]\uli[1][1]\cro\GFbul[0][1][x][left]\GFbul[1][1][y]\dli[0][0]\dli[1][0]}
    +
    \tikzpic{\uli[0][1]\uli[1][1]\cro\GFbul[0][1][y][left]\GFbul[1][1][x]\dli[0][0]\dli[1][0]}
    &\overset{*}{\longrightarrow}_{\affT}&
    \tikzpic{\uli[0][1]\uli[1][1]\cro\GFbul[0][0][y][left]\GFbul[1][0][x]\dli[0][0]\dli[1][0]}
    +
    \tikzpic{\uli[0][1]\uli[1][1]\cro\GFbul[0][0][x][left]\GFbul[1][0][y]\dli[0][0]\dli[1][0]}
    +
    \frac{B(x)+B(y)}{x+y}.
  \end{IEEEeqnarray*}
  Note that $B(x)$ is odd.

  Writing $B(x)=\sum_{k\geq 0}b_kx^{2k+1}$, we have:
  \begin{gather*}
    \frac{B(x)+B(y)}{x+y}
    =
    \sum_{k\geq 0}\sum_{r=0}^{2k} b_k (-1)^rx^{r}y^{2k-r}.
  \end{gather*}
  Choose two pairs $(i,j)$ and $(a,b)$ such that $i+j=a+b$ and extract the coefficients of $x^iy^j$ and $x^ay^b$ in this identity.
  These coefficients are always equal except if $i+j=a+b$ is even and $i$ and $a$ have distinct parities, in which case they differ by a sign.
  In the special case $(a,b)=(i+j,0)$, the latter situation happens precisely when $i$ and $j$ are odd.
  The statement follows.
\end{proof}

\begin{lemma}\label{alt:lem:property_P_merging}
\renewcommand{\sca}{.3}
  Let $(\affS,\affT)$ be a rewriting system of affine Brauer type
such that \cref{alt:lem:diagrams_property_P} holds. 
For all power series $g(x;\bft)$ and $h(x;\bft)$, the element
\renewcommand{\sca}{.5}
\begin{align*}
\tikzpic{\bigbli[0][0][g(x;\bft)][left]\bigbli[.5][0][h(x;\bft)]}
+
\tikzpic{\bigbli[0][0][h(x;\bft)][left]\bigbli[.5][0][g(x;\bft)]}
-
\Big(
\tikzpic{\bigbli[0][0][\psi(x;\bft)][left]\li[.5][0]}
+
\tikzpic{\li\bigbli[.5][0][\psi(x;\bft)]}
\Big)
\end{align*}
commutes with the crossing, where
\begin{align*}
\psi(x;\bft)\coloneqq\tfrac12\big(g(x;\bft)+g(x;-\bft)\big)\big(h(x;\bft)+h(x;-\bft)\big)-g(x;\bft)\,h(x;\bft).
\end{align*}
\end{lemma}
\begin{proof}
Both sides are bilinear in $(g,h)$, so it suffices to check the congruence for monomials $g=\bft^i$ and $h=\bft^j$. In that case the left-hand side is $\tikzpic{\bli[0][0][][i][left]\bli[0.5][0][][j][right]}
+
\tikzpic{\bli[0][0][][j][left]\bli[0.5][0][][i][right]}$, and
\begin{align*}
\psi(x;\bft)=\Big(\tfrac12\big(1+(-1)^i\big)\big(1+(-1)^j\big)-1\Big)\,\bft^{i+j}
=(-1)^{\evenindicator{i}\evenindicator{j}+1}\,\bft^{i+j}.
\end{align*}
The difference of the two sides is exactly the element from \cref{alt:lem:diagrams_property_P}.
\end{proof}

\begin{lemma}\label{alt:lem:involutive_double_strand}
\renewcommand{\sca}{.3}
  Let $(\affS,\affT)$ be a rewriting system of affine Brauer type
such that \cref{alt:lem:diagrams_property_P} holds.
For every power series $g(x;\bft)$ such that $g(x;\bft)\,g(x;-\bft)=1$, the element
\renewcommand{\sca}{.5}
\begin{align*}
\tikzpic{\bigbli[0][0][g(x;\bft)][left]\bigbli[.5][0][g(x;\bft)]}
\end{align*}
commutes with the crossing.
\end{lemma}
\begin{proof}
Apply \cref{alt:lem:property_P_merging} with $h=g$: we wish to show that $\tikzpic{\bigbli[0][0][\psi(x;\bft)][left]\li[.5][0]}
+
\tikzpic{\li\bigbli[.5][0][\psi(x;\bft)]}
$ commutes with the crossing. Using $g(x;\bft)\,g(x;-\bft)=1$:
\begin{align*}
\psi(x;\bft)=\tfrac12\big(g(x;\bft)+g(x;-\bft)\big)^2-g(x;\bft)^2
=1+\varphi(x;\bft)
\end{align*}
where $\varphi(x;\bft)\coloneqq\tfrac12\big(g(x;-\bft)^2-g(x;\bft)^2\big)$.
It also follows from \cref{alt:lem:property_P_merging}, applied to $h=1$, that for any formal power series $g(x;\bft)$ the element
\begin{align*}
\tikzpic{\bigbli[0][0][g(x;\bft)][left]\li[.5][0]}
+
\tikzpic{\li[0][0]\bigbli[.5][0][g(x;\bft)]}
-
\Big(
\tikzpic{\bigbli[0][0][g(x;-\bft)][left]\li[.5][0]}
+
\tikzpic{\li\bigbli[.5][0][g(x;-\bft)]}
\Big)
\end{align*}
commutes with the crossing.
In particular, if $g(x;\bft)=-g(x;-\bft)$ then $\tikzpic{\bigbli[0][0][g(x;\bft)][left]\li[.5][0]}
+
\tikzpic{\li[0][0]\bigbli[.5][0][g(x;\bft)]}$ commutes with the crossing.
Since $\varphi(x;-\bft)=-\varphi(x;\bft)$, this concludes.
\end{proof}

\begin{lemma}
  \label{alt:lem:property_P}
  \renewcommand{\sca}{.3}
  Let $(\affS,\affT)$ be a rewriting system of affine Brauer type such that \cref{alt:lem:GF_dot_pushing} holds, such that $\VARbbsl=1$ or $p(x)=0$
  and such that
  \cref{alt:lem:diagrams_property_P} holds.
  The following bubble-slide critical branching
  \renewcommand{\sca}{.5}
  \begin{gather*}
    \tikzpic{
      \cro[-2][0]
      \bul[0][.5]\ca[0][.5]\custop[0][.5]
      \node[left=-1pt] at (0,.5) {\scriptsize $k$};
      \draw[very thick,red,->] (-2+.9,.6) to (-2+.5,1);
      \draw[very thick,red,->] (-2+.9,.4) to (-2+.5,0);
    }
  \end{gather*}
  is $\rhd$-tamely $\affT$-preconfluent.
\end{lemma}

\begin{proof}
Following \cref{lem:bubble-slide-as-tame}, the statement is equivalent to saying that
$
\tikzpic{\GFbbsl\li[.5][0]}\,(x)
+
\VARbbsl\;\tikzpic{\GFbbsl\li[-.5][0]}\,(x)
$
commutes with the crossing.
Thanks to \cref{alt:lem:double_slide}, it suffices to show that 
\begin{gather*}
\tikzpic{\bigbli[0][0][\GFbbsltemp(x;\bft)][left]\li[.5][0]}
+
\tikzpic{\li[-.5][0]\bigbli[0][0][\GFbbsltemp(x;\bft)]}
+x\,
\tikzpic{\bigbli[0][0][\GFbbsltemp(x;\bft)][left]\bigbli[.5][0][\GFbbsltemp(x;\bft)]}  \end{gather*}
commutes with the crossing.
Applying \cref{alt:lem:involutive_double_strand} to $x\,\GFbbsltemp(x;\bft)+1$ thanks to \cref{alt:lem:properties_of_F}\ref{alt:item:involution} shows that
\begin{gather*}
\tikzpic{\li[0][0]\li[.5][0]}
+x \big(\tikzpic{\bigbli[0][0][\GFbbsltemp(x;\bft)][left]\li[.5][0]}
+
\tikzpic{\li[-.5][0]\bigbli[0][0][\GFbbsltemp(x;\bft)]}\big)
+x\,^2
\tikzpic{\bigbli[0][0][\GFbbsltemp(x;\bft)][left]\bigbli[.5][0][\GFbbsltemp(x;\bft)]}  \end{gather*}
commutes with the crossing.
This concludes.
\end{proof}

%% file: sections/analysis_nilBrauer_new.tex
\subsection{The nil-Brauer category}
\label{subsec:nilBrauer}

\begin{definition}\label{nbalt:def:nilBr}
  \renewcommand{\sca}{.4}
  \newcommand{\vspc}{.5ex}
  \newcommand{\hspc}{60mu}
  Let $\mK$ be an integral domain in which $2$ is invertible and let $\VARev\in\{0,1\}$.
  The \emph{nil-Brauer category $\nilBr[\VARev]$} \cite[Definition~2.1]{BWW_NilBrauerCategory_2024}\footnote{Compared to \cite[Definition~2.1]{BWW_NilBrauerCategory_2024}, we renormalized the crossing by $-1$.} is the $\mK$-linear monoidal category generated by a single object $\tikzpic{\li}$ and morphisms $\tikzpic{\cu}$, $\tikzpic{\ca}$, $\tikzpic{\cro}$, $\tikzpic{\bli}$ subject to the following relations:
  \begin{gather*}
    \tikzpic{
      \cro[0][1]
      \cro 
    } 
    \;=\;
    0
    \mspace{\hspc}
    %
    \tikzpic{
      \cro[0][2] \li[2][2]
      \li[0][1] \cro[1][1] 
      \cro \li[2][0] 
    }[scale=.8]
    \;=\;
    \tikzpic{
      \cro[1][2] \li[0][2]
      \li[2][1] \cro[0][1] 
      \cro[1][0] \li[0][0] 
    }[scale=.8]
    \mspace{\hspc}
    \tikzpic{
      \cu\ca[0][0]
    } 
    \;=\;\VARev
    \\[\vspc]
    \tikzpic{
      \cu[1][0] \uli[2][0] 
      \dli \ca 
    }
    \;=\;
    \tikzpic{\li}
    \mspace{\hspc}
    %
    \tikzpic{
      \ca[0][1]
      \cro 
    }
    \;=\;
    0
    \mspace{\hspc}
    %
    \tikzpic{
      \uli[0][1] \ca[1][1] 
      \cro \li[2][0] 
    }
    \;=\;
    \tikzpic{ 
      \ca[0][1] \uli[2][1]
      \li \cro[1][0] 
    }
    \\[\vspc]
    \tikzpic{
      \ulcr
    }
    \;=\;
    \tikzpic{\drcr}
    \;-\;
    \tikzpic{\lili} 
    \;+\;
    \tikzpic{\cc} 
    \mspace{\hspc}
    \tikzpic{\rca}
    \;=\;-\;
    \tikzpic{\lca}
  \end{gather*}
\end{definition}

In this section, we use the Main Theorem (\cref{thm:affine_Grobner_system_iff_critical_branchings}) to prove the following:

\begin{proposition}\label{nbalt:prop:nilBr_main}
  The nil-Brauer category $\nilBr[\VARev]$ is a category of affine Brauer type with $(2\bN+1)$-bubbles (\cref{defn:category_affine_brauer_type}).
\end{proposition}

We follow the blueprint given by the analogous statement (\cref{alt:prop:affBr_main}) for the affine Brauer category.

\subsubsection{The nil-Brauer rewriting system}

\begin{lemma}\label{nbalt:lem:relations}
  \renewcommand{\sca}{.4}
  \newcommand{\vspc}{.5ex}
  \newcommand{\hspc}{60mu}
  In $\nilBr[\VARev]$, the following relations hold:
  \begin{gather*}
    \tikzpic{
      \dli[2][0] \ca[1][0] 
      \uli \cu 
    }
    \;=\;
    \tikzpic{\li}
    \mspace{\hspc}
    %
    \tikzpic{
      \ca \uli[2][0]
      \li[0][-1] \cro[1][-1]
      \cro[0][-2] \li[2][-2]
    }[scale=.8]
    =
    0
    \mspace{\hspc}
    %
    \tikzpic{
      \cro \li[2][0]
      \dli \cu[1][0]
    }
    =
    \tikzpic{
      \li[0][0] \cro[1][0]
      \cu[0][0] \dli[2][0]
    }
    \mspace{\hspc}
    %
    \tikzpic{
      \cro[0][1] \li[2][1]
      \li[0][0] \cro[1][0]
      \cu \dli[2][0]
    }[scale=.8]
    =
    0
    \\[\vspc]
    \tikzpic{
      \cro
      \cu
    }
    \;=\;
    0
    \mspace{\hspc}
    %
    \tikzpic{
      \ca[0][1] \uli[2][1]
      \li[0][0] \cro[1][0]
      \cu \dli[2][0]
    }[scale=.8]
    \;=\;
    0
    \mspace{\hspc}
    \tikzpic{\urcr}
    \;=\;
    \tikzpic{\dlcr}
    +\tikzpic{\lili}
    -\tikzpic{\cc}
    \mspace{\hspc}
    \tikzpic{\rcu \cu[0][0]}
    \;=\;-\;
    \tikzpic{\lcu \cu[0][0]}.
  \end{gather*}
\end{lemma}

To define a rewriting system of affine Brauer type as in \cref{fig:generic_grobner_of_affine_Brauer_type} associated to $\nilBr[\VARev]$, it remains to define the bubble evaluation and the bubble slide.
Note that $1-2\VARev\in\{1,-1\}$ is invertible in $\mK$.

We denote
\begin{gather*}
  \Delta_>(x)=\sum_{k>0}x^k\Delta_k.
\end{gather*}
In particular, $\Delta(x) = \Delta_>(x) + \VARev$.

\begin{lemma}\label{nbalt:lem:bubble_identities}
In $\nilBr[\VARev]$, we have the following identities:
\begin{gather*}
\Delta_\subeven(x) = \VARev + \frac{1}{1-2\VARev}\Delta_>(x)\Delta_>(-x)
\\
\tikzpic{\li}\;\Delta(x)
=
\Delta(x)\;\tikzpic{\li}
+
\frac{1}{2}
\left(
\tikzpic{
  \li[0][0]\li[0][1]
  \bigbul[0][1][\frac{1}{(1-x\bft)^2}-1]
}
-
\tikzpic{
  \li[0][0]\li[0][1]
  \bigbul[0][1][\frac{1}{(1+x\bft)^2}-1]
}
\right)
-
\left(
\tikzpic{
  \li[0][0]\li[0][1]
  \bigbul[0][1][\frac{1}{(1-x\bft)^2}-1]
}\;\Delta(x)
-
\Delta(x)\;
\tikzpic{
  \li[0][0]\li[0][1]
  \bigbul[0][1][\frac{1}{(1+x\bft)^2}-1]
}
\right)
\end{gather*}
\end{lemma}

\begin{proof}Proceeding similarly to the affine Brauer category, we find:
  \renewcommand{\sca}{.4}
  \begin{IEEEeqnarray*}{rCl}
  0
  &=&
  \tikzpic{\ca[0][1]\cro\GFbul[1][1]\cu}
  =
  \tikzpic{\ca[0][1]\cro\GFbul\cu}
  +
  x
  \tikzpic{\ca[0][1]\lili\GFbul[1][1]\GFbul\cu}
  -
  x
  \tikzpic{\ca[0][1]\cc\GFbul[1][1]\GFbul\cu}
  =
  \frac{x}{2}(\Delta(x)+\Delta(-x))
  - x\Delta(x)\Delta(-x)
  \\[1ex]
  &=&x\left(\frac{1}{2} - \VARev\right)(\Delta(x)+\Delta(-x))
  - x\Delta_>(x)\Delta_>(-x) + x\VARev^2
  \\[1ex]
  \Leftrightarrow
  0 &=& (1-2\VARev)\Delta_\subeven(x) - \Delta_>(x)\Delta_>(-x) + \VARev^2.
\end{IEEEeqnarray*}
We restricted to even coefficients in the latter equivalence. Using that $\VARev\in\{0,1\}$ and hence $-\frac{\VARev^2}{1-2\VARev}=\VARev$, we find the first part of the statement.

For the second part of the statement, the initial computation is almost identical to the affine Brauer case:
\begin{IEEEeqnarray*}{rCl}
  0
  &=&
  \tikzpic{
    \ca[0][2]\uli[2][2]
    \li[0][1]\cro[1][1]
    \cro \li[2][0]
    \dli \cu[1][0]
    \GFbul[0][1]
  }[scale=.8]
  =
  \tikzpic{
    \ca[0][2]\uli[2][2]
    \li[0][1]\cro[1][1]
    \cro \li[2][0]
    \dli \cu[1][0]
    \GFbul[1][0]
  }[scale=.8]
  -x\;
  \tikzpic{
    \ca[0][2]\uli[2][2]
    \li[0][1] \cro[1][1]
    \lili \li[2][0]
    \dli\cu[1][0]
    \GFbul[0][1]
    \GFbul[1][0]
  }[scale=.8]
  +x\;
  \tikzpic{
    \ca[0][2]\uli[2][2]
    \li[0][1] \cro[1][1]
    \cc \li[2][0]
    \dli\cu[1][0]
    \GFbul[0][1]
    \GFbul[1][0]
  }[scale=.8]
  =
  -x
  \tikzpic{
    \ca[0][2]\uli[2][2]
    \li[0][1] \cro[1][1]
    \lili \li[2][0]
    \dli\cu[1][0]
    \GFbul[0][0]
    \GFbul[1][1]
  }[scale=.8]
  +x
  \tikzpic{
    \ca[0][2]\uli[2][2]
    \cro[0][1] \li[2][1]
    \li \cu[1][1]
    \dli
    \GFbul[0][2]
    \GFbulrev[0][0]
  }[scale=.8]
  \\
  &=&
  -x
  \left(
  \tikzpic{
    \ca[0][2]\uli[2][2]
    \li[0][1] \cro[1][1]
    \lili \li[2][0]
    \dli\cu[1][0]
    \GFbul[0][0]
    \GFbul[2][2]
  }[scale=.8]
  -x
  \tikzpic{
    \ca[0][2]\uli[2][2]
    \li[0][1] \lili[1][1]
    \lili \li[2][0]
    \dli\cu[1][0]
    \GFbul[0][0]
    \GFbul[1][1]
    \GFbul[2][2]
  }[scale=.8]
  +x
  \tikzpic{
    \ca[0][2]\uli[2][2]
    \li[0][1] \cc[1][1]
    \lili \li[2][0]
    \dli\cu[1][0]
    \GFbul[0][0]
    \GFbul[1][1]
    \GFbul[2][2]
  }[scale=.8]
  \right)
  +x
  \left(
  \tikzpic{
    \ca[0][2]\uli[2][2]
    \cro[0][1] \li[2][1]
    \li \cu[1][1]
    \dli
    \GFbul[1][1]
    \GFbulrev[0][0]
  }[scale=.8]
  -x
  \tikzpic{
    \ca[0][2]\uli[2][2]
    \lili[0][1] \li[2][1]
    \li \cu[1][1]
    \dli
    \GFbul[0][2]
    \GFbul[1][1]
    \GFbulrev[0][0]
  }[scale=.8]
  +x
  \tikzpic{
    \ca[0][2]\uli[2][2]
    \cc[0][1] \li[2][1]
    \li \cu[1][1]
    \dli
    \GFbul[0][2]
    \GFbul[1][1]
    \GFbulrev[0][0]
  }[scale=.8]
  \right)
  \\
  &=&
  x^2
  \left(
  \tikzpic{
    \GFbli[0][0]
    \GFblirev[0][1]
    \GFbli[0][2]
  }
  -
  \tikzpic{
    \GFblirev[0][0]
    \GFbli[0][1]
    \GFblirev[0][2]
  }
  \right)
  -x^2
  \left(
  \tikzpic{
    \GFbli[0][0]
    \GFbli[0][1]
  }\;\Delta(x)
  -
  \Delta(x)\;
  \tikzpic{
    \GFblirev[0][0]
    \GFblirev[0][1]
  }
  \right)
  \\
  &=&
  \frac{x^2}{2}
  \left(
  \tikzpic{
    \GFbli[0][0]
    \GFbli[0][1]
  }
  -
  \tikzpic{
    \GFblirev[0][0]
    \GFblirev[0][1]
  }
  \right)
  -x^2
  \left(
  \tikzpic{
    \GFbli[0][0]
    \GFbli[0][1]
  }\;\Delta(x)
  -
  \Delta(x)\;
  \tikzpic{
    \GFblirev[0][0]
    \GFblirev[0][1]
  }
  \right)
  \\
  \Leftrightarrow
  \tikzpic{\li}\;\Delta(x)
  -
  \Delta(x)\;\tikzpic{\li}
  &=&
  \frac{1}{2}
  \left(
  \tikzpic{
    \li[0][0]\li[0][1]
    \bigbul[0][1][\frac{1}{(1-x\bft)^2}-1]
  }
  -
  \tikzpic{
    \li[0][0]\li[0][1]
    \bigbul[0][1][\frac{1}{(1+x\bft)^2}-1]
  }
  \right)
  -
  \left(
  \tikzpic{
    \li[0][0]\li[0][1]
    \bigbul[0][1][\frac{1}{(1-x\bft)^2}-1]
  }\;\Delta(x)
  -
  \Delta(x)\;
  \tikzpic{
    \li[0][0]\li[0][1]
    \bigbul[0][1][\frac{1}{(1+x\bft)^2}-1]
  }
  \right)
\end{IEEEeqnarray*}
\end{proof}

\begin{definition}\label{nbalt:def:S_nilBr}
Let $\affS(\nilBr)$ be the rewriting system of affine Brauer type (\cref{fig:generic_grobner_of_affine_Brauer_type}) with set of bubbles $\bubbleset=2\bN+1$ and the following data: all scalars equal $1$, except the dot slides through cups and caps, whose scalar is $-1$; we have
\begin{align*}
\tikzpic{\LOTdcro}\;=\;-\,\tikzpic{\lili}+\tikzpic{\cc}\;=\;-\,\tikzpic{\LOTdcrobis},
\end{align*}
the lower-order term of the bubble evaluation is given by
\begin{gather}\label{nbalt:eq:bbev_formula}
  \evgenfun(x)
  =
  \VARev + \frac{1}{1-2\VARev}\Delta_>(x)\Delta_>(-x)\;;
\end{gather}
the lower-order term of the bubble slide is recursively defined by
\begin{gather}\label{nbalt:eq:bbsl_formula}
\bbslgenfun(x)
=
(\frac{1}{2}-\Delta(x))
\left(
\tikzpic{
  \li[0][0]\li[0][1]
  \bigbul[0][1][\frac{1}{(1-x\bft)^2}-1]
}
-
\tikzpic{
  \li[0][0]\li[0][1]
  \bigbul[0][1][\frac{1}{(1+x\bft)^2}-1]
}
\right)
-
\tikzpic{
  \bigbli[0][0][\frac{1}{(1-x\bft)^2}-1]\GFbbsl[0][1]
}\;;
\end{gather}
and the remaining lower-order terms vanish.
\end{definition}

Being a scalar, respectively a polynomial in bubbles carrying at most $2n-1$ dots, $\tikzpic{\LOTev[0][0][2n]}$ has the shape prescribed in \cref{fig:generic_grobner_of_affine_Brauer_type} for the lower-order term of the bubble evaluation.

Unfolding the recursion, $\tikzpic{\LOTbbsl[0][0][k]}$ is a linear combination of straight lines carrying at most $k$ dots together with bubbles carrying at most $k-1$ dots, each summand carrying exactly $k$ dots in total: it has the shape prescribed in \cref{fig:generic_grobner_of_affine_Brauer_type} for the lower-order term of the bubble slide. In contrast to the affine Brauer case, the lower-order term of the nil-Brauer bubble slide preserves the total number of dots.

\begin{corollary}\label{nbalt:cor:presentation}
The linear monoidal rewriting system $\affS(\nilBr)$ presents $\nilBr[\VARev]$.
\end{corollary}
\begin{proof}
This follows from \cref{nbalt:def:nilBr} and \cref{nbalt:lem:bubble_identities}.
\end{proof}

We define the following vertical symmetry on the generators of $\nilBr[\VARev]$:
  \begin{gather}\label{nbalt:eq:vertsymm}
    \left(\tikzpic{\cu}\right)^{\op}
    \;=\;
    \tikzpic{\ca}
    \qquad
    \left(\tikzpic{\ca}\right)^{\op}
    \;=\;
    \tikzpic{\cu}
    \qquad
    \left(\tikzpic{\cro}\right)^{\op}
    \;=\;
    \tikzpic{\cro}
    \qquad
    \left(\tikzpic{\bli}\right)^{\op}
    \;=\;
    \tikzpic{\bli}
  \end{gather}
Just as for the affine Brauer case, one can easily verify that it extends to a vertical symmetry on $\affS(\nilBr)$ in the sense of \cref{defn:vertical_symmetries}.

\begin{corollary}\label{nbalt:cor:reduction_vertsym}
  If the critical branchings of \cref{fig:critical_branchings_vertsym} are $\rhd$-tamely $\affT$-preconfluent, then $\nilBr[\VARev]$ is a category of affine Brauer type with $(2\bN+1)$-bubbles.
\end{corollary}
\begin{proof}
Identical to the proof in the affine Brauer case (\cref{alt:cor:reduction_vertsym}).
\end{proof}

In the following sections, we will show that the four types of critical branching in \cref{fig:critical_branchings_vertsym} are indeed $\rhd$-tamely $\affT$-preconfluent, following the same strategy as the one used in the affine Brauer case.

\subsubsection{Non-affine critical branchings}
\label{nbalt:subsec:nilBr_nonaffine_critical_branchings}

\begin{lemma}\label{nbalt:lem:nonaffine_branchings}
  For $\affS(\nilBr)$, the non-affine critical branchings from \cref{subfig:critical_branchings_non_affine_vertsym} are $\rhd$-tamely $\affT$-congruent.
\end{lemma}

\begin{proof}
The non-affine rewriting steps of $\affS(\nilBr)$ are those of $\affS(\affBr)$, except that some targets are replaced by zero. In particular, each branch of a non-affine critical branching rewrites into a scalar multiple of the unique normal Brauer diagram associated to the pairing induced by the source, and both branches produce the same scalar; the argument of \cref{alt:lem:nonaffine_branchings} applies verbatim.
Hence, every critical branching is even positively $\affT$-confluent.
\end{proof}
\subsubsection{Dot-slide critical branchings}
\label{nbalt:subsec:nilBr_dotslide_without_stop}

\begin{lemma}\label{nbalt:lem:dotslide_branchings}
  For $\affS(\nilBr)$, the dot-slide critical branchings from \cref{subfig:critical_branchings_dot_slides_vertsym} are $\rhd$-tamely $\affT^+$-preconfluent.
\end{lemma}
\begin{proof}
The dot-slide rewriting steps of $\affS(\nilBr)$ coincide with those of $\affS(\affBr)$, and the non-affine rewriting steps coincide or have target zero. The proof of \cref{alt:lem:dotslide_branchings} applies verbatim: whenever a branch of the affine Brauer case rewrites a term through a non-affine rewriting step whose target is zero in $\affS(\nilBr)$, the corresponding term vanishes in both branches.
\end{proof}

\subsubsection{Exceptional critical branchings}
\label{nbalt:subsec:nilBr_dotslide_with_stop}

\begin{lemma}
  \label{nbalt:lem:nilBr_dotslide_with_stop}
  For $\affS(\nilBr)$, the exceptional critical branchings from \cref{subfig:critical_branchings_exceptional_vertsym} are $\rhd$-tamely $\affT^+$-preconfluent.
\end{lemma}
\begin{proof}
  Analogous to the affine Brauer case (\cref{alt:lem:affBr_dotslide_with_stop}), using the computation in the proof of \cref{nbalt:lem:bubble_identities}.
\end{proof}

\subsubsection{Bubble-slide critical branchings}
\label{nbalt:subsec:nilBr_bubble_slides}

It remains to treat the critical branchings involving the bubble-slide rewriting steps (\cref{subfig:critical_branchings_bubbles_vertsym}).
The proof proceeds similarly to the affine Brauer case.
Since the dot slides are identical, we only need to check that \cref{alt:lem:GF_dot_pushing} still holds in the nil-Brauer case, and give an analogue of \cref{alt:lem:general_bubble_evaluation} for categories of affine Brauer type with $2\bN+1$-bubbles.

\subsubsection*{The bubble-slide as a formal power series}
\begin{lemma}
  \label{nbalt:lem:GF_dot_pushing}
  \Cref{alt:lem:GF_dot_pushing} holds in $\affS(\nilBr)$, with
\begin{align*}
\GFbil(x;\bft)\coloneqq-\sum_{m\geq1}(m+1)\,x^{m-1}\bft^m
=-\frac{1}{x}\left(\frac{1}{(1-x\bft)^2}-1\right)
\qquad\text{and}\qquad
\GFnaba(x)\coloneqq-\frac{x}{2}.
\end{align*}
\end{lemma}
\begin{proof}
  Directly follows from the definition.
\end{proof}

\subsubsection*{Sliding two bubbles past a strand}

\begin{lemma}\label{nbalt:lem:sliding_two_bubbles}
Let $(\affS,\affT)$ be a rewriting system of affine Brauer type such that \cref{alt:lem:GF_dot_pushing} holds and $\VARbbsl=1$.
There exists a $\rhd$-tamed congruence
\begin{gather*}
  \tikzpic{\li}\;\Delta_>(x)\Delta_>(-x)-\Delta_>(x)\Delta_>(-x)\;\tikzpic{\li}
  \;\overset{\Delta(x)}{\tamearrow}\;
\left(\frac{\GFnaba(-x)}{x}-\VARev\right)\,\bbslgenfun(x)\ -\ \left(\frac{\GFnaba(x)}{x}+\VARev\right)\,\bbslgenfun(-x)
.
\end{gather*}
\end{lemma}
\begin{proof}
Sliding the product $\Delta_i\,\Delta_j$ costs
\begin{align*}
\Delta_i\,\tikzpic{\LOTbbsl[0][0][j]}
+\Delta_j\,\tikzpic{\LOTbbsl[0][0][i]}
+\tikzpic{\LOTbbsl[0][0][j] \LOTbbsl[0][1][i]}.
\end{align*}
Summing these against $(-1)^i x^{i+j}$ for $i>0$ and $j>0$, the factor $\Delta_i$ coming from $\Delta(-x)$, gives (using that $\tikzpic{\LOTbbsl[0][0][0]}=0$):
\begin{IEEEeqnarray*}{rCl}
\Delta_>(-x)\,\bbslgenfun(x)\ +\ \Delta_>(x)\,\bbslgenfun(-x)\ +\ \bbslgenfun(-x)\circ\bbslgenfun(x).
\end{IEEEeqnarray*}
Using \cref{alt:lem:properties_of_F,alt:lem:single_slide}, the last term simplifies as
\[\frac{1}{x}\GFnabla(-x)\,\bbslgenfun(x)\ -\ \frac{1}{x}\GFnabla(x)\,\bbslgenfun(-x).\]
Recalling that $\GFnabla(x)\coloneqq\GFnaba(x)+x\Delta(x)$ concludes.
\end{proof}

\subsubsection*{Bubble slide and bubble evaluation}

Recall that in $\nilBr$, we have $p(x)=-\frac{x}{2}$ and $a(x)=\VARev$, and $c(x)=\frac{1}{1-2\VARev}$

\begin{lemma}\label{nbalt:lem:general_bubble_evaluation}
Let $(\affS,\affT)$ be a rewriting system of affine Brauer type with bubble set $\bubbleset=2\bN+1$ and such that \cref{alt:lem:GF_dot_pushing} holds and $\VARbbsl=1$.
Assume that $p(x)$ has the shape $p(x)=p_1x$ and that $\evgenfun(x)$ has the shape
\begin{align}\label{nbalt:eq:shape_of_evgenfun}
\evgenfun(x)\ =\ a(x) + c(x)\,\Delta_>(x)\,\Delta_>(-x)
\end{align}
for power series $a(x),c(x)\in\ring[[x]]$.
If $c(x) = \frac{1}{-2p_1-2\VARev}$, then the branching
\begin{gather*}
  \tikzpic{
    \lorli[-1.4][0] \bul[0][.5]\ca[0][.5]\custop[0][.5]
    \node[left=-1pt] at (0,.5) {\scriptsize $k$};
  }
\end{gather*}
is $\rhd$-tamely $\affT$-preconfluent.
\end{lemma}
\begin{proof}
  Following \cref{nbalt:lem:sliding_two_bubbles}, we have a tamed congruence
  \begin{IEEEeqnarray*}{rCl}
  \tikzpic{\li}\;\evgenfun(x)\ -\ \evgenfun(x)\;\tikzpic{\li}
  &\overset{\Delta(x)}{\tamearrow}&
   c(x)\left(\frac{\GFnaba(-x)}{x}-\VARev\right)\,\bbslgenfun(x)\ -\ c(x)\left(\frac{\GFnaba(x)}{x}+\VARev\right)\,\bbslgenfun(-x)
  \\
  &=&
  -c(x)\left(p_1+\VARev\right)\left(\bbslgenfun(x)\ +\ \bbslgenfun(-x)\right).
  \end{IEEEeqnarray*}
  By assumption, the leading factor is $\frac{1}{2}$.
Restricting to even coefficients, we find $\bbslgenfun(x)$ and conclude with \cref{lem:bubble-slide-as-tame}.
\end{proof}

%% file: sections/analysis_affine_VW_supercategory_new.tex
\subsection{The affine VW supercategory}
\label{subsec:affine_VW_supercategory}

\begin{definition}\label{svalt:def:sVW}
  \renewcommand{\sca}{.4}
  \newcommand{\vspc}{.5ex}
  \newcommand{\hspc}{60mu}
  The \emph{affine VW supercategory $\sVW$} \cite[Section~1.4]{BDE+_AffineVWSupercategory_2020}\footnote{Compared to \cite[Section~1.4]{BDE+_AffineVWSupercategory_2020}, we renormalized the cup by $-1$.} is the supermonoidal category generated by a single object $\di{\sca}{\li}$ and morphisms
  \begin{IEEEeqnarray*}{CcCcCcC}
    \di{0.8*\sca}{\cro}
    &\qquad&
    \di{0.8*\sca}{\ca}
    &\qquad&
    \di{0.8*\sca}{\cu}
    &\qquad&
    \di{0.8*\sca}{\bli}
    \\
    \text{\scriptsize even}
    &&
    \text{\scriptsize odd}
    &&
    \text{\scriptsize odd}
    &&
    \text{\scriptsize even}
  \end{IEEEeqnarray*}
  subject to the super interchange law and the following relations:
  \begin{gather*}
    \di{\sca}{
      \cro[0][\h]
      \cro
    }
    =
    \di{\sca}{\lili\lili[0][\h]}
    \mspace{\hspc}
    %
    \di{\sca*.8}{
      \cro[0][2*\h] \li[2*\br][2*\h]
      \li[0][\h] \cro[\br][\h]
      \cro \li[2*\br][0]
    }
    =
    \di{\sca*.8}{
      \cro[\br][2*\h] \li[0][2*\h]
      \li[2*\br][\h] \cro[0][\h]
      \cro[\br][0] \li[0][0]
    }
    \mspace{\hspc}
    %
    \di{\sca}{
      \cu[\br][0] \uli[2*\br][0]
      \dli \ca
    }
    =
    \di{\sca}{\li}
    \mspace{\hspc}
    %
    \di{\sca}{
      \dli[2*\br][0] \ca[\br][0]
      \uli \cu
    }
    =\;-\;
    \di{\sca}{\li}
    \mspace{\hspc}
    %
    \di{\sca}{
      \uli[0][\h] \ca[\br][\h]
      \cro \li[2*\br][0]
    }
    =
    \di{\sca}{
      \ca[0][\h] \uli[2*\br][\h]
      \li \cro[\br][0]
    }
    \\[\vspc]
    \mspace{\hspc}
    \di{\sca}{
      \ca[0][\h]
      \cro
    }
    \;=\;
    \di{\sca}{\ca}
    \mspace{\hspc}
    \di{\sca}{
      \ulcr
    }
    \;=\;
    \di{\sca}{\drcr}
    -\di{\sca}{\lili}
    -\di{\sca}{\cc}
    \mspace{\hspc}
    \di{\sca}{\rca \ca[0][0]}
    \;=\;
    \di{\sca}{\lca \ca[0][0]}
    +
    \di{\sca}{\ca \ca[0][0]}.
  \end{gather*}
\end{definition}
The following relations are easy to derive, see  \cite[Section~2.1]{BDE+_AffineVWSupercategory_2020}.
\begin{lemma}\label{svalt:lem:relations}
  \renewcommand{\sca}{.4}
  \newcommand{\vspc}{.5ex}
  \newcommand{\hspc}{60mu}
  In $\sVW$, the following relations hold:
  \begin{gather*}
    \di{\sca}{
      \ca \uli[2*\br][0]
      \li[0][-\h] \cro[\br][-\h]
      \cro[0][-2*\h] \li[2*\br][-2*\h]
    }
    =
    \di{\sca}{\li\ca[\br][0]}
    \mspace{\hspc}
    %
    \di{\sca}{
      \cro \li[2*\br][0]
      \dli \cu[\br][0]
    }
    =
    \di{\sca}{
      \li[0][0] \cro[\br][0]
      \cu[0][0] \dli[2*\br][0]
    }
    \mspace{\hspc}
    %
    \di{\sca}{
      \cro[0][\h] \li[2*\br][\h]
      \li[0][0] \cro[\br][0]
      \cu \dli[2*\br][0]
    }
    =
    \di{\sca}{\li\cu[\br][\h]}
    \mspace{\hspc}
    %
    \di{\sca}{\lcu\ca\node[left=-1pt] at (0,0) {\footnotesize $k$};}
    \;=\;
    0
    \\[\vspc]
    \di{\sca}{
      \cro
      \cu
    }
    \;=\;-\;
    \di{\sca}{\cu}
    \mspace{\hspc}
    %
    \di{\sca}{
      \ca[0][\h] \uli[2*\br][\h]
      \li[0][0] \cro[\br][0]
      \cu \dli[2*\br][0]
    }
    \;=\;
    \di{\sca}{\li}
    \mspace{\hspc}
    \di{\sca}{\urcr}
    \;=\;
    \di{\sca}{\dlcr}
    +\di{\sca}{\lili}
    -\di{\sca}{\cc}
    \mspace{\hspc}
    \di{\sca}{\rcu}
    \;=\;
    \di{\sca}{\lcu}
    -
    \di{\sca}{\cu}.
  \end{gather*}
\end{lemma}

In this section, we use the Main Theorem (\cref{thm:affine_Grobner_system_iff_critical_branchings}) to prove the following:

\begin{proposition}\label{svalt:prop:sVW_main}
  The affine VW supercategory $\sVW$ is a category of affine Brauer type with $\emptyset$-bubbles (\cref{defn:category_affine_brauer_type}).
\end{proposition}

We follow the blueprint given by the analogous statement (\cref{alt:prop:affBr_main}) for the affine Brauer category.

\subsubsection{The affine VW rewriting system}
\label{svalt:subsec:sVW_rewriting_system}

Since $\bubbleset=\emptyset$ and $\VARev = 0$ in $\sVW$, we can define the bubble evaluation and the bubble slide to be trivial.

\begin{definition}\label{svalt:def:S_sVW}
Let $\affS(\sVW)$ be the rewriting system of affine Brauer type (\cref{fig:generic_grobner_of_affine_Brauer_type}) with set of bubbles $\bubbleset=\emptyset$, with the scalars and lower-order terms of the non-affine and dot-slide rewriting steps as prescribed by the relations of \cref{svalt:def:sVW} and of \cref{svalt:lem:relations}, and with trivial lower-order terms for the bubble evaluation and the bubble slide:
\begin{gather*}
  \evgenfun(x)
  =
  0
  \qquad\an\qquad
  \bbslgenfun(x)
  =
  0\;.
\end{gather*}
\end{definition}

\begin{corollary}
The linear monoidal rewriting system $\affS(\sVW)$ presents $\sVW$.
\end{corollary}
\begin{proof}
This follows from \cref{svalt:def:sVW} and \cref{svalt:lem:relations}.
\end{proof}

We define the following vertical symmetry $*$ on the generators of $\sVW$ \cite[Proposition~13]{BDE+_AffineVWSupercategory_2020}, playing the role of $\op$ in the affine Brauer and nil-Brauer cases:
  \begin{gather}\label{svalt:eq:vertsymm}
    \left(\di{\sca}{\cu}\right)^*
    \;=\;
    \di{\sca}{\ca}
    \qquad
    \left(\di{\sca}{\ca}\right)^*
    \;=\;
    -\;\di{\sca}{\cu}
    \qquad
    \left(\di{\sca}{\cro}\right)^*
    \;=\;
    -\;\di{\sca}{\cro}
    \qquad
    \left(\di{\sca}{\bli}\right)^*
    \;=\;
    -\;\di{\sca}{\bli}.
  \end{gather}
Just as for the affine Brauer case, one can easily verify that it extends to a vertical symmetry on $\affS(\sVW)$ in the sense of \cref{defn:vertical_symmetries}.

\begin{corollary}\label{svalt:cor:reduction_vertsym}
  If the critical branchings of \cref{fig:critical_branchings_vertsym} are $\rhd$-tamely $\affT$-preconfluent, then $\sVW$ is a category of affine Brauer type with $\emptyset$-bubbles.
\end{corollary}
\begin{proof}
Identical to the proof in the affine Brauer case (\cref{alt:cor:reduction_vertsym}).
\end{proof}

In the following sections, we will show that the four types of critical branching in \cref{fig:critical_branchings_vertsym} are indeed $\rhd$-tamely $\affT$-preconfluent, following the same strategy as the one used in the affine Brauer case.
By \cref{svalt:cor:reduction_vertsym}, this gives a proof of \cref{svalt:prop:sVW_main}.

\subsubsection{Non-affine critical branchings}
\label{svalt:subsec:sVW_nonaffine_critical_branchings}

\begin{lemma}\label{svalt:lem:nonaffine_branchings}
  For $\affS(\sVW)$, the non-affine critical branchings from \cref{subfig:critical_branchings_non_affine_vertsym} are $\rhd$-tamely $\affT$-congruent.
\end{lemma}
\begin{proof}
The non-affine rewriting steps of $\affS(\sVW)$ have the same shape as those of $\affS(\affBr)$, up to the scalars recorded in \cref{svalt:def:S_sVW} and \cref{svalt:lem:relations}.
While we cannot use the bypass argument from the proof of \cref{alt:lem:nonaffine_branchings}, checking that each branching is tamely congruent is straightforward.
\end{proof}

\subsubsection{Dot-slide critical branchings}
\label{svalt:subsec:sVW_dotslide_without_stop}

\begin{lemma}\label{svalt:lem:dotslide_branchings}
  For $\affS(\sVW)$, the dot-slide critical branchings from \cref{subfig:critical_branchings_dot_slides_vertsym} are $\rhd$-tamely $\affT^+$-preconfluent.
\end{lemma}
\begin{proof}
The argument is similar to the affine Brauer case (\cref{alt:lem:dotslide_branchings}): the branch where we first push a dot one way and then back, and the branch where we push it the other way first, agree once the non-affine rewriting steps and dot-slide scalars of \cref{svalt:def:S_sVW} are substituted in.
\end{proof}

\subsubsection{Exceptional critical branchings}
\label{svalt:subsec:sVW_dotslide_with_stop}

\begin{lemma}
  \label{svalt:lem:sVW_dotslide_with_stop}
  For $\affS(\sVW)$, the exceptional critical branchings from \cref{subfig:critical_branchings_exceptional_vertsym} are $\rhd$-tamely $\affT^+$-preconfluent.
\end{lemma}
\begin{proof}
  Since we defined the bubble evaluation and bubble slide to be zero without relying on the exceptional branchings to ``find'' their expression, the argument slightly differs from the affine Brauer and nil-Brauer cases.
  The first exceptional critical branching is trivial, since both sides rewrite to zero; and similarly for the third exceptional critical branching. 
  To compute with the second one, we use generating functions, setting notations similar to the affine Brauer and nil-Brauer cases:
  We write
\begin{gather*}
  \tikzpic{\GFbli} = \tikzpic{\bigbli[0][0][\frac{1}{1-x\bft}]}
  \quad\an\quad
  \tikzpic{\GFblirev[0][0]} = \tikzpic{\bigbli[0][0][\frac{1}{1-x(\bft+1)}]}.
\end{gather*}
Note that this convention differs from the one in the affine Brauer and nil-Brauer cases. In particular, in contrast to the previous two cases:
\[\tikzpic{\ca \GFbul[1][-.5]\lili[0][-1]}= \tikzpic{\ca \GFbulrev[0][-.5]\lili[0][-1]}
\quad\text{while}\quad
\tikzpic{\ca \GFbulrev[1][-.5]\lili[0][-1]}\neq \tikzpic{\ca \GFbul[0][-.5]\lili[0][-1]}\;.
\]
The analogue of the identity \eqref{alt:eq:identity_GF_dot} is the identity
\begin{gather*}
  \tikzpic{
    \GFbli[0][0]
    \GFblirev[0][1]
  }
  =
  -\frac{1}{x}\left(
  \tikzpic{
    \GFbli[0][0]
  }
  -\;
  \tikzpic{
    \GFblirev[0][1]
  }\right).
\end{gather*}
Then, we compute:
\begin{gather*}
  \Delta(x)\;
  \tikzpic{\li}
  \;\overset{*}{\longleftarrow}_{\affT}\;
  \tikzpic{
    \ca[0][2]\uli[2][2]
    \li[0][1]\cro[1][1]
    \cro \li[2][0]
    \dli \cu[1][0]
    \GFbul[0][1][][left]
  }[scale=.8]
  =
  \tikzpic{
    \ca[0][2]\uli[2][2]
    \li[0][1]\cro[1][1]
    \cro \li[2][0]
    \dli \cu[1][0]
    \GFbul[1][0]
  }[scale=.8]
  -x\;
  \tikzpic{
    \ca[0][2]\uli[2][2]
    \li[0][1] \cro[1][1]
    \lili \li[2][0]
    \dli\cu[1][0]
    \GFbul[0][1][][left]
    \GFbul[1][0]
  }[scale=.8]
  -x\;
  \tikzpic{
    \ca[0][2]\uli[2][2]
    \li[0][1] \cro[1][1]
    \cc \li[2][0]
    \dli\cu[1][0]
    \GFbul[0][1][][left]
    \GFbul[1][0]
  }[scale=.8]
  \;\overset{*}{\longrightarrow}_{\affT}\;
  \tikzpic{\li}\;\Delta(x)
  -x\Bigg(
  \tikzpic{
    \ca[0][2]\uli[2][2]
    \li[0][1] \cro[1][1]
    \lili \li[2][0]
    \dli\cu[1][0]
    \GFbul[0][0][][left]
    \GFbul[1][1]
  }[scale=.8]
  +
  \tikzpic{
    \ca[0][2]\uli[2][2]
    \cro[0][1] \li[2][1]
    \li \cu[1][1]
    \dli
    \GFbul[0][2][][left]
    \GFbulrev[0][0]
  }[scale=.8]
  \Bigg)
\end{gather*}
The bubbles are zero, and the term inside parenthesis rewrites as follows:
\begin{gather*}
  \left(
  \tikzpic{
    \ca[0][2]\uli[2][2]
    \li[0][1] \cro[1][1]
    \lili \li[2][0]
    \dli\cu[1][0]
    \GFbul[0][0][][left]
    \GFbul[2][2][][right]
  }[scale=.8]
  -x
  \tikzpic{
    \ca[0][2]\uli[2][2]
    \li[0][1] \lili[1][1]
    \lili \li[2][0]
    \dli\cu[1][0]
    \GFbul[0][0][][left]
    \GFbul[1][1]
    \GFbul[2][2][][right]
  }[scale=.8]
  +x\;
  \tikzpic{
    \ca[0][2]\uli[2][2]
    \li[0][1] \cc[1][1]
    \lili \li[2][0]
    \dli\cu[1][0]
    \GFbul[0][0][][left]
    \GFbul[1][1]
    \GFbul[2][2][][right]
  }[scale=.8]
  \right)
  +
  \left(
  \tikzpic{
    \ca[0][2]\uli[2][2]
    \cro[0][1] \li[2][1]
    \li \cu[1][1]
    \dli
    \GFbul[1][1]
    \GFbulrev[0][0]
  }[scale=.8]
  -x\;
  \tikzpic{
    \ca[0][2]\uli[2][2]
    \lili[0][1] \li[2][1]
    \li \cu[1][1]
    \dli
    \GFbul[0][2][][left]
    \GFbul[1][1]
    \GFbulrev[0][0]
  }[scale=.8]
  -x\;
  \tikzpic{
    \ca[0][2]\uli[2][2]
    \cc[0][1] \li[2][1]
    \li \cu[1][1]
    \dli
    \GFbul[0][2][][left]
    \GFbul[1][1]
    \GFbulrev[0][0]
  }[scale=.8]
  \right)
  \mspace{50mu}
  \\
  \mspace{50mu}
  \;\overset{*}{\longrightarrow}_{\affT}\;
  -
  \left(
  \tikzpic{
    \GFbli[0][0]
    \GFbli[0][1]
  }
  -
  \tikzpic{
    \GFblirev[0][0]
    \GFblirev[0][1]
  }
  \right)
  -x
  \left(
  \tikzpic{
    \GFbli[0][0]
    \GFblirev[0][1]
    \GFbli[0][2]
  }
  +
  \tikzpic{
    \GFblirev[0][0]
    \GFbli[0][1]
    \GFblirev[0][2]
  }
  \right)
  +x
  \left(
  \tikzpic{
    \GFbli[0][0]
    \GFbli[0][1]
  }\;\Delta(x)
  +
  \Delta(x)\;
  \tikzpic{
    \GFblirev[0][0]
    \GFblirev[0][1]
  }
  \right)
  \;\overset{*}{\longrightarrow}_{\affT}\;
  0.
\end{gather*}
This concludes.
\end{proof}

\subsubsection{Bubble-slide critical branchings}
\label{svalt:subsec:sVW_bubble_slides}

\begin{lemma}
  \label{svalt:lem:sVW_bubble_slides}
  For $\affS(\sVW)$, the bubble-slide critical branchings from \cref{subfig:critical_branchings_bubbles_vertsym} are $\rhd$-tamely $\affT^+$-preconfluent.
\end{lemma}

\begin{proof}
  Trivial, since each one of these rewriting steps rewrites to zero.
\end{proof}

%% file: sections/oddnilBrauer.tex
\subsection{The odd nil-Brauer category}

\defnoddnilBr

\oddnilBraffineBrauerthm

Note that here, the bubbles with an odd number of dots have an odd parity, so that the algebra of bubbles in $\bubbleset$ is an exterior algebra rather than a symmetric algebra.
The proof follows the analogous proof for the nil-Brauer category given in \cref{subsec:nilBrauer}, which itself follows the analogous proof for the affine Brauer category in \cref{subsec:affine_brauer_category}.

\subsubsection{Generating functions}
It will be convenient to introduce a new formal variable $y=x^2$ and an evaluation variable $\bfs=\bft^2$. This means that $\bfs$ should be evaluated as a strand with two dots.
We use the following notations:
\begin{gather*}
  \tikzpic{\GFbliev} = \sum_{j\geq 0}x^{2j}\tikzpic{\bli[0][0][][2j]}
  = \tikzpic{\bigbli[0][0][\frac{1}{1-x^2\bft^2}]}
  \quad\an\quad
  \tikzpic{\GFblievrev} = \sum_{j\geq 0}(-1)^jx^{2j}\tikzpic{\bli[0][0][][2j]}
  = \tikzpic{\bigbli[0][0][\frac{1}{1+x^2\bft^2}]}.
\end{gather*}
In the new variable $y$, this translates to:
\begin{gather*}
  \tikzpic{\GFbliev} = \sum_{j\geq 0}y^{j}\tikzpic{\bli[0][0][][2j]}
  = \tikzpic{\bigbli[0][0][\frac{1}{1-y\bfs}]}
  \quad\an\quad
  \tikzpic{\GFblievrev} = \sum_{j\geq 0}(-1)^jy^{j}\tikzpic{\bli[0][0][][2j]}
  = \tikzpic{\bigbli[0][0][\frac{1}{1+y\bfs}]}.
\end{gather*}
We also write $\epsilon =(-1)^{p}$ where $p\coloneqq p(\cup)\in\{0,1\}$ throughout. Note that $\VARev=-\frac{1}{2}(\epsilon+1)$.
We gather below some useful identities.

\begin{lemma}
  The following identities hold in $\oddnilBr[\VARev]$:
  \begingroup
  \renewcommand{\sca}{.4}
  \begin{gather*}
    \tikzpic{\GFbli}
    =
    \tikzpic{\GFbliev}
    +x\;
    \tikzpic{\GFbliev\bli[0][1]}
    \qquad
    \tikzpic{\GFblirev}
    =
    \tikzpic{\GFbliev}
    -x\;
    \tikzpic{\GFbliev\bli[0][1]}
    \qquad
    \tikzpic{\GFbliev\GFblievrev[0][1]}
    =
    \frac{1}{2}
    (\tikzpic{\GFbliev}+\tikzpic{\GFblievrev})
    \qquad
    x^2\;
    \tikzpic{\GFbliev\bli[0][1][][2]}
    =
    \tikzpic{\GFbliev}
    -
    \tikzpic{\li}
    \\
    \tikzpic{\ca\dli\dli[1][0]\GFbulev[1][0]}
    =
    \tikzpic{\ca\dli\dli[1][0]\GFbulevrev[0][0]}
    \qquad
    \tikzpic{\ca\dli\dli[1][0]\GFbulevrev[1][0]}
    =
    \tikzpic{\ca\dli\dli[1][0]\GFbulev[0][0]}
    \qquad
    \tikzpic{\cu\uli\uli[1][0]\GFbulev[1][0]}
    =
    \tikzpic{\cu\uli\uli[1][0]\GFbulevrev[0][0]}
    \qquad
    \tikzpic{\cu\uli\uli[1][0]\GFbulevrev[1][0]}
    =
    \tikzpic{\cu\uli\uli[1][0]\GFbulev[0][0]}
    \\
      \tikzpic{
        \GFbul[0][1] \cro
      }
      \;=\;
      \tikzpic{\cro \GFbulrev[1][0]}
      \;-x\;
      \tikzpic{\li \GFbul[0][1] \li[1][0] \GFbulrev[1][0]} 
      \;-x\epsilon\;
      \tikzpic{\cc \GFbul[0][1] \GFbulrev[1][0] } 
    \qquad   
      \tikzpic{
        \GFbul[1][1] \cro
      }
      \;=\;
      \tikzpic{\cro \GFbulrev[0][0]}
      \;-x\;
      \tikzpic{\li \GFbul[1][1] \li[1][0] \GFbulrev[0][0]} 
      \;-x\;
      \tikzpic{\cc \GFbul[1][1] \GFbulrev[0][0] } 
      \\
      \tikzpic{
        \cro\GFbulev 
      }
      \;=\;
      \tikzpic{\cro \GFbulev[1][1]}
      \;+x^2\;\left(
      \tikzpic{\GFbliev \li[0][1] \GFbliev[1] \bli[1][1]} 
      \;-\;
      \tikzpic{\GFbliev\bli[0][1] \GFbliev[1]  \li[1][1]} 
      \;+\;
      \tikzpic{\cc  \GFbulev[0][0] \bli[1][1] \li[0][1]\GFbulev[1][1] }
      \;-\;
      \tikzpic{\cc \GFbulev[1][1] \bli[0][-1]  \GFbulev[0][0]\li[1][-1]} 
    \right)
      \\
      \tikzpic{
        \cro\GFbulev[0][1] 
      }
      \;=\;
      \tikzpic{\cro \GFbulev[1][0]}
      \;+x^2\;\left(
      \tikzpic{\GFbliev \li[0][1] \GFbliev[1] \bli[1][1]} 
      \;-\;
      \tikzpic{\GFbliev\bli[0][1] \GFbliev[1]  \li[1][1]} 
      \;-\epsilon\;
      \tikzpic{\cc  \GFbulev[1][0] \li[1][1] \bli[0][1]\GFbulev[0][1] }
      \;+\epsilon\;
      \tikzpic{\cc \GFbulev[0][1] \li[0][-1]  \bli[1][-1]\GFbulev[1][0]} 
    \right)
  \end{gather*}
  \endgroup
\end{lemma}

\subsubsection{The odd nil-Brauer rewriting system}

\begin{lemma}\label{nbalt:lem:relations}
  \renewcommand{\sca}{.4}
  \newcommand{\vspc}{.5ex}
  \newcommand{\hspc}{60mu}
  In $\oddnilBr[\VARev]$, the following relations hold:
  \begin{gather*}
    \tikzpic{
      \dli[2][0] \ca[1][0] 
      \uli \cu 
    }
    \;=\;
    (-1)^{p(\cup)}\;
    \tikzpic{\li}
    \mspace{\hspc}
    %
    \tikzpic{
      \ca \uli[2][0]
      \li[0][-1] \cro[1][-1]
      \cro[0][-2] \li[2][-2]
    }[scale=.8]
    =
    0
    \mspace{\hspc}
    %
    \tikzpic{
      \cro \li[2][0]
      \dli \cu[1][0]
    }
    =
    \tikzpic{
      \li[0][0] \cro[1][0]
      \cu[0][0] \dli[2][0]
    }
    \mspace{\hspc}
    %
    \tikzpic{
      \cro[0][1] \li[2][1]
      \li[0][0] \cro[1][0]
      \cu \dli[2][0]
    }[scale=.8]
    =
    0
    \\[\vspc]
    \tikzpic{
      \cro
      \cu
    }
    \;=\;
    0
    \mspace{\hspc}
    %
    \tikzpic{
      \ca[0][1] \uli[2][1]
      \li[0][0] \cro[1][0]
      \cu \dli[2][0]
    }[scale=.8]
    \;=\;
    0
    \mspace{\hspc}
    \tikzpic{\urcr}
    \;=\;
    -\tikzpic{\dlcr}
    -\tikzpic{\lili}
    -\tikzpic{\cc}
    \mspace{\hspc}
    \tikzpic{\rcu \cu[0][0]}
    \;=\;-(-1)^{p(\cup)}\;
    \tikzpic{\lcu \cu[0][0]}.
  \end{gather*}
\end{lemma}

To define a rewriting system of affine Brauer type as in \cref{fig:generic_grobner_of_affine_Brauer_type} associated to $\oddnilBr[\VARev]$, it remains to define the bubble evaluation and the bubble slide.

\begin{lemma}\label{nboddalt:lem:bubble_identities}
In $\oddnilBr[\VARev]$, we have the following identities:
\begin{gather*}
  \Delta_{2k}
  =
  0
  \qquad\text{ for all $k>0$}
  \\
 \tikzpic{
    \li[0][1]  
  }\; \Delta(y)   + \Delta(y) \;    \tikzpic{    \li[0][1]}   \;=\;   \tikzpic{
    \bigbli[3][0][yA(y;\bfs)]
  }[scale=.8] \Delta(y) +\Delta(y)\;     \tikzpic{
    \bigbli[3][0][yA(y;-\bfs)]
  }[scale=.8] 
\end{gather*}
where
\begin{gather*}
\Delta(y)\coloneqq\sum_{k\geq0}y^k\,\Delta_{2k+1}
\quad\an\quad
A(y;\bfs) \coloneqq\bfs \frac{y\bfs-3}{(1-y\bfs)^2}.  
\end{gather*}
\end{lemma}

\begin{proof}
  By moving a dot along a bubble, one sees that a bubble with an even and non-zero number of dots must be zero.
  The second identity will follow from the computation of the second exceptional branching in \cref{lem:odd_nilBr_except_B}.
\end{proof}

We leave to the reader to make explicit the underlying rewriting system of affine Brauer type $\affS(\oddnilBr)$.
As for the nil-Brauer category, its odd counterpart admits vertical symmetry:
\begin{gather}\label{nbalt:eq:vertsymm}
  \left(\tikzpic{\cu}\right)^{*}
  \;=\;
  \tikzpic{\ca}
  \qquad
  \left(\tikzpic{\ca}\right)^{*}
  \;=\;\epsilon\;
  \tikzpic{\cu}
  \qquad
  \left(\tikzpic{\cro}\right)^{*}
  \;=\;
  \tikzpic{\cro}
  \qquad
  \left(\tikzpic{\bli}\right)^{*}
  \;=\;
  \tikzpic{\bli}
\end{gather}
which extends to a vertical symmetry on $\affS(\nilBr[\VARev])$ in the sense of \cref{defn:vertical_symmetries}.

\subsubsection{Critical branchings}

To conclude the proof of \cref{thm:oddnilBr_is_affine_Brauer_type}, we must check that every critical branching is $\rhd$-tamely $\affT^+$-preconf\-luent.
We omit the details for non-affine and dot-slide critical branchings.

\begin{lemma}
  \label{lem:odd_nilBr_except_A}
  In $\affS(\oddnilBr[\VARev])$, the first and third exceptional critical branchings are $\rhd$-tamely $\affT^+$-preconfluent.
\end{lemma}

\begin{proof}
  \renewcommand{\sca}{.4}
  Tamed preconfluence of the third exceptional branching follows from the fact that bubbles with a non-zero even number of dots evaluate to zero.
  For the first exceptional branching, we compute:
\begin{gather*}
  0
  \;\longleftarrow\;
  \tikzpic{\ca[0][1]\cro\GFbul[1][1]\cu}
  \;\longrightarrow\;
  \tikzpic{\ca[0][1]\cro\GFbulrev\cu}
  -
  x\;
  \tikzpic{\ca[0][1]\lili\GFbul[1][1]\GFbulrev\cu}
  -
  x\;
  \tikzpic{\ca[0][1]\cc\GFbul[1][1]\GFbulrev\cu}
  \;\longrightarrow\;
  -
  x\;
  \tikzpic{\ca[0][1]\lili\GFbul[1][1]\GFbulrev\cu}
  -
  x\;
  \tikzpic{\ca[0][1]\cc\GFbul[1][1]\GFbulrev\cu}.
\end{gather*}
We compute each term separately, decomposition the generating function of dots in its even and odd parts, and evaluating bubbles with a non-zero even number of dots to zero:
\begin{IEEEeqnarray*}{rCl}
  \tikzpic{\cu\ca[0][1]\lili\GFbul[1][1]\GFbulrev}
  &=&
  \tikzpic{\cu\ca[0][1]\lili\GFbulev[1][1]\GFbulev}
  +x\;
  \tikzpic{\cu\ca[0][1]\lili\GFbulev[1][1]\GFbulev\bul[1][.7]}
  -x\;
  \tikzpic{\cu\ca[0][1]\lili\GFbulev[1][1]\GFbulev\bul[0][.4]}
  -x^2\;
  \tikzpic{\cu\ca[0][1]\lili\GFbulev[1][1]\GFbulev\bul[1][.7]\bul[0][.4]}
  \;\longleftrightarrow\;
  \VARev-x(1+\epsilon)\;
  \tikzpic{\cu\ca[0][1]\lili\GFbulevrev[0][1]\GFbulev\bul[0][.5]}
  \;\longleftrightarrow\;
  \VARev-x\frac{1}{2}(1+\epsilon)
  \left(\;
  \tikzpic{\cu\ca[0][1]\lili\GFbulev\bul[0][1]}
  +
  \tikzpic{\cu\ca[0][1]\lili\GFbulevrev\bul[0][1]}
  \;\right)
  \\
  \tikzpic{\cu\ca[0][1]\cc\GFbul[1][1]\GFbulrev}
  &=&
  \tikzpic{\cu\ca[0][1]\cc\GFbulev[1][1]\GFbulev}
  +x\;
  \tikzpic{\cu\rca[0][1]\cc\GFbulev[1][1]\GFbulev}
  -x\;
  \tikzpic{\lcu\ca[0][1]\cc\GFbulev[1][1]\GFbulev}
  -x^2\;
  \tikzpic{\lcu\rca[0][1]\cc\GFbulev[1][1]\GFbulev}
  \;\longleftrightarrow\;
  \VARev^2
  -x\VARev\epsilon\;
  \tikzpic{\cu\ca[0][1]\lili\GFbulevrev\bul[0][1]}
  -x\VARev\;
  \tikzpic{\cu\ca[0][1]\lili\GFbulev\bul[0][1]}
\end{IEEEeqnarray*}
Using that $\epsilon=1$ if $\VARev=-1$ and $\epsilon=-1$ if $\VARev=0$, the two terms cancel.
\end{proof}

\begin{lemma}
  \label{lem:odd_nilBr_except_B}
  In $\affS(\oddnilBr[\VARev])$, the second exceptional critical branchings is $\rhd$-tamely $\affT^+$-preconfluent.
\end{lemma}

\begin{proof}
  \renewcommand{\sca}{.4}
We compute the residue of the second exceptional branching:
\begin{IEEEeqnarray*}{rCl}
  0
  \;\longleftarrow\;
  \tikzpic{
    \ca[0][2]\uli[2][2]
    \li[0][1]\cro[1][1]
    \cro \li[2][0]
    \dli \cu[1][0]
    \GFbul[0][1]
  }[scale=.8]
  &
  \;\longrightarrow\;
  &
  \tikzpic{
    \ca[0][2]\uli[2][2]
    \li[0][1]\cro[1][1]
    \cro \li[2][0]
    \dli \cu[1][0]
    \GFbulrev[1][0]
  }[scale=.8]
  -x\;
  \tikzpic{
    \ca[0][2]\uli[2][2]
    \li[0][1] \cro[1][1]
    \lili \li[2][0]
    \dli\cu[1][0]
    \GFbul[0][1]
    \GFbulrev[1][0]
  }[scale=.8]
  -
  \epsilon x\;
  \tikzpic{
    \ca[0][2]\uli[2][2]
    \li[0][1] \cro[1][1]
    \cc \li[2][0]
    \dli\cu[1][0]
    \GFbul[0][1]
    \GFbulrev[1][0]
  }[scale=.8]
  \\
  &=&
  -x
  \left(
  \tikzpic{
    \ca[0][2]\uli[2][2]
    \li[0][1] \cro[1][1]
    \lili \li[2][0]
    \dli\cu[1][0]
    \GFbulev[0][1]
    \GFbulev[1][0]
  }[scale=.8]
  -x^2\;
  \tikzpic{
    \ca[0][2]\uli[2][2]
    \li[0][1] \cro[1][1]
    \lili \li[2][0]
    \dli\lcu[1][0]
    \GFbulev[0][1]\bul[0][.5]
    \GFbulev[1][0]
  }[scale=.8]
  +
  \epsilon\;
  \tikzpic{
    \ca[0][2]\uli[2][2]
    \li[0][1] \cro[1][1]
    \cc \li[2][0]
    \dli\cu[1][0]
    \GFbulev[0][1]
    \GFbulev[1][0]
  }[scale=.8]
  -
  \epsilon x^2\;
  \tikzpic{
    \ca[0][2]\uli[2][2]
    \li[0][1] \cro[1][1]
    \ca\lcu[0][1] \li[2][0]
    \dli\lcu[1][0]
    \GFbulev[0][1]
    \GFbulev[1][0]
  }[scale=.8]
  \right)
  \\
  &&
  {}-x^2
  \left(
  \tikzpic{
    \ca[0][2]\uli[2][2]
    \li[0][1] \cro[1][1]
    \lili \li[2][0]
    \dli\cu[1][0]
    \GFbulev[0][1]\bul[0][.5]
    \GFbulev[1][0]
  }[scale=.8]
  -\;
  \tikzpic{
    \ca[0][2]\uli[2][2]
    \li[0][1] \cro[1][1]
    \lili \li[2][0]
    \dli\lcu[1][0]
    \GFbulev[0][1]
    \GFbulev[1][0]
  }[scale=.8]
  +
  \epsilon\;
  \tikzpic{
    \ca[0][2]\uli[2][2]
    \li[0][1] \cro[1][1]
    \ca\lcu[0][1] \li[2][0]
    \dli\cu[1][0]
    \GFbulev[0][1]
    \GFbulev[1][0]
  }[scale=.8]
  -
  \epsilon\;
  \tikzpic{
    \ca[0][2]\uli[2][2]
    \li[0][1] \cro[1][1]
    \cc \li[2][0]
    \dli\lcu[1][0]
    \GFbulev[0][1]
    \GFbulev[1][0]
  }[scale=.8]
  \right)
\end{IEEEeqnarray*}
On the one hand:
\begin{IEEEeqnarray*}{rCl}
  \tikzpic{
    \ca[0][2]\uli[2][2]
    \li[0][1] \cro[1][1]
    \lili \li[2][0]
    \dli\cu[1][0]
    \GFbulev[0][1]\bul[0][.5]
    \GFbulev[1][0]
  }[scale=.8]
  -\;
  \tikzpic{
    \ca[0][2]\uli[2][2]
    \li[0][1] \cro[1][1]
    \lili \li[2][0]
    \dli\lcu[1][0]
    \GFbulev[0][1]
    \GFbulev[1][0]
  }[scale=.8]
  +
  \epsilon\;
  \tikzpic{
    \ca[0][2]\uli[2][2]
    \li[0][1] \cro[1][1]
    \ca\lcu[0][1] \li[2][0]
    \dli\cu[1][0]
    \GFbulev[0][1]
    \GFbulev[1][0]
  }[scale=.8]
  -
  \epsilon\;
  \tikzpic{
    \ca[0][2]\uli[2][2]
    \li[0][1] \cro[1][1]
    \cc \li[2][0]
    \dli\lcu[1][0]
    \GFbulev[0][1]
    \GFbulev[1][0]
  }[scale=.8]
  &
  \;\longleftrightarrow\;
  &
  -\epsilon\;
  \tikzpic{
    \ca[0][2]\uli[2][2]
    \li[0][1] \lili[1][1]
    \lili \li[2][0]
    \dli\cu[1][0]
    \GFbulev[0][1]
    \GFbulev[1][0]
  }[scale=.8]
  -
  \tikzpic{
    \ca[0][2]\uli[2][2]
    \li[0][1] \cc[1][1]
    \lili \li[2][0]
    \dli\cu[1][0]
    \GFbulev[0][1]
    \GFbulev[1][0]
  }[scale=.8]
  -
  \tikzpic{
    \ca[0][2]\uli[2][2]
    \li[0][1] \lili[1][1]
    \cc \li[2][0]
    \dli\cu[1][0]
    \GFbulev[0][1]
    \GFbulev[1][0]
  }[scale=.8]
  -\epsilon\;
  \tikzpic{
    \ca[0][2]\uli[2][2]
    \li[0][1] \cc[1][1]
    \cc \li[2][0]
    \dli\cu[1][0]
    \GFbulev[0][1]
    \GFbulev[1][0]
  }[scale=.8]
  \\
  &
  \;\longleftrightarrow\;
  &
  -(\epsilon+1)\;
  \tikzpic{\GFbliev\GFblievrev[0][1]}
  -\VARev\left(
    \tikzpic{\GFbliev}
    +
    \tikzpic{\GFblievrev}
  \right)
  \\
  &=&
  -(\epsilon+1)\frac{1}{2}\left(
    \tikzpic{\GFbliev}
    +
    \tikzpic{\GFblievrev}
  \right)
  +\frac{1}{2}(\epsilon+1)\left(
    \tikzpic{\GFbliev}
    +
    \tikzpic{\GFblievrev}
  \right)
  =0.
\end{IEEEeqnarray*}
On the other hand:
\begin{IEEEeqnarray*}{rCl}
  \tikzpic{
    \ca[0][2]\uli[2][2]
    \li[0][1] \cro[1][1]
    \lili \li[2][0]
    \dli\cu[1][0]
    \GFbulev[0][1]
    \GFbulev[1][0]
  }[scale=.8]
  -x^2\;
  \tikzpic{
    \ca[0][2]\uli[2][2]
    \li[0][1] \cro[1][1]
    \lili \li[2][0]
    \dli\lcu[1][0]
    \GFbulev[0][1]\bul[0][.5]
    \GFbulev[1][0]
  }[scale=.8]
  +
  \epsilon\;
  \tikzpic{
    \ca[0][2]\uli[2][2]
    \li[0][1] \cro[1][1]
    \cc \li[2][0]
    \dli\cu[1][0]
    \GFbulev[0][1]
    \GFbulev[1][0]
  }[scale=.8]
  -
  \epsilon x^2\;
  \tikzpic{
    \ca[0][2]\uli[2][2]
    \li[0][1] \cro[1][1]
    \ca\lcu[0][1] \li[2][0]
    \dli\lcu[1][0]
    \GFbulev[0][1]
    \GFbulev[1][0]
  }[scale=.8]
  &
  \;\longleftrightarrow\;
  &
  2\;
  \tikzpic{
    \ca[0][2]\uli[2][2]
    \li[0][1] \cro[1][1]
    \lili \li[2][0]
    \dli\cu[1][0]
    \GFbulev[0][1]
    \GFbulev[1][0]
  }[scale=.8]
  +
  2\epsilon\;
  \tikzpic{
    \ca[0][2]\uli[2][2]
    \li[0][1] \cro[1][1]
    \cc \li[2][0]
    \dli\cu[1][0]
    \GFbulev[0][1]
    \GFbulev[1][0]
  }[scale=.8]
  +
  x^2\epsilon\;
  \tikzpic{
    \ca[0][2]\uli[2][2]
    \li[0][1] \lili[1][1]
    \lili \li[2][0]
    \dli\lcu[1][0]
    \GFbulev[0][1]
    \GFbulev[1][0]
  }[scale=.8]
  +x^2\;
  \tikzpic{
    \ca[0][2]\uli[2][2]
    \li[0][1] \cc[1][1]
    \lili \li[2][0]
    \dli\lcu[1][0]
    \GFbulev[0][1]
    \GFbulev[1][0]
  }[scale=.8]
  +x^2
  \tikzpic{
    \ca[0][2]\uli[2][2]
    \li[0][1] \lili[1][1]
    \ca\lcu[0][1] \li[2][0]
    \dli\cu[1][0]
    \GFbulev[0][1]
    \GFbulev[1][0]
  }[scale=.8]
  +x^2\epsilon\;
  \tikzpic{
    \ca[0][2]\uli[2][2]
    \li[0][1] \cc[1][1]
    \ca\lcu[0][1] \li[2][0]
    \dli\cu[1][0]
    \GFbulev[0][1]
    \GFbulev[1][0]
  }[scale=.8]
  \\
  &
  \;\longleftrightarrow\;
  &
  -
2x^2\left(
  \tikzpic{
    \uli[0][1] \ca[1][1]
    \bli[1][0] \li[2][0]
    \dli\cu[1][0]
    \GFbliev[0][0]     \GFbliev[0][1]
    \GFbulev[1][0] 
  }[scale=.8]
  +
  \tikzpic{
    \uli[3][1] \ca[1][1]
    \bli[1][0] \li[2][0]
    \dli[3][0]\cu[1][0]
    \GFblievrev[3][0]     \GFblievrev[3][1]
    \GFbulev[1][0]
  }[scale=.8]
  \right)
  +
  x^2\;
  \left(
  \tikzpic{
    \uli[0][1] \ca[1][1]
    \bli[1][0] \li[2][0]
    \dli\cu[1][0]
    \GFbliev[0][0]
    \GFbulev[1][0]
  }[scale=.8]
  +
  \tikzpic{
    \uli[3][1] \ca[1][1]
    \bli[1][0] \li[2][0]
    \dli[3][0]\cu[1][0]
    \GFblievrev[3][0]
    \GFbulev[1][0]
  }[scale=.8]
  \right)
\end{IEEEeqnarray*}
where we used the following intermediate computation:
\begin{IEEEeqnarray*}{rCl}
\tikzpic{
\ca[0][2]\uli[2][2]
\li[0][1] \cro[1][1]
\lili \li[2][0]
\dli\cu[1][0]
\GFbulev[0][1]
\GFbulev[1][0]
}[scale=.8]
&\;\overset{}{\longleftrightarrow}\;&
\tikzpic{
\ca[0][2]\uli[2][2]
\li[0][1] \cro[1][1]
\lili \li[2][0]
\dli\cu[1][0]
\GFbulev[0][1]
\GFbulev[2][2]
}[scale=.8]
+x^2
\left(
\tikzpic{
\ca[0][2]\uli[2][2]
\li[0][1] \lili[1][1]
\lili \li[2][0]
\dli\cu[1][0]
\GFbulev[0][1]
\GFbulev[1][0]\GFbulev[2][2]
\bul[2][1]
}[scale=.8]
-
\tikzpic{
\ca[0][2]\uli[2][2]
\li[0][1] \lili[1][1]
\lili \li[2][0]
\dli\cu[1][0]
\GFbulev[0][1]
\GFbulev[1][0]\GFbulev[2][2]
\bul[1][1]
}[scale=.8]
+
\tikzpic{
\ca[0][2]\uli[2][2]
\li[0][1] \cc[1][1]
\lili \li[2][0]
\dli\cu[1][0]
\GFbulev[0][1]
\GFbulev[1][0]\GFbulev[2][2]
\bul[2][1.75]
}[scale=.8]
-
\tikzpic{
\ca[0][2]\uli[2][2]
\li[0][1] \cc[1][1]
\lili \li[2][0]
\dli\cu[1][0]
\GFbulev[0][1]
\GFbulev[1][0]\GFbulev[2][2]
\bul[1][1]
}[scale=.8]
\right)
\\
&\;\overset{}{\longleftrightarrow}\;&
x^2\frac{1}{2}\left(
  \epsilon\tikzpic{\GFbliev\GFbliev[0][1]\bul[0][1]}
  +
  \epsilon\tikzpic{\GFblievrev\GFbliev[0][1]\bul[0][1]}
  +
  \tikzpic{\GFbliev\GFbliev[0][1]\bul[0][1]}
  +
  \tikzpic{\GFblievrev\GFbliev[0][1]\bul[0][1]}
\right)
+x^2\epsilon\VARev
  \tikzpic{\GFbliev\GFbliev[0][1]\bul[0][1]}
-x^2
  \tikzpic{
    \uli[0][1] \ca[1][1]
    \bli[1][0] \li[2][0]
    \dli\cu[1][0]
    \GFbliev[0][0]     \GFbliev[0][1]
    \GFbulev[1][0] 
  }[scale=.8]
\\[2ex]
\tikzpic{
\ca[0][2]\uli[2][2]
\li[0][1] \cro[1][1]
\cc \li[2][0]
\dli\cu[1][0]
\GFbulev[0][1]
\GFbulev[1][0]
}[scale=.8]
&\;\overset{}{\longleftrightarrow}\;&
\tikzpic{
\ca[0][2]\uli[2][2]
\cro[0][1] \li[2][1]
\li \cu[1][1]
\GFbulev[0][2]
\GFbulevrev[0][.5]
}[scale=.8]
\;\overset{}{\longleftrightarrow}\;
\tikzpic{
\ca[0][2]\uli[2][2]
\cro[0][1] \li[2][1]
\li \cu[1][1]
\GFbulev[1][1]
\GFbulevrev[0][.5]
}[scale=.8]
+
x^2\left(
\tikzpic{
\ca[0][2]\uli[2][2]
\lili[0][1] \li[2][1]
\li \cu[1][1]
\GFbulev[0][2]
\GFbulev[1][1]
\GFbulevrev[0][.5]
\bul[1][1.5]
}[scale=.8]
-
\tikzpic{
\ca[0][2]\uli[2][2]
\lili[0][1] \li[2][1]
\li \cu[1][1]
\GFbulev[0][2]
\GFbulev[1][1]
\GFbulevrev[0][.5]
\bul[0][1.5]
}[scale=.8]
-\epsilon
\tikzpic{
\ca[0][2]\uli[2][2]
\ca[0][1]\lcu[0][2] \li[2][1]
\li \cu[1][1]
\GFbulev[0][2]
\GFbulev[1][1]
\GFbulevrev[0][.5]
}[scale=.8]
+\epsilon
\tikzpic{
\ca[0][2]\uli[2][2]
\rca[0][1]\cu[0][2] \li[2][1]
\li \cu[1][1]
\GFbulev[0][2]
\GFbulev[1][1]
\GFbulevrev[0][.5]
}[scale=.8]
\right)
\\
&\;\overset{}{\longleftrightarrow}\;&
x^2\frac{1}{2}\left(
  -\tikzpic{\GFbliev\GFblievrev[0][1]\bul[0][1]}
  -\tikzpic{\GFblievrev\GFblievrev[0][1]\bul[0][1]}
  -\epsilon
  \tikzpic{\GFbliev\GFblievrev[0][1]\bul[0][1]}
  -\epsilon
  \tikzpic{\GFblievrev\GFblievrev[0][1]\bul[0][1]}
\right)
-x^2\epsilon\VARev
  \tikzpic{\GFblievrev\GFblievrev[0][1]\bul[0][1]}
-x^2\epsilon
  \tikzpic{
    \uli[3][1] \ca[1][1]
    \bli[1][0] \li[2][0]
    \dli[3][0]\cu[1][0]
    \GFblievrev[3][0]     \GFblievrev[3][1]
    \GFbulev[1][0]
  }[scale=.8]
\\[2ex]
2\; \tikzpic{\ca[0][2]\uli[2][2]
  \li[0][1] \cro[1][1]
  \lili \li[2][0]
  \dli\cu[1][0]
  \GFbulev[0][1]
  \GFbulev[1][0]
}[scale=.8]
+2
\epsilon\;
\tikzpic{
  \ca[0][2]\uli[2][2]
  \li[0][1] \cro[1][1]
  \cc \li[2][0]
  \dli\cu[1][0]
  \GFbulev[0][1]
  \GFbulev[1][0]
}[scale=.8]
&\;\overset{}{\longleftrightarrow}\;&
x^2(1+\epsilon)\left(
  \tikzpic{\GFbliev\GFbliev[0][1]\bul[0][1]}
  -
  \tikzpic{\GFblievrev\GFblievrev[0][1]\bul[0][1]}
\right)
+2x^2\VARev\left(
  \epsilon\tikzpic{\GFbliev\GFbliev[0][1]\bul[0][1]}
  -
  \tikzpic{\GFblievrev\GFblievrev[0][1]\bul[0][1]}
\right)
-2x^2\left(
  \tikzpic{
    \uli[0][1] \ca[1][1]
    \bli[1][0] \li[2][0]
    \dli\cu[1][0]
    \GFbliev[0][0]     \GFbliev[0][1]
    \GFbulev[1][0] 
  }[scale=.8]
  +
  \tikzpic{
    \uli[3][1] \ca[1][1]
    \bli[1][0] \li[2][0]
    \dli[3][0]\cu[1][0]
    \GFblievrev[3][0]     \GFblievrev[3][1]
    \GFbulev[1][0]
  }[scale=.8]
  \right)
  \\
  &\;\overset{}{\longleftrightarrow}\;&
-
2x^2\left(
  \tikzpic{
    \uli[0][1] \ca[1][1]
    \bli[1][0] \li[2][0]
    \dli\cu[1][0]
    \GFbliev[0][0]     \GFbliev[0][1]
    \GFbulev[1][0] 
  }[scale=.8]
  +
  \tikzpic{
    \uli[3][1] \ca[1][1]
    \bli[1][0] \li[2][0]
    \dli[3][0]\cu[1][0]
    \GFblievrev[3][0]     \GFblievrev[3][1]
    \GFbulev[1][0]
  }[scale=.8]
  \right)
\end{IEEEeqnarray*}
This leads to the following:
\begin{gather*}
  \tikzpic{
    \uli[0][1] \ca[1][1]
    \bli[1][0] \li[2][0]
    \dli\cu[1][0]
    \li[0][0]
    \GFbulev[1][0]
  }[scale=.8]
  +
  \tikzpic{
    \uli[3][1] \ca[1][1]
    \bli[1][0] \li[2][0]
    \dli[3][0]\cu[1][0]
    \li[3][0]
    \GFbulev[1][0]
  }[scale=.8]
  \;\overset{?}{\longleftrightarrow}\;
 -\; \tikzpic{
    \uli[0][1] \ca[1][1]
    \bli[1][0] \li[2][0]
    \dli\cu[1][0]
    \bigbli[0][0][\frac{1+x^2\bft^2}{(1-x^2\bft^2)^2}-1][left]
    \GFbulev[1][0]
  }[scale=.8]
  -
  \tikzpic{
    \uli[3][1] \ca[1][1]
    \bli[1][0] \li[2][0]
    \dli[3][0]\cu[1][0]
    \bigbli[3][0][\frac{1-x^2\bft^2}{(1+x^2\bft^2)^2}-1]
    \GFbulev[1][0]
  }[scale=.8]
\end{gather*}
This gives the expression given in \cref{nboddalt:lem:bubble_identities}, and concludes.
\end{proof}

\begin{lemma}
  \label{lem:odd_nilBr_bubble-slide}
  In $\affS(\oddnilBr[\VARev])$, the bubble-slide critical branchings are $\rhd$-tamely $\affT^+$-preconfluent.
\end{lemma}

\begin{proof}
We have the rewriting step
\begin{gather*}
  \tikzpic{
    \li[0][1]  
  }\; \Delta(y)    \;\longrightarrow\; \VARbbsl \Delta(y) \;    \tikzpic{    \li[0][1]} +  \tikzpic{\GFbbsl}, 
\end{gather*}
with $\VARbbsl=-1$ and where we define inductively 
\begin{gather*}
  \tikzpic{\GFbbsl} = \VARbbsl\; y\;\Delta(y) \; (  \tikzpic{
    \bigbli[3][0][A(y;\bfs)]
  }[scale=.8]-  \tikzpic{
    \bigbli[3][0][A(y;-\bfs)]
  }[scale=.8])+y  \tikzpic{
    \bigbli[0][0][A(y;\bfs)]
    \GFbbsl[0][1]}
\end{gather*}
We conclude that \cref{alt:lem:GF_dot_pushing } holds for the variables $y$ and $\bfs$ in $\oddnilBr[\VARev]$ with $p(y)=0$ and $\VARbbsl=-1$.
 The cap bubble slide then follows from \cref{ alt:lem:cap_slide} if we note that sliding two dots through  a cap gives a minus sign. 
For the bubble slide through crossings, note that \cref{alt:lem:property_P} holds whenever \cref{alt:lem:GF_dot_pushing } and \cref{alt:lem:diagrams_property_P} holds. One can check that a version  of  \cref{alt:lem:diagrams_property_P} with two dots replacing one dot holds in $\oddnilBr[\VARev]$.
We leave the details to the reader.
\end{proof}
\newpage

%% file: sections/quantization.tex
\subsection{The quantized affine VW supercategory}
\label{sec:quantization}

\NewDocumentCommand{\qVWbul}{O{0}O{0}O{}O{right}}{
  \node[diamond,fill=white,inner sep=1.2pt,draw=black] at (#1,#2) {};
  \node[#4] at (#1,#2) {$\scriptstyle{#3}$};
}
\NewDocumentCommand{\qVWbli}{O{0}O{0}O{}O{right}}{
  \li[#1][#2]\qVWbul[#1][#2+.5][#3][#4]
}

\NewDocumentCommand{\bigqVWbul}{O{0}O{0}O{}O{right}}{
  \IfNoValueOrEmptyTF{#3}{}{\node[#4] at (#1,#2) {$\scriptstyle{#3}$};}
  \begin{scope}[scale=3.5,transform shape]
    \node[diamond,fill=white,inner sep=1.2pt,draw=black] at ({(#1)/3.5},{(#2)/3.5}) {};
  \end{scope}
}
\NewDocumentCommand{\bigqVWbli}{O{0}O{0}O{}O{right}}{
  \li[{#1}][{#2}]
  \bigqVWbul[{#1}][{#2+.5}][#3][#4]
}

\NewDocumentCommand{\qVWulcr}{O{0}O{0}}{%
  \cro[{#1}][{#2}]
  \qVWbul[{#1+0.18}][{#2+0.75}]
}
\NewDocumentCommand{\qVWurcr}{O{0}O{0}}{%
  \cro[{#1}][{#2}]
  \qVWbul[{#1+0.82}][{#2+0.75}]
}
\NewDocumentCommand{\qVWdrcr}{O{0}O{0}}{%
  \cro[{#1}][{#2}]
  \qVWbul[{#1+0.82}][{#2+0.25}]
}
\NewDocumentCommand{\qVWdlcr}{O{0}O{0}}{%
  \cro[{#1}][{#2}]
  \qVWbul[{#1+0.18}][{#2+0.25}]
}
\NewDocumentCommand{\qVWlcu}{O{0}O{0}O{}O{left}}{%
  \cu[{#1}][{#2}]
  \qVWbul[{#1+0.18}][{#2-0.3}][#3][#4]
}
\NewDocumentCommand{\qVWrcu}{O{0}O{0}O{}O{right}}{%
  \cu[{#1}][{#2}]
  \qVWbul[{#1+0.82}][{#2-0.3}][#3][#4]
}
\NewDocumentCommand{\qVWlca}{O{0}O{0}O{}O{left}}{%
  \ca[{#1}][{#2}]
  \qVWbul[{#1+0.18}][{#2+0.3}][#3][#4]
}
\NewDocumentCommand{\qVWrca}{O{0}O{0}O{}O{right}}{%
  \ca[{#1}][{#2}]
  \qVWbul[{#1+0.82}][{#2+0.3}][#3][#4]
}

Recall the definition of the quantized affine VW supercategory from the introduction. \getkeytheorem{qsVWformdefn}


This section is devoted to the proof that quantized affine VW supercategory $\quantsVW$ is a category of affine Brauer type (\cref{thm:quantsVW_is_affine_Brauer_type}).
The proof ressembles the analogue proof for the affine VW supercategory in \cref{subsec:affine_VW_supercategory}.

\getkeytheorem{quantsVWaffineBrauerthm}

\subsubsection{Preliminaries}
The associated rewriting system of affine Brauer type is defined using the following secondary relations:

\begin{lemma}
  \label{lem:quantVW_extra_relations}
  \renewcommand{\sca}{.4}
  \newcommand{\vspc}{.8ex}
  \newcommand{\hspc}{40mu}
  The following relations hold in $\quantsVW$:
  \begin{gather*}
    \tikzpic{
        \dli[2][0] \ca[1][0]
        \uli \cu
      }
      \;=\;-\;
      \tikzpic{\li}
      \mspace{\hspc}
      \mspace{\hspc}
      %
      \tikzpic{
        \cu
        \cro
      }
      \;=\;q^{-1}\;
      \tikzpic{\cu}
      \mspace{\hspc}
      %
      \tikzpic{
        \cro \li[2][0]
        \dli \cu[1][0]
      }
      \;=\;
      \tikzpic{
        \li[0][0] \cro[1][0]
        \cu[0][0] \dli[2][0]
      }
      \;+\;(q-q^{-1})\;
      \tikzpic{\cu[0][1]\li[2][0]}
      \\[\vspc]
      \tikzpic{\urcr}
      \;=\;
      \tikzpic{\dlcr}
      \;-\;q\;
      \tikzpic{\lili}
      \;+\;(q-q^{-1})\;
      \tikzpic{\bli\li[1][0]}
      \;+\;q\;
      \tikzpic{\cc}
      \;-\;(q-q^{-1})\;
      \tikzpic{\lca\cu[0][1]}
      \mspace{\hspc}
      \tikzpic{\rcu}
      \;=\;q^{2}\;
      \tikzpic{\lcu}
      \;-q^2\;
      \tikzpic{\cu}
  \end{gather*}
\end{lemma}

To simplify notations, we set:
\[\tikzpic{\qVWbli[0][0]} \coloneqq (q-q^{-1})\;\tikzpic{\bli} - q\;\tikzpic{\li}\;.\]
As usual, we use generating functions to simplify the proof.
We set:
\begin{gather*}
  \renewcommand{\sca}{.4}
  \tikzpic{\GFbli} \coloneqq \tikzpic{\bigbli[0][0][\frac{1}{1-x\bft}]}
  \quad\an\quad
  \tikzpic{\GFblirev[0][0]} \coloneqq \tikzpic{\bigbli[0][0][\frac{1}{1-x(q^{-2}\bft+1)}]},
\end{gather*}
We these notations in place, the dot-slide relations can be written as follows:
\begin{gather*}
  \tikzpic{\ca\GFbulrev\dli\dli[1][0]}
  =
  \tikzpic{\ca\dli\GFbul[1][0]\dli[1][0]}
  \;,\quad
  \tikzpic{\cu\GFbul\uli\uli[1][0]}
  =
  \tikzpic{\cu\uli\GFbulrev[1][0]\uli[1][0]}
  %
  \quad\an\quad
  \tikzpic{\uli[0][1]\uli[1][1]\cro\GFbul[0][1]\dli[0][0]\dli[1][0]}[scale=.8]
  =
  \tikzpic{\uli[0][1]\uli[1][1]\cro\GFbul[1][0]\dli[0][0]\dli[1][0]}[scale=.8]
  -x\;
  \tikzpic{\uli[0][1]\uli[1][1]\GFbul[0][1]\lili\dli[0][0]\dli[1][0]\GFbul[1][0]\qVWbul[0][.5]}[scale=.8]
  -xq^{-2}\;
  \tikzpic{\uli[0][1]\uli[1][1]\GFbul[0][1] \cu[0][1] \qVWlca[0][0]\dli[0][0]\dli[1][0]\GFbul[1][0]}[scale=.8]
  %
  \;.
\end{gather*}

We have the following lemma:
\begin{lemma}
  \label{lem:qsVW_form_derived}
  \renewcommand{\sca}{.4}
  \newcommand{\vspc}{1ex}
  \newcommand{\hspc}{50mu}
  In the rewriting system associated to $\quantsVW$, we have the following congruence, tamed by any diagram with at least two crossings:
  \begin{gather*}
    \tikzpic{\cu\cro\GFbul[0][0][][left]}
    \;\overset{*}{\longleftrightarrow}\;q^{-1}\;
    \tikzpic{\cu\GFbul\uli\uli[1][0]}
    \mspace{\hspc}
    \an
    \mspace{\hspc}
    \tikzpic{
      \ca[0][1]\uli[2][1]
      \GFbli[0][0]\cro[1][0]
      \cu\dli[2][0]
    }[scale=.8]
    \;\overset{*}{\longleftrightarrow}\;-q\;
    \tikzpic{\GFbli}
    \;.
  \end{gather*}
\end{lemma}
\begin{proof}
  The first relation is shown by induction on the coefficients $x^k$, noting that the case $x^0$ follows from \cref{lem:quantVW_extra_relations}:
  \begin{IEEEeqnarray*}{rCl}
    \tikzpic{\cu\cro\bul[0][0][k][left]}[scale=.8]
    &\;\longleftrightarrow\;&
    \tikzpic{\cu\urcr\bul[0][0][k-1][left]}[scale=.8]
    \;+\;q\tikzpic{\cu[0][0]\bul[0][0][k-1][left]\lili}[scale=.8]
    \;-\;(q-q^{-1})\tikzpic{\cu[0][0]\bul[0][0][k-1][left]\bli\li[1][0]}[scale=.8]
    \\
    &\;\overset{*}{\longleftrightarrow}\;&
    q^{-1}\tikzpic{\cu[0][0]\bul[0][0][k-1][left]\li\bli[1][0]}[scale=.8]
    \;+\;q\tikzpic{\cu[0][0]\bul[0][0][k-1][left]\lili}[scale=.8]
    \;-\;(q-q^{-1})\tikzpic{\cu[0][0]\bul[0][0][k-1][left]\bli\li[1][0]}[scale=.8]
    \quad\text{(induction)}
    \\
    &\;\overset{}{\longleftrightarrow}\;&
    q^{-1}
    \left(
      q^2\tikzpic{\cu[0][0]\bul[0][0][k][left]\li\li[1][0]}[scale=.8]
      -q^2 
      \tikzpic{\cu[0][0]\bul[0][0][k-1][left]\li\li[1][0]}[scale=.8]
    \right)
    \;+\;q^{-1}\tikzpic{\cu[0][0]\bul[0][0][k-1][left]\lili}[scale=.8]
    \;-\;(q-q^{-1})\tikzpic{\cu[0][0]\bul[0][0][k-1][left]\bli\li[1][0]}[scale=.8]
    =
    q^{-1}\tikzpic{\cu[0][0]\bul[0][0][k][left]\li\li[1][0]}[scale=.8]\;.
  \end{IEEEeqnarray*}
  The second relation follows from the first, starting the congruence with an (inverse) cap-slide congruence.
\end{proof}

With these preliminaries in place, we can now proceed and apply \cref{thm:affine_Grobner_system_iff_critical_branchings}, checking the tamed congruence of every critical branching.

\subsubsection{Non-affine and dot-slide critical branchings}

\begingroup
\renewcommand{\sca}{.35}
Checking the confluence of the non-affine and dot-slide rewriting steps is cumbersome but straightforward.
They have been checked using our companion computer program (see \cref{subsec:implementation}).
We give only one example, namely the dot-slide critical branching $\tikzpic{\ulcr[0][1]\cro}[scale=.7]$.
As a preliminary, we compute (each $\longleftrightarrow$ is a congruence, suitably tamed):
\begin{gather*}
  \tikzpic{\lca[0][1]\cu[0][2]\cro}
  \;\longleftrightarrow\;
  \tikzpic{\drcr\cc[0][1]}
  \;+\;q\;
  \tikzpic{\lili\cc[0][1]}
  \;-\;(q-q^{-1})\;
  \tikzpic{\bli\li[1][0]\cc[0][1]}
  \;+\;q^{-1}\;
  \tikzpic{\cc\cc[0][1]}
  \;+\;(q^{-3}-q^{-1})\;
  \tikzpic{\lca\cu[0][1]\cc[0][1]}
  \;\longleftrightarrow\;
  -q\;\tikzpic{\lca\cu[0][1]}
  \;,
\end{gather*}
the first term reducing by an upward kink and rewriting $-q\;\tikzpic{\rca\cu[0][1]}$ as $-q^{-1}\;\tikzpic{\lca\cu[0][1]}-q\;\tikzpic{\cc}$.
Using twice the downward kink of \cref{lem:quantVW_extra_relations} and $\stRtwo$, we also get
\begin{IEEEeqnarray*}{rCl}
  \tikzpic{\drcr[0][1]\cro}
  \;&\longrightarrow&\;
  \tikzpic{\dlcr\cro[0][1]}
  \;-\;q\;
  \tikzpic{\lili\cro[0][1]}
  \;+\;(q-q^{-1})\;
  \tikzpic{\bli\li[1][0]\cro[0][1]}
  \;+\;q\;
  \tikzpic{\cc\cro[0][1]}
  \;-\;(q-q^{-1})\;
  \tikzpic{\lca\cu[0][1]\cro[0][1]}
  \\[1ex]
  &\longrightarrow&
  \tikzpic{\bli\li[1][0]}
  \;-\;q\;
  \tikzpic{\cro[0][1]}
  \;+\;
  \tikzpic{\cc}
  \;-\;q^{-1}(q-q^{-1})\;
  \tikzpic{\lca\cu[0][1]}
\end{IEEEeqnarray*}
We can now derive the ``harder'' branch of the critical branching $\tikzpic{\ulcr[0][1]\cro}[scale=.7]$:
\begin{IEEEeqnarray*}{rCl}
  \tikzpic{\ulcr[0][1]\cro}
  \;&\longrightarrow&\;
  \tikzpic{\drcr[0][1]\cro}
  \;+\;q\;\tikzpic{\lili[0][1]\cro}
  \;-\;(q-q^{-1})\;\tikzpic{\bli[0][1]\li[1][1]\cro}
  \;+\;q^{-1}\;\tikzpic{\cc[0][1]\cro}
  \;+\;(q^{-3}-q^{-1})\;\tikzpic{\lca[0][1]\cu[0][2]\cro}
  \;\longrightarrow\;
  \tikzpic{\bli\li[1][0]}
  \;-\;(q-q^{-1})\;\tikzpic{\ulcr}\;.
\end{IEEEeqnarray*}
The other branch, which begins with $\stRtwo$, reduces to the same normal form.
\endgroup

\subsubsection{Bubble-slide and exceptional critical branchings}

Since bubbles all evaluate to zero, one only needs to check the tamed preconfluence of the second exceptional critical branching.
A direct computation gives the following rewriting sequence:
\begin{gather*}
  \tikzpic{
    \li[0][1]\cro[1][1]
    \cro \li[2][0]
    \dli \cu[1][0]
  }[scale=.8]
  \;\overset{*}{\longrightarrow}\;
  \tikzpic{\cu[0][1]\li[2][0]}
\end{gather*}
Since bubbles are zero, the branch that starts by the above steps gives zero:
\begin{gather*}
  \tikzpic{
    \ca[0][2]\uli[2][2]
    \GFbli[0][1]\cro[1][1]
    \cro \li[2][0]
    \dli \custop[1][0]
  }[scale=.8]
  \;\overset{*}{\longrightarrow}\;
  \Delta(x)\;\tikzpic{\li}
  \;\overset{*}{\longrightarrow}\;
  0.
\end{gather*}

The second branch slides the generating function down through the lower crossing.
With the relations above and evaluating all bubbles to zero, only the two terms below remain (ignoring the factor $-x$ in both of them).
Finally, we simplify these terms with \cref{lem:qsVW_form_derived} and zigzag relations:
\begin{gather*}
  \tikzpic{
    \ca[0][2]\uli[2][2]
    \GFbli[0][1]\cro[1][1]
    \lili \li[2][0]
    \dli \cu[1][0]
    \qVWbul[0][.5]
    \GFbul[1][0]
  }
  \;=\;q^{-1}\;
  \tikzpic{\GFbli[0][0]\qVWbli[0][1]\GFblirev[0][2]}
  \quad\an\quad
  q^{-2}\;\tikzpic{
    \ca[0][2]\uli[2][2]
    \GFbli[0][1]\cro[1][1]
    \cc \li[2][0]
    \dli \cu[1][0]
    \qVWbul[.18][.3]
    \GFbul[1][0]
  }
  \;=\;-q^{-2}q\;
  \tikzpic{\GFbli[0][0]\qVWbli[0][1]\GFblirev[0][2]}\;.
\end{gather*}
The two terms cancel. This concludes the proof.

%% file: sections/classification_nonquantum.tex
\subsection{Classification results}
\label{subsec:classification}

In this section, we prove the following classification result:

\begin{restatable}[label=thm:thmclassificationsimpleminded]{theorem}{thmclassificationsimpleminded} 
  \label{classification simple-minded}
  \renewcommand{\sca}{.4}
  \newcommand{\vspc}{.8ex}
  \newcommand{\hspc}{50mu}
  Fix a field $\Bbbk$ where $2$ is invertible.
  Let $\cC$ be a supermonoidal category over $\Bbbk$ presented by an airy presentation of affine Brauer type (\cref{defn:airy_presentation_affine_Brauer_type}).
  Then $\cC$ is a category of affine Brauer type if and only if, up to isomorphism of filtered supermonoidal categories (where the filtration is the filtration by the number of dots), the category $\cC$ admits a presentation of the form
  \begin{gather*}
    \tikzpic{
      \cro[0][1]
      \cro 
    } 
    \;=\;\VARuk\;
    \tikzpic{
      \lili[0][1] 
      \lili[0][0] 
    } 
    \mspace{\hspc}
    %
    \tikzpic{
      \cro[0][2] \li[2][2]
      \li[0][1] \cro[1][1] 
      \cro \li[2][0] 
    }[scale=.8]
    \;=\;\VARRthree\;
    \tikzpic{
      \cro[1][2] \li[0][2]
      \li[2][1] \cro[0][1] 
      \cro[1][0] \li[0][0] 
    }[scale=.8]
    \mspace{\hspc}
    \tikzpic{
      \cu\ca[0][0]
    } 
    \;=\;\VARev
    \\[\vspc]
    \tikzpic{
      \cu[1][0] \uli[2][0] 
      \dli \ca 
    }
    \;=\;
    \tikzpic{\li}
    %
    \mspace{\hspc}
    %
    \tikzpic{
      \ca[0][1]
      \cro 
    }
    \;=\;\VARuk\;
    \tikzpic{\ca}
    \mspace{\hspc}
    %
    \tikzpic{
      \uli[0][1] \ca[1][1] 
      \cro \li[2][0] 
    }
    \;=\;\VARcasl\;
    \tikzpic{ 
      \ca[0][1] \uli[2][1]
      \li \cro[1][0] 
    }
    \\[\vspc]
    \tikzpic{
      \ulcr
    }
    \;=\;
    \VARdcro\;
    \tikzpic{\drcr}
    \;-\;
    \tikzpic{\lili} 
    \;+\;\VARdcroCC\;
    \tikzpic{\cc} 
    \mspace{\hspc}
    \tikzpic{\rca}
    \;=\;\VARdcc\;
    \tikzpic{\lca}
    \;+\;\VARdccLOT\;
    \tikzpic{\ca}
  \end{gather*}
  such that $p(\cup)=p(\cap)$ and the remaining parities and parameters are as follows:
  \begin{center}
    \renewcommand{\sca}{.4}
    \begin{tabular}{@{}c@{\hskip 3ex}cc@{\hskip 3ex}ccc@{\hskip 3ex}cccc@{\hskip 3ex}c@{}}
      $\tikzpic{\cro}$ \& $\tikzpic{\bli}$ &
      $\bubbleset$
      & $\VARev$ & $\VARuk$
      & $\VARRthree$ & $\VARcasl^2$ & $\VARdcro$ & $\VARdcroCC$ & $\VARdcc$ & $\VARdccLOT$ &
      \eg
      \\
      \midrule
      even & $2\bN_>$ & $\VARev$ & $1$ & $1$ & $(-1)^{p(\cup)}$ & $1$ & $\VARcasl^{-1}$ & $-1$ & $0$ & $\affBr$
      \\[2ex]
      even & $\emptyset$ & $0$ & $1$ & $1$ & $-(-1)^{p(\cup)}$ & $1$ & $-\VARcasl^{-1}$ & $1$ & $1$ & $\sVW$
      \\[2ex]
      even & $2\bN+1$ & $\VARev\in\{0,\VARRthree\VARcasl\}$ & $0$ & $\VARRthree^2=1$ & $(-1)^{p(\cup)}$ & $\VARRthree$ & $\VARRthree\VARcasl^{-1}$ & $-\VARRthree$ & $0$ & $\nilBr[\delta]$
      \\[2ex]
      odd & $2\bN+1$ & $\VARev=-\VARRthree\VARcasl$ & $0$ & $\VARRthree^2=1$ & $(-1)^{p(\cup)}$ & $-\VARRthree$ & $-\VARRthree\VARcasl^{-1}$ & $-\VARRthree$ & $0$ & $\oddnilBr[-1]$
      \\
      odd & $2\bN+1$ & $\VARev=0$ & $0$ & $\VARRthree^2=1$ & $-(-1)^{p(\cup)}$ & $-\VARRthree$ & $\VARRthree\VARcasl^{-1}$ & $\VARRthree$ & $0$ & $\oddnilBr[0]$
    \end{tabular}
  \end{center}
\end{restatable}

This theorem implies \cref{thm:classification_as_lin_categories} thanks to \cref{lem:isomorphism-simple-minded-not-monoidal}.
The strategy is simple: derive the system of polynomial equations on the parameters encoded by critical branchings, and solve this system.
One could, a posteriori, give a more concise proof by computing only the necessary critical branchings to arrive at the desired conclusion.
We refrain from doing this, as our purpose is also to highlight the technique of proof: thanks to the Main Theorem (\cref{thm:affine_Grobner_system_iff_critical_branchings}), critical branchings encode every necessary condition, so computing them is garanteed to give a classification.
In fact, this is how the statement \cref{classification simple-minded} was found: see what examples come out of applying the Main Theorem, and then make sense of them.

\medbreak

Throughout, write $p(v)$ the parity of a morphism $v$ and we use the notation $\inter{v}{w}\coloneqq\bil(p(v),p(w))$ from \cref{not:scalar_interchange} for the scalar appearing in the interchange.

\subsubsection{Non-affine conditions}

This subsubsection gives the classification in the non-affine case, obtained by computing the conditions coming from non-affine critical branchings; see \cref{fig:critical_branchings_non_affine_non_quantum}.
\Cref{defn:airy_presentation_Brauer_type} gives the analogous non-affine definition, while \cref{lem:simple_presentation_simple_minded} derives the relevant relations to define the rewriting system.

\begin{definition}\label{defn:airy_presentation_Brauer_type}
  \renewcommand{\sca}{.4}
  \newcommand{\vspc}{.8ex}
  \newcommand{\hspc}{40mu}
  Fix a commutative ring $\Bbbk$.
  An \defnemph{airy presentation of Brauer type} is a presentation of a supermonoidal category with generators crossing, cup, and cap and generating relations
  \begin{gather*}
    \tikzpic{\cro \cro[0][1]} 
    \;=\;
    \VARRtwo\;
    \tikzpic{\lili \lili[0][1]}
    \mspace{\hspc}
    %
    \tikzpic{
      \cro[0][2] \li[2][2]
      \li[0][1] \cro[1][1] 
      \cro \li[2][0] 
    }[scale=.8] 
    \;=\;\VARRthree\;
    \tikzpic{
      \cro[1][2] \li[0][2]
      \li[2][1] \cro[0][1] 
      \cro[1][0] \li[0][0] 
    }[scale=.8] 
    \mspace{\hspc}
    \tikzpic{\cu\ca[0][0]} 
    \;=\;\VARev
    \\[\vspc]
    \tikzpic{
      \cu[1][0] \uli[2][0] 
      \dli \ca 
    }
    \;=\;\VARzz\;
    \tikzpic{\li}
    %
    %
    \mspace{\hspc}
    %
    \tikzpic{
      \ca[0][1]
      \cro 
    }
    \;=\;\VARuk\;
    \tikzpic{\ca}
    \mspace{\hspc}
    %
    \tikzpic{
      \uli[0][1] \ca[1][1] 
      \cro \li[2][0] 
    }
    \;=\;\VARcasl\;
    \tikzpic{ 
      \ca[0][1] \uli[2][1]
      \li \cro[1][0] 
    }
  \end{gather*}
  where each symbol is a scalar in $\Bbbk$ and the scalars $\VARRthree$, $\VARzz$ and $\VARcasl$ are invertible.
\end{definition}

\begin{lemma}\label{lem:simple_presentation_simple_minded}
  Let $\cC$ be a category presented by an airy presentation of Brauer type.
  Up to isomorphism, $\cC$ admits a presentation of the form given in \cref{fig:simple_presentation_simple_minded}.
\end{lemma}

\begin{proof}
  We can normalize the cup (or the cap) such that $\VARzz=1$.
  Moreover, computing the branching
  $\tikzpic{
    \ca[0][2]
    \cro[0][1]
    \cro[0][0]
  }[scale=.5]$
  shows that $\VARRtwo=\VARuk^2$.
  The rest are direct computations.
\end{proof}


\begin{proposition}\label{classification simple minded non-affine}
  \renewcommand{\sca}{.4}
  \newcommand{\vspc}{.8ex}
  \newcommand{\hspc}{40mu}
  Let $\cC$ be a category presented by an airy presentation of Brauer type.
  Up to isomorphism of supermonoidal categories, $\cC$ admits a presentation of the form
  \begin{gather*}
      \tikzpic{
        \cro[0][1]
        \cro 
      } 
      \;=\;
      \VARuk^2\;
      \tikzpic{
        \li \li[1][0] 
      } 
      \mspace{\hspc}
      %
      \tikzpic{
        \cro[0][2] \li[2][2]
        \li[0][1] \cro[1][1] 
        \cro \li[2][0] 
      }[scale=.8]  
      \;=\;\VARRthree\;
      \tikzpic{
        \cro[1][2] \li[0][2]
        \li[2][1] \cro[0][1] 
        \cro[1][0] \li[0][0] 
      }[scale=.8] 
      \mspace{\hspc}
      \tikzpic{
        \cu\ca[0][0]
      } 
      \;=\;\VARev
      \\[\vspc]
      \tikzpic{
        \cu[1][0] \uli[2][0] 
        \dli \ca 
      }
      \;=\;
      \tikzpic{\li}
      %
      %
      \mspace{\hspc}
      %
      \tikzpic{
        \ca[0][1]
        \cro 
      }
      \;=\;\VARuk\;
      \tikzpic{\ca}
      \mspace{\hspc}
      %
      \tikzpic{
        \uli[0][1] \ca[1][1] 
        \cro \li[2][0] 
      }
      \;=\;\VARcasl\;
      \tikzpic{ 
        \ca[0][1] \uli[2][1]
        \li \cro[1][0] 
      }
  \end{gather*}
  where $\VARRthree$ and $\VARcasl$ are invertible, and with the following possible parities and value of parameters:
  \begin{itemize}
    \item If $\VARuk=0$, then $\VARRthree^2 =1$;
    \item if $\VARuk\neq 0$ and $\VARev=0$, then $\VARRthree = \inter{\times}{\times}=1$ and $\VARcasl^4=1$;
    \item if $\VARuk\neq 0$ and $\VARev\neq 0$, then $\VARRthree = \inter{\times}{\times}=1$ and $\VARcasl^2 = \inter{\cup}{\cup}$.
  \end{itemize}
\end{proposition}
\begin{proof}
  We go over the conditions imposed by non-affine critical branchings (see \cref{fig:critical_branchings_non_affine_non_quantum}).
  We find the following conditions on the parameters:
  \begin{gather*}
    \VARRthree^2 =1
    ,\quad
    \VARRthree\VARuk = \VARuk
    ,\quad
    \VARcasl^4\inter{\times}{\times}\VARuk^2 = \VARuk^2
    ,\quad
    \VARcasl^4\VARuk = \VARuk
    ,\quad
    \VARRthree\VARuk^3 = \inter{\times}{\times}\VARuk^3
    \quad\an\quad
    \VARuk\VARev = \VARuk\VARcasl^2\inter{\cup}{\cup}\VARev.
  \end{gather*} 
  The conditions simplify to the ones in the statement.
\end{proof}

\begin{figure}[p]
  \begingroup
  \renewcommand{\sca}{.45}
  \newcommand{\vspc}{1.5ex}
  \newcommand{\hspc}{40mu}
  \begin{gather*}
    \tikzpic{
      \cro[0][1]
      \cro
    }
    =
    \VARuk^2\;
    \tikzpic{
      \lili\lili[0][1]
    }
    \mspace{\hspc}
    %
    \tikzpic{
      \cro[0][2] \li[2][2]
      \li[0][1] \cro[1][1]
      \cro \li[2][0]
    }[scale=.8]
    =\VARRthree\;
    \tikzpic{
      \cro[1][2] \li[0][2]
      \li[2][1] \cro[0][1]
      \cro[1][0] \li[0][0]
    }[scale=.8]
    \mspace{\hspc}
    \tikzpic{
      \cu\ca[0][0]
    }
    =\VARev
    \mspace{\hspc}
    \tikzpic{
      \cu[1][0] \uli[2][0]
      \dli \ca
    }
    =
    \tikzpic{\li}
    \mspace{\hspc}
    %
    \tikzpic{
      \dli[2][0] \ca[1][0]
      \uli \cu
    }
    \;=\inter{\cup}{\cup}\;
    \tikzpic{\li}
    \\[\vspc]
    \tikzpic{
      \uli[0][1] \ca[1][1]
      \cro \li[2][0]
    }
    =\VARcasl\;
    \tikzpic{
      \ca[0][1] \uli[2][1]
      \li \cro[1][0]
    }
    \mspace{\hspc}
    %
    \tikzpic{
      \ca \uli[2][0]
      \li[0][-1] \cro[1][-1]
      \cro[0][-2] \li[2][-2]
    }[scale=.8]
    =
    \VARuk^2\VARcasl^{-1}\;
    \tikzpic{
      \li[0][0] \ca[1][0]
    }
    \mspace{\hspc}
    \tikzpic{
      \cro \li[2][0]
      \dli \cu[1][0]
    }
    =\VARcasl^{-1}\;
    \tikzpic{
      \li[0][0] \cro[1][0]
      \cu[0][0] \dli[2][0]
    }
    \mspace{\hspc}
    \tikzpic{
      \cro[0][1] \li[2][1]
      \li[0][0] \cro[1][0]
      \cu \dli[2][0]
    }[scale=.8]
    =
    \VARuk^2\,\VARcasl\;
    \tikzpic{
      \li[0][0] \cu[1][1]
    }
    \\[\vspc]
    \tikzpic{
      \ca[0][1]
      \cro
    }
    =\VARuk\;
    \tikzpic{\ca}
    \mspace{\hspc}
    %
    \tikzpic{
      \cro
      \cu
    }
    =
      \VARuk\,\VARcasl^2\inter{\cup}{\cup}\;
    \tikzpic{
      \cu
    }
    \mspace{\hspc}
    \tikzpic{
      \ca[0][1] \uli[2][1]
      \li[0][0] \cro[1][0]
      \cu \dli[2][0]
    }[scale=.8]
    =
      \VARuk\,\VARcasl\;
    \tikzpic{
      \li
    }
  \end{gather*}
  \caption{Relations in a category presented by an airy presentation of Brauer type. These relations underpine the associated rewriting system of Brauer type.}
  \label{fig:simple_presentation_simple_minded}
  \endgroup
\end{figure}

\begin{figure}[p]
  \def\tempsp{30mu}
  \renewcommand{\sca}{.35}
  \newcommand{\vertsp}{.5ex}
  \newcommand{\arOrNot}[1]{}
  \begin{gather*}
    \begin{IEEEeqnarraybox}{CcCcCcC}
      \tikzpic{
        \cro[0][2]
        \cro[0][1]
        \cro
      }
      &\mspace{\tempsp}&
      \tikzpic{
        \cro[0][3]\li[2][3]
        \cro[0][2]\li[2][2]
        \li[0][1]\cro[1][1]
        \cro[0][0]\li[2][0]
      }
      \mid
      \tikzpic{
        \cro[0][3]\li[2][3]
        \li[0][2]\cro[1][2]
        \cro[0][1]\li[2][1]
        \cro[0][0]\li[2][0]
      }^{\;\vertsym}
      &\mspace{\tempsp}&
      \tikzpic{
        \cro[0][4]\li[2][4]
        \li[0][3]\cro[1][3]
        \cro[0][2]\cro[2][2]
        \li[0][1]\cro[1][1]
        \cro\li[2][0]
        \draw[] (3,0) to (3,2);\draw[] (3,3) to (3,5);
      }
      &\mspace{\tempsp}&
      \tikzpic{
        \ca[2][1]
        \ca\lili[2][0]
        \dli\cu[1][0]\dli[3][0]
      }
      \mid
      \tikzpic{
        \uli[0][1]\ca[1][1]\uli[3][1]
        \cu[0][1]\lili[2][0]
        \cu[2][0]
      }^{\;\vertsym}
      \\
      &&
      \VARRthree^2a^2=a^2
      &&
      &&
    \end{IEEEeqnarraybox}
    \\[\vertsp]
    \begin{IEEEeqnarraybox}{CcCcCcCcC}
    \tikzpic{
      \uli[0][2]\ca[1][2]
      \cro[0][1]\li[2][1]
      \cro\li[2][0]
    }
    &\mspace{\tempsp}&
    \tikzpic{
      \uli[0][4]\ca[1][4]
      \cro[0][3]\li[2][3]
      \lili[0][2]\li[2][2]
      \li[0][1]\cro[1][1]
      \cro\li[2][0]
    }
    &\mspace{\tempsp}&
    \tikzpic{
      \cro[0][2]\ca[2][2]
      \li[0][1]\cro[1][1]\li[3][1]
      \cro\lili[2][0]
    }
    &\mspace{\tempsp}&
    \tikzpic{
      \uli[0][1]\ca[1][1]\uli[3][1]
      \cro\lili[2][0]
      \dli[0][0]\dli[1][0]\cu[2][0]
    }
    &\mspace{\tempsp}&
    \tikzpic{
      \ca[2][3]
      \ca[0][2]\lili[2][2]
      \li[0][1]\cro[1][1]\li[3][1]
      \cro\lili[2][0]
    }
    \\
    && \VARRthree\VARuk=\VARuk
    && \VARRthree^2 = 1 &&
    && \VARcasl^4\inter{\times}{\times}\VARuk^2=\VARuk^2
    \end{IEEEeqnarraybox}
    \\[\vertsp]
    \begin{IEEEeqnarraybox}{CcCcCcCcC}
    \tikzpic{
      \cro[0][1]\li[2][1]
      \cro\li[2][0]
      \dli[0][0]\cu[1][0]
    }^{\;\vertsym}
    &\mspace{\tempsp}&
    \tikzpic{
      \cro[0][3]\li[2][3]
      \li[0][2]\cro[1][2]
      \lili[0][1]\li[2][1]
      \cro\li[2][0]
      \dli\cu[1][0]
    }^{\;\vertsym}
    &\mspace{\tempsp}&
    \tikzpic{
      \cro[0][2]\lili[2][2]
      \li[0][1]\cro[1][1]\li[3][1]
      \cro\cu[2][1]
    }^{\;\vertsym}
    &\mspace{\tempsp}&
    \tikzpic{
      \uli[0][1]\uli[1][1]\ca[2][1]
      \cro\lili[2][0]
      \dli[0][0]\cu[1][0]\dli[3][0]
    }^{\;\vertsym}
    &\mspace{\tempsp}&
    \tikzpic{
      \cro[0][2]\lili[2][2]
      \li[0][1]\cro[1][1]\li[3][1]
      \cu[0][1]\lili[2][0]
      \cu[2][0]
    }^{\;\vertsym}
    \\
    && \VARRthree\VARuk=\VARuk
    && \VARRthree^2 = 1 &&
    && \VARcasl^4\inter{\times}{\times}\VARuk^2=\VARuk^2
    \end{IEEEeqnarraybox}
    \\[\vertsp]
    \begin{IEEEeqnarraybox}{CcCcCcCcCcC}
    \tikzpic{
      \uli[0][1]\ca[1][1]
      \cro\li[2][0]
      \dli\cu[1][0]
    }
    &\mspace{\tempsp}&
    \tikzpic{
      \uli[0][1]\ca[1][1]
      \cro\cro[2][0]
      \dli\cu[1][0]
      \draw[] (3,-.4) to (3,0);\draw[] (3,1) to (3,1.4);
    }
    &\mspace{\tempsp}&
    \tikzpic{
      \ca[0][3]\uli[2][3]
      \li[0][2]\cro[1][2]
      \lili[0][1]\li[2][1]
      \cro\li[2][0]
      \dli\cu[1][0]
    }
    &\mspace{\tempsp}&
    \tikzpic{
      \ca[0][3]\uli[2][3]
      \li[0][2]\cro[1][2]
      \lili[0][1]\cro[2][1]
      \cro\li[2][0]
      \dli\cu[1][0]
      \draw[] (3,-.4) to (3,1);\draw[] (3,2) to (3,3.4);
    }
    &\mspace{\tempsp}&
    \tikzpic{
      \uli[0][3]\ca[1][3]
      \cro[0][2]\li[2][2]
      \lili[0][1]\li[2][1]
      \li[0][0]\cro[1][0]
      \cu \dli[2][0]
    }^{\;\vertsym}
    &\mspace{\tempsp}&
    \tikzpic{
      \uli[0][3]\ca[1][3]
      \cro[0][2]\li[2][2]
      \lili[0][1]\cro[2][1]
      \li[0][0]\cro[1][0]
      \cu \dli[2][0]
      \draw[] (3,-.4) to (3,1);\draw[] (3,2) to (3,3.4);
    }^{\;\vertsym}
    \\
    \VARcasl^4 a = a
    && \VARcasl^4 a^2 = \inter{\times}{\times}a^2 &&
    && \VARRthree\VARuk^3 = \inter{\times}{\times}\VARuk^3
    && && \VARRthree\VARuk^3 = \inter{\times}{\times}\VARuk^3
    \end{IEEEeqnarraybox}
    \\[\vertsp]
    \begin{IEEEeqnarraybox}{CcCcCcCcCcCcC}
    \tikzpic{
      \ca[0][2]
      \cro[0][1]
      \cro[0][0]
    }
    &\mspace{\tempsp}&
    \tikzpic{
      \ca[0][3]\uli[2][3]
      \cro[0][2]\li[2][2]
      \cro[1][1]\li[0][1]
      \cro[0][0]\li[2][0]
    }
    &\mspace{\tempsp}&
    \tikzpic{
      \ca[1][1]\uli[0][1]
      \cro[0][0]\li[2][0]
      \cu\dli[2][0]
    }
    &\mspace{\tempsp}&
    \tikzpic{
      \ca[0][1]
      \cro[0][0]
      \cu\
    }
    &\mspace{\tempsp}&
    \tikzpic{
      \cro[0][1]
      \cro[0][0]
      \cu[0][0]
    }^{\;\vertsym}
    &\mspace{\tempsp}&
    \tikzpic{
      \cro[0][2]\li[2][2]
      \cro[1][1]\li[0][1]
      \cro[0][0]\li[2][0]
      \cu[0][0]\dli[2][0]
    }^{\;\vertsym}
    &\mspace{\tempsp}&
    \tikzpic{
      \ca[0][1]\uli[2][1]
      \cro[0][0]\li[2][0]
      \dli \cu[1][0]
    }^{\;\vertsym}
    \\
    \VARRtwo = \VARuk^2 
    && \VARRthree\VARuk^3 = \VARuk^3
    && \VARuk\VARcasl^2 = \VARuk\VARcasl^2
    && \VARuk\VARev = \VARuk\VARcasl^2\inter{\cup}{\cup}\VARev
    && \VARcasl^4\VARuk^2=\VARuk^2
    && \VARRthree\VARuk^3 = \VARuk^3
    &&
    \end{IEEEeqnarraybox}
  \end{gather*}
  \caption{Conditions derived from non-affine critical branchings for categories presented by an airy presentation of Brauer type.}
  \label{fig:critical_branchings_non_affine_non_quantum}
\end{figure}

\subsubsection{Dot-slide conditions}

Next, we derive the conditions coming from dot-slide critical branchings.

\begin{lemma}\label{affine-simple-minded-basic-relations}
  If $\cC$ is a supermonoidal category presented by an airy presentation of affine Brauer type, it admits a presentation of the form given in \cref{fig:affine_simple_minded_basic_relations}, where $\tilde{\VARdcroID}=\VARdcroCC\VARcasl^{-1}\inter{\cup}{\cup}$ and $\tilde{\VARdcroCC}=\VARdcroID\VARcasl^{-1}$.
\end{lemma}

\begin{proof}
  Direct computation.
  We use that if $\VARdcroID\neq 0$, then we have $p(\times)+p(\bullet)=0$, so that $\tilde{\VARdcroID}=\VARdcroCC\VARcasl^{-1}\inter{\cup}{\times}\inter{\cup}{\bullet}\inter{\cup}{\cup}=\VARdcroCC\VARcasl^{-1}\inter{\cup}{\cup}$ and $\tilde{\VARdcroCC}=\VARdcroID\VARcasl^{-1}\inter{\cup}{\times}\inter{\cup}{\bullet}=\VARdcroID\VARcasl^{-1}$.
\end{proof}

\begin{proof}[Proof of \cref{classification simple-minded}]
  We have already shown that $\affBr$ (\cref{alt:prop:affBr_main}), $\sVW$ (\cref{svalt:prop:sVW_main}), $\nilBr[\VARev]$ (\cref{nbalt:prop:nilBr_main}) and $\oddnilBr$ (\cref{thm:oddnilBr_is_affine_Brauer_type}) are categories of affine Brauer type for their respective set of bubbles.
  The supermonoidal categories appearing in the statement of \cref{classification simple-minded} are all obtained from $\affBr$, $\sVW$, $\nilBr[\VARev]$ or $\oddnilBr$ by applying \cref{lem:isomorphism-simple-minded-not-monoidal}: it follows that they are categories of affine Brauer type with their respective set of bubbles.
  This shows direction $\Leftarrow$ of the statement. It remains to show that if $\cC$ is a category of affine Brauer type with a presentation as stated, then, up to isomorphism of filtered supermonoidal categories, it must be one of the categories given in \cref{classification simple-minded}.
  To prove this, we go over the conditions imposed by dot-slide critical branchings; see \cref{fig:critical_branchings_dot_slides_non_quantum}.

  Recall that $\tilde{\VARdcroID}=\VARdcroCC\VARcasl^{-1}\inter{\cup}{\cup}$ and $\tilde{\VARdcroCC}=\VARdcroID\VARcasl^{-1}$.
As a preliminary, one can compute the leading term of
$\tikzpic{
  \dlcr \uli[0][1] \ca[1][1] \li[2][0]
  \AntiDiag[0][0]
}[scale=.7]$
and find the condition $\VARdcro^2=1$. We use this equality when computing the remaining conditions.
We mark with $ \;\checkmark$ relations that are readily redundant or follow directly by definition of $\tilde{\VARdcroCC}$ and $\tilde{\VARdcroID}$.
The first Reidemeister-2-move gives $\VARdcroCC=-\VARcasl\VARdcro\inter{\cup}{\cup}\VARdcroID$, which in particular implies that $\VARdcroID\neq 0$ if and only if $\VARdcroCC\neq 0$.
Together with the definitions of $\tilde{\VARdcroCC}$ and $\tilde{\VARdcroID}$ the conditions on parameters become:
\begin{gather*}
  \begin{cases}
    \VARdcro\VARuk\VARcasl^2\inter{\cup}{\cup}{\VARdcroCC} + \VARuk\tilde{\VARdcroCC} &= 0 \\
    \VARdcro\VARuk\VARcasl^2\inter{\cup}{\cup}\tilde{\VARdcroCC} + \VARuk\VARdcroCC &= 0 \\
  \VARdcro\VARdcroID - \VARRthree\inter{\bullet}{\times}\VARdcroID &=0\\
    \VARdcroCC - \VARdcc\inter{\bullet}{\bullet}\tilde{\VARdcroCC} &= 0\\
    (\VARdcroCC - \VARdcro\VARdcc\tilde{\VARdcroCC})\VARdccLOT &= 0 \\
    (\VARdcc - \VARdcc^{-1})\VARuk &=0 \\
    (\inter{\times}{\bullet} - \VARdcro) \VARdccLOT &= 0 \\
    \VARcasl\VARdcroCC + \VARdcro\VARdcc\inter{\times}{\bullet}\tilde{\VARdcroID} &= 0 \\
    \VARcasl\VARdcroID + \VARdcro\VARdcc\inter{\times}{\bullet}\inter{\cup}{\cup}\tilde{\VARdcroCC} &= 0
  \end{cases}
  \quad\Leftrightarrow\quad
  \begin{cases}
    (\VARcasl^4 - 1)\VARdcroID\VARuk &= 0\;\checkmark \\
    \VARdcro\VARuk\VARcasl\inter{\cup}{\cup}\VARdcroID - \VARuk\VARcasl\VARdcro\inter{\cup}{\cup}\VARdcroID &= 0 \;\checkmark \\
  (\VARdcro - \VARRthree\inter{\bullet}{\times})\VARdcroID &=0\\
    (\VARcasl^2\VARdcro\inter{\cup}{\cup} + \VARdcc\inter{\bullet}{\bullet})\VARdcroID &= 0\\
    (\VARcasl^2\inter{\cup}{\cup}+ \VARdcc)\VARdcro\VARdcroID \VARdccLOT &= 0 \\
    (\VARdcc - \VARdcc^{-1})\VARuk &=0 \\
    (\inter{\times}{\bullet} - \VARdcro) \VARdccLOT &= 0 \\
    (\VARcasl^2\VARdcro\inter{\cup}{\cup} + \VARdcc\inter{\times}{\bullet})\VARdcroID &= 0 \;\checkmark  \\
    (\VARcasl^2 + \VARdcro\VARdcc\inter{\times}{\bullet}\inter{\cup}{\cup})\VARdcroID &= 0\;\checkmark 
  \end{cases}
\end{gather*} 
and
\begin{gather*}
  \begin{cases}
    (\VARdcc + \VARdcro)\VARdcroID
      + (\VARdcro + \VARdcc^{-1}) \VARdcro\VARdccLOT\VARuk
      - (\VARdcc\VARcasl^2\VARdcro\inter{\cup}{\cup} + 1)\VARcasl^{-1}\VARdcroID\VARev &= 0 \\
    (\VARdcc + \VARdcro)\VARdcroID
      - (\VARdcro + \VARdcc^{-1})\VARcasl^2\inter{\cup}{\cup} \VARdccLOT\VARuk
      - (\VARdcc\VARcasl^2\VARdcro\inter{\cup}{\cup} + 1)\VARcasl^{-1}\VARdcroID\VARev &= 0 
  \end{cases}
\end{gather*}

\paragraph{Symmetric-like case.}
Assume $\VARuk\neq 0$.
By \cref{classification simple minded non-affine}, this forces $p(\times)=0$, $\VARRthree = 1$ and $\VARcasl^4=1$; moreover, if $\VARev\neq 0$, then $\VARcasl^2 = \inter{\cup}{\cup}$.
Also, since $\VARdcroID\neq 0$, we have that $p(\bullet)=p(\times)=0$.
Using what we know of parities and non-affine parameters, we find that the conditions above are equivalent to $\VARdcro = 1$, $\VARdcc = -\VARcasl^{2}\inter{\cup}{\cup}$ and
\begin{gather*}
  \begin{cases}
    (\VARdcc + \VARdcro)\VARdcroID
      + (\VARdcro + \VARdcc^{-1}) \VARdcro\VARdccLOT\VARuk
      - (\VARdcc\VARcasl^2\VARdcro\inter{\cup}{\cup} + 1)\VARcasl^{-1}\VARdcroID\VARev &= 0 \\
    (\VARdcc + \VARdcro)\VARdcroID
      - (\VARdcro + \VARdcc^{-1})\VARcasl^2\inter{\cup}{\cup} \VARdccLOT\VARuk
      - (\VARdcc\VARcasl^2\VARdcro\inter{\cup}{\cup} + 1)\VARcasl^{-1}\VARdcroID\VARev &= 0 
  \end{cases}
  \\
  \Leftrightarrow
  \begin{cases}
    (1-\VARcasl^{2}\inter{\cup}{\cup})\VARdcroID
      + (1-\VARcasl^{2}\inter{\cup}{\cup}) \VARdccLOT\VARuk &= 0 \\
    (1-\VARcasl^{2}\inter{\cup}{\cup})\VARdcroID
      - (1-\VARcasl^{2}\inter{\cup}{\cup})\VARcasl^2\inter{\cup}{\cup} \VARdccLOT\VARuk &= 0 
  \end{cases}
\end{gather*}
which boils down to the condition that if $\VARcasl^{2}\inter{\cup}{\cup}=-1$, then $\VARdcroID+\VARdccLOT\VARuk=0$.

We can normalize the crossing so that $\VARuk$ is $1$ and the dot so that $\VARdcroID$ is $-1$.
If $\VARcasl^{2}\inter{\cup}{\cup}\neq -1$ then $\VARcasl^{2}\inter{\cup}{\cup} =1$ and $\VARdcc = -1$, so that with the change of generator $\tikzpic{\bli}\mapsto\tikzpic{\bli} + \VARdccLOT/2\;\tikzpic{\li}$, we can assume that $\VARdccLOT = 0$.
Moreover, $\VARdcroCC = \VARcasl^{-1}$.
This gives the first row in \cref{thm:thmclassificationsimpleminded}.
If $\VARcasl^{2}\inter{\cup}{\cup}=-1$, then $\VARdccLOT=1$ by the additional condition above.
Moreover, $\VARdcroCC = -\VARcasl^{-1}$.
This gives the second row in \cref{thm:thmclassificationsimpleminded}.

\paragraph{Nil-like case.}
  Assume $\VARuk = 0$.
  By \cref{classification simple minded non-affine}, $\VARRthree^2=1$.
  The dot-slide conditions boil down to:
  \begin{gather*}
    \VARdcro-\VARRthree\inter{\bullet}{\times}=0,
    \quad
    \VARcasl^2\VARdcro\inter{\cup}{\cup}+\VARdcc\inter{\bullet}{\bullet} = 0
    \quad\an\quad
    (\VARdcc + \VARdcro)
    - (\VARdcc\VARcasl^2\VARdcro\inter{\cup}{\cup} + 1)\VARcasl^{-1}\VARev = 0,
    \\
    (\textbf{if }\VARdccLOT\neq 0:
    \quad
    \VARcasl^2\inter{\cup}{\cup}+\VARdcc = 0
    \quad\an\quad \inter{\times}{\bullet} - \VARdcro = 0).
  \end{gather*}
  Assume $\VARdccLOT\neq 0$.
  This forces $p(\bullet)=0$.
  Then $\VARdcro = \inter{\times}{\bullet} = 1$ and $\VARdcc = -\VARcasl^2\inter{\cup}{\cup}$, so that:
  \begin{align*}
    (\VARdcc + \VARdcro)
      - (\VARdcc\VARcasl^2\VARdcro\inter{\cup}{\cup} + 1)\VARcasl^{-1}\VARev = 0
    &\quad\Leftrightarrow\quad
    (\VARdcc + 1)
      - (-\VARdcc^2 + 1)\VARcasl^{-1}\VARev = 0
    \\
    &\quad\Leftrightarrow\quad
    (\VARdcc+1)(1+(\VARdcc-1)\VARcasl^{-1}\VARev)=0.
  \end{align*}
  It follows that $\VARdcc\neq 1$; with the isomorphism $\tikzpic{\bli}\mapsto\tikzpic{\bli} - \VARdccLOT/(\VARdcc - 1)\;\tikzpic{\li}$, we can assume that $\VARdccLOT = 0$.

  Assume then that $\VARdccLOT = 0$.
  The dot-slide conditions boil down to $\VARdcro=\VARRthree\inter{\bullet}{\bullet}$, $\VARdcc = -\VARRthree\VARcasl^2\inter{\cup}{\cup}$ and
  \begin{equation}
    \label{eq:nil-simple-minded-special-condition}
    (\VARcasl^4\inter{\bullet}{\bullet} - 1)\VARcasl^{-1}\VARev = (\VARcasl^2\inter{\cup}{\cup} - \inter{\bullet}{\bullet})\VARRthree.
  \end{equation}
  Moreover, the first exceptional critical branching with one dot gives that either $\VARev=0$ or $\VARev=\VARdcro\VARcasl$:
  \begin{IEEEeqnarray*}{rCl}
    \renewcommand{\sca}{.3}
  0
  &=&
  \tikzpic{\ca[0][1]\cro\bul[1][1]\cu}
  =
  \VARdcro\;
  \tikzpic{\ca[0][1]\cro\bul\cu}
  +
  x\tilde{\VARdcroID}\;
  \tikzpic{\ca[0][1]\lili\cu}
  +
  x\tilde{\VARdcroCC}\;
  \tikzpic{\ca[0][1]\cc\cu}
  =-x\VARdcro \VARdcroID \VARev + x\VARcasl^{-1}\VARdcroID\VARev^2 \Leftrightarrow \VARev(\VARcasl^{-1} - \VARdcro\VARev)=0.
\end{IEEEeqnarray*}
If $\VARcasl^2\inter{\cup}{\cup}\neq \inter{\bullet}{\bullet}$, then \eqref{eq:nil-simple-minded-special-condition} forces $\VARev\neq 0$ and hence $\VARev=\VARdcro\VARcasl$ by the above, so that \eqref{eq:nil-simple-minded-special-condition} becomes
$\VARcasl^2=\inter{\cup}{\cup}$. In particular, $p(\bullet)=1$.
If instead $\VARcasl^2\inter{\cup}{\cup}= \inter{\bullet}{\bullet}$, then $\VARcasl^4=1$ and \eqref{eq:nil-simple-minded-special-condition} gives that either $p(\bullet)=0$ or $\VARev=0$.

We can normalize the dot so that $\VARdcroID$ is $-1$.
Assume $p(\bullet)=0$.
Then $\VARcasl^2\inter{\cup}{\cup}=1$ by the above, so that $\VARdcro=\VARRthree$, $\VARdcroCC=\VARRthree\VARcasl^{-1}$ and $\VARdcc = -\VARRthree$.
Assume $p(\bullet)=1$.
If $\VARev\neq 0$,
then $\VARev=-\VARRthree\VARcasl$ and $\VARcasl^2=\inter{\cup}{\cup}$;
if instead $\VARev= 0$, then $\VARcasl^2=-\inter{\cup}{\cup}$.
These cases respectively give the third, fourth and fifth rows in \cref{thm:thmclassificationsimpleminded}.
\end{proof}

\begin{figure}[p]
  \begingroup
  \renewcommand{\sca}{.4}
  \newcommand{\vspc}{1ex}
  \newcommand{\hspc}{40mu}
  \begin{gather*}
    \tikzpic{
      \ulcr
    }
    \;=\;\VARdcro\;
    \tikzpic{\drcr}
    \;+\;\VARdcroID\;
    \tikzpic{\lili}
    \;+\;\VARdcroCC\;
    \tikzpic{\cc}
    \mspace{\hspc}
    \tikzpic{
      \drcr
    }
    \;=\;\VARdcro^{-1}\;
    \tikzpic{\ulcr}
    \;-\;\VARdcro^{-1}\VARdcroID\;
    \tikzpic{\lili}
    \;-\;\VARdcro^{-1}\VARdcroCC\;
    \tikzpic{\cc}
    \\[2ex]
    \tikzpic{
      \urcr
    }
    \;=\;\VARdcro\;
    \tikzpic{\dlcr}
    \;+\;\tilde{\VARdcroID}\;
    \tikzpic{\lili}
    \;+\;\tilde{\VARdcroCC}\;
    \tikzpic{\cc}
    \mspace{\hspc}
    \tikzpic{
      \dlcr
    }
    \;=\;\VARdcro^{-1}\;
    \tikzpic{\urcr}
    \;-\;\VARdcro^{-1}\tilde{\VARdcroID}\;
    \tikzpic{\lili}
    \;-\;\VARdcro^{-1}\tilde{\VARdcroCC}\;
    \tikzpic{\cc}
    \\[2ex]
    \tikzpic{\rca}
    \;=\;\VARdcc\;
    \tikzpic{\lca}
    \;+\;\VARdccLOT\;
    \tikzpic{\ca}
    \mspace{\hspc}
    \tikzpic{\lca}
    \;=\;\VARdcc^{-1}\;
    \tikzpic{\rca}
    \;-\;\VARdcc^{-1}\VARdccLOT\;
    \tikzpic{\ca}
    \\[\vspc]
    \tikzpic{\rcu}
    \;=\;\VARdcc^{-1}\;
    \tikzpic{\lcu}
    \;-\;
    \VARdcc^{-1}\VARdccLOT\;
    \tikzpic{\cu}
    \mspace{\hspc}
    \tikzpic{\lcu}
    \;=\;\VARdcc\;
    \tikzpic{\rcu}
    \;+\;
    \VARdccLOT\;
    \tikzpic{\cu}
  \end{gather*}
  \caption{Dot-slide relations in a category presented by an airy presentation of affine Brauer type. These relations, together with \cref{fig:simple_presentation_simple_minded}, underpine the associated rewriting system of affine Brauer type. Here $\tilde{\VARdcroID}=\VARdcroCC\VARcasl^{-1}\inter{\cup}{\cup}$ and $\tilde{\VARdcroCC}=\VARdcroID\VARcasl^{-1}$.}
  \label{fig:affine_simple_minded_basic_relations}
  \endgroup

  \begingroup
    \def\tempsp{20mu}
    \renewcommand{\sca}{.35}
    \newcommand{\vertsp}{.5ex}
    \newcommand{\arOrNot}[1]{}

    \begin{gather*}
      \begin{IEEEeqnarraybox}{CcCcC}
      \tikzpic{
        \drcr \cro[0][1]
        \Diag[0][0][\arOrNot{->}]
      }
      &\mspace{\tempsp}&
      \tikzpic{
        \dlcr \cro[0][1]
        \AntiDiag[0][0][\arOrNot{->}]
      }
      &\mspace{\tempsp}\mspace{\tempsp}&
      \tikzpic{
        \dlcr \li[2][0] \li[0][1] \cro[1][1] \cro[0][2] \li[2][2]
        \AntiDiag[0][0][\arOrNot{->}]
      }
      \mid
      \tikzpic{
        \drcr \li[2][0] \li[0][1] \cro[1][1] \cro[0][2] \li[2][2]
        \Diag[0][0][\arOrNot{->}]
      }
      \mid
      \tikzpic{
        \cro \li[2][0] \li[0][1] \drcr[1][1] \cro[0][2] \li[2][2]
        \Diag[1][1][\arOrNot{->}]
      }
      \\
      \begin{cases}
        \VARdcro{\VARdcroID} + \tilde{\VARdcroID} &= 0 \\
        \VARdcro\VARuk\VARcasl^2\inter{\cup}{\cup}{\VARdcroCC} + \VARuk\tilde{\VARdcroCC} &= 0
      \end{cases}
      &&
      \begin{cases}
        \VARdcro\tilde{\VARdcroID} + \VARdcroID &= 0 \;\checkmark \\
        \VARdcro\VARuk\VARcasl^2\inter{\cup}{\cup}\tilde{\VARdcroCC} + \VARuk\VARdcroCC &= 0
      \end{cases}
      && 
      \VARdcro\VARdcroID = \VARRthree\inter{\bullet}{\times}\VARdcroID
      \end{IEEEeqnarraybox}
      \\[\vertsp]
      \begin{IEEEeqnarraybox}{CcCcCcC}
      \tikzpic{
        \dli \lca[0][0]\ca[0][0][\arOrNot{->}] \cu[1][0] \uli[2][0]
      }
      \;\mid\;
      \tikzpic{
        \uli \lcu[0][0]\cu[0][0][\arOrNot{->}] \ca[1][0] \dli[2][0]
      }^{\;\vertsym}
      &\mspace{\tempsp}\mspace{\tempsp}&
      \tikzpic{
        \ulcr \dlcr
        \Diag[0][0][\arOrNot{<-}]
        \AntiDiag[0][0][\arOrNot{->}]
      }
      &\mspace{\tempsp}\mspace{\tempsp}&
      \tikzpic{
        \dlcr \uli[0][1] \ca[1][1] \li[2][0]
        \AntiDiag[0][0][\arOrNot{->}]
      }
      \mid
      \tikzpic{
        \ulcr \dli \cu[1][0] \li[2][0]
        \Diag[0][0][\arOrNot{<-}]
      }^{\;\vertsym}
      &\mspace{\tempsp}&
      \tikzpic{
        \drcr \uli[0][1] \ca[1][1] \li[2][0]
        \Diag[0][0][\arOrNot{->}]
      }
      \mid
      \tikzpic{
        \dlcr \dli \cu[1][0] \li[2][0]
        \AntiDiag[0][0][\arOrNot{->}]
      }^{\;\vertsym}
      \\
      \\
      &&
      \begin{cases}
        \VARdcroID + \VARdcro\tilde{\VARdcroID} &= 0 \;\checkmark \\
        \VARdcroCC - \VARdcc\inter{\bullet}{\bullet}\tilde{\VARdcroCC} &= 0\\
        (\VARdcroCC - \VARdcro\VARdcc\tilde{\VARdcroCC})\VARdccLOT &= 0
      \end{cases}
      &&
      \begin{cases}
        \VARdcro^2 &= 1\\
        (\inter{\times}{\bullet} - \VARdcro) \VARdccLOT &= 0 \\
        \VARcasl\VARdcroCC + \VARdcro\VARdcc\inter{\times}{\bullet}\tilde{\VARdcroID} &= 0 \\
        \VARcasl\VARdcroID + \VARdcro\VARdcc\inter{\times}{\bullet}\inter{\cup}{\cup}\tilde{\VARdcroCC} &= 0
      \end{cases}
      &&
      \begin{cases}
        \VARdcroID - \VARcasl\tilde{\VARdcroCC} &= 0 \;\checkmark \\
        \VARdcroCC - \VARcasl\inter{\cup}{\cup}\tilde{\VARdcroID} &=0 \;\checkmark
      \end{cases} 
      \end{IEEEeqnarraybox}
      \\[\vertsp]
      \begin{IEEEeqnarraybox}{CcC}
      \tikzpic{
        \dlcr \ca[0][1]
        \AntiDiag[0][0][\arOrNot{->}]
      }
      &\mspace{\tempsp}\mspace{\tempsp}&
      \tikzpic{
        \ulcr
        \cu
        \Diag[0][0][\arOrNot{<-}]
      }^{\;\vertsym}
      \\
      \begin{cases}
        \VARdcc\VARuk = \VARdcc^{-1}\VARuk \\
        (\VARdcc\VARdcroID - \tilde{\VARdcroID})
        + (\VARdcro + \VARdcc^{-1}) \VARdcro\VARdccLOT\VARuk
        + (\VARdcc\VARdcroCC - \tilde{\VARdcroCC})\VARev = 0
      \end{cases}
      &&
      \begin{cases}
        \VARdcc\VARuk = \VARdcc^{-1}\VARuk \;\checkmark \\
        (\VARdcc\VARdcroID - \tilde{\VARdcroID})
        - (\VARdcro + \VARdcc^{-1})\VARcasl^2\inter{\cup}{\cup} \VARdccLOT\VARuk
        + (\VARdcc\VARdcroCC - \tilde{\VARdcroCC})\VARev = 0
      \end{cases}
      \end{IEEEeqnarraybox}
    \end{gather*}
    \caption{Conditions derived from dot-slide critical branchings for categories presented by an airy presentation of affine Brauer type.}
    \label{fig:critical_branchings_dot_slides_non_quantum}
  \endgroup
\end{figure}

%% file: sections/proof_normal_forms.tex
\section{Description of normal forms}
\label{sec:proof_normal_forms}

In this appendix we prove \cref{thm:normal_form_are_standard_diagrams} (as well as \cref{lem:basis_induced_matching_works_for_all}, see below):

\thmnormalformstandarddiagrams*


\begin{lemma}
  \label{lem:reduced_additional_pattern_properties}
  Let $\trafficrule=(\trafficend,\bubbleset)$ be a choice of traffic rules.
  A diagram is a left-cc-justified and braid-justified $(\trafficend,\bubbleset)$-reduced matching if and only if it is a $\affT_\trafficrule$-normal form that avoids, up to interchange, the following additional patterns:
  \begin{gather*}
    \tikzpic{
      \cro[-1][1]\li[1][1]
      \bul[0][1]
      \node[below left=-3pt] at (0,1) {\scriptsize $k$};
      \li[-1][0]\custop[0][1]
    }[scale=.8]
    \qquad
    \tikzpic{
      \ca[-1][1]\li[1][1]
      \bul[0][1]
      \node[below left=-3pt] at (0,1) {\scriptsize $k$};
      \li[-1][0]\custop[0][1]
    }[scale=.8]
    \qquad
    \tikzpic{
      \cro[0][2]
      \bli[0][1]\indexNF[1][1][]
      \node[left] at (0,1.5) {\small $k$};
      \cro
    }[scale=.8]
    \qquad
    \tikzpic{
      \ca[0][2]
      \bli[0][1]\indexNF[1][1][]
      \node[left] at (0,1.5) {\small $k$};
      \cro
    }[scale=.8]
    \qquad
    \tikzpic{
      \cro[0][2]
      \bli[0][1]\indexNF[1][1][]
      \node[left] at (0,1.5) {\small $k$};
      \cu[0][1]
    }[scale=.8]
    \qquad
    \tikzpic{
      \ca[0][2]
      \bli[0][1]\indexNF[1][1][]
      \node[left] at (0,1.5) {\small $k$};
      \cu[0][1]
    }[scale=.8]
  \end{gather*}
  Here $k\geq 0$ and $\tikzpic{\indexNF[1][1][]}[scale=.7]$ is a diagram with at least one input and one output, assumed to be distinct from the identity in the last case.
\end{lemma}

\begin{proof}
  The direction $\Rightarrow$ is clear.
  Consider then a $\affT_\trafficrule$-normal diagram that avoids, up to interchange, the given additional patterns.
  We show that it has the following properties:
  \begin{enumerate}[(a)]
    \setlength{\itemsep}{0em}
    \item it is left-cc-justified;
    \item except for bubbles, a strand contains at most one cup or cap;
    \item a strand does not self-intersect;
    \item two strands intersect at most once.
  \end{enumerate}
  Property (a) follows from the (reduced or additional) pattern-avoidance.
  Indeed, the patterns cover all the ways the left strand of a cup or cap can meet another generator (cup, cap or crossing); sometimes, one needs to reduce to a pattern without dots, arguing that dots must lie close to their $\trafficend$.

  We prove property (b).
  Take a cup (or cap) in $D$ which does not belong to a bubble.
  By property (a), its left strand is a straight line and ends without encountering any cup or cap. Follow instead the right strand of the cup and assume it meets a cap; it must meet that cap from the cap's right strand, since otherwise the cap's left strand is not straight, contradicting property (a).
  But in this situation, the region below the left strand of the cap is enclosed by the cup and its straight left strand:
  \begin{gather*}
    \tikzpic{\draw (0,2) to (0,1) to[out=down,in=left] (1,-.5) to[out=right,in=down] (2,1) to[out=up,in=right] (1.5,1.5) to[out=left,in=up] (1,1);\draw[dotted] (1,1) to (1,0);\draw[very thick] (-.5,2) to (.5,2);}
  \end{gather*}
  The left strand of the cap cannot be a straight line, contradicting property (a).

  Property (c) follows from property (b), as a strand that self-intersects has at least one cup and one cap.

  We prove property (d).
  It follows from properties (a) and (b) that a piece of strand going from one crossing to another does not have any critical point.
  Hence, if two strands intersect more than once, they must intersect under the following schematic:
  \begin{gather*}
    \tikzpic{
      \cro[0][3]\cro
      \draw[dotted] (0,1) to[out=135,in=-135] (0,3);
      \draw[dotted] (1,1) to[out=45,in=-45] (1,3);
    }[scale=.8]\;,
  \end{gather*}
  where a dotted line indicates a piece of strand possibly crossing other strands, but without critical point.
  Consider the left dotted line; if it is a straight line, we find one of the additional forbidden patterns.
  Else, arguing by induction on the number of crossings along the left dotted line, we find two crossings on it that exhibit the same forbidden pattern.
\end{proof}

\begin{lemma}
  \label{lem:standard_is_reduced_plus_patterns}
  Let $\trafficrule=(\trafficend,\bubbleset)$ be a choice of traffic rules.
  A diagram is $(\trafficend,\bubbleset)$-standard if and only if it is a $\affT_\trafficrule$-normal form which, up to interchange, avoids the additional patterns of \cref{lem:reduced_additional_pattern_properties}.
\end{lemma}

\begin{proof}
  A direct check shows that a $(\trafficend,\bubbleset)$-standard diagram is a $\affT_\trafficrule$-normal form that avoids the given patterns.
  Consider then a $\affT_\trafficrule$-normal diagram $D$ that further avoids the given patterns.
  Throughout we use \cref{lem:reduced_additional_pattern_properties}.
  Property (b) shows that $D$ does not have closed strands apart from the bubbles on its leftmost region, and all remaining dots sit close to their dead ends. Going forward, we can assume $D$ does not have any closed strands nor dots.

  By property (a), the left strand of any cap in $D$ is a straight line; reading its bottom endpoint from left to right determines an order on caps.
  Consider the first cap in that order.
  Up to interchange, we are in the following situation:
  \begin{gather*}
    \newcommand{\vsh}{.1}
    \newcommand{\hsh}{.15}
    \tikzpic{
      \draw (0-\hsh,4+\vsh) rectangle (7+\hsh,5-\vsh);
      \draw (0,0) to (0,4+\vsh);
      \draw (1,0) to (1,4+\vsh);
      \draw (2,1-\vsh) to (2,3);
      \draw (3,1-\vsh) to (3,2) to (4,3) to (4,4+\vsh);
      \draw (4,1-\vsh) to (5,2) to (5,4+\vsh);
      \draw (6,1-\vsh) to (6,4+\vsh);
      \draw (7,1-\vsh) to (7,4+\vsh);
      \ca[2][3]\draw (3,3) to (5,1) to (5,1-\vsh);
      \draw (2-\hsh,0+\vsh) rectangle (7+\hsh,1-\vsh);
      \node[] at (.5,0) {\scriptsize\ldots};
      \node[] at (.5,4) {\scriptsize\ldots};
      \node[] at (3.5,1) {\scriptsize\ldots};
      \node[] at (4.5,4) {\scriptsize\ldots};
      \node[] at (6.5,1) {\scriptsize\ldots};
      \node[] at (6.5,4) {\scriptsize\ldots};
      \node at (3.5,4.5) {\scriptsize $T$};
      \draw[<->,thick] (3.5-.2,2.5-.2) to node[draw=none,circle,fill=white,inner sep=.2pt] {\scriptsize $n$} (4.5-.2,1.5-.2);
      \node at (4.5,.5-\vsh) {\scriptsize $B$};
      \draw[line width=2pt] (2,1-2*\vsh) to (5,1-2*\vsh);
    }
  \end{gather*}
  Here $n$ denotes the number of crossings; take $n$ to be maximal.
  Consider the second cap in the order given above; up to interchange, we can assume it lies in $B$. We can push it to the top of $B$; indeed, since cups and crossings cannot have endpoints on the thick line, they can be pushed up to $T$.

  Continuing in this way, we can write the bottom part of $D$ as a cap-standard diagram, and such that all caps of $D$ lie in that bottom part. Arguing similarly with cups, we reduce the proof to $D$ consisting only of crossings.

  Order crossings by reading their bottom endpoints from left to right (considering the first endpoint that occurs).
  Consider the first crossing in that order.
  Using the third additional forbidden pattern, that crossing is also the first when reading \emph{top} endpoints from left to right.
  This leads to the context defining braid-standard diagrams.
\end{proof}

Combining \cref{lem:reduced_additional_pattern_properties} and \cref{lem:standard_is_reduced_plus_patterns}, which respectively characterize conditions (iii) and (ii) of the theorem via pattern-avoidance,
we already have the directions (i) $\Leftarrow$ (ii) $\Leftrightarrow$ (iii).
The additional claim follows from the fact that there is a one-to-one correspondence between matchings and standard diagrams without dots or bubbles.
It remains to show (i) $\Rightarrow$ (ii), i.e.\ a monomial $\affT_\trafficrule$-normal form $D$ is, up to interchange, $(\trafficend,\bubbleset)$-standard, and that it has properties (a)-(d). We proceed by induction on the number of generators in $D$.
Following \cref{lem:reduced_additional_pattern_properties}, it suffices to show that, up to interchange, $D$ avoids the given additional patterns.
This is done in the following two lemmas.

\begin{lemma}
  \label{lem:defn_normal_form_avoid_stop_sign_patterns}
  Under the induction hypothesis, the diagram $D$ avoids the following patterns:
  \begin{gather*}
    \tikzpic{
      \cro[-1][1]\li[1][1]
      \bul[0][1]
      \node[below left=-3pt] at (0,1) {\scriptsize $k$};
      \li[-1][0]\custop[0][1]
    }[scale=.8]\;,
    \qquad
    \tikzpic{
      \ca[-1][1]\li[1][1]
      \bul[0][1]
      \node[below left=-3pt] at (0,1) {\scriptsize $k$};
      \li[-1][0]\custop[0][1]
    }[scale=.8]\;,
    \qquad
    \tikzpic{
      \cro[0][1]
      \bul[0][1]
      \node[left] at (0,1) {\scriptsize $k$};
      \custop[0][1]
    }[scale=.8]\;,
    \qquad
    \tikzpic{
      \cro[0][2]
      \bli[0][1]\cro[1][1]\li[2][2]
      \node[left] at (0,1.5) {\scriptsize $k$};
      \custop[0][1]\dli[2][1]
    }[scale=.8]
    \qquad\an\qquad
    \tikzpic{
      \ca[0][2]\uli[2][2]
      \bli[0][1]\cro[1][1]
      \node[left] at (0,1.5) {\scriptsize $k$};
      \custop[0][1]\dli[2][1]
    }[scale=.8]\;.
  \end{gather*}
\end{lemma}

\begin{proof}
  Let $P$ be such a pattern and $\Gamma$ a context such that $\Gamma[P]=D$; we assume $\Gamma$ is such that there is a stop sign on the cup, as depicted.
  Write $\Gamma = T\circ(\id_a\otimes -\otimes \id_b)\circ B$ and call $T$ the top diagram and $B$ the bottom diagram.
  
  Consider the top rightmost endpoint of $P$, and follow the strand into the top diagram $T$. Since the cup has a stop sign, the strand is closed in $\Gamma[P]$, so that we must meet a cap in $T$.
  If we meet the cap from its left strand, by the induction hypothesis applied to $T$, the strand is a straight line from the endpoint to the cap; up to interchange, this gives a forbidden pattern.
  If we meet the cap from its right strand, we continue to follow the strand; if it meets a top endpoint of $P$, we can again apply the induction hypothesis to $T$ and find forbidden patterns.
  Finally, consider the situation where the strand continues to the left of $P$, into the bottom diagram $B$.
  Since the strand is closed, it must meet a cup in $B$.
  Schematically, we have the following situation:
  \begin{gather*}
    \tikzpic{
      \node[fill,circle,inner sep=1.5pt] (A) at (1,0) {};
      \node[fill,circle,inner sep=1.5pt] (B) at (0,1.5) {};
      \node[fill,circle,inner sep=1.5pt] (C) at (-1,-1.5) {};
      \draw (A) to[out=90,in=-45] (B) to [out=-135,in=90] (C);
      \draw[dashed] (-1.5,.75) to (1.5,.75);
      \draw[dashed] (-1.5,-.75) to (1.5,-.75);
      \node at (-3,1.5) {$T$};
      \node at (-3,0) {$P$};
      \node at (-3,-1.5) {$B$};
    }
  \end{gather*}
  Any interchange applied to the strand induces an interchange to the diagram above; in particular, the cup in $P$ cannot be $\textsc{lint}$-minimal, and hence cannot have a stop sign. (This kind of argument is an idea from \cite{DV_NormalizationPlanarString_2022}.)
\end{proof}

\begin{lemma}
  Under the induction hypothesis, the diagram $D$ avoids the additional patterns of \cref{lem:reduced_additional_pattern_properties}.
\end{lemma}

\begin{proof}
  For the first two additional patterns, this follows from the above lemma.
  Let $P$ be one of the last four additional patterns and $\Gamma$ a context such that $\Gamma[P]=D$. By the induction hypothesis, we can assume $\tikzpic{\indexNF[0][0][]}[scale=.8]$ avoids the additional patterns of \cref{lem:reduced_additional_pattern_properties} and satisfies properties (a)-(d) from its proof.

  Assume $P$ is the first such pattern; we can ignore the $k$ dots, as there exists an orientation of the strand that allows one to slide them away.
  Follow the bottom right strand of the top crossing. It goes to the bottom of $\tikzpic{\indexNF[0][0][]}[scale=.8]$ without critical points.
  Indeed, if it meets a cap from its left strand, the latter must be a straight line to the top of $\tikzpic{\indexNF[0][0][]}[scale=.8]$, that is, to the crossing; up to interchange, this leads to the forbidden pattern $\tikzpic{\cro[0][1]\li[2][1]\dli[0][1]\cu[1][1]}[scale=.7]$.
  On the other hand, if it meets the cap from its right strand, then its left strand must necessarily cross a strand in $\tikzpic{\indexNF[0][0][]}[scale=.8]$, again contradicting the induction hypothesis.
  As the same argument applies when one follows the top right strand of the bottom crossing, we are in the following situation:
  \begin{gather*}
    \renewcommand{\tempdiag}{
      \draw (0,0) to (3,3);
      \draw (1,0) to (0,1) to (0,6) to (1,7);
      \draw (2,0) to (2,1) to (1,2) to (1,3);
      \draw (1,4) to (1,5) to (2,6) to (2,7);
      \begin{scope}[shift={(1,0)}]
        \draw (2,0) to (2,1) to (1,2) to (1,3);
        \draw (1,4) to (1,5) to (2,6) to (2,7);
      \end{scope}
      \indexNF[1][3][]
      \draw (3,4) to (0,7);
      \node at (2.5,0) {\scriptsize $\ldots$};
      \node at (2.5,7) {\scriptsize $\ldots$};
      \draw[line width=1pt,<->] (2,.5) to node[above=-1pt] {\scriptsize $m$}(3,.5);
      \draw[line width=1pt,<->] (2,6.5) to node[below=-1pt] {\scriptsize $n$}(3,6.5);
    }
    \tikzpic{
      \tempdiag
      \li[3][3]
    }[scale=.9]
    \qquad\text{ or }\qquad
    \tikzpic{
      \tempdiag
      \cro[3][3]
      \draw (4,4) to (4,7);\draw (4,0) to (4,3);
    }[scale=.9]
    \;.
  \end{gather*}
  Using again the left-justified property, we deduce that $\tikzpic{\indexNF[0][0][]}[scale=.8]$ in the diagram above consists only of crossings (bubbles can be slid away). If $n=m=0$, this is the R2 or R3 forbidden pattern.
  Otherwise, we have $n=m\neq 0$ and viewing $\tikzpic{\indexNF[0][0][]}[scale=.8]$ as a braid-standard diagram, we find one of the forbidden additional patterns.
  
  A similar argument applies to the other patterns.
  For the second pattern, we can slide away the $k$ dots, and the case $n=m=0$ corresponds to the forbidden patterns $\tikzpic{\cro\ca[0][1]}[scale=.7]$
  and $\tikzpic{\cro\li[0][1]\cro[1][1]\ca[0][2]}[scale=.7]$ from \eqref{eq:forbidden_patterns} (\cref{sec:definitions}).
  For the third pattern, the case $n=m=0$ corresponds to the third and fourth patterns of the previous lemma.
  Finally, for the fourth pattern, the case $n=m=0$ corresponds to the last pattern of the previous lemma.
\end{proof}


\begin{lemma}
  \label{lem:reduced_curves}
  Let $(\trafficend,\bubbleset)$ be a choice of traffic rules.
  Any two equivalent $(\trafficend,\bubbleset)$-reduced matchings are related by interchange, R3 moves, and dot, cup and cap slides.
\end{lemma}

\begin{proof}
  As a consequence of \cref{thm:normal_form_are_standard_diagrams}, any diagram without closed strands can be written as a standard diagram by successively applying interchange, R2 and R3 moves, dot, cup and cap slides, cup and cap pullings, kinks and zigzags, and such that one never increases the number of generators.
  In particular, if the diagram was a reduced matching, this only uses interchange, R3 moves, and dot, cup and cap slides.
  It also follows from \cref{thm:normal_form_are_standard_diagrams} that two equivalent reduced matchings rewrite as the same standard diagram. The result follows.
  (One could also give topological arguments following \cite{HS_ShorteningCurvesSurfaces_1994}\footnote{We thank David Freund for pointing us to that reference.}).
\end{proof}

We can now prove \cref{lem:basis_induced_matching_works_for_all}:

\lembasisinduced*

\begin{proof}
  Pairings are partitioned by their number of inversions, which amounts to the number of crossings in a reduced matching.
  Let $A=\sqcup A_k$ and $B=\sqcup B_k$ be two choices of representative sets for $(\trafficend,\bubbleset)$-reduced matchings, such that $A$ is a hom-basis. (Here $k$ is the number of crossings.)
  Thanks to \cref{lem:reduced_curves} and the invertibility conditions (ii) (see \cref{subsubsec:defn_rewriting_system}), any element $a\in A_k$ can be expressed as $a=\alpha b + v$ where $\alpha$ is an invertible scalar, $b\in B_k$ and $v$ is a linear combination of elements in $B_{\leq k-1}\coloneqq B_0\sqcup\ldots\sqcup B_{k-1}$. (And vice versa.)
  In particular, the set $B$ is spanning. We show linear independence by induction on $k$. When $k=0$, the vector $v$ is zero, and going from a linear combination in $B_0$ to a linear combination in $A_0$ is just a renormalization by invertible scalars.
  Let then $k$ be generic and assume that $B_{\leq k}$ is a linearly independent set.
  Consider a linear combination $\sum_i\lambda_ib_i + w$ in $B_{\leq k+1}$, with $b_i\in B_{k+1}$ and $w\in\langle B_{\leq k}\rangle$, such that $\sum_i\lambda_ib_i=0$. It rewrites as $\sum_i\lambda_ib_i=\sum_i\lambda_i\alpha_ia_i+v$, with $\alpha_i$ invertible scalars, $a_i$ elements in $A_k$ and $v\in\langle B_{\leq k}\rangle$.
  Using filtration and invertibility, we have $\lambda_i\alpha_i=0$ and hence $\lambda_i=0$, and induction shows that the coefficients in $w$ are zero.
\end{proof}

%% file: sections/appendix_symmetries.tex
\section{Schematics for horizontal symmetries}
\label{sec:horizontal_symmetries}

This appendix contains the schematics for horizontal symmetries; see \cref{subsubsec:horizontal_symmetry}.

\begin{figure}[H]
  \centering
  \begin{gather*}
  \begingroup
  \savebox{\tempboxa}{\tikzpic{
    \uli[0][4]\uli[1][4]\ca[2][4]
    \cro[0][3]\lili[2][3]
    \li[0][2]\cro[1][2]\li[3][2]
    \cro[0][1]\lili[2][1]
    \cro[0][0]\lili[2][0]
    \draw[branch_overlay] (-.3,1) rectangle (3+.3,4.4);
  }[scale=.4]}
  \savebox{\tempboxb}{\tikzpic{
    \uli[0][4]\uli[1][4]\ca[2][4]
    \cro[0][3]\lili[2][3]
    \li[0][2]\cro[1][2]\li[3][2]
    \cro[0][1]\lili[2][1]
    \cro[0][0]\lili[2][0]
    \draw[branch_overlay] (-.3,0) rectangle (2+.3,4);
  }[scale=.4]}
  \def\scl{.6}
  \begin{tikzcd}[ampersand replacement=\&,row sep=large]
  	\tikzpic{
      \uli[0][4]\ca[1][4]\uli[3][4]
      \lili[0][3]\cro[2][3]
      \li[0][2]\cro[1][2]\li[3][2]
      \cro[0][1]\lili[2][1]
      \li[0][0]\cro[1][0]\li[3][0]
    }[scale=\scl]
    \&
    \tikzpic{
      \uli[0][4]\ca[1][4]
      \cro[0][3]\li[2][3]
      \lili[0][2]\cro[2][2]
      \li[0][1]\cro[1][1]
      \cro\li[2][0]
      \draw[] (3,0) to (3,2);\draw[] (3,3) to (3,4.4);
    }[scale=\scl*1.2]
    \&
    \tikzpic{
      \ca[0][4]\uli[2][4]
      \li[0][3]\cro[1][3]
      \lili[0][2]\cro[2][2]
      \li[0][1]\cro[1][1]
      \cro\li[2][0]
      \draw[] (3,0) to (3,2);\draw[] (3,3) to (3,4.4);
    }[scale=\scl]
    \&
    \tikzpic{
      \cro[2][4]
      \ca[0][3]\lili[2][3]
      \li[0][2]\cro[1][2]\li[3][2]
      \lili[0][1]\cro[2][1]
      \cro\lili[2][0]
    }[scale=\scl]
    \\
  	\tikzpic{
      \uli[0][4]\uli[1][4]\ca[2][4]
      \li[0][3]\cro[1][3]\li[3][3]
      \li[0][2]\cro[1][2]\li[3][2]
      \cro[0][1]\lili[2][1]
      \li[0][0]\cro[1][0]\li[3][0]
    }[scale=\scl]
    \&
    \tikzpic{
      \uli[0][4]\uli[1][4]\ca[2][4]
      \li[0][3]\cro[1][3]\li[3][3]
      \cro[0][2]\lili[2][2]
      \li[0][1]\cro[1][1]\li[3][1]
      \cro[0][0]\lili[2][0]
    }[scale=\scl]
    \&
    \tikzpic{
      \uli[0][4]\uli[1][4]\ca[2][4]
      \cro[0][3]\lili[2][3]
      \li[0][2]\cro[1][2]\li[3][2]
      \cro[0][1]\lili[2][1]
      \cro[0][0]\lili[2][0]
    }[scale=\scl]
    \&
    \tikzpic{
      \cro[1][4]
      \AntiDiag[0][3]\ca[1][3]\Diag[2][3]
      \lili[0][2]\cro[2][2]
      \cro[0][1]\lili[2][1]
      \cro[0][0]\lili[2][0]
    }[scale=\scl]
    \\
  	\tikzpic{
      \uli[0][4]\uli[1][4]\ca[2][4]
      \li[0][3]\li[3][3]
      \li[0][2]\LOTRtwo[1][2]\li[3][2]
      \cro[0][1]\lili[2][1]
      \li[0][0]\cro[1][0]\li[3][0]
    }[scale=\scl]
  	\&\&
    \tikzpic{
      \uli[0][4]\uli[1][4]\ca[2][4]
      \cro[0][3]\lili[2][3]
      \li[0][2]\cro[1][2]\li[3][2]
      \lili[2][1]
      \LOTRtwo[0][0]\lili[2][0]
    }[scale=\scl]
    \&
    \tikzpic{
      \cro[1][4]
      \AntiDiag[0][3]\ca[1][3]\Diag[2][3]
      \lili[0][2]\cro[2][2]
      \lili[2][1]
      \LOTRtwo[0][0]\lili[2][0]
    }[scale=\scl]
  	\arrow[curve={height=40pt}, from=1-1, to=3-1]
  	\arrow[from=1-2, to=1-1]
  	\arrow[from=1-2, to=1-3]
  	\arrow[from=1-3, to=1-4]
  	\arrow[draw=none, from=1-3, to=2-3,"\usebox{\tempboxa}"{description}]
  	\arrow[curve={height=-40pt}, from=1-4, to=3-4]
  	\arrow[from=2-1, to=1-1]
  	\arrow[from=2-1, to=3-1]
  	\arrow["\ulcorner"{description, pos=0}, draw=none, from=2-2, to=1-1]
  	\arrow[from=2-2, to=1-2]
  	\arrow[from=2-2, to=2-1]
  	\arrow[from=2-3, to=2-2]
  	\arrow[from=2-3, to=2-4]
  	\arrow[from=2-3, to=3-3]
  	\arrow["\lrcorner"{description, pos=0}, draw=none, from=2-3, to=3-4]
  	\arrow[from=2-4, to=1-4]
  	\arrow[from=2-4, to=3-4]
  	\arrow[""{name=0, anchor=center, inner sep=0}, tail reversed, dottedcd, from=3-1, to=3-3]
  	\arrow[tail reversed, dottedcd, from=3-3, to=3-4]
  	\arrow[draw=none, from=2-2, to=0,"\usebox{\tempboxb}"{description}]
  \end{tikzcd}
  \endgroup
  \\[2ex]
  \begingroup
  \savebox{\tempboxa}{\tikzpic{
    \uli[0][2]\ca[1][2]\uli[3][2]
    \cro[0][1]\cu[2][2]
    \cro[0][0]
    \draw[branch_overlay] (-.3,1) rectangle (3.3,2.4);
  }[scale=.4]}
  \def\scl{.6}
  \begin{tikzcd}[ampersand replacement=\&]
  	\tikzpic{ 
      \ca[0][2]\uli[2][2]\uli[3][2]
      \li[0][1]\cro[1][1]\li[3][1]
      \cro[0][0]\cu[2][1]
    }[scale=\scl*1.2]
    \&
  	\tikzpic{ 
      \uli[0][2]\ca[1][2]\uli[3][2]
      \cro[0][1]\cu[2][2]
      \cro[0][0]
    }[scale=\scl]
    \&
  	\tikzpic{ 
      \uli[0][2]\ca[1][2]\uli[3][2]
      \cu[2][2]
      \LOTRtwo[0][0]
    }[scale=\scl]
    \\
  	\tikzpic{ 
      \cro[2][4]
      \ca[0][3]\lili[2][3]
      \li[0][2]\cu[1][3]\li[3][2]
      \Diag[0][1]\AntiDiag[2][1]
      \cro[1][0]
    }[scale=\scl*.8]
  	\&
    \tikzpic{ 
      \cro[0][1]
      \cro[0][0]
    }[scale=\scl]
    \& 
    \tikzpic{ 
      \LOTRtwo[0][0]
    }[scale=\scl]
  	\arrow[curve={height=-40pt}, from=1-1, to=1-3]
  	\arrow[from=1-1, to=2-1]
  	\arrow[from=1-2, to=1-1]
  	\arrow[from=1-2, to=1-3]
  	\arrow[draw=none, from=1-2, to=2-1, "\usebox{\tempboxa}"{description}]
  	\arrow[from=1-2, to=2-2]
  	\arrow["\lrcorner"{description, pos=0}, draw=none, from=1-2, to=2-3]
  	\arrow[tail reversed, dottedcd, from=1-3, to=2-3]
  	\arrow[from=2-1, to=2-2]
  	\arrow[from=2-2, to=2-3]
  \end{tikzcd}
  \endgroup
  \qquad
  \begingroup
  \savebox{\tempboxa}{\tikzpic{
    \uli[0][3]\ca[1][3]
    \cro[0][2]\li[2][2]
    \cro[0][1]\indexNF[2][1][] 
    \cro[0][0]\li[2][0]
    \dli[0][0]\cu[1][0]
    \draw[branch_overlay] (-.3,0) rectangle (1+.3,3);
  }[scale=.4]}
  \def\scl{.6}
  \begin{tikzcd}[ampersand replacement=\&]
  	\tikzpic{
      \ca[0][3]\uli[2][3]
      \li[0][2]\cro[1][2] 
      \cro[0][1]\indexNF[2][1][]
      \li[0][0]\cro[1][0]
      \dli[2][0]\cu
    }[scale=\scl]
    \&
  	\tikzpic{
      \uli[0][3]\ca[1][3]
      \cro[0][2]\li[2][2]
      \cro[0][1]\indexNF[2][1][] 
      \li[0][0]\cro[1][0]
      \dli[2][0]\cu
    }[scale=\scl]
    \&
  	\tikzpic{
      \uli[0][3]\ca[1][3]
      \li[2][2]
      \LOTRtwo[0][1]\indexNF[2][1][] 
      \li[0][0]\cro[1][0]
      \dli[2][0]\cu
    }[scale=\scl]
    \\
  	\tikzpic{
      \ca[0][3]\uli[2][3]
      \li[0][2]\cro[1][2] 
      \cro[0][1]\indexNF[2][1][] 
      \cro[0][0]\li[2][0]
      \dli[0][0]\cu[1][0]
    }[scale=\scl]
  	\&
  	\tikzpic{
      \uli[0][3]\ca[1][3]
      \cro[0][2]\li[2][2]
      \cro[0][1]\indexNF[2][1][] 
      \cro[0][0]\li[2][0]
      \dli[0][0]\cu[1][0]
    }[scale=\scl]
    \&
  	\tikzpic{
      \uli[0][3]\ca[1][3]
      \li[2][2]
      \LOTRtwo[0][1]\indexNF[2][1][] 
      \cro[0][0]\li[2][0]
      \dli[0][0]\cu[1][0]
    }[scale=\scl]
    \\
  	\tikzpic{
      \ca[0][3]\uli[2][3]
      \li[0][2]\cro[1][2] 
      \indexNF[2][1][] 
      \LOTRtwo[0][0]\li[2][0]
      \dli[0][0]\cu[1][0]
    }[scale=\scl]
    \&
  	\tikzpic{
      \uli[0][3]\ca[1][3]
      \cro[0][2]\li[2][2]
      \indexNF[2][1][] 
      \LOTRtwo[0][0]\li[2][0]
      \dli[0][0]\cu[1][0]
    }[scale=\scl]
  	\arrow[curve={height=-30pt}, from=1-1, to=1-3]
  	\arrow[curve={height=30pt}, from=1-1, to=3-1]
  	\arrow[from=1-2, to=1-1]
  	\arrow[from=1-2, to=1-3]
  	\arrow[draw=none, from=1-2, to=2-1]
  	\arrow[tail reversed, dottedcd, from=1-3, to=2-3]
  	\arrow[from=2-1, to=1-1]
  	\arrow[from=2-1, to=3-1]
  	\arrow["\ulcorner"{description, pos=0}, draw=none, from=2-2, to=1-1]
  	\arrow[from=2-2, to=1-2]
  	\arrow["\urcorner"{description, pos=0}, draw=none, from=2-2, to=1-3]
  	\arrow[from=2-2, to=2-1]
  	\arrow[from=2-2, to=2-3]
  	\arrow["\llcorner"{description, pos=0}, draw=none, from=2-2, to=3-1]
  	\arrow[from=2-2, to=3-2]
  	\arrow[""{name=0, anchor=center, inner sep=0}, curve={height=-30pt}, tail reversed, dottedcd, from=2-3, to=3-2]
  	\arrow[tail reversed, dottedcd, from=3-2, to=3-1]
  	\arrow[draw=none, from=2-2, to=0,"\usebox{\tempboxa}"{description}]
  \end{tikzcd}
  \endgroup
  \\[2ex]
  \begingroup
  \savebox{\tempboxa}{\tikzpic{
    \uli[0][2]\ca[1][2]
    \cro[0][1]\li[2][1]
    \cro[0][0]\li[2][0]
    \cu\dli[2][0]
    \draw[branch_overlay] (-.3,-1+.6) rectangle (1+.3,2);
  }[scale=.4]}
  \def\scl{.6}
  \begin{tikzcd}[ampersand replacement=\&]
  	\tikzpic{ 
      \ca[0][2]\uli[2][2]
      \li[0][1]\cro[1][1]
      \cro[0][0]\li[2][0]
      \cu\dli[2][0]
    }[scale=\scl*1.2]
    \&
  	\tikzpic{ 
      \uli[0][2]\ca[1][2]
      \cro[0][1]\li[2][1]
      \cro[0][0]\li[2][0]
      \cu\dli[2][0]
    }[scale=\scl]
    \&
  	\tikzpic{ 
      \uli[0][2]\ca[1][2]
      \li[2][1]
      \LOTRtwo[0][0]\li[2]
    }[scale=\scl]
    \\
  	\tikzpic{ 
      \ca[0][2]\uli[2][2]
      \li[0][1]\cro[1][1]
      \LOTdk[0][0]\li[2][0]
      \dli[2][0]
    }[scale=\scl]
  	\&
    \tikzpic{ 
      \uli[0][2]\ca[1][2]
      \cro[0][1]\li[2][1]
      \LOTdk[0][0]\li[2][0]
      \cu\dli[2][0]
    }[scale=\scl]
  	\arrow[curve={height=-30pt}, from=1-1, to=1-3]
  	\arrow[from=1-1, to=2-1]
  	\arrow[from=1-2, to=1-1]
  	\arrow[from=1-2, to=1-3]
  	\arrow[draw=none, from=1-2, to=2-1, "\llcorner"{description, pos=0},]
  	\arrow[from=1-2, to=2-2]
  	\arrow[from=2-1, to=2-2, tail reversed, dottedcd]
  	\arrow[""{name=0},curve={height=30pt}, tail reversed, dottedcd, from=2-2, to=1-3]
  	\arrow["\usebox{\tempboxa}"{description}, draw=none, from=1-2, to=0]
  \end{tikzcd}
  \endgroup
  \qquad
  \begingroup
  \savebox{\tempboxa}{\tikzpic{
    \ca[0][2]
    \li[0][1]\li[1][1]\ca[2][1]
    \cro[0][0]\li[2][0]\li[3][0]
    \dli[0][0]\cu[1][0]\dli[3][0]
    \draw[branch_overlay] (-.3,-1+.6) rectangle (3+.3,1.4);
  }[scale=.4]}
  \def\scl{.6}
  \begin{tikzcd}[ampersand replacement=\&]
  	\tikzpic{ 
      \ca[0][2]
      \li[0][1]\li[1][1]\ca[2][1]
      \li\cro[1][0]\li[3][0]
      \cu\dli[2][0]\dli[3][0]
    }[scale=\scl*1.2]
    \&
  	\tikzpic{ 
      \ca[0][2]
      \li[0][1]\li[1][1]\ca[2][1]
      \cro[0][0]\li[2][0]\li[3][0]
      \dli[0][0]\cu[1][0]\dli[3][0]
    }[scale=\scl]
    \&
  	\tikzpic{ 
      \LOTuk[0][1]
      \ca[2][0]\li[0][0]\li[1][0]
      \dli[3][0]\cu[1][0]\dli[0][0]
    }[scale=\scl]
    \\
  	\tikzpic{ 
      \ca[1][2]
      \AntiDiag[0][1]\ca[1][1]\Diag[2][1]
      \cu[0][1]\cro[2][0]
    }[scale=\scl]
  	\&
    \tikzpic{ 
      \cro\ca[0][1]
    }[scale=\scl]
  	\&
    \tikzpic{ 
      \LOTuk
    }[scale=\scl]
  	\arrow[curve={height=-30pt}, from=1-1, to=1-3]
  	\arrow[from=1-1, to=2-1]
  	\arrow[from=1-2, to=1-1]
  	\arrow[from=1-2, to=1-3]
  	\arrow[draw=none, from=1-2, to=2-1, "\usebox{\tempboxa}"{description}, ]
  	\arrow[from=1-2, to=2-2]
  	\arrow["\lrcorner"{description, pos=0},draw=none, from=1-2, to=0]
  	\arrow[from=2-1, to=2-2]
  	\arrow[from=2-2, to=2-3]
  	\arrow[tail reversed, dottedcd, from=1-3, to=2-3]
  \end{tikzcd}
  \endgroup
  \end{gather*}
  \caption{Reduction of critical branchings using horizontal symmetries, part I. One of them has a box: it denotes either a straight line or a crossing.}
  \label{fig:horizontal_symmetries_A}
\end{figure}
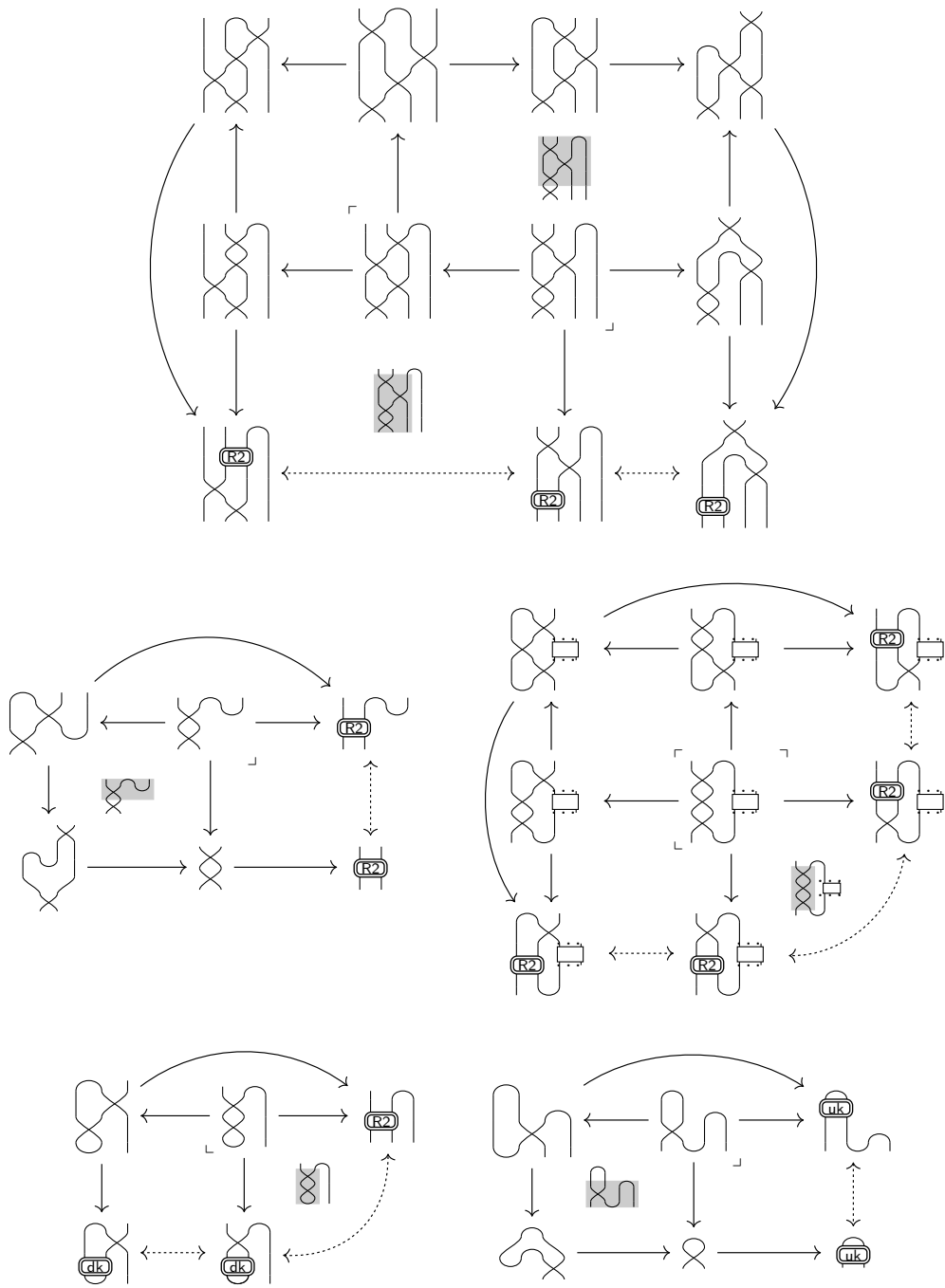

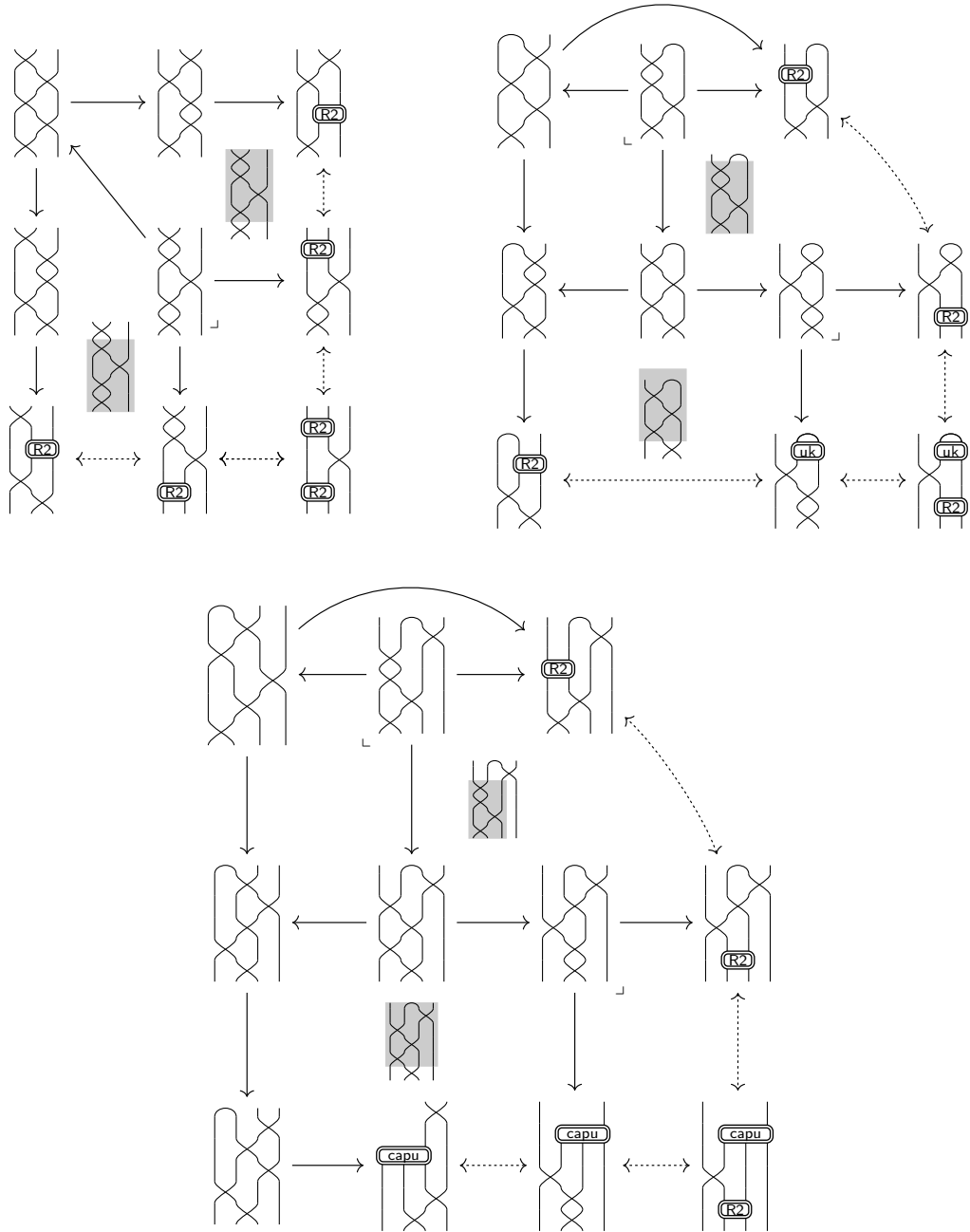
\begin{figure}[p]
  \centering
  \begin{gather*}
  \begingroup
  \savebox{\tempboxa}{\tikzpic{
	\cro[0][4]\li[2][4]
    \cro[0][3]\li[2][3]
    \li[0][2]\cro[1][2]
    \cro[0][1]\li[2][1]
    \cro[0][0]\li[2][0]
    \draw[branch_overlay] (-.3,1) rectangle (2+.3,5);
  }[scale=.5]}
  \savebox{\tempboxb}{\tikzpic{
    \cro[0][4]\li[2][4]
    \cro[0][3]\li[2][3]
    \li[0][2]\cro[1][2]
    \cro[0][1]\li[2][1]
    \cro[0][0]\li[2][0]
    \draw[branch_overlay] (-.3,0) rectangle (2+.3,4);
  }[scale=.5]}
  \def\scl{.6}
  \begin{tikzcd}[ampersand replacement=\&]
    \tikzpic{
      \cro[0][4]\li[2][4]
      \li[0][3]\cro[1][3]
      \cro[0][2]\li[2][2]
      \li[0][1]\cro[1][1]
      \cro[0][0]\li[2][0]
    }[scale=\scl]
    \&
    \tikzpic{
      \li[0][4]\cro[1][4]
      \cro[0][3]\li[2][3]
      \li[0][2]\cro[1][2]
      \li[0][1]\cro[1][1]
      \cro[0][0]\li[2][0]
    }[scale=\scl]
    \&
    \tikzpic{
      \li[0][4]\cro[1][4]
      \cro[0][3]\li[2][3]
      \li[0][2]
      \li[0][1]\LOTRtwo[1][1]
      \cro[0][0]\li[2][0]
    }[scale=\scl]
    \\
    \tikzpic{
      \cro[0][4]\li[2][4]
      \li[0][3]\cro[1][3]
      \li[0][2]\cro[1][2]
      \cro[0][1]\li[2][1]
      \li[0][0]\cro[1][0]
    }[scale=\scl]
    \&
    \tikzpic{
      \cro[0][4]\li[2][4]
      \cro[0][3]\li[2][3]
      \li[0][2]\cro[1][2]
      \cro[0][1]\li[2][1]
      \cro[0][0]\li[2][0]
    }[scale=\scl]
    \&
    \tikzpic{
      \li[2][4]
      \LOTRtwo[0][3]\li[2][3]
      \li[0][2]\cro[1][2]
      \cro[0][1]\li[2][1]
      \cro[0][0]\li[2][0]
    }[scale=\scl]
    \\
    \tikzpic{
      \cro[0][4]\li[2][4]
      \li[0][3]
      \li[0][2]\LOTRtwo[1][2]
      \cro[0][1]\li[2][1]
      \li[0][0]\cro[1][0]
    }[scale=\scl]
    \&
    \tikzpic{
      \cro[0][4]\li[2][4]
      \cro[0][3]\li[2][3]
      \li[0][2]\cro[1][2]
      \li[2][1]
      \LOTRtwo[0][0]\li[2][0]
    }[scale=\scl]
    \&
    \tikzpic{
      \li[2][4]
      \LOTRtwo[0][3]\li[2][3]
      \li[0][2]\cro[1][2]
      \li[2][1]
      \LOTRtwo[0][0]\li[2][0]
    }[scale=\scl]
    \arrow[from=1-1, to=1-2]
    \arrow[from=1-1, to=2-1]
    \arrow[from=1-2, to=1-3]
    \arrow[tail reversed, dottedcd,from=1-3, to=2-3]
    \arrow[from=2-1, to=3-1]
    \arrow[from=2-2, to=1-1]
    \arrow[draw=none, from=2-2, to=1-3,"{\usebox{\tempboxa}}"{description}]
    \arrow[from=2-2, to=2-3]
    \arrow[draw=none, from=2-2, to=3-1,"{\usebox{\tempboxb}}"{description}]
    \arrow[from=2-2, to=3-2]
    \arrow["\lrcorner"{description, pos=0}, draw=none, from=2-2, to=3-3]
    \arrow[tail reversed, dottedcd, from=2-3, to=3-3]
    \arrow[tail reversed, dottedcd, from=3-1, to=3-2]
    \arrow[tail reversed, dottedcd, from=3-2, to=3-3]
    \arrow[tail reversed, dottedcd, from=3-2, to=3-3]
  \end{tikzcd}
  \endgroup
  \mspace{80mu}
  \begingroup
  \savebox{\tempboxa}{\tikzpic{ 
      \uli[0][4]\ca[1][4]
      \cro[0][3]\li[2][3]
      \li[0][2]\cro[1][2]
      \cro[0][1]\li[2][1]
      \li[0][0]\cro[1][0]
      \draw[branch_overlay] (-.3,1) rectangle (2+.3,5);
    }[scale=.5]}
  \savebox{\tempboxb}{\tikzpic{ 
      \uli[0][4]\ca[1][4]
      \cro[0][3]\li[2][3]
      \cro[0][2]\li[2][2]
      \li[0][1]\cro[1][1]
      \cro[0][0]\li[2][0]
      \draw[branch_overlay] (-.3,0) rectangle (2+.3,4);
    }[scale=.5]}
    \def\scl{.6}
  \begin{tikzcd}[ampersand replacement=\&,row sep=large]
  	\tikzpic{ 
      \ca[0][4]\uli[2][4]
      \li[0][3]\cro[1][3]
      \cro[0][2]\li[2][2]
      \li[0][1]\cro[1][1]
      \cro[0][0]\li[2][0]
    }[scale=\scl*1.2]
    \&
    \tikzpic{ 
      \uli[0][4]\ca[1][4]
      \cro[0][3]\li[2][3]
      \cro[0][2]\li[2][2]
      \li[0][1]\cro[1][1]
      \cro[0][0]\li[2][0]
    }[scale=\scl]
    \&
    \tikzpic{ 
      \uli[0][4]\ca[1][4]
      \LOTRtwo[0][2]\li[2][3]\li[2][2]
      \li[0][1]\cro[1][1]
      \cro[0][0]\li[2][0]
    }[scale=\scl] \& 
    \\
  	\tikzpic{ 
      \ca[0][4]\uli[2][4]
      \li[0][3]\cro[1][3]
      \li[0][2]\cro[1][2]
      \cro[0][1]\li[2][1]
      \li[0][0]\cro[1][0]
    }[scale=\scl]
    \&
  	\tikzpic{ 
      \uli[0][4]\ca[1][4]
      \cro[0][3]\li[2][3]
      \li[0][2]\cro[1][2]
      \cro[0][1]\li[2][1]
      \li[0][0]\cro[1][0]
    }[scale=\scl]
    \&
    \tikzpic{ 
      \uli[0][4]\ca[1][4]
      \li[0][3]\cro[1][3]
      \cro[0][2]\li[2][2]
      \li[0][1]\cro[1][1]
      \li[0][0]\cro[1][0]
    }[scale=\scl]
    \&
    \tikzpic{ 
      \uli[0][4]\ca[1][4]
      \li[0][3]\cro[1][3]
      \cro[0][2]\li[2][2]
      \li[0][1]\li[0][0]\LOTRtwo[1][0]
    }[scale=\scl]
    \\
  	\tikzpic{ 
      \ca[0][4]\uli[2][4]
      \li[0][3]\li[0][2]\LOTRtwo[1][2]
      \cro[0][1]\li[2][1]
      \li[0][0]\cro[1][0]
    }[scale=\scl]
    \&\&
    \tikzpic{ 
      \uli[0][4]\ca[1][4]
      \li[0][3]\LOTuk[1][3]
      \cro[0][2]\li[2][2]
      \li[0][1]\cro[1][1]
      \li[0][0]\cro[1][0]
    }[scale=\scl]
    \&
    \tikzpic{ 
      \uli[0][4]\ca[1][4]
      \li[0][3]\LOTuk[1][3]
      \cro[0][2]\li[2][2]
      \li[0][1]\li[0][0]\LOTRtwo[1][0]
    }[scale=\scl]
  	\arrow[from=1-1, to=2-1]
  	\arrow[from=1-2, to=1-1]
  	\arrow[from=1-2, to=1-3]
  	\arrow[from=1-1, to=1-3, curve={height=-40pt}]
  	\arrow["\llcorner"{description, pos=0}, draw=none, from=1-2, to=2-1]
  	\arrow[from=1-2, to=2-2]
  	\arrow[draw=none, from=1-3, to=2-2,"\usebox{\tempboxb}"{description}]
  	\arrow[curve={height=-12pt}, from=1-3, to=2-4, tail reversed, dottedcd]
  	\arrow[from=2-1, to=3-1]
  	\arrow[from=2-2, to=2-1]
  	\arrow[from=2-2, to=2-3]
  	\arrow[from=2-3, to=2-4]
  	\arrow[from=2-3, to=3-3]
  	\arrow["\lrcorner"{description, pos=0}, draw=none, from=2-3, to=3-4]
  	\arrow[from=2-4, to=3-4, tail reversed, dottedcd]
  	\arrow[""{name=0, anchor=center, inner sep=0}, tail reversed, from=3-1, to=3-3,dottedcd]
  	\arrow[from=3-3, to=3-4, tail reversed, dottedcd]
  	\arrow[draw=none, from=2-2, to=0,"\usebox{\tempboxa}"{description}]
  \end{tikzcd}
  \endgroup
  \\[2ex]
  \begingroup
  \savebox{\tempboxa}{\tikzpic{
    \uli[0][5]\ca[1][5]\uli[3][5]
      \lili[0][4]\cro[2][4]
      \cro[0][3]\lili[2][3]
      \cro[0][2]\lili[2][2]
      \li[0][1]\cro[1][1]\li[3][1]
      \cro\lili[2][0]
      \draw[branch_overlay] (-.3,0) rectangle (2+.3,4);
  }[scale=.4]}
  \savebox{\tempboxb}{\tikzpic{
    \uli[0][5]\ca[1][5]\uli[3][5]
    \lili[0][4]\cro[2][4]
    \cro[0][3]\lili[2][3]
    \li[0][2]\cro[1][2]\li[3][2]
    \cro[0][1]\lili[2][1]
    \li[0][0]\cro[1][0]\li[3][0]
    \draw[branch_overlay] (-.3,1) rectangle (3+.3,5.4);
  }[scale=.4]}
  \def\scl{.6}
  \begin{tikzcd}[ampersand replacement=\&,row sep=huge]
  	\tikzpic{
      \ca[0][5]\uli[2][5]
      \li[0][4]\cro[1][4]
      \cro[0][3]\li[2][3]
      \lili[0][2]\cro[2][2]
      \li[0][1]\cro[1][1]
      \cro\li[2][0]
      \draw[] (3,0) to (3,2);\draw[] (3,3) to (3,5.4);
    }[scale=\scl*1.2]
    \&
    \tikzpic{
      \uli[0][5]\ca[1][5]\uli[3][5]
      \lili[0][4]\cro[2][4]
      \cro[0][3]\lili[2][3]
      \cro[0][2]\lili[2][2]
      \li[0][1]\cro[1][1]\li[3][1]
      \cro\lili[2][0]
      %
    }[scale=\scl]
    \&
    \tikzpic{
      \uli[0][5]\ca[1][5]\uli[3][5]
      \lili[0][4]\cro[2][4]
      \LOTRtwo[0][2]\lili[2][2]\lili[2][3]
      \li[0][1]\cro[1][1]\li[3][1]
      \cro\lili[2][0]
    }[scale=\scl] \& 
    \\
  	\tikzpic{
      \ca[0][5]\uli[2][5]
      \li[0][4]\cro[1][4]
      \lili[0][3]\cro[2][3]
      \li[0][2]\cro[1][2]
      \cro[0][1]\li[2][1]
      \li[0][0]\cro[1][0]
      \draw[] (3,0) to (3,3);\draw[] (3,4) to (3,5.4);
    }[scale=\scl]
    \&
    \tikzpic{
      \uli[0][5]\ca[1][5]\uli[3][5]
      \lili[0][4]\cro[2][4]
      \cro[0][3]\lili[2][3]
      \li[0][2]\cro[1][2]\li[3][2]
      \cro[0][1]\lili[2][1]
      \li[0][0]\cro[1][0]\li[3][0]
      %
    }[scale=\scl]
    \&
    \tikzpic{
      \uli[0][5]\ca[1][5]
      \lili[0][4]\cro[2][4]
      \li[0][3]\cro[1][3]
      \cro[0][2]\li[2][2]
      \li[0][1]\cro[1][1]
      \li[0][0]\cro[1][0]
      \draw[] (3,0) to (3,4);\draw[] (3,5) to (3,5.4);
    }[scale=\scl]
    \&
    \tikzpic{
      \uli[0][5]\ca[1][5]
      \lili[0][4]\cro[2][4]
      \li[0][3]\cro[1][3]
      \cro[0][2]\li[2][2]
      \li[0][1]\li[0][0]\LOTRtwo[1][0]
      \draw[] (3,0) to (3,4);\draw[] (3,5) to (3,5.4);
    }[scale=\scl]
    \\
  	\tikzpic{
      \ca[0][5]\uli[2][5]\uli[3][5]
      \lili[0][4]\cro[2][4]
      \li[0][3]\cro[1][3]\li[3][3]
      \lili[0][2]\cro[2][2]
      \cro[0][1]\li[2][1]
      \li[0][0]\cro[1][0]
      \draw[] (3,0) to (3,2);
    }[scale=\scl]
    \&
    \tikzpic{
      \cro[2][5]
      \LOTcapu[0][2]\li[3][2]\li[3][3]\li[3][4]
      \lili[0][1]\cro[2][1]
      \li[0][0]\cro[1][0]\li[3][0]
    }[scale=\scl]
    \&
    \tikzpic{
      \li[0][5]
      \li[0][4]
      \li[0][3]\LOTcapu[1][3]
      \cro[0][2]\li[2][2]
      \li[0][1]\cro[1][1]
      \li[0][0]\cro[1][0]
      \draw[] (3,0) to (3,3);
    }[scale=\scl]
    \&\tikzpic{
      \li[0][5]
      \li[0][4]
      \li[0][3]\LOTcapu[1][3]
      \cro[0][2]\li[2][2]
      \li[0][1]\li[0][0]\LOTRtwo[1][0]
      \draw[] (3,0) to (3,3);
    }[scale=\scl]
  	\arrow[from=1-1, to=2-1]
  	\arrow[from=1-2, to=1-1]
  	\arrow[from=1-2, to=1-3]
  	\arrow[from=1-1, to=1-3, curve={height=-40pt}]
  	\arrow["\llcorner"{description, pos=0}, draw=none, from=1-2, to=2-1]
  	\arrow[from=1-2, to=2-2]
  	\arrow[draw=none, from=1-3, to=2-2,"\usebox{\tempboxa}"{description}]
  	\arrow[dottedcd, tail reversed, curve={height=-12pt}, from=1-3, to=2-4]
  	\arrow[from=2-1, to=3-1]
  	\arrow[from=2-2, to=2-1]
  	\arrow[from=2-2, to=2-3]
  	\arrow[draw=none, from=2-2, to=3-2,"\usebox{\tempboxb}"{description}]
  	\arrow[from=2-3, to=2-4]
  	\arrow[from=2-3, to=3-3]
  	\arrow["\lrcorner"{description, pos=0}, draw=none, from=2-3, to=3-4]
  	\arrow[dottedcd, tail reversed, from=2-4, to=3-4]
  	\arrow[from=3-1, to=3-2]
  	\arrow[dottedcd, tail reversed, from=3-2, to=3-3]
  	\arrow[dottedcd, tail reversed, from=3-3, to=3-4]
  \end{tikzcd}
  \endgroup
  \end{gather*}
  \caption{Reduction of critical branchings using horizontal symmetries, part II.}
  \label{fig:horizontal_symmetries_B}
\end{figure}

\begin{figure}[p]
  \centering
  \begin{gather*}
  \begingroup
  \savebox{\tempboxa}{\tikzpic{
      \ca[0][1]\li[2][1]
      \cro[0][0] \bli[2][0] 
      \dli\cu[1][0]
    \draw[branch_overlay] (-.3,-.4) rectangle (2+.3,1);
  }[scale=.4]}
  \savebox{\tempboxb}{\tikzpic{
      \ca[0][1]\uli[2][1]
      \drcr[0][0] \li[2][0] 
      \dli\cu[1][0]
    \draw[branch_overlay] (-.3,0) rectangle (1+.3,1.4);
  }[scale=.4]}
  \savebox{\tempboxc}{\tikzpic{
      \ca[0][1]\uli[2][1]
      \cro[0][0] \li[2][0] 
      \li[0][-1]\cu[1][0]
      \bli[0][-2]
    \draw[branch_overlay] (-.3,-.4) rectangle (2+.3,1.4);
  }[scale=.4]}
  \savebox{\tempboxd}{\tikzpic{
      \ca[0][1]\uli[2][1]
      \urcr[0][0] \li[2][0] 
      \dli\cu[1][0]
    \draw[branch_overlay] (-.3,-.4) rectangle (2+.3,1);
  }[scale=.4]}
  \def\scl{.6}
  \begin{tikzcd}[ampersand replacement=\&]
  	\tikzpic{
      \ca[0][1]\uli[2][1]
      \li[0][0] \dlcr[1][0] 
      \cu\dli[2][0]
    }[scale=\scl]
    \&
  	\tikzpic{
      \ca[0][1]\uli[2][1]
      \li[0][0] \urcr[1][0] 
      \cu\dli[2][0]
    }[scale=\scl*1.2]
    \&
  	\tikzpic{
      \ca[0][1]\uli[2][1]
      \cro[0][0] \bli[2][0] 
      \dli\cu[1][0]
    }[scale=\scl] 
    \&
  	\tikzpic{
      \uli[2][1]
      \LOTuk[0][0] \bli[2][0] 
      \dli\cu[1][0]
    }[scale=\scl] 
    \\
  	\tikzpic{
      \ca[0][1]\uli[2][1]
      \bli[0][0] \cro[1][0] 
      \cu\dli[2][0]
    }[scale=\scl] 
  	\&
  	\tikzpic{
      \ca[0][1]\uli[2][1]
      \ulcr[0][0] \li[2][0] 
      \dli\cu[1][0]
    }[scale=\scl]
    \&
  	\tikzpic{
      \ca[0][1]\uli[2][1]
      \drcr[0][0] \li[2][0] 
      \dli\cu[1][0]
    }[scale=\scl] 
    \&
  	\tikzpic{
      \uli[2][1]
      \LOTuk[0][0] \li[2][0] 
      \li[0][-1]\lcu[1][0]
    }[scale=\scl] 
    \\
  	\tikzpic{
      \ca[0][1]\uli[2][1]
      \li[0][0] \ulcr[1][0] 
      \cu\dli[2][0]
    }[scale=\scl] 
    \&
  	\tikzpic{
      \ca[0][1]\uli[2][1]
      \urcr[0][0] \li[2][0] 
      \dli\cu[1][0]
    }[scale=\scl] 
    \&
  	\tikzpic{
      \ca[0][1]\uli[2][1]
      \dlcr[0][0] \li[2][0] 
      \dli\cu[1][0]
    }[scale=\scl] 
    \&
  	\tikzpic{
      \uli[2][1]
      \LOTuk[0][0] \li[2][0] 
      \bli[0][-1]\cu[1][0]
    }[scale=\scl] 
    \\
  	\&\&
  	\tikzpic{
      \ca[0][1]\uli[2][1]
      \li[0][0] \drcr[1][0] 
      \cu\dli[2][0]
    }[scale=\scl] 
    \&
\tikzpic{
      \uli[2][1]
      \LOTuk[0][0] \li[2][0] 
      \bli[0][-1]\cu[1][0]
    }[scale=\scl] 
  	\arrow[from=1-1, to=2-1]
  	\arrow[from=1-2, to=1-1]
  	\arrow[curve={height=-24pt}, from=1-2, to=1-4]
  	\arrow[draw=none, from=1-2, to=2-2,"\usebox{\tempboxa}"{description}]
  	\arrow[from=1-3, to=1-2]
  	\arrow[from=1-3, to=1-4]
  	\arrow["\lrcorner"{description, pos=0}, draw=none, from=1-3, to=2-4]
  	\arrow[tail reversed, dottedcd, from=1-4, to=2-4]
  	\arrow[from=2-1, to=3-1]
  	\arrow[from=2-2, to=2-1]
  	\arrow["\llcorner"{description, pos=0}, draw=none, from=2-2, to=3-1]
  	\arrow[from=2-2, to=3-2]
  	\arrow[from=1-3, to=2-3]
  	\arrow[from=2-3, to=2-2]
  	\arrow[from=2-3, to=2-4]
  	\arrow[draw=none, from=2-3, to=3-3,"\usebox{\tempboxb}"{description}]
  	\arrow[tail reversed, dottedcd, from=2-4, to=3-4]
  	\arrow[""{name=0, anchor=center, inner sep=0}, curve={height=12pt}, from=3-1, to=4-3]
  	\arrow[from=3-2, to=3-1]
  	\arrow[from=3-2, to=3-3]
  	\arrow[from=3-3, to=3-4]
  	\arrow[from=3-3, to=4-3]
  	\arrow[draw=none, from=3-3, to=4-4,"\usebox{\tempboxc}"{description}]
  	\arrow[tail reversed, dottedcd, from=3-4, to=4-4]
  	\arrow[from=4-3, to=4-4]
  	\arrow[draw=none, from=3-3, to=0,"\usebox{\tempboxd}"{description}]
  \end{tikzcd}
  \endgroup
  \\[2ex]
  \begingroup
  \savebox{\tempboxa}{\tikzpic{
    \uli[0][2]\ca[1][2]
    \cro[0][1]\li[2][1]
    \dlcr[0][0]\li[2][0] 
    \draw[branch_overlay] (-.3,0) rectangle (1+.3,2);
  }[scale=.4]}
  \savebox{\tempboxb}{\tikzpic{
    \uli[0][2]\ca[1][2]
    \drcr[0][1]\li[2][1]
    \cro[0][0]\li[2][0] 
    \draw[branch_overlay] (-.3,1) rectangle (2+.3,2.4);
  }[scale=.4]}
  \savebox{\tempboxc}{\tikzpic{
    \bli[0][3]
    \li[0][2]\ca[1][2]
    \cro[0][1]\li[2][1]
    \cro[0][0]\li[2][0] 
    \draw[branch_overlay] (-.3,0) rectangle (2+.3,2.4);
  }[scale=.4]}
  \def\scl{.6}
  \begin{tikzcd}[ampersand replacement=\&]
  	\tikzpic{ 
      \ca[0][2]\uli[2][2]
      \li[0][1]\cro[1][1]
      \dlcr[0][0]\li[2][0] 
    }[scale=\scl*1.2]
    \&
  	\tikzpic{ 
      \uli[0][2]\ca[1][2]
      \cro[0][1]\li[2][1]
      \dlcr[0][0]\li[2][0] 
    }[scale=\scl]
    \&
  	\tikzpic{ 
      \uli[0][2]\ca[1][2]
      \li[2][1]
      \LOTRtwo[0][0]\li[2][0] 
      \bli[0][-1]\lili[1][-1]
    }[scale=\scl]
    \\
  	\tikzpic{ 
      \ca[0][2]\uli[2][2]
      \li[0][1]\cro[1][1]
      \bul[1][1]
      \cro[0][0]\li[2][0] 
    }[scale=\scl]
  	\&
    \tikzpic{ 
      \uli[0][2]\ca[1][2]
      \cro[0][1]\li[2][1]
      \bul[1][1]
      \cro[0][0]\li[2][0] 
    }[scale=\scl]
    \\
  	\tikzpic{ 
      \ca[0][2]\uli[2][2]
      \li[0][1]\urcr[1][1]
      \cro[0][0]\li[2][0] 
    }[scale=\scl]
  	\&
    \tikzpic{ 
      \uli[0][2]\ca[1][2]
      \ulcr[0][1]\li[2][1]
      \cro[0][0]\li[2][0] 
    }[scale=\scl]
  	\&
    \tikzpic{ 
      \bli[0][2]\ca[1][2]
      \li[2][1]
      \LOTRtwo[0][0]\li[2][0] 
    }[scale=\scl]
    \\
    \&\&
  	\tikzpic{ 
      \bli[2][3]
      \LOTcapu[0][0]
    }[scale=\scl]
  	\arrow[curve={height=-30pt}, from=1-1, to=1-3]
  	\arrow[from=1-2, to=1-1]
  	\arrow[from=1-2, to=1-3]
  	\arrow[from=1-1, to=2-1]
  	\arrow[from=1-2, to=2-2]
  	\arrow[from=2-2, to=2-1]
  	\arrow[from=1-3, to=3-3, tail reversed, dottedcd]
  	\arrow[from=2-1, to=3-1]
  	\arrow[from=2-2, to=3-2]
  	\arrow[from=3-2, to=3-1]
  	\arrow[from=3-2, to=3-3]
  	\arrow[from=3-3, to=4-3, tail reversed, dottedcd]
  	\arrow[draw=none, from=1-2, to=2-1, "\llcorner"{description, pos=0},]
  	\arrow["\usebox{\tempboxa}"{description}, draw=none, from=1-2, to=3-3]
  	\arrow["\usebox{\tempboxb}"{description}, draw=none, from=2-2, to=3-1]
  	\arrow[""{name=0,pos=.6},curve={height=30pt}, from=3-1, to=4-3]
  	\arrow["\usebox{\tempboxc}"{description}, draw=none, from=3-3, to=0]
  \end{tikzcd}
  \endgroup
  \qquad
  \begingroup
  \savebox{\tempboxa}{\tikzpic{
    \uli[0][2]\ca[1][2]
    \cro[0][1]\li[2][1]
    \drcr[0][0]\li[2][0] 
    \draw[branch_overlay] (-.3,0) rectangle (1+.3,2);
  }[scale=.4]}
  \savebox{\tempboxb}{\tikzpic{
    \uli[0][2]\ca[1][2]
    \dlcr[0][1]\li[2][1]
    \cro[0][0]\li[2][0] 
    \draw[branch_overlay] (-.3,1) rectangle (2+.3,2.4);
  }[scale=.4]}
  \savebox{\tempboxc}{\tikzpic{
    \li[0][2]\ca[1][2]
    \cro[0][1]\li[2][1]
    \cro[0][0]\li[2][0] 
    \lili[0][-1]\bli[2][-1]
    \draw[branch_overlay] (-.3,0) rectangle (2+.3,2.4);
  }[scale=.4]}
  \def\scl{.6}
  \begin{tikzcd}[ampersand replacement=\&]
  	\tikzpic{ 
      \ca[0][2]\uli[2][2]
      \li[0][1]\cro[1][1]
      \drcr[0][0]\li[2][0] 
    }[scale=\scl*1.2]
    \&
  	\tikzpic{ 
      \uli[0][2]\ca[1][2]
      \cro[0][1]\li[2][1]
      \drcr[0][0]\li[2][0] 
    }[scale=\scl]
    \&
  	\tikzpic{ 
      \uli[0][2]\ca[1][2]
      \li[2][1]
      \LOTRtwo[0][0]\li[2][0] 
      \li[0][-1]\bli[1][-1]\li[2][-1]
    }[scale=\scl]
    \\
  	\tikzpic{ 
      \ca[0][2]\uli[2][2]
      \li[0][1]\cro[1][1]
      \bul[0][1]
      \cro[0][0]\li[2][0] 
    }[scale=\scl]
  	\&
    \tikzpic{ 
      \uli[0][2]\ca[1][2]
      \cro[0][1]\li[2][1]
      \bul[0][1]
      \cro[0][0]\li[2][0] 
    }[scale=\scl]
    \\
  	\tikzpic{ 
      \ca[0][2]\uli[2][2]
      \li[0][1]\ulcr[1][1]
      \cro[0][0]\li[2][0] 
    }[scale=\scl]
  	\&
    \tikzpic{ 
      \uli[0][2]\ca[1][2]
      \urcr[0][1]\li[2][1]
      \cro[0][0]\li[2][0] 
    }[scale=\scl]
  	\&
    \tikzpic{ 
      \li[0][2]\lca[1][2]
      \li[2][1]
      \LOTRtwo[0][0]\li[2][0] 
    }[scale=\scl]
    \\
  	\tikzpic{ 
      \ca[0][2]\uli[2][2]
      \li[0][1]\drcr[1][1]
      \cro[0][0]\li[2][0] 
    }[scale=\scl]
  	\&
    \tikzpic{ 
      \uli[0][2]\rca[1][2]
      \cro[0][1]\li[2][1]
      \cro[0][0]\li[2][0] 
    }[scale=\scl]
  	\&
    \tikzpic{ 
      \li[0][2]\rca[1][2]
      \li[2][1]
      \LOTRtwo[0][0]\li[2][0] 
    }[scale=\scl]
    \\
    \&\&
  	\tikzpic{ 
      \LOTcapu[0][0]
      \lili[0][-1]\bli[2][-1]
    }[scale=\scl]
  	\arrow[curve={height=-30pt}, from=1-1, to=1-3]
  	\arrow[from=1-2, to=1-1]
  	\arrow[from=1-2, to=1-3]
  	\arrow[from=1-1, to=2-1]
  	\arrow[from=1-2, to=2-2]
  	\arrow[from=2-2, to=2-1]
  	\arrow[from=1-3, to=3-3, tail reversed, dottedcd]
  	\arrow[from=2-1, to=3-1]
  	\arrow[from=2-2, to=3-2]
  	\arrow[from=3-2, to=3-3]
  	\arrow[from=3-1, to=4-1]
  	\arrow[from=3-2, to=4-2]
  	\arrow[from=4-2, to=4-1]
  	\arrow[from=4-2, to=4-3]
  	\arrow[from=4-3, to=5-3, tail reversed, dottedcd]
  	\arrow[from=3-3, to=4-3, tail reversed, dottedcd]
  	\arrow[draw=none, from=1-2, to=2-1, "\llcorner"{description, pos=0},]
  	\arrow["\usebox{\tempboxa}"{description}, draw=none, from=1-2, to=3-3]
  	\arrow["\usebox{\tempboxb}"{description}, draw=none, from=3-2, to=3-1]
  	\arrow["\lrcorner"{description, pos=0}, draw=none, from=3-2, to=4-3]
  	\arrow[""{name=0,pos=.6},curve={height=30pt}, from=4-1, to=5-3]
  	\arrow["\usebox{\tempboxc}"{description}, draw=none, from=4-3, to=0]
  \end{tikzcd}
  \endgroup
  \end{gather*}
  \caption{Reduction of critical branchings using horizontal symmetries, part III.}
  \label{fig:horizontal_symmetries_C}
\end{figure}
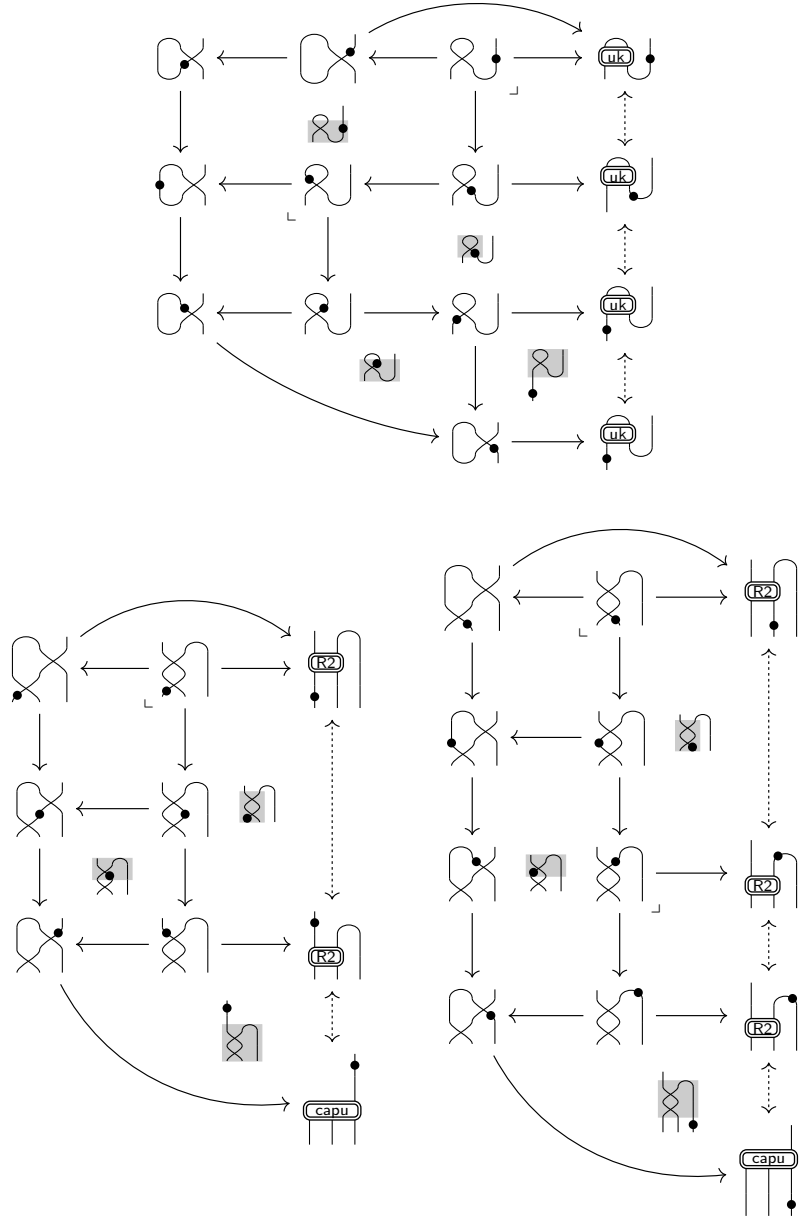

\begin{figure}
  \centering
\begin{gather*}
  \begingroup
    \savebox{\tempboxa}{\tikzpic{
      \uli[0][2]\ca[1][2]
      \cro[0][1]\li[2][1]\bul[1][2]
      \cro\li[2][0]
      \dli\cu[1][0]
      \draw[branch_overlay] (-.3,0) rectangle (1+.3,2);
    }[scale=.4]}
    \savebox{\tempboxb}{\tikzpic{
      \ca[0][2]\uli[2][2]
      \li[0][1]\cro[1][1]\bul[0][1]
      \cro\li[2][0]
      \dli\cu[1][0]
      \draw[branch_overlay] (-.3,-.4) rectangle (1+1.3,2.4);
    }[scale=.4]}
    \def\scl{.6}
  \begin{tikzcd}[ampersand replacement=\&]
  	\tikzpic{
      \uli[0][2]\ca[1][2]
      \cro[0][1]\li[2][1]\bigbul[1][2]
      \li\cro[1][0]
      \cu\dli[2][0]
    }[scale=\scl] 
    \&
  	\tikzpic{
      \uli[0][2]\ca[1][2]
      \cro[0][1]\li[2][1]\bul[1][2]
      \cro\li[2][0]
      \dli\cu[1][0]
    }[scale=\scl] 
    \&
  	\tikzpic{
      \uli[0][2]\ca[1][2]
      \li[2][1]\bul[1][2]
      \LOTRtwo\li[2][0]
      \dli\cu[1][0]
    }[scale=\scl] 
    \\
  	\tikzpic{
      \uli[0][2]\ca[1][2]
      \cro[0][1]\li[2][1]\bul[0][1]
      \li\cro[1][0]
      \cu\dli[2][0]
    }[scale=\scl] 
    \&
  	\tikzpic{
      \uli[0][2]\ca[1][2]
      \cro[0][1]\li[2][1]\bul[0][1]
      \cro\li[2][0]
      \dli\cu[1][0]
    }[scale=\scl] 
    \&
  	\tikzpic{
      \uli[0][2]\ca[1][2]
      \cro[0][1]\li[2][1]\bul[2][0]
      \cro\li[2][0]
      \dli\cu[1][0]
    }[scale=\scl] 
  	\\
  	\tikzpic{
      \ca[0][2]\uli[2][2]
      \li[0][1]\cro[1][1]\bul[0][1]
      \li\cro[1][0]
      \cu\dli[2][0]
    }[scale=\scl] 
  	\&
  	\tikzpic{
      \ca[0][2]\uli[2][2]
      \li[0][1]\cro[1][1]\bul[0][1]
      \cro\li[2][0]
      \dli\cu[1][0]
    }[scale=\scl]
    \&
  	\tikzpic{
      \ca[0][2]\uli[2][2]
      \li[0][1]\cro[1][1]\bul[1][0]
      \cro\li[2][0]
      \dli\cu[1][0]
    }[scale=\scl]
    \\
  	\&\&\&
  	\tikzpic{
      \ca[0][2]\uli[2][2]
      \li[0][1]\bul[0][1]
      \li\LOTRtwo[1][0]
      \cu\dli[2][0]
    }[scale=\scl]
  	\arrow[curve={height=-20pt}, from=1-1, to=1-3]
  	\arrow[from=1-1, to=2-1]
  	\arrow[from=1-2, to=1-1]
  	\arrow[from=1-2, to=1-3]
  	\arrow["\llcorner"{description, pos=0}, draw=none, from=1-2, to=2-1]
  	\arrow[from=1-2, to=2-2]
  	\arrow[draw=none, from=1-2, to=2-3, "\usebox{\tempboxa}"{description}]
  	\arrow[curve={height=-30pt}, from=1-3, to=4-4]
  	\arrow[from=2-1, to=3-1]
  	\arrow[from=2-2, to=2-1]
  	\arrow[from=2-2, to=2-3]
  	\arrow["\llcorner"{description, pos=0}, draw=none, from=2-2, to=3-1]
  	\arrow[from=2-2, to=3-2]
  	\arrow["\lrcorner"{description, pos=0}, draw=none, from=2-2, to=3-3]
  	\arrow[from=2-3, to=1-3]
  	\arrow[from=2-3, to=3-3]
  	\arrow[curve={height=30pt}, from=3-1, to=4-4]
  	\arrow[from=3-2, to=3-1]
  	\arrow[from=3-2, to=3-3]
  	\arrow[curve={height=20pt}, from=3-3, to=1-3]
  	\arrow[draw=none, from=3-3, to=4-4, "\usebox{\tempboxb}"{description}]
  \end{tikzcd}
  \endgroup
\end{gather*}  
\caption{Horizontal symmetry for exceptional critical branchings; see the proof of \cref{lem:dot-slide_branching_tamed_congruent}.}
  \label{fig:horizontal_symmetries_D}
\end{figure}
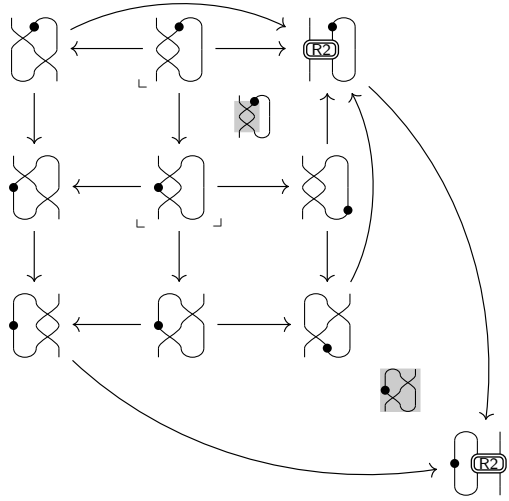

%% file: rewriting_affine_brauer.bib
@book{ABG+_PolygraphsRewritingHigher_2025,
  title = {Polygraphs: From Rewriting to Higher Categories},
  shorttitle = {Polygraphs},
  author = {Ara, Dimitri and Burroni, Albert and Guiraud, Yves and Malbos, Philippe and Métayer, François and Mimram, Samuel},
  date = {2025},
  series = {London {{Mathematical Society Lecture Note Series}}},
  volume = {495},
  publisher = {Cambridge University Press, Cambridge},
  isbn = {978-1-009-49898-2},
  mrnumber = {4866320},
  pagetotal = {xx+648}
}

@article{AGG_QuantizedEnvelopingSuperalgebra_2021,
  title = {Quantized Enveloping Superalgebra of Type {{P}}},
  author = {Ahmed, Saber and Grantcharov, Dimitar and Guay, Nicolas},
  date = {2021},
  journaltitle = {Letters in Mathematical Physics},
  shortjournal = {Lett. Math. Phys.},
  volume = {111},
  number = {3},
  pages = {Paper No. 84, 17},
  issn = {0377-9017,1573-0530},
  doi = {10.1007/s11005-021-01424-y},
  mrnumber = {4277315}
}

@article{Alleaume_RewritingHigherDimensional_2018,
  title = {Rewriting in Higher Dimensional Linear Categories and Application to the Affine Oriented {{Brauer}} Category},
  author = {Alleaume, Clément},
  date = {2018},
  journaltitle = {Journal of Pure and Applied Algebra},
  shortjournal = {J. Pure Appl. Algebra},
  volume = {222},
  number = {3},
  pages = {636--673},
  issn = {0022-4049,1873-1376},
  doi = {10.1016/j.jpaa.2017.05.002},
  mrnumber = {3710720}
}

@unpublished{Barbier_DiagramCategoriesBrauer_2024,
  title = {Diagram Categories of {{Brauer}} Type},
  author = {Barbier, Sigiswald},
  date = {2024-06-26},
  eprint = {2406.18436},
  eprinttype = {arXiv},
  eprintclass = {math},
  url = {http://arxiv.org/abs/2406.18436},
  urldate = {2024-07-08},
  langid = {english},
  pubstate = {prepublished}
}

@book{BCK+_GrobnerShirshovBases_2020,
  title = {Gröbner-{{Shirshov}} Bases},
  author = {Bokut, Leonid and Chen, Yuqun and Kalorkoti, Kyriakos and Kolesnikov, Pavel and Lopatkin, Viktor},
  date = {2020},
  publisher = {World Scientific Publishing Co. Pte. Ltd., Hackensack, NJ},
  doi = {10.1142/9287},
  isbn = {978-981-4619-48-6},
  mrnumber = {4211773},
  pagetotal = {xxii+285}
}

@book{BD_AlgebraicOperads_2016,
  title = {Algebraic Operads},
  author = {Bremner, Murray R. and Dotsenko, Vladimir},
  date = {2016},
  publisher = {CRC Press, Boca Raton, FL},
  isbn = {978-1-4822-4856-2},
  mrnumber = {3642294},
  pagetotal = {xvii+365}
}

@article{BDE+_AffineVWSupercategory_2020,
  title = {The Affine {{VW}} Supercategory},
  author = {Balagović, M. and Daugherty, Z. and Entova-Aizenbud, I. and Halacheva, I. and Hennig, J. and Im, M. S. and Letzter, G. and Norton, E. and Serganova, V. and Stroppel, C.},
  date = {2020},
  journaltitle = {Selecta Mathematica. New Series},
  shortjournal = {Sel. Math. New Ser.},
  volume = {26},
  number = {2},
  pages = {Paper No. 20, 42},
  issn = {1022-1824,1420-9020},
  doi = {10.1007/s00029-020-0541-4},
  mrnumber = {4073971}
}

@article{BE_MonoidalSupercategories_2017,
  title = {Monoidal Supercategories},
  author = {Brundan, Jonathan and Ellis, Alexander P.},
  date = {2017},
  journaltitle = {Communications in Mathematical Physics},
  shortjournal = {Commun. Math. Phys.},
  volume = {351},
  number = {3},
  pages = {1045--1089},
  issn = {0010-3616,1432-0916},
  doi = {10.1007/s00220-017-2850-9},
  mrnumber = {3623246}
}

@article{BE_SuperKacMoody_2017,
  title = {Super {{Kac}}–{{Moody}} 2‐categories},
  author = {Brundan, Jonathan and Ellis, Alexander P.},
  date = {2017-11},
  journaltitle = {Proceedings of the London Mathematical Society},
  shortjournal = {Proc. Lond. Math. Soc.},
  volume = {115},
  number = {5},
  pages = {925--973},
  issn = {0024-6115, 1460-244X},
  doi = {10.1112/plms.12055},
  langid = {english}
}

@article{Bergman_DiamondLemmaRing_1978,
  title = {The Diamond Lemma for Ring Theory},
  author = {Bergman, George M.},
  date = {1978},
  journaltitle = {Advances in Mathematics},
  shortjournal = {Adv. Math.},
  volume = {29},
  number = {2},
  pages = {178--218},
  issn = {0001-8708},
  doi = {10.1016/0001-8708(78)90010-5},
  mrnumber = {506890}
}

@article{BK_OddGrassmannianBimodules_2026,
  title = {Odd {{Grassmannian}} Bimodules and Derived Equivalences for Spin Symmetric Groups},
  author = {Brundan, Jonathan and Kleshchev, Alexander},
  date = {2026-02},
  journaltitle = {Selecta Mathematica},
  shortjournal = {Sel. Math.},
  volume = {32},
  number = {1},
  pages = {9},
  issn = {1022-1824, 1420-9020},
  doi = {10.1007/s00029-025-01102-0},
  langid = {english}
}

@article{Bokut_ImbeddingsSimpleAssociative_1976,
  title = {Imbeddings into Simple Associative Algebras},
  author = {Bokut', L. A.},
  date = {1976},
  journaltitle = {Akademiya Nauk SSSR. Sibirskoe Otdelenie. Institut Matematiki. Algebra i Logika},
  shortjournal = {Akad. Nauk SSSR Sib. Otd. Inst. Mat. Algebra Log.},
  volume = {15},
  number = {2},
  pages = {117--142, 245},
  issn = {0373-9252},
  mrnumber = {506423}
}

@software{BS_BarbierSKAffineBrauerReductionAlgorithm_2026,
  title = {{{BarbierSK}}/{{AffineBrauerReductionAlgorithm}}},
  author = {Barbier, Sigiswald and Schelstraete, Léo},
  date = {2026-09-17T10:18:59Z},
  origdate = {2026-09-17T08:20:39Z},
  url = {https://github.com/BarbierSK/AffineBrauerReductionAlgorithm},
  urldate = {2026-09-17}
}

@book{Buchberger_AlgorithmusAuffindenBasiselemente_1965,
  title = {Ein Algorithmus zum Auffinden der Basiselemente des Restklassenringes nach einem nulldimensionalen Polynomideal},
  author = {Buchberger, Bruno},
  date = {1965},
  publisher = {Innsbruck: Univ. Innsbruck, Mathematisches Institut (Diss.)},
  url = {www.risc.jku.at/people/buchberg/papers/1965-00-00-A.pdf},
  langid = {german},
  zbmath = {6070338},
  zmnumber = {1245.13020}
}

@incollection{Buchberger_CriterionDetectingUnnecessary_1979,
  title = {A Criterion for Detecting Unnecessary Reductions in the Construction of {{Gröbner-bases}}},
  booktitle = {Symbolic and Algebraic Computation ({{EUROSAM}} '79, {{Internat}}. {{Sympos}}., {{Marseille}}, 1979)},
  author = {Buchberger, B.},
  date = {1979},
  series = {Lecture {{Notes}} in {{Comput}}. {{Sci}}.},
  volume = {72},
  pages = {3--21},
  publisher = {Springer, Berlin-New York},
  isbn = {978-3-540-09519-4},
  mrnumber = {575678}
}

@article{BW_CanonicalBasesArising_2018,
  title = {Canonical Bases Arising from Quantum Symmetric Pairs},
  author = {Bao, Huanchen and Wang, Weiqiang},
  date = {2018-09},
  journaltitle = {Inventiones mathematicae},
  shortjournal = {Invent. Math.},
  volume = {213},
  number = {3},
  pages = {1099--1177},
  issn = {0020-9910, 1432-1297},
  doi = {10.1007/s00222-018-0801-5},
  langid = {english}
}

@unpublished{BWW_CategorificationQuasisplitIquantum_2025,
  title = {Categorification of Quasi-Split Iquantum Groups},
  author = {Brundan, Jonathan and Wang, Weiqiang and Webster, Ben},
  date = {2025-05-28},
  eprint = {2505.22929},
  eprinttype = {arXiv},
  eprintclass = {math},
  doi = {10.48550/arXiv.2505.22929},
  langid = {english},
  pubstate = {prepublished}
}

@unpublished{BWW_NilBrauerCategorifiesSplit_2025,
  title = {Nil-{{Brauer}} Categorifies the Split Iquantum Group of Rank One},
  author = {Brundan, Jonathan and Wang, Weiqiang and Webster, Ben},
  date = {2025-09-01},
  eprint = {2305.05877},
  eprinttype = {arXiv},
  eprintclass = {math.QA},
  doi = {10.48550/arXiv.2305.05877},
  pubstate = {prepublished}
}

@article{BWW_NilBrauerCategory_2024,
  title = {The Nil-{{Brauer}} Category},
  author = {Brundan, Jonathan and Wang, Weiqiang and Webster, Benjamin},
  date = {2024-05-23},
  journaltitle = {Annals of Representation Theory},
  shortjournal = {Ann. Represent. Theory},
  volume = {1},
  number = {1},
  pages = {21--58},
  issn = {2704-2081},
  doi = {10.5802/art.2},
  langid = {english}
}

@article{CFL+_QuantumSupergroupsIII_2014,
  title = {Quantum Supergroups {{III}}. {{Twistors}}},
  author = {Clark, Sean and Fan, Zhaobing and Li, Yiqiang and Wang, Weiqiang},
  date = {2014},
  journaltitle = {Communications in Mathematical Physics},
  shortjournal = {Commun. Math. Phys.},
  volume = {332},
  number = {1},
  pages = {415--436},
  issn = {0010-3616,1432-0916},
  doi = {10.1007/s00220-014-2071-4},
  mrnumber = {3253707}
}

@article{CH_QuantumSupergroupsBraid_2016,
  title = {Quantum Supergroups {{V}}. {{Braid}} Group Action},
  author = {Clark, Sean and Hill, David},
  date = {2016},
  journaltitle = {Communications in Mathematical Physics},
  shortjournal = {Commun. Math. Phys.},
  volume = {344},
  number = {1},
  pages = {25--65},
  issn = {0010-3616,1432-0916},
  doi = {10.1007/s00220-016-2630-y},
  mrnumber = {3493137}
}

@unpublished{Chung_CanonicalBasesArising_2021,
  author = {Chung, Christopher},
  date = {2021-07-13},
  eprint = {2107.06322},
  eprinttype = {arXiv},
  eprintclass = {math.QA},
  doi = {10.48550/arXiv.2107.06322},
  pubstate = {prepublished},
  title = {Canonical Bases Arising from {$\imath$}quantum Covering Groups}
}

@unpublished{Chung_SerrePresentation$imath$quantum_2019,
  author = {Chung, Christopher},
  date = {2019-12-17},
  eprint = {1912.09281},
  eprinttype = {arXiv},
  eprintclass = {math.RT},
  doi = {10.48550/arXiv.1912.09281},
  pubstate = {prepublished},
  title = {A {{Serre}} Presentation for the {$\imath$}quantum Covering Groups}
}

@article{CHW_QuantumSupergroupsFoundations_2013,
  title = {Quantum Supergroups {{I}}. {{Foundations}}},
  author = {Clark, Sean and Hill, David and Wang, Weiqiang},
  date = {2013-12},
  journaltitle = {Transformation Groups},
  shortjournal = {Transform. Groups},
  volume = {18},
  number = {4},
  pages = {1019--1053},
  issn = {1083-4362, 1531-586X},
  doi = {10/gpdzwd},
  langid = {english}
}

@article{CHW_QuantumSupergroupsII_2014,
  title = {Quantum Supergroups {{II}}. {{Canonical}} Basis},
  author = {Clark, Sean and Hill, David and Wang, Weiqiang},
  date = {2014},
  journaltitle = {Representation Theory of the American Mathematical Society},
  shortjournal = {Represent. Theory Am. Math. Soc.},
  volume = {18},
  number = {9},
  pages = {278--309},
  issn = {1088-4165},
  doi = {10.1090/S1088-4165-2014-00453-9},
  langid = {english}
}

@article{Clark_QuantumSupergroupsIV_2014,
  title = {Quantum Supergroups {{IV}}: The Modified Form},
  shorttitle = {Quantum Supergroups {{IV}}},
  author = {Clark, Sean},
  date = {2014},
  journaltitle = {Mathematische Zeitschrift},
  shortjournal = {Math. Z.},
  volume = {278},
  number = {1--2},
  pages = {493--528},
  issn = {0025-5874,1432-1823},
  doi = {10.1007/s00209-014-1324-4},
  mrnumber = {3267588}
}

@article{Coulembier_PeriplecticBrauerAlgebra_2018,
  title = {The Periplectic {{Brauer}} Algebra},
  author = {Coulembier, Kevin},
  date = {2018-09},
  journaltitle = {Proceedings of the London Mathematical Society},
  shortjournal = {Proc. Lond. Math. Soc.},
  volume = {117},
  number = {3},
  pages = {441--482},
  issn = {0024-6115, 1460-244X},
  doi = {10.1112/plms.12137},
  langid = {english}
}

@article{CSW_QuantumSupergroupsVI_2019,
  title = {Quantum Supergroups {{VI}}: Roots of 1},
  shorttitle = {Quantum Supergroups {{VI}}},
  author = {Chung, Christopher and Sale, Thomas and Wang, Weiqiang},
  date = {2019-12-01},
  journaltitle = {Letters in Mathematical Physics},
  shortjournal = {Lett. Math. Phys.},
  volume = {109},
  number = {12},
  pages = {2753--2777},
  issn = {1573-0530},
  doi = {10.1007/s11005-019-01209-4},
  langid = {english}
}

@article{DEL_SuperRewritingTheory_2024,
  author = {Dupont, Benjamin and Ebert, Mark and Lauda, Aaron D.},
  date = {2024-07-03},
  journaltitle = {Quantum Topology},
  shortjournal = {Quantum Topol.},
  volume = {16},
  number = {3},
  pages = {569--653},
  issn = {1663-487X, 1664-073X},
  doi = {10.4171/qt/185},
  title = {Super Rewriting Theory and Nondegeneracy of Odd Categorified {$\mathfrak{sl}_{2}$}}
}

@article{DK_GrobnerBasesOperads_2010,
  title = {Gröbner Bases for Operads},
  author = {Dotsenko, Vladimir and Khoroshkin, Anton},
  date = {2010},
  journaltitle = {Duke Mathematical Journal},
  shortjournal = {Duke Math. J.},
  volume = {153},
  number = {2},
  pages = {363--396},
  issn = {0012-7094,1547-7398},
  doi = {10.1215/00127094-2010-026},
  mrnumber = {2667136}
}

@unpublished{DN_CategorificationInfinitedimensional$mathfraksl_2$modules_2021,
  shorttitle = {Categorification of Infinite-Dimensional \$\textbackslash mathfrak\{sl\}\_2\$-Modules and Braid Group 2-Actions {{I}}},
  author = {Dupont, Benjamin and Naisse, Grégoire},
  date = {2021-03-26},
  eprint = {2103.14760},
  eprinttype = {arXiv},
  eprintclass = {math},
  url = {http://arxiv.org/abs/2103.14760},
  urldate = {2022-10-21},
  langid = {english},
  pubstate = {prepublished},
  title = {Categorification of Infinite-Dimensional {$\mathfrak{sl}_2$}-Modules and Braid Group 2-Actions {{I}}: Tensor Products}
}

@article{Dupont_RewritingModuloIsotopies_2021,
  title = {Rewriting modulo Isotopies in {{Khovanov-Lauda-Rouquier}}'s Categorification of Quantum Groups},
  author = {Dupont, Benjamin},
  date = {2021},
  journaltitle = {Advances in Mathematics},
  shortjournal = {Adv. Math.},
  volume = {378},
  pages = {Paper No. 107524, 75},
  issn = {0001-8708,1090-2082},
  doi = {10.1016/j.aim.2020.107524},
  mrnumber = {4186573}
}

@article{Dupont_RewritingModuloIsotopies_2022,
  title = {Rewriting modulo Isotopies in Pivotal Linear (2,2)-Categories},
  author = {Dupont, Benjamin},
  date = {2022},
  journaltitle = {Journal of Algebra},
  shortjournal = {J. Algebra},
  volume = {601},
  pages = {1--53},
  issn = {0021-8693,1090-266X},
  doi = {10.1016/j.jalgebra.2022.02.006},
  mrnumber = {4397096}
}

@article{DV_NormalizationPlanarString_2022,
  title = {Normalization for Planar String Diagrams and a Quadratic Equivalence Algorithm},
  author = {Delpeuch, Antonin and Vicary, Jamie},
  date = {2022},
  journaltitle = {Logical Methods in Computer Science},
  shortjournal = {Log. Methods Comput. Sci.},
  volume = {18},
  number = {1},
  pages = {Paper No. 10, 38},
  issn = {1860-5974},
  doi = {10.46298/lmcs-18(1:10)2022},
  mrnumber = {4369873}
}

@article{EL_OddCategorification$U_qmathfraksl_2$_2016,
  author = {Ellis, Alexander P. and Lauda, Aaron D.},
  date = {2016},
  journaltitle = {Quantum Topology},
  shortjournal = {Quantum Topol.},
  volume = {7},
  number = {2},
  pages = {329--433},
  issn = {1663-487X,1664-073X},
  doi = {10.4171/QT/78},
  mrnumber = {3459963},
  title = {An Odd Categorification of {$U_q(\mathfrak{sl}_2)$}}
}

@article{Elias_DiamondLemmaHecketype_2022,
  title = {A Diamond Lemma for {{Hecke-type}} Algebras},
  author = {Elias, Ben},
  date = {2022},
  journaltitle = {Transactions of the American Mathematical Society},
  shortjournal = {Trans. Am. Math. Soc.},
  volume = {375},
  number = {3},
  pages = {1883--1915},
  issn = {0002-9947,1088-6850},
  doi = {10.1090/tran/8554},
  mrnumber = {4378083}
}

@article{FORM4,
  title = {{{FORM}} Version 4.0},
  author = {Kuipers, J. and Ueda, T. and Vermaseren, J.A.M. and Vollinga, J.},
  date = {2013},
  journaltitle = {Computer Physics Communications},
  shortjournal = {Comput. Phys. Commun.},
  volume = {184},
  number = {5},
  pages = {1453--1467},
  issn = {0010-4655},
  doi = {10.1016/j.cpc.2012.12.028}
}

@unpublished{FORM5,
  title = {{{FORM}} Version 5.0},
  author = {Davies, J. and Kaneko, T. and Marinissen, C. and Ueda, T. and Vermaseren, J.A.M.},
  date = {2026-01-27},
  eprint = {2601.19982},
  eprinttype = {arXiv},
  url = {http://https://arxiv.org/abs/2601.19982},
  langid = {english},
  pubstate = {prepublished}
}

@article{GHM_ConvergentPresentationsPolygraphic_2019,
  title = {Convergent Presentations and Polygraphic Resolutions of Associative Algebras},
  author = {Guiraud, Yves and Hoffbeck, Eric and Malbos, Philippe},
  date = {2019},
  journaltitle = {Mathematische Zeitschrift},
  shortjournal = {Math. Z.},
  volume = {293},
  number = {1--2},
  pages = {113--179},
  issn = {0025-5874,1432-1823},
  doi = {10.1007/s00209-018-2185-z},
  mrnumber = {4002273}
}

@article{GM_HigherdimensionalCategoriesFinite_2009,
  title = {Higher-Dimensional Categories with Finite Derivation Type},
  author = {Guiraud, Yves and Malbos, Philippe},
  date = {2009},
  journaltitle = {Theory and Applications of Categories},
  shortjournal = {Theory Appl. Categ.},
  volume = {22},
  pages = {No. 18, 420--478},
  issn = {1201-561X},
  mrnumber = {2559651}
}

@article{GM_PolygraphsFiniteDerivation_2018,
  title = {Polygraphs of Finite Derivation Type},
  author = {Guiraud, Yves and Malbos, Philippe},
  date = {2018},
  journaltitle = {Mathematical Structures in Computer Science. A Journal in the Applications of Categorical, Algebraic and Geometric Methods in Computer Science},
  shortjournal = {Math. Struct. Comput. Sci. J. Appl. Categ. Algebr. Geom. Methods Comput. Sci.},
  volume = {28},
  number = {2},
  pages = {155--201},
  issn = {0960-1295,1469-8072},
  doi = {10.1017/S0960129516000220},
  mrnumber = {3742562}
}

@article{Guetta_HomologyCategoriesPolygraphic_2021,
  title = {Homology of Categories via Polygraphic Resolutions},
  author = {Guetta, Léonard},
  date = {2021},
  journaltitle = {Journal of Pure and Applied Algebra},
  shortjournal = {J. Pure Appl. Algebra},
  volume = {225},
  number = {10},
  pages = {Paper No. 106688, 33},
  issn = {0022-4049,1873-1376},
  doi = {10.1016/j.jpaa.2021.106688},
  mrnumber = {4207332}
}

@article{Guiraud_TerminationOrdersThreedimensional_2006,
  title = {Termination Orders for Three-Dimensional Rewriting},
  author = {Guiraud, Yves},
  date = {2006},
  journaltitle = {Journal of Pure and Applied Algebra},
  shortjournal = {J. Pure Appl. Algebra},
  volume = {207},
  number = {2},
  pages = {341--371},
  issn = {0022-4049,1873-1376},
  doi = {10.1016/j.jpaa.2005.10.011},
  mrnumber = {2254890}
}

@article{HS_ShorteningCurvesSurfaces_1994,
  title = {Shortening Curves on Surfaces},
  author = {Hass, Joel and Scott, Peter},
  date = {1994},
  journaltitle = {Topology. An International Journal of Mathematics},
  shortjournal = {Topol. Int. J. Math.},
  volume = {33},
  number = {1},
  pages = {25--43},
  issn = {0040-9383},
  doi = {10.1016/0040-9383(94)90033-7},
  mrnumber = {1259513}
}

@unpublished{JLS+_BasisSchurWeylDuality_2025,
  title = {A Basis and {{Schur-Weyl}} Duality for the Loop {{Hecke}} Algebra},
  author = {Janssens, Geoffrey and Lacabanne, Abel and Schelstraete, Léo and Vaz, Pedro},
  date = {2025-07-17},
  eprint = {2507.12839},
  eprinttype = {arXiv},
  eprintclass = {math},
  doi = {10.48550/arXiv.2507.12839},
  pubstate = {prepublished}
}

@article{KKO_SupercategorificationQuantumKacMoody_2013,
  title = {Supercategorification of Quantum {{Kac-Moody}} Algebras},
  author = {Kang, Seok-Jin and Kashiwara, Masaki and Oh, Se-jin},
  date = {2013},
  journaltitle = {Advances in Mathematics},
  shortjournal = {Adv. Math.},
  volume = {242},
  pages = {116--162},
  issn = {0001-8708,1090-2082},
  doi = {10.1016/j.aim.2013.04.008},
  mrnumber = {3055990}
}

@article{KL_CategorificationQuantum$sln$_2010,
  author = {Khovanov, Mikhail and Lauda, Aaron},
  date = {2010},
  journaltitle = {Quantum Topology},
  shortjournal = {Quantum Topol.},
  volume = {1},
  pages = {1--92},
  issn = {1663-487X},
  doi = {10.4171/QT/1},
  langid = {english},
  title = {A Categorification of Quantum {$sl(n)$}}
}

@article{KL_DiagrammaticApproachCategorification_2009,
  title = {A Diagrammatic Approach to Categorification of Quantum Groups {{I}}},
  author = {Khovanov, Mikhail and Lauda, Aaron},
  date = {2009-07-28},
  journaltitle = {Representation Theory of the American Mathematical Society},
  shortjournal = {Represent. Theory Am. Math. Soc.},
  volume = {13},
  number = {14},
  pages = {309--347},
  issn = {1088-4165},
  doi = {10.1090/S1088-4165-09-00346-X},
  langid = {english}
}

@article{Kolb_QuantumSymmetricKac_2014,
  title = {Quantum Symmetric {{Kac}}–{{Moody}} Pairs},
  author = {Kolb, Stefan},
  date = {2014-12},
  journaltitle = {Advances in Mathematics},
  shortjournal = {Adv. Math.},
  volume = {267},
  pages = {395--469},
  issn = {00018708},
  doi = {10.1016/j.aim.2014.08.010},
  langid = {english}
}

@article{Lafont_AlgebraicTheoryBoolean_2003,
  title = {Towards an Algebraic Theory of {{Boolean}} Circuits},
  author = {Lafont, Yves},
  date = {2003},
  journaltitle = {Journal of Pure and Applied Algebra},
  shortjournal = {J. Pure Appl. Algebra},
  volume = {184},
  number = {2--3},
  pages = {257--310},
  issn = {0022-4049,1873-1376},
  doi = {10.1016/S0022-4049(03)00069-0},
  mrnumber = {2004976}
}

@article{Lauda_CategorificationQuantum$sl2$_2010,
  author = {Lauda, Aaron D.},
  date = {2010-12},
  journaltitle = {Advances in Mathematics},
  shortjournal = {Adv. Math.},
  volume = {225},
  number = {6},
  pages = {3327--3424},
  issn = {00018708},
  doi = {10.1016/j.aim.2010.06.003},
  langid = {english},
  title = {A Categorification of Quantum {$sl(2)$}}
}

@article{Letzter_SymmetricPairsQuantized_1999,
  title = {Symmetric {{Pairs}} for {{Quantized Enveloping Algebras}}},
  author = {Letzter, Gail},
  date = {1999-10},
  journaltitle = {Journal of Algebra},
  shortjournal = {J. Algebra},
  volume = {220},
  number = {2},
  pages = {729--767},
  issn = {00218693},
  doi = {10.1006/jabr.1999.8015},
  langid = {english}
}

@article{Metayer_ResolutionsPolygraphs_2003,
  title = {Resolutions by Polygraphs},
  author = {Métayer, François},
  date = {2003},
  journaltitle = {Theory and Applications of Categories},
  shortjournal = {Theory Appl. Categ.},
  volume = {11},
  pages = {No. 7, 148--184},
  issn = {1201-561X},
  mrnumber = {1988395}
}

@article{Mimram_3dimensionalRewritingTheory_2014,
  title = {Towards 3-Dimensional Rewriting Theory},
  author = {Mimram, Samuel},
  date = {2014},
  journaltitle = {Logical Methods in Computer Science},
  shortjournal = {Log. Methods Comput. Sci.},
  volume = {10},
  number = {2},
  pages = {2:1, 47},
  issn = {1860-5974},
  doi = {10.2168/LMCS-10(2:1)2014},
  mrnumber = {3194020}
}

@unpublished{Rouquier_2KacMoodyAlgebras_2008,
  title = {2-{{Kac--Moody}} Algebras},
  author = {Rouquier, Raphael},
  date = {2008-12-30},
  eprint = {0812.5023},
  eprinttype = {arXiv},
  doi = {10.48550/arXiv.0812.5023},
  pubstate = {prepublished}
}

@article{RS_AffineBrauerCategory_2019,
  author = {Rui, Hebing and Song, Linliang},
  date = {2019},
  journaltitle = {Mathematische Zeitschrift},
  shortjournal = {Math. Z.},
  volume = {293},
  number = {1--2},
  pages = {503--550},
  issn = {0025-5874,1432-1823},
  doi = {10.1007/s00209-018-2207-x},
  mrnumber = {4002288},
  title = {Affine {{Brauer}} Category and Parabolic Category {$\mathcal{O}$} in Types {{B}}, {{C}}, {{D}}}
}

@article{RS_Periplectic$q$BrauerCategory_2025,
  author = {Rui, Hebing and Song, Linliang},
  date = {2025},
  journaltitle = {Journal of Algebra},
  shortjournal = {J. Algebra},
  volume = {661},
  pages = {82--122},
  issn = {0021-8693,1090-266X},
  doi = {10.1016/j.jalgebra.2024.06.043},
  mrnumber = {4788650},
  title = {The Periplectic {$q$}-{{Brauer}} Category}
}

@thesis{Schelstraete_OddKhovanovHomology_2024,
  type = {phdthesis},
  title = {Odd {{Khovanov}} Homology, Higher Representation Theory and Higher Rewriting Theory},
  author = {Schelstraete, Léo},
  date = {2024},
  institution = {UCLouvain},
  url = {https://doi.org/10.48550/arXiv.2410.11405},
  addendum = {available at arXiv:2410.11405},
  langid = {english},
  pagetotal = {264}
}

@article{Schelstraete_RewritingModuloDiagrammatic_2026,
  title = {Rewriting modulo in Diagrammatic Algebras and Application to Categorification},
  author = {Schelstraete, Léo},
  date = {2026-10-01},
  journaltitle = {Advances in Mathematics},
  shortjournal = {Adv. Math.},
  volume = {502},
  pages = {111198},
  issn = {0001-8708},
  doi = {10.1016/j.aim.2026.111198}
}

@article{Shirshov_AlgorithmicProblemsLie_1962,
  title = {Some algorithmic problems for Lie algebras},
  author = {Shirshov, Anatoliı̆ Illarionovich},
  date = {1962},
  volume = {3},
  pages = {292--296},
  issn = {0037-4474},
  fjournal = {Sibirskiı̆ Matematicheskiı̆ Zhurnal},
  langid = {russian},
  zbmath = {3170376},
  zmnumber = {0104.26004}
}

@unpublished{SV_OddKhovanovHomology_2023,
  title = {Odd {{Khovanov}} Homology and Higher Representation Theory},
  author = {Schelstraete, Léo and Vaz, Pedro},
  date = {2023-11-24},
  eprint = {2311.14394},
  eprinttype = {arXiv},
  eprintclass = {math},
  doi = {10.48550/arXiv.2311.14394},
  pubstate = {prepublished}
}

@unpublished{SW_BubblesAffineBrauer_2024,
  title = {Bubbles in the Affine {{Brauer}} and {{Kauffman}} Categories},
  author = {Savage, Alistair and Webster, Ben},
  date = {2024-08-13},
  eprint = {2408.07000},
  eprinttype = {arXiv},
  eprintclass = {math},
  url = {http://arxiv.org/abs/2408.07000},
  urldate = {2024-08-14},
  langid = {english},
  pubstate = {prepublished}
}

@incollection{WB_CriterionEliminatingUnnecessary_1986,
  title = {A Criterion for Eliminating Unnecessary Reductions in the {{Knuth-Bendix}} Algorithm},
  booktitle = {Algebra, Combinatorics and Logic in Computer Science, {{Vol}}.\textbackslash{} {{I}}, {{II}} ({{Győr}}, 1983)},
  author = {Winkler, F. and Buchberger, B.},
  date = {1986},
  series = {Colloq. {{Math}}. {{Soc}}. {{János Bolyai}}},
  volume = {42},
  pages = {849--869},
  publisher = {North-Holland, Amsterdam},
  isbn = {978-0-444-87869-4},
  mrnumber = {875913}
}

@unpublished{Webster_UnfurlingKhovanovLaudaRouquierAlgebras_2024,
  title = {Unfurling {{Khovanov-Lauda-Rouquier}} Algebras},
  author = {Webster, Ben},
  date = {2024-08-06},
  eprint = {1603.06311},
  eprinttype = {arXiv},
  eprintclass = {math},
  doi = {10.48550/arXiv.1603.06311},
  langid = {english},
  pubstate = {prepublished}
}
